%% file: main.tex
\def\isarxiv{1} 

\ifdefined\isarxiv
\documentclass[11pt]{article}

\usepackage[numbers]{natbib}

\else
\documentclass{article}
\usepackage{iclr2025_conference,times}
\iclrfinalcopy

\usepackage{hyperref}
\usepackage{url}
\usepackage{booktabs}

\title{Faster Algorithms for Structured Linear and Kernel Support Vector Machines}
\author{Yuzhou Gu \\
New York University \\
\texttt{yuzhougu@nyu.edu}
\And
Zhao Song \\
Simons Institute for the Theory of Computing, UC Berkeley\\
\texttt{magic.linuxkde@gmail.com}
\And
Lichen Zhang \\
MIT CSAIL \\
\texttt{lichenz@csail.mit.edu}
}

\newcommand{\new}{\marginpar{NEW}}

\fi

\usepackage{amsmath}
\usepackage{amsthm}
\usepackage{amssymb}
\usepackage{algorithm}
\usepackage{algpseudocode}
\usepackage{grffile}
\usepackage{wrapfig,epsfig}
\usepackage{url}
\usepackage{xcolor}
\usepackage{epstopdf}
\usepackage{accents}
\usepackage{enumerate}
\usepackage{enumitem}

\usepackage{bbm}
\usepackage{dsfont}

\allowdisplaybreaks

\ifdefined\isarxiv

\usepackage{tikz}
\usepackage{hyperref}  
\hypersetup{colorlinks=true,citecolor=blue,linkcolor=blue} 
\usetikzlibrary{arrows}
\usepackage[margin=1in]{geometry}

\else

\usepackage{microtype}
\usepackage{hyperref}
\usepackage{url}
\definecolor{mydarkblue}{rgb}{0,0.08,0.45}
\hypersetup{colorlinks=true, citecolor=mydarkblue,linkcolor=mydarkblue}

\fi

\theoremstyle{plain}
\newtheorem{theorem}{Theorem}[section]
\newtheorem{lemma}[theorem]{Lemma}

\newtheorem{proposition}[theorem]{Proposition}
\newtheorem{corollary}[theorem]{Corollary}

\newtheorem{claim}[theorem]{Claim}

\theoremstyle{definition}
\newtheorem{definition}[theorem]{Definition}
\newtheorem{assumption}[theorem]{Assumption}

\newtheorem{remark}[theorem]{Remark}

\DeclareFontFamily{U}{mathx}{}
\DeclareFontShape{U}{mathx}{m}{n}{<-> mathx10}{}
\DeclareSymbolFont{mathx}{U}{mathx}{m}{n}
\DeclareMathAccent{\widecheck}{0}{mathx}{"71}

\newcommand{\wh}{\widehat}
\newcommand{\wt}{\widetilde}

\newcommand{\ov}{\overline}

\newcommand{\R}{\mathbb{R}}

\renewcommand{\d}{\mathrm{d}}

\newcommand{\seth}{\textsf{SETH}}

\renewcommand{\hat}{\wh}
\newcommand{\Tmat}{{\cal T}_{\mathrm{mat}}}
\newcommand{\triangleup}{\vartriangle}

\DeclareMathOperator{\OPT}{OPT}

\DeclareMathOperator{\poly}{poly}

\DeclareMathOperator{\nnz}{nnz}

\DeclareMathOperator{\tr}{tr}

\def\new{\mathrm{new}}

\def\cB{\mathcal{B}}

\def\cI{\mathcal{I}}
\def\cK{\mathcal{K}}
\def\cN{\mathcal{N}}

\def\cS{\mathcal{S}}
\def\cT{\mathcal{T}}
\def\cW{\mathcal{W}}

\def\bN{\mathbb{N}}

\def\tot{\mathrm{tot}}
\def\apx{\mathrm{apx}}

\makeatletter
\newcommand*{\RN}[1]{\expandafter\@slowromancap\romannumeral #1@}
\makeatother

\usepackage{lineno}

\begin{document}

\ifdefined\isarxiv

\date{}

\title{Faster Algorithms for Structured Linear and Kernel Support Vector Machines}

\author{Yuzhou Gu\thanks{\texttt{yuzhougu@nyu.edu}. New York University.}
\and
Zhao Song\thanks{\texttt{magic.linuxkde@gmail.com}. Simons Institute for the Theory of Computing, UC Berkeley.}
\and 
Lichen Zhang\thanks{\texttt{lichenz@mit.edu}. Massachusetts Institute of Technology.}
}

\else

\fi

\ifdefined\isarxiv
\begin{titlepage}
  \maketitle
  \begin{abstract}
\input{abstract}

  \end{abstract}
  \thispagestyle{empty}
\end{titlepage}

{\hypersetup{linkcolor=black}
\tableofcontents
}
\newpage

\else
\maketitle
\begin{abstract}
\input{abstract}
\end{abstract}

\fi

\input{intro} 
\input{tech} 


\input{ack}
\ifdefined\isarxiv
\else
\bibliography{ref}
\bibliographystyle{iclr2025_conference}

\fi

\ifdefined\isarxiv
\else
\newpage
\onecolumn
\appendix
\fi
\paragraph{Roadmap.} In Section~\ref{sec:preli}, we present some basic definitions and tools that will be used in the remainder of the paper. In Section~\ref{sec:svm-formulations}, we introduce a few more SVM formulations, including classification and distribution estimation. In Section~\ref{sec:general_qp}, we show convex quadratic programming can be reduced to convex empirical risk minimization. Therefore, under boundedness assumptions and given an eigendecomposition of the objective matrix, LCQP on $n$ variables can be solved approximately in expected time $(n^{\omega+o(1)}+n^{2.5-\alpha/2+o(1)}+n^{2+1/6+o(1)})\cdot\wt O(\log(nLR/\epsilon))$ by the algorithm of~\cite{lsz19}, with the equality constraints holding only approximately, and QCQP with $m$ quadratic constraints in the same time with $(m+1)n$ in place of $n$ (Propositions~\ref{prop:lcqp} and~\ref{prop:qcqp}). Here $L$ and $R$ bound the coefficients of the program and the entries of some optimal solution, and $\alpha$ denotes the dual matrix multiplication exponent, not the dual vector of the SVM. In Section~\ref{sec:treewidth} and~\ref{sec:rank}, we prove results on low-treewidth and low-rank QPs, respectively. In Section~\ref{sec:robust_ipm}, we present a robust IPM framework for QPs, generalize beyond LPs and convex ERMs with linear objective. In Section~\ref{sec:gaussian_kernel}, we present our algorithms for Gaussian kernel SVMs, with complementary lower bound.
\input{preli}

\input{svm_form}

\input{reduction}

\input{ds-treewidth}

\input{ds-rank}

\input{ipm}

\input{kernel}

\ifdefined\isarxiv
\newpage
\bibliographystyle{alpha}
\bibliography{ref}
\else
\fi



\end{document}

%% file: abstract.tex
Quadratic programming is a ubiquitous prototype in convex programming. Many machine learning problems can be formulated as quadratic programming, including the famous Support Vector Machines (SVMs). Linear and kernel SVMs have been among the most popular models in machine learning over the past three decades, prior to the deep learning era.

Generally, a quadratic program has an input size of $\Theta(n^2)$, where $n$ is the number of variables. Assuming the Strong Exponential Time Hypothesis ($\textsf{SETH}$), it is known that no $O(n^{2-o(1)})$ time algorithm exists when the quadratic objective matrix is positive semidefinite (Backurs, Indyk, and Schmidt, NeurIPS'17). However, problems such as SVMs usually admit much smaller input sizes: one is given $n$ data points, each of dimension $d$, and $d$ is oftentimes much smaller than $n$. Furthermore, the SVM program has only $O(1)$ equality linear constraints. This suggests that faster algorithms are feasible, provided the program exhibits certain structures.

In this work, we design the first nearly-linear time algorithm for solving quadratic programs whenever a low-rank factorization of the quadratic objective matrix is given, the number of linear constraints is small, the coefficients of the objective and some optimal solution are bounded by $\poly(n)$, and some feasible point has all entries between $1/\poly(n)$ and $\poly(n)$. Consequently, we obtain results for SVMs:
\begin{itemize}
    \item For linear SVM when the input data is $d$-dimensional, our randomized algorithm runs in expected time $\widetilde O(nd^{(\omega+1)/2}\log(n/\epsilon))$ when $d$ is much smaller than $n$, precisely $(d+1)^{\omega-1}\le n$, and the data points have norms at most $\poly(n)$ and margin at least $1/\poly(n)$, where $\omega$ is the matrix multiplication exponent, which is less than $2.258$. Its output satisfies the equality constraint only approximately.
   
    \item For Gaussian kernel SVM, when the data dimension $d = {\color{black}O(\log n)}$ and the squared dataset radius is $o(\log n/\log\log n)$, our randomized algorithm runs in expected time $n^{1+o(1)}$ for any accuracy $\epsilon\ge 1/\poly(n)$, when the dual variables are restricted to a box $[0,R]^n$ with $R\le\poly(n)$, as in soft-margin SVMs or when a hard-margin optimum is known to lie in such a box. Its output satisfies the equality constraint up to error $3\epsilon$. We also prove that, assuming \seth{}, when the squared dataset radius is at least $\Omega(\log^2 n)$, then $\Omega(n^{2-o(1)})$ time is required on some inputs to reach error $\exp(-\Theta(\log^2 n))$ for the program without the bias term. This improves upon the prior best lower bound in both the dimension $d$ and the squared dataset radius.
\end{itemize}

%% file: intro.tex
\section{Introduction}

Quadratic programming (QP) represents a class of convex optimization problems that optimize a quadratic objective over the intersection of an affine subspace and the non-negative orthant\footnote{There are classes of QPs with quadratic constraints as well. However, in this paper, we focus on cases where the constraints are linear.}. QPs naturally extend linear programming by incorporating a quadratic objective, and they find extensive applications in operational research, theoretical computer science, and machine learning~\citep{ktk79,w99,gt00,ghn01,pu04,ct06}. The quadratic objective introduces challenges: QPs with a general (not necessarily positive semidefinite) symmetric quadratic objective matrix are NP-hard to solve~\citep{sahni1974computationally,pardalos1991quadratic}. When the quadratic objective matrix is positive semidefinite, the problem becomes weakly polynomial-time solvable, as it can be reduced to convex empirical risk minimization~\citep{lsz19} (refer to Section~\ref{sec:general_qp} for further discussion).

Formally, the QP problem is defined as follows:
\begin{definition}[Quadratic Programming] \label{def:lcqp}
Given an $n \times n$ symmetric, positive semidefinite objective matrix $Q$, a vector $c \in \R^n$, and a polytope described by a pair $(A \in \R^{m \times n}, b \in \R^m)$, the linearly constrained quadratic programming (LCQP) or simply quadratic programming (QP) problem seeks to solve the following optimization problem:
\begin{align} \label{eqn:qp}
\min_{x \in \R^n} ~ & ~ \frac{1}{2} x^\top Q x + c^\top x \\
\mathrm{s.t.~} & ~ A x = b \nonumber\\
& ~ x \geq 0. \nonumber
\end{align}
\end{definition}

A classic application of QP is the Support Vector Machine (SVM) problem~\citep{bgv92,cv95}. In SVMs, a dataset ${x_1,\ldots,x_n} \in \R^d$ is provided, along with corresponding labels ${y_1,\ldots,y_n} \in \{\pm 1\}$. The objective is to identify a hyperplane that separates the two groups of points with opposite labels, while maintaining a large margin from both. Remarkably, this popular machine learning problem can be formulated as a QP and subsequently solved using specialized QP solvers~\citep{mmrts01}. Thus, advancements in QP algorithms could potentially lead to runtime improvements for SVMs.

Despite its practical and theoretical significance, algorithmic quadratic programming has garnered relatively less attention compared to its close relatives in convex programming, such as linear programming~\citep{cls19,jswz21,b20,sy21}, convex empirical risk minimization~\citep{lsz19,qszz23}, and semidefinite programming~\citep{jkl+20,hjs+22,gs22}. In this work, we take a pioneering step in developing fast and robust interior point-type algorithms for solving QPs. We particularly focus on improving the runtime for high-precision hard- and soft-margin SVMs. For the purposes of this discussion, we will concentrate on hard-margin SVMs, with the understanding that our results naturally extend to soft-margin variants. We begin by introducing the hard-margin linear SVMs:

\begin{definition}[Linear SVM]
\label{def:linear-svm}
Given a dataset $X \in \R^{n \times d}$ and a collection of labels $y_1, \ldots, y_n$ each in ${\pm 1}$, the \emph{linear SVM problem} requires solving the following quadratic program:
\begin{align} \label{eqn:linear-svm}
\max_{\alpha \in \R^n}&~{\bf 1}_n^\top \alpha - \frac{1}{2} \alpha^\top (yy^\top \circ XX^\top) \alpha, \\
{\rm s.t.~} & ~ \alpha^\top y = 0, \nonumber \\
& ~ \alpha \geq 0. \nonumber
\end{align}
where $\circ$ denotes the Hadamard product.
\end{definition}
It should be noted that this formulation is actually the \emph{dual} of the SVM optimization problem. The primal program seeks a vector $w \in \R^d$ such that
\begin{align*}
\min_{w \in \R^d} & ~ \frac{1}{2}\|w\|_2^2, \\
\text{s.t.} & ~ y_i (w^\top x_i - b) \geq 1, \qquad \forall i \in [n],
\end{align*}
where $b \in \R$ is the bias term. Given the solution $\alpha \in \R^n$, one can conveniently convert it to a primal solution: $w^* = \sum_{i=1}^n \alpha_i^* y_i x_i$. At first glance, one might be inclined to solve the primal problem directly, especially in cases where $d \ll n$, as it presents a lower-dimensional optimization problem compared to the dual. The dual formulation becomes particularly advantageous when solving the kernel SVM, which maps features to a high or potentially infinite-dimensional space.

\begin{definition}[Kernel SVM]
\label{def:kernel_svm}
Given a dataset $X \in \R^{n \times d}$ and a positive definite kernel function \({\sf K}:\R^d \times \R^d \rightarrow \R\), let \(K \in \R^{n \times n}\) denote the kernel matrix, where \(K_{i,j} = {\sf K}(x_i, x_j)\). With a collection of labels \(y_1, \ldots, y_n\) each in \(\{\pm 1\}\), the \emph{kernel SVM problem} requires solving the following quadratic program: 
\begin{align}\label{eqn:kernel_svm}
    \max_{\alpha \in \R^n}~&~{\bf 1}_n^\top \alpha - \frac{1}{2} \alpha^\top (yy^\top \circ K) \alpha, \\
    {\rm s.t.~} & ~ \alpha^\top y = 0,\notag \\
    & ~ \alpha \geq 0.\notag
\end{align}
\end{definition}
The positive definite kernel function \({\sf K}\) corresponds to a feature mapping, implying that \({\sf K}(x_i, x_j) = \phi(x_i)^\top \phi(x_j)\) for some \(\phi: \R^d \rightarrow \R^s\). Thus, solving the primal SVM can be viewed as solving the optimization problem on the transformed dataset. However, the primal program's dimension depends on the (transformed) data's dimension \(s\), which can be infinite. Conversely, the dual program, with dimension \(n\), is typically easier to solve. Throughout this paper, when discussing the SVM program, we implicitly refer to the dual quadratic program, not the primal.

One key aspect of the SVM program is its minimal equality constraints. Specifically, for both linear and kernel SVMs, there is only a single equality constraint of the form \(\alpha^\top y = 0\). This constraint arises naturally from the bias term in the primal SVM formulation and its Lagrangian. The limited number of constraints enables us to design QP solvers with favorable dependence on the number of data points \(n\), albeit with a higher dependence on the number of constraints \(m\), thus offering effective end-to-end guarantees for SVMs.

Previous efforts to solve the SVM program typically involve breaking down the large QP into smaller, constant-sized QPs. These algorithms, while easy to implement and well-suited to modern hardware architectures, oftentimes lack tight theoretical analysis and the estimation of iteration count is usually pessimistic~\citep{cl11}. Theoretically,~\cite{j06} systematically analyzed this class of algorithms, demonstrating that to achieve an $\epsilon$-approximation solution, $\wt O(\epsilon^{-2}B \cdot \nnz(A))$ time is sufficient, where $B$ is the squared-radius of the dataset and $\nnz(\cdot)$ denotes the number of nonzero entries. {\color{black}This is subsequently improved in~\cite{sssc11} with a subgradient-based method that runs in $\wt O(\epsilon^{-1}d)$ time.} Unfortunately, the polynomial dependence on the precision $\epsilon^{-1}$ makes them hard to be adapted for even moderately small $\epsilon$. For example, when $\epsilon$ is set to be $10^{-3}$ to account for the usual machine precision errors, these algorithms would require at least {\color{black}$10^3$} iterations to converge.

To develop a high-precision algorithm with $\poly\log(\epsilon^{-1})$ dependence instead of $\poly(\epsilon^{-1})$, we focus on second-order methods for QPs. A variety of approaches have been explored in previous works, including the interior point method~\citep{k84}, active set methods~\citep{m88}, augmented Lagrangian techniques~\citep{dg03}, conjugate gradient, gradient projection, and extensions of the simplex algorithm~\citep{d55,w59,m88}. Our interest is particularly piqued by the interior point method (IPM). Recent advances in the robust IPM framework have led to significant successes for convex programming problems~\citep{cls19,lsz19,b20,jkl+20,blss20,jswz21,sy21,jnw22,hjs+22,gs22,qszz23,syyz23_dp}. These successes are a result of combining robust analysis of IPM with dedicated data structure design.

Applying IPM to solve QPs with a constant number of constraints is not entirely novel; existing work~\citep{fm02} has already adapted IPM to solve the linear SVM problem. However, the runtime of their algorithm is sub-optimal. Each iteration of their algorithm requires multiplying a $d\times n$ matrix with an $n\times d$ matrix in $O(nd^{\omega-1})$ time, where $\omega$ is the matrix multiplication exponent, which is less than $2.258$ over every field of characteristic zero~\citep[Theorem~1.1 and Corollaries~14.1--14.2]{OpenAIMM26}. Moreover, the IPM requires $O(\sqrt n \log(1/\epsilon))$ iterations to converge. This ends up with an overall runtime $O(n^{1.5}d^{\omega-1}\log(1/\epsilon))$, which is super-linear in the dataset size even when the dimension $d$ is small. In practical scenarios where $n$ is usually large, the $n^{1.5}$ dependence becomes prohibitive. Therefore, it is crucial to develop an algorithm with almost- or nearly-linear dependence on $n$ and logarithmic dependence on $\epsilon^{-1}$.

For linear SVM, we propose a nearly-linear time algorithm with high-precision guarantees, applicable when the dimension of the dataset is much smaller than the number of points, under the polynomial bounds stated in the theorem:

\begin{theorem}[Low-rank QP and Linear SVM, informal version of Theorem~\ref{thm:rank-formal}]
\label{thm:main-rank-informal}
Given a quadratic program as defined in Definition~\ref{def:lcqp}, and assuming a low-rank factorization of the quadratic objective matrix $Q = UV^\top$, where $U, V \in \R^{n \times k}$, there exists an algorithm that can solve the program~\eqref{eqn:qp} up to $\epsilon$-error\footnote{We say an algorithm that solves the program up to $\epsilon$-error if it returns an approximate solution vector $\wt \alpha$ whose objective value is at most $\epsilon$ more than the optimal objective value.} in $\wt O(n(k+m)^{(\omega+1)/2}\log(n/\epsilon))$ time, provided $k$ and $m$ are much smaller than $n$, precisely $(k+m)^{\omega-1}\le n$, some optimal solution satisfies $\|x^*\|_\infty\le R$, some $z$ with $Az=b$ satisfies $r\le z_i\le R-r$ for all $i\in[n]$, and $R, 1/r, \|Q\|_{2\to 2}, \|c\|_2\le \poly(n)$ (so that Theorem~\ref{thm:rank-formal} applies after adding the constraints $x\le R\cdot{\bf 1}_n$). The algorithm is randomized and its time bound holds in expectation. The returned $\wt \alpha\ge 0$ has objective value at most $\epsilon$ above the optimum and satisfies $\|A\wt\alpha-b\|_1\le \epsilon(\|A\|_1+\|b\|_1)$. In particular, the equality constraints hold only approximately.

Specifically, for linear SVM (as per Definition~\ref{def:linear-svm}) with $(d+1)^{\omega-1} \leq n$, one can solve program~\eqref{eqn:linear-svm} up to $\epsilon$-error in $\wt O(nd^{(\omega+1)/2}\log(n/\epsilon))$ time under the same conditions. They hold, for example, for a hard-margin SVM whose data points have norms at most $\poly(n)$ and are separable with margin $\gamma\ge 1/\poly(n)$. Then $\|\alpha^*\|_1=1/\gamma^2$, so one may take $R=1/\gamma^2$, and $z_i=R/(2|\{j: y_j=y_i\}|)$ gives $r=R/(2n)$.
\end{theorem}

The assumption $(k+m)^{\omega-1}\le n$ says that the rank $k$ and the number of equality constraints $m$ are much smaller than $n$. We use it to bound the cost of inverting a $k\times k$ matrix and an $m\times m$ matrix in every iteration (Section~\ref{sec:tech:rank}). It holds in our applications, where $k$ and $m$ are much smaller than $n$. The SVM programs have at most two equality constraints, the linear SVM has $k=d$ and $m=1$, and our Gaussian kernel application has $k=n^{o(1)}$.

While a nearly-linear time algorithm for linear SVMs is appealing, most applications look at kernel SVMs as they provide more expressive power to the linear classifier. This poses significant challenge in algorithm design, as forming the kernel matrix exactly would require $\Omega(n^2)$ time. Moreover, the kernel matrix could be full-rank without any structural assumptions, rendering our low-rank QP solver inapplicable. In fact, it has been shown that for data dimension $d=\omega(\log n)$, no algorithm can approximately solve kernel SVM within an error $\exp(-\omega(\log^2 n))$ in time $O(n^{2-o(1)})$, assuming the famous Strong Exponential Time Hypothesis ({\sf SETH})\footnote{{\sf SETH} is a standard complexity theoretical assumption~\citep{ipz98,ip01}. Informally, it states that for every $c<2$ there is a $k$ such that $k$-SAT on $n$ variables and $m$ clauses cannot be solved in time $O(c^n\cdot\poly(m))$.}~\citep{backurs2017fine}.


Conversely, a long line of works aim to speed up computation with the kernel matrix faster than quadratic, especially when the kernel has certain smooth and Lipschitz properties~\citep{acss20,aa22,biksz23,ckw24}. For instance, when kernel functions are sufficiently smooth, efficient approximation using low-degree polynomials is feasible, leading to an approximate low-rank factorization. A prime example is the Gaussian RBF kernel, where~\cite{aa22} showed that for dimension $d = \Theta(\log n)$ and squared dataset radius (defined as $\max_{i, j\in [n]} \|x_i-x_j\|_2^2$) $B = o(\log n)$, for every accuracy $\epsilon\ge 1/\poly(n)$ there exist low-rank matrices $U, V \in \R^{n \times n^{o(1)}}$ such that for any vector $x \in \R^n$, $\|(K - UV^\top)x\|_\infty \leq \epsilon \|x\|_1$. For much smaller $\epsilon$, such as $\exp(-\Theta(\log^2 n))$, the rank of their factorization becomes superpolynomial (Remark~\ref{rem:gaussian_kernel_svm_parameters}). They subsequently develop an algorithm to solve the Batch Gaussian KDE problem in $O(n^{1+o(1)})$ time.

Based on this dichotomy in fast kernel matrix algebra, we establish two results: 1) Solving Gaussian kernel SVM to accuracy $\epsilon\ge 1/\poly(n)$ in expected time $n^{1+o(1)}$, with the equality constraint holding up to $3\epsilon$, is feasible when $d=O(\log n)$, $B = o(\frac{\log n}{\log \log n})$, and either the program is soft-margin with $C\le\poly(n)$ or some optimal solution has entries at most a known $R\le\poly(n)$, and 2) Assuming \seth, no sub-quadratic time algorithm solves the program to error $\exp(-\Theta(\log^2 n))$ for $B = \Omega(\log^2 n)$ in SVMs without bias and $B = \Omega(\log^6 n)$ in SVMs with bias (the latter with a proof that is only sketched). This improves the lower bound established by~\cite{backurs2017fine} in terms of dimension $d$.

\begin{theorem}[Gaussian Kernel SVM, informal version of Theorem~\ref{thm:gaussian_kernel}, Corollary~\ref{cor:gaussian_kernel_almost_linear}, Theorem~\ref{thm:svm_hard_no_bias} and Corollary~\ref{cor:svm_hard_bias} (whose proof is only sketched, following~\cite{backurs2017fine})]
\label{thm:gaussian_kernel_informal}
Given a dataset $X \in \R^{n \times d}$ with dimension $d$ and squared radius denoted by $B$, let ${\sf K}(x_i, x_j) = \exp(-\|x_i - x_j\|_2^2)$ be the Gaussian kernel function. Then, for the kernel SVM problem defined in Definition~\ref{def:kernel_svm},
\begin{itemize}
    \item If ${\color{black}d=O(\log n)}$, $B = o(\frac{\log n}{\log \log n})$, both labels occur, $\epsilon\in[n^{-O(1)},1)$, and some optimal solution $\alpha^*$ satisfies $\|\alpha^*\|_\infty\le R$ for a known $R\le\poly(n)$, then a randomized algorithm outputs, in expected time $n^{1+o(1)}$, a vector $\wh\alpha\in[0,R]^n$ whose objective value is at least the optimum minus $\epsilon$ and which satisfies $|\wh\alpha^\top y|\le 3\epsilon$ (Theorem~\ref{thm:gaussian_kernel} and Corollary~\ref{cor:gaussian_kernel_almost_linear}). The same holds for the soft-margin program (Definition~\ref{def:C_SVC}) with $C\le\poly(n)$, where $R=C$ and no assumption on $\alpha^*$ is needed. For much smaller $\epsilon$ the time bound can be superpolynomial, for example $n^{\Theta(\log\log n)}$ at $d=\Theta(\log n)$, $B=\Theta(1)$ and $\epsilon=\exp(-\Theta(\log^2 n))$ (Remark~\ref{rem:gaussian_kernel_svm_parameters});
    \item If ${\color{black}d=\Omega(\log n)}, B = \Omega(\log^2 n)$, then assuming \seth, any algorithm that solves the Gaussian kernel SVM \emph{without} a bias term up to $\exp(-\Theta(\log^2 n))$ error would require $\Omega(n^{2-o(1)})$ time;
    \item If ${\color{black}d=\Omega(\log n)}, B = \Omega(\log^6 n)$ (Corollary~\ref{cor:svm_hard_bias}, whose proof is only sketched), then assuming \seth, any algorithm that solves the Gaussian kernel SVM \emph{with} a bias term up to $\exp(-\Theta(\log^2 n))$ error would require $\Omega(n^{2-o(1)})$ time.
\end{itemize}
\end{theorem}

To our knowledge, this is the first almost-linear time algorithm for Gaussian kernel SVM even when $d = \log n$ and the radius is small. Our algorithm effectively utilizes the rank-$n^{o(1)}$ factorization of the Gaussian kernel matrix alongside our low-rank QP solver. 

\paragraph{Note on hardness assumptions.} On 26 October 2026, OpenAI~\cite{openai26eth} and Anthropic~\cite{anthropic26eth} independently announced deterministic algorithms, with Lean formalizations, that solve $k$-SAT in time $2^{o(n)}$ for every fixed $k$. As of 27 October 2026, neither proof had been peer reviewed. The OpenAI Lean formalization checks a slightly weaker time bound than its paper states, but this bound is still $2^{o(n)}$. The Anthropic paper noted that no person was known to have checked its Lean development. As of that date, R.~Impagliazzo, R.~Paturi, and R.~Williams were organizing an independent examination of the OpenAI proof. If the proofs are correct, they refute the Exponential Time Hypothesis (\textsf{ETH}), and hence also the Strong Exponential Time Hypothesis (\textsf{SETH}), which implies \textsf{ETH}~\cite{ipz98,ip01}, as both papers say. Our lower bounds (Theorem~\ref{thm:svm_hard_no_bias} and Corollaries~\ref{cor:svm_hard_bias} and~\ref{cor:svm_soft_bias}) assume \textsf{SETH}, and so do the lower bounds of~\cite{backurs2017fine} and~\cite{aa22} that we cite. All of them use \textsf{SETH} only through problems that are hard under the Orthogonal Vectors hypothesis, ours and that of~\cite{aa22} through the approximate Hamming nearest neighbor problem of~\cite{r18}. So they still hold under this hypothesis, which neither announcement claims to refute (Remark~\ref{rem:eth-2026}). Our algorithms do not depend on any of these hypotheses.

\subsection{Related Work}

\paragraph{Support Vector Machines.} SVM, one of the most prominent machine learning models before the rise of deep learning, has a rich literature dedicated to its algorithmic speedup. For linear SVM, \cite{j06} offers a first-order algorithm that solves its QP in nearly-linear time, but with a runtime dependence of $\epsilon^{-2}$, limiting its use in high precision settings. {\color{black} This runtime is later significantly improved by~\cite{sssc11} to $\wt O(\epsilon^{-1}d)$ via a stochastic subgradient descent algorithm.} For SVM classification, existing algorithms such as SVM-Light~\citep{j99}, SMO~\citep{p98}, LIBSVM~\citep{cl11}, and SVM-Torch~\citep{cb01} perform well in high-dimensional data settings. However, their runtime scales super-linearly with $n$, making them less viable for large datasets. Previous investigations into solving linear SVM via interior point methods~\citep{fm02} have been somewhat basic, leading to an overall runtime of $O(n^{1.5}d^{\omega-1}\log(1/\epsilon))$. For a more comprehensive survey on efficient algorithms for SVM, refer to~\cite{cgf+20}. On the hardness side, \cite{backurs2017fine} provides an efficient reduction from the Bichromatic Closest Pair problem to Gaussian kernel SVM, establishing an almost-quadratic lower bound assuming \seth.

\paragraph{Interior Point Method.} The interior point method, a well-established approach for solving convex programs under constraints, was first proposed by \cite{k84} as a (weakly) polynomial-time algorithm for linear programs, later improved by \cite{v89_lp} in terms of runtime. Recent work by~\cite{cls19} used a robust IPM framework to solve linear programs in $\wt O(n^\omega+n^{2.5-\alpha/2}+n^{2+1/6})$ time, where $\alpha$ is the dual matrix multiplication exponent. Under the bounds on $\omega$ and $\alpha$ known when it appeared, this bound was $\wt O(n^\omega)$, the time of matrix multiplication. Subsequent studies~\citep{lsz19,b20,jswz21,sy21,hjs+22,jnw22,qszz23} have further refined their algorithm or applied it to different optimization problems.

\paragraph{Kernel Matrix Algebra.}
Kernel methods, fundamental in machine learning, enable feature mappings to potentially infinite dimensions for $n$ data points in $d$ dimensions. The kernel matrix, a crucial component of kernel methods, often has a prohibitive quadratic size for explicit formation. Recent active research focuses on computing and approximating kernel matrices and related tasks in sub-quadratic time, such as kernel matrix-vector multiplication, spectral sparsification, and Laplacian system solving. The study by \cite{acss20} introduces a comprehensive toolkit for solving these problems in almost-linear time for small dimensions, leveraging techniques like polynomial methods and ultra Johnson-Lindenstrauss transforms. Alternatively, \cite{bimw21,biksz23} reduce various kernel matrix algebra tasks to kernel density estimation (KDE), which recent advancements in KDE data structures~\citep{cs17,bcis18,ckns20} have made more efficient. A recent contribution by \cite{aa22} provides a tighter characterization of the low-degree polynomial approximation for the $e^{-x}$ function, leading to more efficient algorithms for the Batch Gaussian KDE problem. Another line of works aim to apply \emph{oblivious sketching} to the kernel matrix, so that a spectral approximation can be computed in time nearly-linear in the dataset size and a low-rank factorization with rank being the \emph{statistical dimension} of the problem can be computed~\citep{anw14,akk+20,zha+21,wz22}. Unfortunately, for the statistical dimension to be smaller than $n$, one has to add a regularization parameter $\lambda$, which usually appears when one wants to solve a kernel ridge regression. The work~\citep{swyz21} does not require the parameter $\lambda$, but it obtains runtime improvements for $d\geq n$, and the sketching dimension is super-linear in $n$.

\paragraph{Subsequent work.} After the first version of this paper, the recent preprint~\cite{lptwy25} gave an interior point method for empirical risk minimization, $\min c^\top x$ subject to $A^\top x=b$ and $x\in K_1\times\cdots\times K_m$, where $A\in\mathbb{R}^{n\times d}$ and each $K_i$ is a compact convex set of dimension $O(1)$. It solves this problem to high accuracy in time $\widetilde O(nd+d^6\sqrt n)$, and its output satisfies $A^\top x=b$ exactly, while ours satisfies the equality constraints only approximately. Our linear SVM program~\eqref{eqn:linear-svm} fits this form with $d+1$ equality constraints once we add the box constraint $\alpha\le R\cdot{\bf 1}_n$ of Theorem~\ref{thm:main-rank-informal}, the variables $w=\sum_{i=1}^n\alpha_iy_ix_i$ and the epigraph variables $w_j^2/2\le t_j\le T$, where the cap $T\le\operatorname{poly}(n)$ follows from the box and the bound on the data and makes each set $K_i$ compact. Under the assumptions of that theorem, the boundedness and conditioning assumptions of~\cite{lptwy25} hold with polynomially bounded parameters. Both algorithms are randomized interior point methods with logarithmic dependence on $1/\epsilon$. Hence the bound of~\cite{lptwy25} improves upon our $\widetilde O(nd^{(\omega+1)/2}\log(n/\epsilon))$ time when $n\ge d^{11-\omega}$, and it is nearly linear in the input size when $n\ge d^{10}$.

%% file: tech.tex
\section{Technical Overview} \label{sec:tech}

In this section, we provide an overview of the techniques employed in our development of two nearly-linear time algorithms for structured QPs. In Section~\ref{sec:tech:general}, we detail the robust IPM framework, which forms the foundation of our algorithms. Subsequent section, namely Section~\ref{sec:tech:treewidth} and~\ref{sec:tech:rank}, delves into dedicated data structures designed for efficiently solving low-treewidth and low-rank QPs, respectively. Finally, in Section~\ref{sec:tech:kernel}, we discuss the adaptation of these advanced QP solvers for both linear and kernel SVMs.

{\color{black} Due to the heavily-technical nature, we recommend that in the first read, the audience can skip Section~\ref{sec:tech:treewidth} and~\ref{sec:tech:rank}.}

\subsection{General Strategy} \label{sec:tech:general}

Our algorithm is built upon the robust IPM framework, an efficient variant of the primal-dual central path method~\citep{r88}. This framework maintains a primal-dual solution pair $(x,s) \in \R^n \times \R^n$. To understand the central path for QPs, we first consider the central path equations for linear programming (see~\cite{cls19,lsz19} for reference):
\begin{align*}
s/t + \nabla \phi(x) = & ~ \mu, \\
A x = & ~ b, \\
A^\top y + s = & ~ c,
\end{align*}
where $x$ is the primal variable, $s$ is the slack variable, $y$ is the dual variable, $\phi(x)$ is a self-concordant barrier function, and $\mu$ denotes the error. The central path is defined by the trajectory of $(x,s)$ as $t$ approaches $0$. 

In quadratic programming, we modify these equations:
\begin{align*}
s/t + \nabla \phi(x) = & ~ \mu, \\
A x = & ~ b, \\
-Q x + A^\top y + s = & ~ c,
\end{align*}
where $Q$ is the positive semidefinite objective matrix. The key difference in the central path equations for LP and QP is the inclusion of the $-Q x$ term in the third equation, significantly affecting algorithm design.

Fundamentally, IPM is a Newton's method in which we update the variables $x, y$ and $s$ through the second-order information from the self-concordant barrier function. We derive the update rules for QP (detailed derivation in Section~\ref{sec:robust_ipm:derive}):
\begin{align*}
\delta_x &= t {\color{black}M^{-1/2}} (I-P) {\color{black}M^{-1/2}} \delta_\mu,\\
\delta_y &= -t (A{\color{black}M^{-1}} A^\top)^{-1} A {\color{black}M^{-1}} \delta_\mu,\\
\delta_s &= t \delta_\mu - t^2 H {\color{black}M^{-1/2}} (I-P) {\color{black}M^{-1/2}} \delta_\mu,\\
\text{where} \qquad H &= \nabla^2 \phi(x), \qquad {\color{black}M} = Q + t H,\\
P &= {\color{black}M^{-1/2}}  A^\top  (A{\color{black}M^{-1}} A^\top)^{-1} A {\color{black}M^{-1/2}},
\end{align*}
where $\delta_x$, $\delta_y$, $\delta_s$, and $\delta_\mu$ are the incremental steps for $x$, $y$, $s$, and $\mu$, respectively.

The robust IPM approximates these updates rather than computing them exactly. It maintains an approximate primal-dual solution pair $(\ov x, \ov s) \in \R^n \times \R^n$ and computes the steps using this approximation. Provided the approximation is sufficiently accurate, it can be shown (see Section~\ref{sec:robust_ipm} for more details) that the algorithm converges efficiently to the optimal solution along the robust central path.

Therefore, the critical challenge lies in efficiently maintaining $(\ov x, \ov s)$, an approximation to $(x,s)$, when $(x,s)$ evolves following the robust central path steps. The primary difficulty is that explicitly managing the primal-dual solution pair $(x,s)$ is inefficient due to potential dense changes. Such changes can lead to dense updates in $H$, slowing down the computation of steps. The innovative aspect of robust IPM is recognizing that $(x,s)$ are only required at the algorithm's conclusion, not during its execution. Instead, we can identify entries with significant changes and update the approximation $(\ov x, \ov s)$ correspondingly. With IPM's lazy updates, only a nearly-linear number of entries are adjusted throughout the algorithm:
\begin{align*}
    \sum_{t=1}^T \|\ov x^{(t)}-\ov x^{(t-1)}\|_0 + \|\ov s^{(t)}-\ov s^{(t-1)}\|_0 = & ~ \wt O(n\log(1/\epsilon))
\end{align*}
where $T = \wt O(\sqrt{n}\log(1/\epsilon))$ is the number of iterations for IPM convergence. This indicates that, on average, each entry of $\ov x$ and $\ov s$ is updated $\log(1/\epsilon)$ times, facilitating rapid updates to these quantities and, consequently, to $H$.

In the special case where $Q = 0$, the path reverts to the LP case, with ${\color{black}M} = tH$ being a diagonal matrix, allowing for efficient computation and updates of ${\color{black}M^{-1}}$. This simplifies maintaining $A{\color{black}M^{-1}}A^\top$, as updates to ${\color{black}M^{-1}}$ correspond to row and column scaling of $A$. However, in the QP scenario, where ${\color{black}M}$ is symmetric positive semidefinite, maintaining the term $A{\color{black}M^{-1}}A^\top$ becomes more complex. Nevertheless, when the number of constraints is small, as in SVMs, this issue is less problematic. Yet, even with this simplification, the challenge is far from trivial, given the presence of terms like ${\color{black}M^{-1/2}}$ in the robust central path steps. While the matrix Woodbury identity could be considered, it falls short when maintaining a square root term. Despite these hurdles, we construct efficient data structures for ${\color{black}M^{-1/2}}$ maintenance when $Q$ possesses succinct representations, such as low-rank. 

Before diving into the particular techniques for low-rank QPs, we start by exploring the \emph{low-treewidth QPs}, which could be viewed as a structured sparsity condition. It provides valuable insights for the low-rank scenario.


\subsection{Low-Treewidth Setting: How to Leverage Sparsity} \label{sec:tech:treewidth}

Treewidth is parameter for graphs that captures the sparsity pattern. Given a graph $G=(V, E)$ with $n$ vertices and $m$ edges, a \emph{tree decomposition} of $G$ arranges its vertices into bags, which collectively form a tree structure. For any two bags $X_i, X_j$, if a vertex $v$ is present in both, it must also be included in all bags along the path between $X_i$ and $X_j$. Additionally, each pair of adjacent vertices in the graph must be present together in at least one bag. The treewidth $\tau$ is defined as the maximum size of a bag minus one. Intuitively, a graph $G$ with a small treewidth $\tau$ implies a structure akin to a tree. For a formal definition, see Definition~\ref{def:treewidth}. When relating this combinatorial structure to linear algebra, we could treat the quadratic objective matrix $Q$ as a generalized adjacency matrix, where we put a vertex $v_i$ on $i$-th row of $Q$, and put an edge $\{v_i, v_j\}$ whenever the entry $Q_{i,j}$ is nonzero. The low-treewidth structure of the graph corresponds to a sparsity pattern that allows one to compute a column-sparse Cholesky factorization of $Q$. Since ${\color{black}M}=Q+tH$ and $H$ is diagonal, we can decompose ${\color{black}M}=LL^\top$ into sparse Cholesky factors\footnote{Note that adding a non-negative diagonal matrix to $Q$ does not change its sparsity pattern, hence ${\color{black}M}$ also retains the treewidth $\tau$.}.


Under any coordinate update to $\ov x$, ${\color{black}M}$ is updated on only one diagonal entry, enabling efficient updates to $L$. The remaining task is to use this Cholesky decomposition to maintain the central path step.

By expanding the central path equations and substituting ${\color{black}M} = LL^\top$, we derive
\begin{align*}
\delta_x &= t {\color{black}M^{-1/2}} (I-P) {\color{black}M^{-1/2}} \delta_\mu \\
&= t {\color{black}M^{-1}} \delta_\mu - t {\color{black}M^{-1}} A^\top (A {\color{black}M^{-1}} A^\top)^{-1} A {\color{black}M^{-1}} \delta_\mu \\
&= t L^{-\top} L^{-1} \delta_\mu - t L^{-\top} L^{-1} A^\top (A L^{-\top} L^{-1} A^\top)^{-1} A L^{-\top} L^{-1} \delta_\mu,\\
\delta_s &= t\delta_\mu - t^2 H{\color{black}M^{-1/2}}(I-P){\color{black}M^{-1/2}}\delta_\mu \\
&= t\delta_\mu - t^2 H L^{-\top} L^{-1} \delta_\mu + t^2 H L^{-\top} L^{-1} A^\top (A L^{-\top} L^{-1} A^\top)^{-1} A L^{-\top} L^{-1} \delta_\mu.
\end{align*}
Updates to the diagonal of ${\color{black}M}$ do not change $L$'s nonzero pattern, allowing for efficient utilization of the sparse factor and maintenance of $L^{-1}A^\top\in \R^{n\times m}$ and $L^{-1} \delta_\mu\in \R^n$. Terms like $(A L^{-\top} L^{-1} A^\top)^{-1} A L^{-\top} L^{-1} \delta_\mu \in \R^m$ can also be explicitly maintained.

With this approach, we propose the following implicit representation for maintaining $(x,s)$:
\begin{align}
x &= \wh x + H^{-1/2} \cW^\top ( h \beta_x - \wt h \wt \beta_x +\epsilon_x), \label{eqn:treewidth-implicit-x-rep} \\
s &= \wh s + H^{1/2} c_s \beta_{c_s} - H^{1/2} \cW^\top ( h \beta_s - \wt h \wt \beta_s +\epsilon_s), \label{eqn:treewidth-implicit-s-rep}
\end{align}
where $\wh x,\wh s \in \R^{n}$, $\cW = L^{-1} H^{1/2} \in \R^{n \times n}$, $h = L^{-1} \ov\delta_\mu \in \R^{n}$, $c_s = H^{-1/2} \ov\delta_\mu \in \R^{n}$, $\beta_x,\beta_s, \beta_{c_s} \in \R$, $\wt h = L^{-1} A^\top \in \R^{n \times m}$, $\wt \beta_x, \wt \beta_s \in \R^{m}$, $\epsilon_x, \epsilon_s \in \R^{n}$. All quantities except for $\cW$ can be explicitly maintained. For linear programming, the implicit representation is as follows:
\begin{align*}
x &= \wh x + H^{-1/2} \beta_x c_x - H^{-1/2} \cW^\top (\beta_x h + \epsilon_x) \\
s &= \wh s + H^{1/2} \cW^\top(\beta_s h + \epsilon_s),
\end{align*}
with $\cW = L^{-1} A H^{-1/2}$ maintained implicitly and the other terms explicitly.

The representation in \eqref{eqn:treewidth-implicit-x-rep} and \eqref{eqn:treewidth-implicit-s-rep} enables us to maintain the central path step using a combination of ``coefficients'' $h+\wt h\wt \beta_x$ and ``basis'' $\cW^\top$. We need to detect entries of $\ov x$ that deviate significantly from $x$ and capture these changes with $\|H^{1/2}(\ov x-x)\|_2$. We maintain this vector using $x_0+\cW^\top (h+\wt h\wt \beta_x)$. Here, $\cW^\top$ acts as a wavelet basis and the vector $h+\wt h\wt \beta_x$ as its multiscale coefficients. While computing and maintaining $\cW^\top h$ seems challenging, leveraging column-sparsity of $L^{-1}$ is possible through contraction with a vector $v$:
\begin{align*}
    v^\top \cW^\top = & ~ (\cW v)^\top \\
    = & ~ (L^{-1} H^{1/2} v)^\top.
\end{align*}
By applying the Johnson-Lindenstrauss transform (JL) in place of $v$, we can quickly approximate $\|\cW^\top h\|_2$ by maintaining $\Phi \cW^\top$ for a JL matrix $\Phi$. Similarly, we handle $\cW^\top \wt h\wt \beta_x$ by explicitly computing $A^\top \wt \beta_x$ and using the sparsity of $L^{-1}$ for $\wt h\wt \beta_x$.

We focus on entries significantly deviating from $x_0$, the heavy entries of $\cW^\top (h+\wt h\wt \beta_x)$. Here, the treewidth-$\tau$ decomposition enables quick computation of an elimination tree based on $L^{-1}$'s sparsity, facilitating efficient estimation of $\|(\cW^\top (h+\wt h\wt \beta_x))_{\chi(v)}\|_2$ for any subtree $\chi(v)$\footnote{Given any tree node $v$, we use $\chi(v)$ to denote the subtree rooted at $v$.}. With an elimination tree of height $\wt O(\tau)$, we can employ heavy-light decomposition~\citep{sleator1981data} for an $O(\log n)$-height tree.

Using these data structures, convergence is established using the robust IPM framework~\citep{y20,lv21}. While the framework is generally applicable to QPs, computing an initial point remains a challenge. We propose a simpler objective $x_0=\arg\min_{x\in \R^n}\sum_{i=1}^n \phi_i(x_i)$ with $\phi_i$ as the log-barrier function, resembling the initial point reduction in~\cite{lsz19}. This initial point enables us to solve an augmented quadratic program that increases dimension by 1.

\subsection{Low-Rank Setting: How to Utilize Small Factorization} \label{sec:tech:rank}
The low-treewidth structure can be considered a form of sparsity, allowing for a sparse factorization ${\color{black}M}=LL^\top$. Another significant structure arises when the matrix $Q$ admits a low-rank factorization. Let $Q = UV^\top$ where $U, V \in \R^{n \times k}$ and $k \ll n$, then ${\color{black}M} = Q + tH = UV^\top + tH$. Although $Q$ has a low-rank structure, ${\color{black}M}$ may not be low-rank due to the diagonal matrix being dense. However, in the central path equations, we need only handle ${\color{black}M^{-1}}$, which can be efficiently maintained using the matrix Woodbury identity:
\begin{align*}
    {\color{black}M^{-1}} = & ~ t^{-1}H^{-1} - t^{-2} H^{-1} U (I + t^{-1} V^\top H^{-1} U)^{-1} V^\top H^{-1},
\end{align*}
Given that $H$ is diagonal, the complex term $(I + t^{-1} V^\top H^{-1} U)^{-1}$ can be quickly updated under sparse changes to $H^{-1}$ by simply scaling rows of $U$ and $V$. With only a nearly-linear number of updates to $H^{-1}$, the $k\times k$ matrix $V^\top H^{-1} U$ is updated in $\wt O(k^2)$ time per changed entry of $H$ and recomputed in $\wt O(nk^{\omega-1})$ time at each restart of the data structure. The inverse $(I + \ov t^{-1} V^\top H^{-1} U)^{-1}$, with $t$ replaced by an approximation $\ov t$, is recomputed in $O(k^\omega)$ time per iteration, and the inverse of the $m\times m$ matrix $A{\color{black}M^{-1}}A^\top$ is recomputed in $O(m^\omega)$ time per iteration (Section~\ref{sec:rank}). Over the $\wt O(\sqrt n\log(n/\epsilon))$ iterations, these inversions cost $\wt O(\sqrt n(k^\omega+m^\omega)\log(n/\epsilon))$ time in total. We assume that $k$ and $m$ are much smaller than $n$, precisely that $(k+m)^{\omega-1}\le n$. Then $(k+m)^{\omega}\le \sqrt n(k+m)^{(\omega+1)/2}$, so this cost is at most $\wt O(n(k+m)^{(\omega+1)/2}\log(n/\epsilon))$. All of these costs are counted in the running time bound of Lemma~\ref{lem:rank:cpm-time}, which yields the overall $\wt O(n(k+m)^{(\omega+1)/2}\log(n/\epsilon))$ bound of Theorem~\ref{thm:main-rank-informal}. We modify the $(x,s)$ implicit representation as follows:
\begin{align}
x &= \wh x + H^{-1/2} h \beta_x + H^{-1/2} \wh h \wh \beta_x + H^{-1/2} \wt h \wt \beta_x, \label{eqn:rank-implicit-x-rep} \\
s &= \wh s + H^{1/2} h \beta_s + H^{1/2} \wh h \wh \beta_s + H^{1/2} \wt h \wt \beta_s, \label{eqn:rank-implicit-s-rep}
\end{align}
where $\wh x, \wh s \in \R^n$, $h = H^{-1/2}\ov\delta_\mu \in \R^n$, $\wh h = H^{-1/2}U \in \R^{n \times k}$, and $\wt h = H^{-1/2}A^\top \in \R^{n \times m}$, with $\wt \beta_x, \wt \beta_s \in \R^m$ and $\beta_x,\beta_s\in \R$. The nontrivial terms to maintain are $\wh h$ and $\wt h$, but both can be managed straightforwardly: updates to $H^{-1/2}$ correspond to scaling rows of $U$ and $A^\top$, and can be performed in total $\wt O(nk)$ and $\wt O(nm)$ time, respectively. The key observation is that we \emph{never explicitly form} ${\color{black}M^{-1/2}}$, hence matrix Woodbury identity suffices for fast updates.

The remaining task is to design a data structure for detecting heavy entries. Instead of starting with an elimination tree and re-balancing it through heavy-light decomposition, we construct a balanced tree on $n$ nodes, hierarchically dividing length-$n$ vectors by their indices. Sampling is then performed by traversing down to the tree's leaves. While a heavy-hitter data structure could lead to improvements in poly-logarithmic and sub-logarithmic factors, we primarily focus on polynomial dependencies on various parameters and leave this enhancement for future exploration.

\subsection{Gaussian Kernel SVM: Algorithm and Hardness}
\label{sec:tech:kernel}
Our specialized QP solvers provide fast implementations for linear SVMs when the data dimension $d$ is much smaller than $n$. However, for kernel SVM, forming the kernel matrix exactly would take $\Theta(n^2)$ time. Fortunately, advancements in kernel matrix algebra~\citep{acss20,bimw21,aa22,biksz23} have enabled sub-quadratic algorithms when the data dimension $d$ is small or the kernel matrix has a relatively large minimum entry. Both~\cite{acss20} and~\cite{biksz23} introduce algorithms for spectral sparsification, generating an approximate matrix $\wt K \in \R^{n \times n}$ such that $(1-\epsilon)\cdot K \preceq \wt K \preceq (1+\epsilon)\cdot K$, with $\wt K$ having only $O(\epsilon^{-2}n\log n)$ nonzero entries. \cite{acss20} achieves this in $O(n^{1+o(1)})$ time for multiplicatively Lipschitz kernels when $d = {\color{black}O(\log n)}$, while \cite{biksz23} overcomes limitations for Gaussian kernels by basing their algorithm on KDE and the magnitude of the minimum entry of the kernel matrix, parameterized by $\tau$. Their algorithm for Gaussian kernels runs in time $\wt O(nd/\tau^{3.173+o(1)})$. Unfortunately, spectral sparsifiers do not aid our primitives since a sparsifier only reduces the number of nonzero entries, but not the rank of the kernel matrix.

Besides spectral sparsification, \cite{acss20,aa22} also demonstrate that with $d = {\color{black}O(\log n)}$ and suitable kernels, there exists an $O(n^{1+o(1)})$ time algorithm to multiply the kernel matrix with an arbitrary vector $v \in \R^n$ to accuracy $1/\poly(n)$. This operation is crucial in Batch KDE as shown in \cite{aa22}. Moreover, \cite{aa22} establishes an almost-quadratic lower bound for this operation when the squared dataset radius $B = \omega(\log n)$, assuming \seth. These results rely on computing a rank-$n^{o(1)}$ factorization for the Gaussian kernel matrix. The function $e^{-x}$ can be approximated by a low-degree polynomial of degree
\begin{align*}
    q := & ~ \Theta(\max\{\sqrt{B\log(1/\epsilon)},\frac{\log(1/\epsilon)}{\log(\log(1/\epsilon)/B)}\})
\end{align*}
for $x \in [0, B]$. Using this polynomial, one can create matrices $U, V$ with rank $\binom{2d+2q}{2q}$ in time $O(nd\binom{2d+2q}{2q})$. For $d=O(\log n)$, $B=o(\log n)$ and $\epsilon\ge 1/\poly(n)$ the degree is $q=o(\log n)$, so the rank is $n^{o(1)}$ and the time is $n^{1+o(1)}$. For superpolynomially small $\epsilon$ the rank can be superpolynomial (Remark~\ref{rem:gaussian_kernel_svm_parameters}). Given this factorization, multiplying it with a vector $v$ as $U(V^\top v)$ takes $O(n^{1+o(1)})$ time. Let $\wt K = UV^\top$ where $\wt K_{i,j} = f(\|x_i-x_j\|_2^2)$, we have for any $(i,j) \in [n] \times [n]$,
\begin{align*}
    |f(\|x_i-x_j\|_2^2) - \exp(-\|x_i-x_j\|_2^2)| \leq & ~ \epsilon,
\end{align*}
and for any row $i \in [n]$,
\begin{align*}
    |(\wt K v)_i - (Kv)_i| = & ~ |\sum_{j=1}^n v_j(f(\|x_i-x_j\|_2^2) - \exp(-\|x_i-x_j\|_2^2))| \\
    \leq & ~ (\max_{j\in [n]} |f(\|x_i-x_j\|_2^2) - \exp(-\|x_i-x_j\|_2^2)|) \|v\|_1 \\
    \leq & ~ \epsilon \|v\|_1,
\end{align*}
using H{\"o}lder's inequality. This provides an $\ell_\infty$-guarantee of the error vector $(\wt K - K)v$, useful for Batch Gaussian KDE. Transforming this $\ell_\infty$-guarantee into a bound on quadratic forms (Lemma~\ref{lem:linf_to_spectral}) yields, for every $v\in\R^n$,
\begin{align*}
    v^\top(K-\epsilon n\cdot I)v \le v^\top\wt Kv \le v^\top(K+\epsilon n\cdot I)v.
\end{align*}
For dual variables in a box $[0,R]^n$, whose squared Euclidean norms are at most $nR^2$, this changes the SVM objective by at most $\epsilon n^2R^2/2$, so a kernel accuracy $\epsilon=\Theta(\epsilon'/(n^2R^2))$ suffices for a target accuracy $\epsilon'$ (Theorem~\ref{thm:gaussian_kernel}). 

Given $\wt K = UV^\top$ for $U, V \in \R^{n \times n^{o(1)}}$, the matrix $(yy^\top)\circ\wt K=(D_yU)(D_yV)^\top$, where $D_y$ is the diagonal matrix of the labels, need not be symmetric or positive semidefinite, while our low-rank QP algorithm needs a symmetric positive semidefinite matrix. We therefore take its symmetric part and add $\epsilon n$ times the orthogonal projection onto the column space of $[D_yU,\ D_yV]$. The result is positive semidefinite, has rank $n^{o(1)}$, and is within $2\epsilon n$ of $(yy^\top)\circ K$ in spectral norm (Step~1 of the proof of Theorem~\ref{thm:gaussian_kernel}). With this matrix in place of $(yy^\top)\circ K$, we can solve program~\eqref{eqn:kernel_svm} using our low-rank QP algorithm in expected time $n^{1+o(1)}$ for any accuracy $\epsilon'\ge 1/\poly(n)$.\footnote{The algorithm also requires $B = o(\frac{\log n}{\log\log n})$ and that both labels occur. In the program it solves, the constraint $\alpha\ge 0$ is replaced by $\alpha\in[0,R]^n$ with $R\le\poly(n)$. The soft-margin program has this box, with $R=C$. For the hard-margin program, the box needs a known bound $R$ on the entries of some optimal solution. The algorithm is randomized, and the equality constraint holds only up to $3\epsilon'$. See Theorem~\ref{thm:gaussian_kernel} and Corollary~\ref{cor:gaussian_kernel_almost_linear}.} This is the first almost-linear time algorithm for Gaussian kernel SVM, even in low-precision settings, as prior works either lack machinery to approximately form the kernel matrix efficiently, or do not possess faster convex optimization solvers for solving a structured quadratic program associated with a kernel SVM.

The requirements $d = {\color{black}O(\log n)}$ and $B = o(\frac{\log n}{\log\log n})$ may seem restrictive. Some bound on $B$ is needed for high-precision solutions, as the following lower bounds show. Assuming \seth, no sub-quadratic time algorithm solves Gaussian kernel SVM to error $\exp(-\Theta(\log^2 n))$ without bias on some inputs with $d = \Theta(\log n)$ and $B = \Theta(\log^2 n)$ (Theorem~\ref{thm:svm_hard_no_bias}), or with bias on some inputs with $d=\Theta(\log n)$ and $B = \Theta(\log^6 n)$ (Corollary~\ref{cor:svm_hard_bias}, whose proof is only sketched). These lower bounds do not show that the condition on $d$ is needed. Their error is also much smaller than the $1/\poly(n)$ accuracy of our algorithm. This is based on a reduction from the $(1+\delta)$-approximate Hamming nearest neighbor (bichromatic closest pair) problem of~\cite{r18} (Theorem~\ref{thm:r18}) to Gaussian kernel SVM, following the framework of~\cite{backurs2017fine}. So, at error $\exp(-\Theta(\log^2 n))$, an almost-linear time algorithm needs some upper bound on $B$. A gap remains between our assumption $B=o(\frac{\log n}{\log\log n})$ and the squared radius $\Theta(\log^2 n)$ of the hard instances.

{\color{black}We note that in other variants of definitions for Gaussian kernels, one requires an additional parameter called the \emph{kernel width}, and the kernel function is defined as $\exp(-\frac{\|x_i-x_j\|_2^2}{2\sigma^2})$. In commonly used heuristics~\citep{rrpsw15}, $\sigma=O(\sqrt d)$, hence we could without loss of generality assuming $\sigma=1$ by requiring the squared radius to be $B/d$.}

\section{Conclusion}
\label{sec:conclusion}
On the algorithmic front, we introduce the first nearly-linear time algorithms for low-rank convex quadratic programming, which are randomized and need the assumptions of Theorem~\ref{thm:main-rank-informal} (rank and number of constraints much smaller than $n$, and polynomial bounds), leading to nearly-linear time algorithms for linear SVMs whose data points have norms at most $\poly(n)$ and margin at least $1/\poly(n)$. For Gaussian kernel SVMs, we utilize a low-rank approximation from \cite{aa22} when ${\color{black}d = O(\log n)}$ and the squared dataset radius is small, enabling a randomized algorithm with expected time $n^{1+o(1)}$ for accuracy $\epsilon\ge 1/\poly(n)$ when the dual variables are restricted to a box $[0,R]^n$ with $R\le\poly(n)$ (soft margin with $R=C$, or a known bound on a hard-margin optimum). On the hardness aspect, we establish that when ${\color{black}d = \Omega(\log n)}$, if the squared dataset radius is sufficiently large ($\Omega(\log^2 n)$ without bias and $\Omega(\log^6 n)$ with bias, the latter with a proof that is only sketched), then assuming \seth, no sub-quadratic algorithm exists for solving the program to error $\exp(-\Theta(\log^2 n))$. As our work is theoretical in nature, we do not foresee any potential negative societal impact. Several open problems arise from our work:


\paragraph{Better dependence on $k$ for low-rank QPs.} Our low-rank QP solver exhibits a dependence of $k^{(\omega+1)/2}$ on the rank $k$. Given the precomputed factorization, can we improve the exponents on $k$? Ideally, an algorithm with nearly-linear dependence on $k$ would align more closely with input size.

\paragraph{Better dependence on $m$ for general QPs.} Focusing on SVMs with a few equality constraints, our QP solvers do not exhibit strong dependence on the number of equality constraints $m$. Without structural assumptions on the constraint matrix $A$, this is expected. However, many QPs, particularly in graph contexts, involve large $m$. Is there a pathway to an algorithm with better dependence on $m$? More broadly, can we achieve a result akin to that of~\cite{ls19}, where the number of iterations depends on the square root of the rank of $A$, with minimal per iteration cost?

\paragraph{High-precision Gaussian kernel SVMs.} We prove almost-linear time for our Gaussian kernel algorithm only for accuracy $\epsilon\ge 1/\poly(n)$. At $d=\Theta(\log n)$ and $\epsilon=\exp(-\Theta(\log^2 n))$, the rank of the kernel factorization is $n^{\Theta(\log\log n)}$ (Remark~\ref{rem:gaussian_kernel_svm_parameters}). Can Gaussian kernel SVMs with $d=O(\log n)$ and small $B$ be solved in time $n^{1+o(1)}\poly\log(1/\epsilon)$ for every $\epsilon$? Our approach would need kernel matrix-vector products with error as small as $\exp(-\Theta(\log^2 n))$ in almost linear time, which the methods surveyed in~\cite{aa22} do not provide when $d=\Omega(\log n)$.

\paragraph{Stronger lower bound in terms of $B$ for Gaussian kernel SVMs.} We establish hardness results for Gaussian kernel SVM when $B = \Omega(\log^2 n)$ without bias and $B = \Omega(\log^6 n)$ with bias (the latter only sketched). This contrasts with our algorithm, which requires $B=o(\frac{\log n}{\log\log n})$. For Batch Gaussian KDE, \cite{aa22} demonstrated that fast algorithms are feasible for $B = o(\log n)$, with no sub-quadratic time algorithms for $B = \omega(\log n)$ assuming \seth. Can a stronger lower bound be shown for SVM programs with a bias term, reflecting a more natural setting?

%% file: ack.tex
\section*{Acknowledgement}
Part of the work was done while Yuzhou Gu was supported by the National Science Foundation under Grant No.~DMS-1926686.
Lichen Zhang was supported in part by NSF CCF-1955217 and NSF DMS-2022448.
The corrections in this version were found and written with the help of Claude Opus 5.5. The authors have verified the corrections and take full responsibility for the content of this paper.

%% file: preli.tex
\section{Preliminary} \label{sec:preli}

\subsection{Notations}

For a positive integer $n$, we use $[n]$ to denote the set $\{1,2,\cdots,n\}$. For a matrix $A$, we use $ A^\top$ to denote its transpose. 
For a matrix $A$, we define $\|A\|_{p \rightarrow q} := \sup_x \| A x \|_q / \| x \|_p$. When $p=q=2$, we recover the spectral norm.

We define the entrywise $\ell_p$-norm of a matrix $A$ as $\|A\|_p := (\sum_{i,j} |A_{i,j}|^p)^{1/p}$.

For any function $f:\mathbb{N}\rightarrow \mathbb{N}$ and $n\in \mathbb{N}$, we use $\wt O(f(n))$ to denote $O(f(n) \poly\log f(n))$. We use $\mathbbm{1}\{E\}$ to denote the indicator for event $E$, i.e., if $E$ happens, $\mathbbm{1}\{E\}=1$ and otherwise it's 0.

\subsection{Treewidth}\label{sec:preli:treewidth}

Treewidth captures the sparsity and tree-like structures of a graph. 

\begin{definition}[Tree Decomposition and Treewidth] \label{def:treewidth}
Let $G=(V,E)$ be a graph, a tree decomposition of $G$ is a tree $T$ with $b$ vertices, and $b$ sets $J_1,\ldots, J_b \subseteq V$ (called bags), satisfying the following properties:
\begin{itemize}
    \item For every edge $(u,v) \in E$, there exists $j\in [b]$ such that $u,v\in J_j$;
    \item For every vertex $v\in V$, $\{j\in [b]: v\in J_j\}$ is a non-empty subtree of $T$.
\end{itemize}

The treewidth of $G$ is defined as the minimum value of $\max\{|J_j| : j\in [b]\}-1$ over all tree decompositions.
\end{definition}
A near-optimal tree decomposition of a graph can be computed in almost linear time.
\begin{theorem}[\cite{bernstein2022deterministic}]
Given a graph $G$, there is an $O(m^{1+o(1)})$ time algorithm that produces a tree decomposition of $G$ of maximum bag size $O(\tau \log^3 n)$, where $\tau$ is the actual (unknown) treewidth of $G$.
\end{theorem}
Therefore, when $\tau = m^{\Theta(1)}$, we can compute an $\wt O(\tau)$-size tree decomposition in time $O(m \tau^{o(1)})$, which is negligible in the final running time of Theorem~\ref{thm:treewidth-formal}.

\subsection{Sparse Cholesky Decomposition}
In this section we state a few results on sparse Cholesky decomposition.
Fast sparse Cholesky decomposition algorithms are based on the concept of elimination tree, introduced in~\cite{schreiber1982new}.
\begin{definition}[Elimination tree] \label{defn:elim_tree}
Let $G$ be an undirected graph on $n$ vertices. An elimination tree $\cT$ is a rooted tree on $V(G)$ together with an ordering $\pi$ of $V(G)$ such that for any vertex $v$, its parent is the smallest (under $\pi$) element $u$ with $\pi(u)>\pi(v)$ such that there exists a path $P$ from $v$ to $u$, such that $\pi(w) \le \pi(v)$ for all $w\in P-u$.
\end{definition}
The following lemma relates the elimination tree and the structure of Cholesky factors.
\begin{lemma}[\cite{schreiber1982new}] \label{lem:elim_tree_imply_cholesky}
    Let $M$ be a PSD matrix and $\cT$ be an elimination tree of the adjacency graph of $M$ (i.e., $(i,j)\in E(G)$ iff $M_{i,j} \ne 0$) together with an elimination ordering $\pi$. Let $P$ be the permutation matrix $P_{i,v} = \mathbbm{1}\{v=\pi(i)\}$.
    Then the Cholesky factor $L$ of $PMP^\top$ (i.e., $PMP^\top = LL^\top$) satisfies $L_{i,j}\ne 0$ only if $\pi(i)$ is an ancestor of $\pi(j)$.
\end{lemma}
The following result is the current best result for computing a sparse Cholesky decomposition.
\begin{lemma}[{\cite[Lemma 8.4]{gs22}}] \label{lem:cholesky_fast}
Let $M\in \R^{n\times n}$ be a PSD matrix whose adjacency graph has treewidth $\tau$.
Then we can compute the Cholesky factorization $M = LL^\top$ in $\wt O(n \tau^{\omega-1})$ time.
\end{lemma}

The following result is the current best result for updating a sparse Cholesky decomposition.
\begin{lemma}[\cite{dh99}]\label{lem:cholesky_update}
Let $M\in \R^{n\times n}$ be a PSD matrix whose adjacency graph has treewidth $\tau$.
Assume that we are given the Cholseky factorization $M = LL^\top$.
Let $w\in \R^n$ and $\sigma\in\{\pm 1\}$ be such that $M+\sigma ww^\top$ is positive definite and has the same adjacency graph as $M$.
Then we can compute $\Delta_L\in \R^{n\times n}$ such that $L+\Delta_L$ is the Cholesky factor of $M+\sigma ww^\top$ in $\wt O(\tau^2)$ time.
\end{lemma}

Throughout our algorithm, we need to compute matrix-vector multiplications involving Cholesky factors. We use the following results from \cite{gs22}.
\begin{lemma}[{\cite[Lemma 4.7]{gs22}}] \label{lem:mat_vec_mult_time}
Let $M\in \R^{n\times n}$ be a PSD matrix whose adjacency graph has treewidth $\tau$.
Assume that we are given the Cholseky factorization $M = LL^\top$.
Then we have the following running time for matrix-vector multiplications.
\begin{enumerate}[label=(\roman*)]
    \item For $v\in \R^{n}$, computing $L v$, $L^\top v$, $L^{-1} v$, $L^{-\top} v$ takes $O(n\tau)$ time.
    \label{item:lem_mat_vec_mult_time_Linvv}
    \item For $v\in \R^{n}$, computing $L v$ takes $O(\|v\|_0 \tau)$ time.
    \label{item:lem_mat_vec_mult_time_Lv_sparse}
    \item For $v\in \R^{n}$, computing $L^{-1} v$ takes $O(\|L^{-1}v\|_0 \tau)$ time.
    \label{item:lem_mat_vec_mult_time_Linvv_sparse}
    \item For $v\in \R^{n}$, if $v$ is supported on a path in the elimination tree, then computing $L^{-1}v$ takes $O(\tau^2)$ time.
    \label{item:lem_mat_vec_mult_time_Linvv_path} 
    \item For $v\in \R^n$, computing $\cW^\top v$ takes $O(n \tau)$ time, where $\cW = L^{-1} H^{1/2}$ with $H\in \R^{n\times n}$ is a non-negative diagonal matrix.
    \label{item:lem_mat_vec_mult_time_WTv}
\end{enumerate}
\end{lemma}

\begin{lemma}[{\cite[Lemma 4.8]{gs22}}] \label{lem:mat_vec_mult_coord_time}
Let $M\in \R^{n\times n}$ be a PSD matrix whose adjacency graph has treewidth $\tau$.
Assume that we are given the Cholseky factorization $M = LL^\top$.
Then we have the following running time for matrix-vector multiplications, when we only need result for a subset of coordinates.
\begin{enumerate}[label=(\roman*)]
    \item Let $S$ be a path in the elimination tree whose one endpoint is the root. For $v\in \R^{n}$, computing $(L^{-\top}v)_S$ takes $O(\tau^2)$ time.
    \label{item:mat_vec_mult_coord_time_LinvtvS}
    \item For $v\in \R^{n}$, for $i\in [n]$, computing $(\cW^\top v)_i$ takes $O(\tau^2)$ time, where $\cW = L^{-1} H^{1/2}$ with $H\in \R^{n\times n}$ be a non-negative diagonal matrix.
    \label{item:mat_vec_mult_coord_time_Wtvi}
\end{enumerate}
\end{lemma}

\subsection{Johnson-Lindenstrauss Lemma}\label{sec:preli:sketching}
We recall the Johnson-Lindenstrauss lemma, a powerful algorithmic primitive that reduces dimension while preserving $\ell_2$ norms.
\begin{lemma}[\cite{jl84}]  
\label{lem:jl_matrix}
    Let $\epsilon \in (0, 1)$ be the precision parameter.  Let $\delta \in (0, 1)$ be the failure probability. Let $A \in \R^{m \times n}$ be a real matrix.
    Let $r = C\epsilon^{-2} \log(mn/\delta)$ for a sufficiently large absolute constant $C$.
    For $R \in \R^{r\times n}$ whose entries are i.i.d~$\cN(0, \frac 1r)$, the following holds with probability at least $1-\delta$:
    \begin{align*}
        (1 - \epsilon)\|a_i\|_2 \le \| R  a_i\|_2 \le (1 + \epsilon)\|a_i\|_2, ~\forall i \in [m],
    \end{align*}
    where for a matrix $A$, $a_i^\top$ denotes the $i$-th row of matrix $A \in \R^{m \times n}$.
\end{lemma}

\subsection{Heavy-Light Decomposition} \label{sec:preli:heavy-light}
Heavy-light decomposition is useful when one wants to re-balance a binary tree with height $O(\log n)$.
\begin{lemma}[\cite{sleator1981data}] \label{lem:heavy_light}
Given a rooted tree $\cT$ with $n$ vertices, we can construct in $O(n)$ time an ordering $\pi$ of the vertices such that (1) every path in $\cT$ can be decomposed into $O(\log n)$ contiguous subseqeuences under $\pi$, and (2) every subtree in $\cT$ is a single contiguous subsequence under $\pi$.
\end{lemma}

%% file: svm_form.tex
\section{SVM Formulations} \label{sec:svm-formulations}

In this section, we review a list of formulations of SVM. These formulations have been implemented in the \texttt{LIBSVM} library~\cite{cl11}.

Throughout this section, we use $\phi:\R^d\rightarrow \R^s$ to denote the feature mapping, ${\sf K}$ to denote the associated kernel function and $K\in \R^{n\times n}$ to denote the kernel matrix. For linear SVM, $\phi$ is just the identity mapping. We will focus on the dual quadratic program formulation as usual. We will also assume for each problem, a dataset $X\in \R^{n\times d}$ is given together with binary labels $y\in \R^n$. Let $Q:=(yy^\top) \circ K$.

\subsection{\texorpdfstring{$C$}{}-Support Vector Classification}

This formulation is also referred as the \emph{soft-margin SVM}. It can be viewed as imposing a regularization on the primal program to allow mis-classification.

\begin{definition}[\texorpdfstring{$C$}{}-Support Vector Classification]\label{def:C_SVC}

 Given a parameter $C>0$, the $C$-support vector classification ($C$-SVC) is defined as
\begin{align*}
    \max_{\alpha\in \R^n}~ & ~ {\bf 1}_n^\top \alpha-\frac{1}{2}\alpha^\top Q\alpha \\
    \textnormal{s.t.}~ & ~ \alpha^\top y = 0, \\
    & ~ 0 \leq \alpha \leq C\cdot {\bf 1}_n.
\end{align*}
\end{definition}

\subsection{\texorpdfstring{$\nu$}{}-Support Vector Classification}

The $C$-SVC (Definition~\ref{def:C_SVC}) penalizes large values of $\alpha$ by limiting the magnitude of it. The $\nu$-SVC (Definition~\ref{def:nu_SVC}) turns ${\bf 1}_n^\top \alpha$ from an objective into a constraint on $\ell_1$ norm.

\begin{definition}[$\nu$-Support Vector Classification]\label{def:nu_SVC}

Given a parameter $\nu>0$, the $\nu$-support vector classification ($\nu$-SVC) is defined as 
\begin{align*}
    \min_{\alpha\in \R^n}~ & ~ \frac{1}{2}\alpha^\top Q \alpha \\
    \textnormal{s.t.}~ & ~ \alpha^\top y = 0, \\
    & ~ {\bf 1}_n^\top \alpha = \nu,\\
    & ~ 0 \leq \alpha \leq \frac{1}{n}\cdot {\bf 1}_n.
\end{align*}
\end{definition}

One can interpret this formulation as to find a vector that lives in the orthogonal complement of $y$ that is non-negative, each entry is at most $\frac{1}{n}$ and its $\ell_1$ norm is $\nu$. Clearly, we must have $\nu\in (0, 1]$. More specifically, let $k_+$ be the number of positive labels and $k_-$ be the number of negative labels. It is shown by~\cite{cl01} that the above problem is feasible if and only if
\begin{align*}
    \nu \leq & ~ \frac{2\min\{k_-,k_+ \}}{n}.
\end{align*}

\subsection{Distribution Estimation}

SVM is widely-used for predicting binary labels. It can also be used to estimate the support of a high-dimensional distribution. The formulation is similar to $\nu$-SVC, except the PSD matrix $Q$ is \emph{label-less}.

\begin{definition}[Distribution Estimation]\label{def:nu_distribution_estimation}

Given a parameter $\nu>0$, the $\nu$-distribution estimation problem is defined as 
\begin{align*}
    \min_{\alpha\in \R^n}~ & ~ \frac{1}{2}\alpha^\top K \alpha \\
    \textnormal{s.t.}~ & ~ 0 \leq \alpha \leq \frac{1}{n}\cdot {\bf 1}_n, \\
    & ~ {\bf 1}_n^\top \alpha = \nu.
\end{align*}
\end{definition}

\subsection{\texorpdfstring{$\epsilon$}{}-Support Vector Regression}

In addition to classification, support vector framework can also be adapted for regression. 

\begin{definition}[$\epsilon$-Support Vector Regression]\label{def:epsilon_SVR}
Given parameters $\epsilon, C>0$, the $\epsilon$-support vector regression ($\epsilon$-SVR) is defined as
\begin{align*}
    \min_{\alpha, \alpha^*\in \R^n} ~ & ~ \frac{1}{2}(\alpha-\alpha^*)^\top K (\alpha-\alpha^*)+\epsilon  {\bf 1}_n^\top (\alpha+\alpha^*)+y^\top (\alpha-\alpha^*) \\
    \textnormal{s.t.}~&~ {\bf 1}_n^\top (\alpha-\alpha^*) = 0, \\
    & ~ 0 \leq \alpha \leq C\cdot {\bf 1}_n, \\
    & ~ 0 \leq \alpha^* \leq C\cdot {\bf 1}_n.
\end{align*}
\end{definition}

\subsection{\texorpdfstring{$\nu$}{}-Support Vector Regression}

One can similar adapt the parameter $\nu$ to control the $\ell_1$ norm of the regression.

\begin{definition}[$\nu$-Support Vector Regression]\label{def:nu_SVR}
    Given parameters $\nu, C>0$, the $\nu$-support vector regression ($\nu$-SVR) is defined as
\begin{align*}
    \min_{\alpha, \alpha^*\in \R^n} ~ & ~ \frac{1}{2}(\alpha-\alpha^*)^\top K (\alpha-\alpha^*)+y^\top (\alpha-\alpha^*) \\
    \textnormal{s.t.}~&~ {\bf 1}_n^\top (\alpha-\alpha^*) = 0, \\
    & ~ {\bf 1}_n^\top (\alpha+\alpha^*) = C \nu, \\
    & ~ 0 \leq \alpha \leq \frac{C}{n}\cdot {\bf 1}_n, \\
    & ~ 0 \leq \alpha^* \leq \frac{C}{n}\cdot {\bf 1}_n.
\end{align*}
\end{definition}

\subsection{One Equality Constraint}

We classify $C$-SVC (Definition~\ref{def:C_SVC}), $\epsilon$-SVR (Definition~\ref{def:epsilon_SVR}) and $\nu$-distribution estimation (Definition~\ref{def:nu_distribution_estimation}) into the following generic form:

\begin{align*}
    \min_{\alpha\in \R^n} ~ & ~ \frac{1}{2} \alpha^\top Q \alpha + p^\top \alpha \\
    \textnormal{s.t.} ~ & ~ \alpha^\top y = \Delta \\
    & ~ 0 \leq \alpha \leq C \cdot{\bf 1}_n.
\end{align*}

Note that $C$-SVC (Definition~\ref{def:C_SVC}) and distribution estimation (Definition~\ref{def:nu_distribution_estimation}) are readily in this form. For $\epsilon$-SVR (Definition~\ref{def:epsilon_SVR}), we need to perform a simple transformation:

Set $\wh \alpha=\begin{bmatrix}
    \alpha \\
    \alpha^*
\end{bmatrix} \in \R^{2n}$, then it can be written as 
\begin{align*}
    \min_{\wh \alpha\in \R^{2n}}~&~ \frac{1}{2} \wh \alpha^\top \begin{bmatrix}
        K & -K \\
        -K & K
    \end{bmatrix} \wh \alpha+\begin{bmatrix}
        \epsilon {\bf 1}_n+y \\ \epsilon {\bf 1}_n-y 
    \end{bmatrix}^\top \wh \alpha \\
    \textnormal{s.t.}~&~ \begin{bmatrix}
        {\bf 1}_n \\
        -{\bf 1}_n
    \end{bmatrix}^\top \wh \alpha = 0\\
    & ~ 0 \leq \wh \alpha \leq C\cdot {\bf 1}_{2n}.
\end{align*}

\subsection{Two Equality Constraints}

Both $\nu$-SVC (Definition~\ref{def:nu_SVC}) and $\nu$-SVR (Definition~\ref{def:nu_SVR}) require one extra constraint. They can be formulated as follows:
\begin{align*}
    \min_{\alpha\in \R^n} ~ & ~ \frac{1}{2} \alpha^\top Q \alpha+p^\top \alpha \\
    \textnormal{s.t.} ~ & ~ {\bf 1}_n^\top \alpha = \Delta_1, \\
    & ~ y^\top \alpha = \Delta_2, \\
    & ~ 0 \leq \alpha \leq C\cdot {\bf 1}_n.
\end{align*}
For $\nu$-SVR (Definition~\ref{def:nu_SVR}), one can leverage a similar transformation as $\epsilon$-SVR (Definition~\ref{def:epsilon_SVR}).

\begin{remark}
All variants of SVM-related quadratic programs can all be solved using our QP solvers for three cases:
\begin{itemize}
    \item Linear SVM with $n\gg d$, precisely $(d+2)^{\omega-1}\le n$, we can solve it in expected time $\wt O(nd^{(\omega+1)/2}\log(n/\epsilon))$;
    \item Linear SVM with a small treewidth decomposition with width $\tau$ on $XX^\top$, we can solve it in expected time $\wt O(n\tau^{(\omega+1)/2}\log(n/\epsilon))$;
    \item Gaussian kernel $C$-SVC (Definition~\ref{def:C_SVC}) with $d=O(\log n)$, $B=o(\frac{\log n}{\log\log n})$, $C\le\poly(n)$ and both labels present, to accuracy $\epsilon\ge 1/\poly(n)$, we can solve it in expected time $n^{1+o(1)}$, with the equality constraint holding up to $3\epsilon$ (Theorem~\ref{thm:gaussian_kernel} with $R=C$, and Corollary~\ref{cor:gaussian_kernel_almost_linear}).
\end{itemize}
The first two items need $C\le\poly(n)$, $1/r\le\poly(n)$ and data points of norm at most $\poly(n)$, where $r$ is the radius of a ball inside the constraint set relative to the equality constraints (Theorems~\ref{thm:rank-formal} and~\ref{thm:treewidth-formal}). For these box constraints, the outer radius in these theorems is at most $C\sqrt{n}$.
Even though our solvers have relatively bad dependence on the number of equality constraints, for all these SVM formulations, at most 2 equality constraints are presented and thus can be solved very fast.
\end{remark}

%% file: reduction.tex
\section{Algorithms for General QP} \label{sec:general_qp}
In this section, we discuss algorithms for general (convex) quadratic programming.
We show that these programs reduce to linear programming with convex constraints, so the algorithm of~\cite{lsz19} solves them approximately, under the boundedness assumptions and in the time of the propositions below.

\subsection{LCQP via Linear Programming with Convex Constraints}
We use the following result of~\cite{lsz19}. Consider a program $\min c'^\top z$ subject to $A'z=b'$ and $z\in\prod_{i}K_i$, with $N$ variables, where the $K_i$ are compact convex sets of dimension $O(1)$ and $A'$ has full row rank. Suppose that for each $K_i$ we are given a $\nu_i$-self-concordant barrier with $\nu_i=O(1)$ and its minimizer, that every $z\in\prod_iK_i$ satisfies $\|z\|_2\le R'$, and that $\|c'\|_2\le L'$. Then for every $\delta>0$ the algorithm of~\cite[Theorems 1.1 and C.3]{lsz19} returns $z\in\prod_iK_i$ with
\begin{align}\label{eq:lsz_guarantee}
c'^\top z\le \min_{A'z=b',\,z\in\prod_iK_i}c'^\top z+L'R'\delta \qquad\text{and}\qquad \|A'z-b'\|_1\le 3\delta(R'\|A'\|_1+\|b'\|_1),
\end{align}
in expected time $(N^{\omega+o(1)} + N^{2.5-\alpha/2+o(1)}+N^{2+1/6+o(1)})\cdot\wt O(\log(N/\delta))$, where $\omega$ is the matrix multiplication exponent and $\alpha$ is the dual matrix multiplication exponent. The sets must be bounded, and the equality constraints are satisfied only approximately. Our propositions inherit both features.

\begin{proposition}\label{prop:lcqp}
Consider LCQP (Definition~\ref{def:lcqp}). Suppose that $A$ has full row rank, that some optimal solution $x^*$ satisfies $\|x^*\|_\infty\le R$, and that $\|Q\|_{2\to 2}\le L$ and $\|c\|_2\le L$, where $R,L\ge 1$. Suppose also that we are given an eigendecomposition $Q=PDP^\top$, where $D$ is diagonal and $P$ is orthogonal. Let $\OPT$ be the optimal value of the program. Then for every $0<\epsilon\le\frac12$ there is a randomized algorithm that returns $x\in[0,R]^n$ with
\begin{align*}
\frac12x^\top Qx+c^\top x\le\OPT+\epsilon \qquad\text{and}\qquad \|Ax-b\|_1\le\epsilon(\|A\|_1+\|b\|_1),
\end{align*}
in expected time $(n^{\omega+o(1)} + n^{2.5-\alpha/2+o(1)}+n^{2+1/6+o(1)})\cdot\wt O(\log(nLR/\epsilon))$. In particular, the equality constraints hold only approximately.
\end{proposition}

\begin{proof}
\textbf{The reduction.}
Let $q_i := D_{i,i}\ge 0$ be the $i$-th element on the diagonal of $D$, let $T:=nR^2/2$, and let $a:=1/(\|A\|_1+\|b\|_1)$. (If there are no equality constraints $Ax=b$, we omit the first block of constraints below.) We introduce variables $\wt x\in\R^n$, which stand for $P^\top x$, and $t\in\R^n$, and we consider the program
\begin{align}\label{eq:lcqp_reduced}
\min ~ & ~ c^\top x + q^\top t \\
\mathrm{s.t.~} & ~ a A x = a b \notag\\
& ~ x - P \wt x = 0 \notag\\
&~ (\wt x_i,t_i)\in K:=\{(u,v)\in\R^2: \tfrac 12 u^2\le v\le T\} \qquad \forall i\in [n] \notag\\
& ~ x_i\in[0,R] \qquad \forall i\in [n]. \notag
\end{align}
Here the squared variables are the unconstrained coordinates $\wt x_i$. We do not split them, since splitting would make the quadratic form block diagonal rather than diagonal. The sets $[0,R]$ and $K$ are compact and convex. The function $-\log x_i-\log(R-x_i)$ is a $2$-self-concordant barrier for $[0,R]$ with minimizer $R/2$. The function $-\log(v-\frac12u^2)-\log(T-v)$ is a $2$-self-concordant barrier for $K$, and its gradient vanishes at $(u,v)=(0,T/2)$, which is its minimizer. Both are sums of two functions $-\log f$ with $f$ affine or concave quadratic, and each such function is a $1$-self-concordant barrier for the set where $f>0$. Indeed, for a direction $h$ let $\theta:=Df[h]/f$ and $\beta:=-D^2f[h,h]/f\ge 0$. The first three derivatives of $-\log f$ in the direction $h$ are $-\theta$, $\theta^2+\beta$ and $-2\theta^3-3\theta\beta$, and we have $\theta^2\le\theta^2+\beta$ and $|\theta|(2\theta^2+3\beta)\le 2(\theta^2+\beta)^{3/2}$.

\textbf{The reduced program has value at most $\OPT$.}
Let $x^*$ be the optimal solution from the hypothesis, let $\wt x^*:=P^\top x^*$ and let $t^*_i:=\frac12(\wt x^*_i)^2$. Then $\|\wt x^*\|_2=\|x^*\|_2\le\sqrt nR$, so $t_i^*\le T$ for all $i$, and $(x^*,\wt x^*,t^*)$ is feasible for~\eqref{eq:lcqp_reduced}. Its objective value is $c^\top x^*+\frac12\sum_{i\in[n]}q_i(\wt x_i^*)^2=c^\top x^*+\frac12x^{*\top}Qx^*=\OPT$.

\textbf{Parameters.}
Program~\eqref{eq:lcqp_reduced} has $N=3n$ variables. Every point of the product set satisfies $\|x\|_2\le\sqrt nR$, $|\wt x_i|\le\sqrt{2T}=\sqrt nR$ and $0\le t_i\le T$, so its Euclidean norm is at most $\sqrt nR+nR+\sqrt nT\le R':=3n^{3/2}R^2$. The cost vector has norm at most $\|c\|_2+\|q\|_2\le L':=2\sqrt nL$, because $0\le q_i\le\|Q\|_{2\to 2}\le L$. The constraint matrix is
\begin{align*}
A'=\begin{bmatrix} aA & 0 & 0\\ I & -P & 0\end{bmatrix}, \qquad\text{with right-hand side}\qquad b'=\begin{bmatrix} ab\\ 0\end{bmatrix}.
\end{align*}
It has full row rank, because $P$ is invertible and $A$ has full row rank. Each column of $P$ has Euclidean norm $1$ and hence $\ell_1$ norm at most $\sqrt n$, so $\|A'\|_1\le a\|A\|_1+n+\|P\|_1\le 1+n+n^{3/2}\le 3n^{3/2}$. Also $\|b'\|_1=a\|b\|_1\le 1$.

\textbf{Applying~\cite{lsz19}.}
We apply the result of~\cite{lsz19} stated above with $\delta:=\epsilon/(100n^4LR^3)$. It returns $(x,\wt x,t)$ in the product set. By~\eqref{eq:lsz_guarantee} and the bounds above,
\begin{align*}
c^\top x+q^\top t & \le\OPT+L'R'\delta\le\OPT+6n^2LR^2\delta, \\
a\|Ax-b\|_1+\|x-P\wt x\|_1 & \le\eta:=30n^3R^2\delta,
\end{align*}
where we used $3(R'\|A'\|_1+\|b'\|_1)\le 3(9n^3R^2+1)\le 30n^3R^2$. Note that $\eta=\frac{3\epsilon}{10nLR}\le 1$.

\textbf{Back to the original program.}
Let $e:=P^\top x-\wt x$, so that $\|e\|_2=\|x-P\wt x\|_2\le\eta$. Since $q\ge 0$ and $t_i\ge\frac12\wt x_i^2$, we have $\frac12\wt x^\top D\wt x\le q^\top t$. Also $\|\wt x\|_2\le nR$ and $\|D\|_{2\to2}\le L$. Hence
\begin{align*}
\frac12x^\top Qx & =\frac12(\wt x+e)^\top D(\wt x+e)
=\frac12\wt x^\top D\wt x+e^\top D\wt x+\frac12e^\top De \\
& \le q^\top t+L\eta\cdot nR+\frac L2\eta^2
\le q^\top t+2nLR\eta .
\end{align*}
Therefore
\begin{align*}
\frac12x^\top Qx+c^\top x & \le c^\top x+q^\top t+2nLR\eta\le\OPT+6n^2LR^2\delta+2nLR\eta \\
& \le\OPT+\frac{6\epsilon}{100}+\frac{6\epsilon}{10}\le\OPT+\epsilon .
\end{align*}
Moreover $\|Ax-b\|_1\le\eta/a\le\epsilon(\|A\|_1+\|b\|_1)$, and $x\in[0,R]^n$. Finally $N=3n$ and $\log(N/\delta)=O(\log(nLR/\epsilon))$, which gives the running time.
\end{proof}

\subsection{Algorithm for QCQP}
Our algorithm for LCQP in the previous section can be generalized to quadratically constrained quadratic programs (QCQP).
QCQP is defined as follows.
\begin{definition}[QCQP] \label{defn:qcqp-wiki}
Let $Q_0,\ldots, Q_m \in \R^{n \times n}$ be PSD matrices.
Let $q_0,\ldots, q_m\in \R^n$.
Let $r\in \R^m$.
Let $A\in \R^{d \times n}$, $b\in \R^{d}$.
The quadratically constrained quadratic programming (QCQP) problem asks the solve the following program.
\begin{align*}
 \min_{x \in \R^n} ~ & ~ \frac{1}{2} x^\top Q_0 x + q_0^\top x \\
\mathrm{s.t.~} & ~ \frac 12 x^\top Q_i x + q_i^\top x +r_i \le 0 \qquad \forall i\in [m]\\
&~ Ax = b \\
&~ x\ge 0
\end{align*}
\end{definition}

\begin{proposition}\label{prop:qcqp}
Consider QCQP (Definition~\ref{defn:qcqp-wiki}). Suppose that $A$ has full row rank, that some optimal solution $x^*$ satisfies $\|x^*\|_\infty\le R$, and that $\|Q_i\|_{2\to 2}\le L$ and $\|q_i\|_2\le L$ for all $0\le i\le m$ and $|r_i|\le L$ for all $i\in[m]$, where $R,L\ge 1$. Suppose also that we are given eigendecompositions $Q_i=P_iD_iP_i^\top$ for all $0\le i\le m$, where $D_i$ is diagonal and $P_i$ is orthogonal. Let $\OPT$ be the optimal value of the program and let $N:=(m+1)n$. Then for every $0<\epsilon\le\frac12$ there is a randomized algorithm that returns $x\in[0,R]^n$ with
\begin{align*}
\frac12x^\top Q_0x+q_0^\top x & \le\OPT+\epsilon, \\
\frac12x^\top Q_ix+q_i^\top x+r_i & \le\epsilon \qquad \forall i\in[m],\\
\|Ax-b\|_1 & \le\epsilon(\|A\|_1+\|b\|_1),
\end{align*}
in expected time $(N^{\omega+o(1)} + N^{2.5-\alpha/2+o(1)}+N^{2+1/6+o(1)})\cdot\wt O(\log(NLR/\epsilon))$, where $\omega$ is the matrix multiplication exponent and $\alpha$ is the dual matrix multiplication exponent. In particular, the equality constraints and the quadratic constraints hold only approximately.
\end{proposition}

\begin{proof}
We follow the proof of Proposition~\ref{prop:lcqp}.

\textbf{The reduction.}
Let $d_i\in\R^n$ be the diagonal of $D_i$, so that $0\le d_{i,j}\le L$. Let $T:=nR^2/2$, $S:=3nLR^2$ and $a:=1/(\|A\|_1+\|b\|_1)$. For $0\le i\le m$ we introduce variables $\wt x_i\in\R^n$, which stand for $P_i^\top x$, and $t_i\in\R^n$. For $i\in[m]$ we introduce a slack variable $\sigma_i\in\R$ for the $i$-th quadratic constraint. We consider the program
\begin{align}\label{eq:qcqp_reduced}
\min ~ & ~ q_0^\top x + d_0^\top t_0 \\
\mathrm{s.t.~} & ~ a A x = a b \notag\\
& ~ x - P_i \wt x_i = 0 \qquad \forall 0\le i\le m \notag\\
& ~ L^{-1}(d_i^\top t_i+q_i^\top x+r_i+\sigma_i)=0 \qquad \forall i\in[m] \notag\\
&~ (\wt x_{i,j},t_{i,j})\in K:=\{(u,v)\in\R^2: \tfrac 12 u^2\le v\le T\} \qquad \forall 0\le i\le m,\ j\in [n] \notag\\
& ~ x_j\in[0,R] \quad \forall j\in [n], \qquad \sigma_i\in[0,S] \quad \forall i\in[m]. \notag
\end{align}
We use the barriers from the proof of Proposition~\ref{prop:lcqp} for $K$ and $[0,R]$, and the barrier $-\log\sigma_i-\log(S-\sigma_i)$ for $[0,S]$, whose minimizer is $S/2$.

\textbf{The reduced program has value at most $\OPT$.}
Let $x^*$ be the optimal solution from the hypothesis, and let $\wt x_i^*:=P_i^\top x^*$, $t^*_{i,j}:=\frac12(\wt x^*_{i,j})^2$ and $\sigma_i^*:=-(\frac12x^{*\top}Q_ix^*+q_i^\top x^*+r_i)$. As before $t^*_{i,j}\le T$. Since $x^*$ is feasible and $\|x^*\|_2\le\sqrt nR$, we have $0\le\sigma_i^*\le\frac12LnR^2+L\sqrt nR+L\le S$. Since $d_i^\top t_i^*=\frac12x^{*\top}Q_ix^*$, this point is feasible for~\eqref{eq:qcqp_reduced}, with objective value $\OPT$.

\textbf{Parameters.}
Program~\eqref{eq:qcqp_reduced} has $n+2(m+1)n+m\le 4N$ variables. Every point of the product set satisfies $\|x\|_2\le\sqrt nR$, $\|\wt x_i\|_2\le nR$, $\|t_i\|_2\le\sqrt nT$ and $\|\sigma\|_2\le\sqrt mS$, so its Euclidean norm is at most
\begin{align*}
\sqrt nR+(m+1)(nR+\sqrt nT)+\sqrt mS\le R':=6(m+1)n^{3/2}LR^2 .
\end{align*}
The cost vector has norm at most $\|q_0\|_2+\|d_0\|_2\le L':=2\sqrt nL$. The constraint matrix $A'$ has full row rank: the rows of the third group are independent of the others because of the columns of $\sigma$, and the remaining rows are independent because the $P_i$ are invertible and $A$ has full row rank. As in the proof of Proposition~\ref{prop:lcqp}, and using $\|d_i\|_1\le nL$ and $\|q_i\|_1\le\sqrt nL$,
\begin{align*}
\|A'\|_1 & \le 1+(m+1)(n+n^{3/2})+m(n+\sqrt n+1)\le 6(m+1)n^{3/2}, \\
\|b'\|_1 & \le a\|b\|_1+L^{-1}\sum_{i\in[m]}|r_i|\le m+1 .
\end{align*}

\textbf{Applying~\cite{lsz19}.}
We apply the result of~\cite{lsz19} stated above with $\delta:=\epsilon/(400(m+1)^2n^4L^2R^3)$. It returns a point of the product set. By~\eqref{eq:lsz_guarantee}, this point satisfies $q_0^\top x+d_0^\top t_0\le\OPT+L'R'\delta$ and
\begin{align*}
a\|Ax-b\|_1+\sum_{0\le i\le m}\|x-P_i\wt x_i\|_1+L^{-1}\sum_{i\in[m]}|d_i^\top t_i+q_i^\top x+r_i+\sigma_i|\le\eta:=111(m+1)^2n^3LR^2\delta ,
\end{align*}
where we used $3(R'\|A'\|_1+\|b'\|_1)\le 3(36(m+1)^2n^3LR^2+m+1)\le 111(m+1)^2n^3LR^2$. Note that $\eta=\frac{111\epsilon}{400nLR}\le 1$.

\textbf{Back to the original program.}
For every $0\le i\le m$, the computation in the proof of Proposition~\ref{prop:lcqp}, with $e:=P_i^\top x-\wt x_i$, gives $\frac12x^\top Q_ix\le d_i^\top t_i+2nLR\eta$. For the objective we get
\begin{align*}
\frac12x^\top Q_0x+q_0^\top x & \le q_0^\top x+d_0^\top t_0+2nLR\eta\le\OPT+L'R'\delta+2nLR\eta \\
& \le\OPT+\frac{12\epsilon}{400}+\frac{222\epsilon}{400}\le\OPT+\epsilon,
\end{align*}
where we used $L'R'\delta=12(m+1)n^2L^2R^2\delta$. For $i\in[m]$ we use $\sigma_i\ge 0$ and $|d_i^\top t_i+q_i^\top x+r_i+\sigma_i|\le L\eta$ and get
\begin{align*}
\frac12x^\top Q_ix+q_i^\top x+r_i\le(d_i^\top t_i+q_i^\top x+r_i+\sigma_i)-\sigma_i+2nLR\eta\le L\eta+2nLR\eta\le 3nLR\eta\le\epsilon .
\end{align*}
Moreover $\|Ax-b\|_1\le\eta/a\le\epsilon(\|A\|_1+\|b\|_1)$, and $x\in[0,R]^n$. Finally the number of variables is at most $4N$ and $\log(4N/\delta)=O(\log(NLR/\epsilon))$, which gives the running time.
\end{proof}

%% file: ds-treewidth.tex
\section{Algorithm for Low-Treewidth QP} \label{sec:treewidth}
In this section we present a nearly-linear time algorithm for solving low-treewidth QP with small number of linear constraints. We briefly describe the outline of this section.
\begin{itemize}
    \item In Section~\ref{sec:treewidth:intro}, we present the main statement of Section~\ref{sec:treewidth}.
    \item In Section~\ref{sec:treewidth:cpm-structure}, we present the main data structure \textsc{CentralPathMaintenance}.
    \item In Section~\ref{sec:treewidth:cpm-ds}, we present several data structures used in \textsc{CentralPathMaintenance}, including \textsc{ExactDS} (Section~\ref{sec:treewidth:cpm-ds:exactds}), \textsc{ApproxDS} (Section~\ref{sec:treewidth:cpm-ds:approxds}), \textsc{BatchSketch} (Section~\ref{sec:treewidth:cpm-ds:batchsketch}), \textsc{VectorSketch} (Section~\ref{sec:treewidth:cpm-ds:vsketch}), \textsc{BalancedSketch} (Section~\ref{sec:treewidth:cpm-ds:bsketch}).
    \item In Section~\ref{sec:treewidth:cpm-analysis}, we prove correctness and running time of \textsc{CentralPathMaintenance} data structure.
    \item In Section~\ref{sec:treewidth:main}, we prove the main result (Theorem~\ref{thm:treewidth-formal}).
\end{itemize}

\subsection{Main Statement} \label{sec:treewidth:intro}
We consider programs of the form~\eqref{eqn:qp-general}, i.e.,
\begin{align*}
 \min_{x \in \R^n} ~ & ~ \frac{1}{2} x^\top Q x + c^\top x \\
\mathrm{s.t.~} & ~ A x = b  \\
& ~ x_i \in \cK_i \qquad \forall i\in [n]
\end{align*}
where $Q\in \cS^{n_\tot}$ is symmetric positive semidefinite, $c\in \R^{n_\tot}$, 
$A\in \R^{m\times n_\tot}$, $b\in \R^m$, $\cK_i \subset \R^{n_i}$ is a convex set.
For simplicity, we assume that $n_i = O(1)$ for all $i\in [n]$.

\begin{theorem} \label{thm:treewidth-formal}
Consider the convex program~\eqref{eqn:qp-general}.
Let $\phi_i: \cK_i \to \R$ be a $\nu_i$-self-concordant barrier for all $i\in [n]$.
Suppose the program satisfies the following properties:
\begin{itemize}
    \item Inner radius $r$: There exists $z\in \R^{n_\tot}$ such that $Az=b$ and $B(z,r) \subseteq \cK$.
    \item Outer radius $R$: $\cK \subseteq B(0,R)$ where $0\in \R^{n_\tot}$.
    \item Lipschitz constant $L$: $\|Q\|_{2\to 2}\le L$, $\|c\|_2 \le L$.
    \item Treewidth $\tau$: Treewidth (Definition~\ref{def:treewidth}) of the block adjacency graph of $Q$ (vertex set $[n]$, with $\{i,j\}$ an edge whenever the block $Q_{i,j}\in\R^{n_i\times n_j}$ is nonzero) is at most $\tau$.
\end{itemize}
Let $(w_i)_{i\in [n]} \in \R_{\ge 1}^n$ and $\kappa = \sum_{i\in [n]} w_i \nu_i$.
Given any $0<\epsilon\le \frac 12$, we can find an approximate solution $x\in \cK$ satisfiying
\begin{align*}
\frac 12 x^\top Q x + c^\top x &\le \min_{Ax=b, x\in \cK} \left(\frac 12 x^\top Q x + c^\top x\right)+ \epsilon L R(R+1),\\
\|Ax-b\|_1 &\le 3\epsilon(R \|A\|_1 + \|b\|_1),
\end{align*}
with probability at least $0.9$, in expected time
\begin{align*}
    \wt O\big( & (\sqrt \kappa n^{-1/2} + \max_{i\in [n]} (\nu_i/w_i+\sqrt{\nu_i})) \log(\frac{n\kappa R}{r\epsilon}) \cdot n (\tau^2 m+\tau m^2)^{1/2} (\tau^{\omega-1} + \tau m + m^{\omega-1})^{1/2} \\
    & + \sqrt \kappa \log(\frac{n\kappa R}{r\epsilon}) \cdot m^{\omega}\big).
\end{align*}
We assume that $m$ is much smaller than $n$, precisely that $m^{\omega-1}\le n\tau$. Under this assumption, which the bound above does not use, the last term is at most the first one, and the expected time is
\begin{align*}
    \wt O((\sqrt \kappa n^{-1/2} + \max_{i\in [n]} (\nu_i/w_i+\sqrt{\nu_i})) \log(\frac{n\kappa R}{r\epsilon}) \cdot n (\tau^2 m+\tau m^2)^{1/2} (\tau^{\omega-1} + \tau m + m^{\omega-1})^{1/2}).
\end{align*}
When $\max_{i\in [n]} \nu_i = \wt O(1)$, $w_i=1$, $m = \wt O(\tau^{\omega-2})$, the running time simplifies to
\begin{align*}
    \wt O(n \tau^{(\omega+1)/2} m^{1/2} \log(R/(r\epsilon))).
\end{align*}
\end{theorem}
The term $\sqrt \kappa \log(\frac{n\kappa R}{r\epsilon}) \cdot m^{\omega}$ is the cost of solving a linear system with an $m\times m$ matrix in every iteration (Lemma~\ref{lem:treewidth:cpm-time}). The assumption $m^{\omega-1}\le n\tau$ holds in our applications, where the number of equality constraints is much smaller than $n$: the SVM programs of Section~\ref{sec:svm-formulations} have $m\le 2$.

\subsection{Algorithm Structure and Central Path Maintenance} \label{sec:treewidth:cpm-structure}
Our algorithm is based on the robust Interior Point Method (robust IPM). Details of the robust IPM will be given in Section~\ref{sec:robust_ipm}.
During the algorithm, we maintain a primal-dual solution pair $(x,s) \in \R^{n_\tot} \times \R^{n_\tot}$ on the robust central path. In addition, we maintain a sparsely-changing approximation $(\ov x, \ov s) \in \R^{n_\tot} \times \R^{n_\tot}$ to $(x,s)$.
In each iteration, we implicitly perform update
\begin{align*}
x&\gets x + \ov t B_{w,\ov x,\ov t}^{-1/2} (I-P_{w,\ov x,\ov t}) B_{w,\ov x,\ov t}^{-1/2} \delta_\mu \\
s&\gets s + \ov t \delta_\mu - \ov t^2 H_{w,\ov x} B_{w,\ov x,\ov t}^{-1/2} (I-P_{w,\ov x,\ov t}) B_{w,\ov x,\ov t}^{-1/2} \delta_\mu
\end{align*}
where
\begin{align*}
    H_{w,\ov x} &= \nabla^2 \phi_w(\ov x) \tag{see Eq.~\eqref{eqn:def:Hwx}} \\
    B_{w,\ov x,\ov t} &= Q + \ov t H_{w,\ov x} \tag{see Eq.~\eqref{eqn:def:B}}\\
    P_{w,\ov x,\ov t} &= B_{w,\ov x,\ov t}^{-1/2} A^\top (A B_{w,\ov x,\ov t}^{-1} A^\top)^{-1} A B_{w,\ov x,\ov t}^{-1/2} \tag{see Eq.~\eqref{eqn:def:P}}
\end{align*}
and explicitly maintain $(\ov x, \ov s)$ such that they remain close to $(x,s)$ in $\ell_\infty$-distance.

This task is handled by the \textsc{CentralPathMaintenance} data structure, which is our main data structure. The robust IPM algorithm (Algorithm~\ref{alg:robust-ipm},~\ref{alg:centering}) directly calls it in every iteration.

The \textsc{CentralPathMaintenance} data structure (Algorithm~\ref{alg:treewidth:cpm}) has two main sub data structures, \textsc{ExactDS} (Algorithm~\ref{alg:treewidth:exactds-part1},~\ref{alg:treewidth:exactds-part2}) and \textsc{ApproxDS} (Algorithm~\ref{alg:treewidth:approxds-part1},~\ref{alg:treewidth:approxds-part2}).
\textsc{ExactDS} is used to maintain $(x,s)$, and \textsc{ApproxDS} is used to maintain $(\ov x, \ov s)$.

\begin{algorithm}[!ht]\caption{Our main data structure for low-treewidth QP solver.}\label{alg:treewidth:cpm}
\begin{algorithmic}[1]
\State {\bf data structure} \textsc{CentralPathMaintenance} \Comment{Theorem~\ref{thm:treewidth:cpm}}
\State {\bf private : member}
\State \hspace{4mm} \textsc{ExactDS} $\mathsf{exact}$ \Comment{Algorithm~\ref{alg:treewidth:exactds-part1}, \ref{alg:treewidth:exactds-part2}}
\State \hspace{4mm}  \textsc{ApproxDS} $\mathsf{approx}$ \Comment{ Algorithm~\ref{alg:treewidth:approxds-part1}}
\State \hspace{4mm} $\ell\in \bN$
\State {\bf end members}
\Procedure{\textsc{Initialize}}{$x, s\in \R^{n_\tot}, t\in \R_+, \ov \epsilon \in (0, 1)$}
    \State $\mathsf{exact}.\textsc{Initialize}(x, s, x, s, t)$ \Comment{Algorithm~\ref{alg:treewidth:exactds-part1}}
    \State $\ell \leftarrow 0$
    \State $w \gets \nu_{\max}$, $N\gets \sqrt \kappa \log n \log \frac{n \kappa R}{\ov \epsilon r}$
    \State $q \gets n^{1/2} (\tau^2 m+\tau m^2)^{-1/2} (\tau^{\omega-1} + \tau m + m^{\omega-1})^{1/2}$
    \State $\epsilon_{\apx,x} \leftarrow \ov \epsilon, \zeta_x \leftarrow 2 \alpha, \delta_\apx \leftarrow \frac 1{20N}$
    \State $\epsilon_{\apx,s} \leftarrow (1-\epsilon_t)\cdot\ov \epsilon \cdot \ov t , \zeta_s \leftarrow 3\alpha \ov t$
    \State \begin{align*}
        & ~ \mathsf{approx}.\textsc{Initialize}(x,s,h, \wt h, \epsilon_x, \epsilon_s, H_{w,\ov x}^{1/2} \wh x, H_{w,\ov x}^{-1/2} \wh s, c_s, \beta_x, \beta_s, \beta_{c_s},\\
        & ~ \wt \beta_x, \wt \beta_s, q, \& \mathsf{exact}, \epsilon_{\apx,x}, \epsilon_{\apx,s}, \delta_{\apx})
    \end{align*}
    \State
    \Comment{Algorithm~\ref{alg:treewidth:approxds-part1}.Parameters from $x$ to $\wt \beta_s$ come from $\mathsf{exact}$. $\& \mathsf{exact}$ is pointer to $\mathsf{exact}$}
\EndProcedure
\Procedure{\textsc{MultiplyAndMove}}{$t\in \R_+$}
    \State $\ell\gets \ell + 1$
    \If{$|\ov t - t| > \ov t \cdot \epsilon_t$ or $\ell \ge q$}
        \State $x, s \gets \mathsf{exact}.\textsc{Output}()$ \Comment{Algorithm~\ref{alg:treewidth:exactds-part1}}
        \State \textsc{Initialize}$(x,s,t,\ov \epsilon)$
    \EndIf
    \State $\beta_x, \beta_s, \beta_{c_s}, \wt \beta_x, \wt \beta_s \gets \mathsf{exact}.\textsc{Move}()$ \Comment{Algorithm~\ref{alg:treewidth:exactds-part1}}
    \State $\delta_{\ov x}, \delta_{\ov s } \gets \mathsf{approx}.\textsc{MoveAndQuery}(\beta_x, \beta_s, \beta_{c_s}, \wt \beta_x, \wt \beta_s)$ \Comment{Algorithm~\ref{alg:treewidth:approxds-part1}}
    \State $\delta_{h}, \delta_{\wt h}, \delta_{\epsilon_x}, \delta_{\epsilon_s}, \delta_{H_{w,\ov x}^{1/2} \wh x}, \delta_{H_{w,\ov x}^{-1/2} \wh s}, \delta_{c_s}\gets \mathsf{exact}.\textsc{Update}(\delta_{\ov x}, \delta_{\ov s})$ \Comment{Algorithm~\ref{alg:treewidth:exactds-part2}}
    \State $\mathsf{approx}.\textsc{Update}(\delta_{\ov x}, \delta_{h}, \delta_{\wt h}, \delta_{\epsilon_x}, \delta_{\epsilon_s}, \delta_{H_{w,\ov x}^{1/2} \wh x}, \delta_{H_{w,\ov x}^{-1/2} \wh s}, \delta_{c_s})$ \Comment{Algorithm~\ref{alg:treewidth:approxds-part1}}
\EndProcedure
\Procedure{\textsc{Output}}{$ $}
    \State \Return $\mathsf{exact}.\textsc{Output}()$ \Comment{Algorithm~\ref{alg:treewidth:exactds-part1}}
\EndProcedure
\State {\bf end data structure}
\end{algorithmic}
\end{algorithm}
\begin{theorem} \label{thm:treewidth:cpm}
Data structure \textsc{CentralPathMaintenance} (Algorithm~\ref{alg:treewidth:cpm}) implicitly maintains the central path primal-dual solution pair $(x,s) \in \R^{n_\tot} \times \R^{n_\tot}$ and explicitly maintains its approximation $(\ov x, \ov s) \in \R^{n_\tot} \times \R^{n_\tot}$ using the following functions:
\begin{itemize}
    \item \textsc{Initialize}$(x\in \R^{n_\tot}, s\in \R^{n_\tot}, t_0\in \R_{>0}, \epsilon\in (0,1))$: Initializes the data structure with initial primal-dual solution pair $(x,s)\in \R^{n_\tot} \times \R^{n_\tot}$, initial central path timestamp $t_0\in \R_{>0}$ in $\wt O(n (\tau^{\omega-1} + \tau m + m^{\omega-1}))$ time.
    \item \textsc{MultiplyAndMove}$(t\in \R_{>0})$: It implicitly maintains
    \begin{align*}
    x&\gets x + \ov t B_{w,\ov x,\ov t}^{-1/2} (I-P_{w,\ov x,\ov t}) B_{w,\ov x,\ov t}^{-1/2} \delta_\mu(\ov x, \ov s, \ov t) \\
    s&\gets s + \ov t \delta_\mu - \ov t^2 H_{w,\ov x} B_{w,\ov x,\ov t}^{-1/2} (I-P_{w,\ov x,\ov t}) B_{w,\ov x,\ov t}^{-1/2} \delta_\mu(\ov x, \ov s, \ov t)
    \end{align*}
    where $H_{w,\ov x}$, $B_{w,\ov x, \ov t}$, $P_{w,\ov x, \ov t}$ are defined in Eq.~\eqref{eqn:def:Hwx}\eqref{eqn:def:B}\eqref{eqn:def:P} respectively, and $\ov t$ is some timestamp satisfying $|\ov t-t| \le \epsilon_t \cdot \ov t$.
    
    It also explicitly maintains $(\ov x, \ov s) \in \R^{n_{\tot} \times n_{\tot}}$ such that $\|\ov x_i-x_i\|_{\ov x_i} \le \ov\epsilon$ and $\|\ov s_i-s_i\|^*_{\ov x_i} \le t \ov\epsilon w_i$ for all $i\in [n]$ with probability at least $0.9$.
    
    Assuming the function is called at most $N$ times and $t$ decreases from $t_{\max}$ to $t_{\min}$, on the same event of probability at least $0.9$ the total running time is 
    \begin{align*}
        \wt O((N n^{-1/2} + \epsilon_t^{-1}\log(t_{\max}/t_{\min})) \cdot n (\tau^2 m+\tau m^2)^{1/2} (\tau^{\omega-1} + \tau m + m^{\omega-1})^{1/2} + N m^{\omega}).
    \end{align*}
    
    \item \textsc{Output}: Computes $(x,s) \in \R^{n_{\tot}} \times \R^{n_{\tot}}$ exactly and outputs them in 
    $\wt O(n \tau m)$
    time.

\end{itemize}
\end{theorem}

\subsection{Data Structures Used in \textsc{CentralPathMaintenance}} \label{sec:treewidth:cpm-ds}
In this section we present several data structures used in \textsc{CentralPathMaintenance}, including:
\begin{itemize}
    \item \textsc{ExactDS} (Section~\ref{sec:treewidth:cpm-ds:exactds}): This data structure maintains an implicit representation of the primal-dual solution pair $(x,s)$. This is directly used by \textsc{CentralPathMaintenance}.
    \item \textsc{ApproxDS} (Section~\ref{sec:treewidth:cpm-ds:approxds}): This data structure explicitly maintains an approximation $(\ov x,\ov s)$ of $(x,s)$. This data structure is directly used by \textsc{CentralPathMaintenance}.
    \item \textsc{BatchSketch} (Section~\ref{sec:treewidth:cpm-ds:batchsketch}): This data structure maintains a sketch of $(x,s)$. This data structure is used by \textsc{ApproxDS}.
    \item \textsc{VectorSketch} (Section~\ref{sec:treewidth:cpm-ds:vsketch}): This data structure maintains a sketch of sparsely-changing vectors. This data structure is used by \textsc{BatchSketch}.
    \item \textsc{BalancedSketch} (Section~\ref{sec:treewidth:cpm-ds:bsketch}): This data structure maintains a sketch of vectors of form $\cW^\top v$, where $v$ is sparsely-changing. This data structure is used by \textsc{BatchSketch}.
\end{itemize}

Notation: In this section, for simplicity, we write $B_{\ov x}$ for $B_{w,\ov x,\ov t}$, and $L_{\ov x}$ for the Cholesky factor of $B_{\ov x}$, i.e., $B_{\ov x} = L_{\ov x}L_{\ov x}^\top$.

\subsubsection{\textsc{ExactDS}} \label{sec:treewidth:cpm-ds:exactds}
In this section we present the data structure \textsc{ExactDS}. It maintains an implicit representation of the primal-dual solution pair $(x,s)$ by maintaining several sparsely-changing vectors (see Eq.~\eqref{eqn:thm:treewidth:exactds:x-rep}\eqref{eqn:thm:treewidth:exactds:s-rep}).
This data structure has a similar spirit as \textsc{ExactDS} in \cite{gs22}, but we have a different representation from the previous works because we are working with quadratic programming rather than linear programming.

\begin{theorem} \label{thm:treewidth:exactds}
Data structure \textsc{ExactDS} (Algorithm~\ref{alg:treewidth:exactds-part1}, \ref{alg:treewidth:exactds-part2}) implicitly maintains the primal-dual pair $(x,s) \in \R^{n_{\tot}} \times \R^{ n_{\tot} }$, computable via the expression
\begin{align}
x &= \wh x + H_{w,\ov x}^{-1/2} \cW^\top ( h \beta_x - \wt h \wt \beta_x +\epsilon_x), \label{eqn:thm:treewidth:exactds:x-rep} \\
s &= \wh s + H_{w,\ov x}^{1/2} c_s \beta_{c_s}- H_{w,\ov x}^{1/2} \cW^\top ( h \beta_s - \wt h \wt \beta_s +\epsilon_s), \label{eqn:thm:treewidth:exactds:s-rep}
\end{align}
where
$\wh x,\wh s \in \R^{n_\tot}$, $\cW = L_{\ov x}^{-1} H_{w,\ov x}^{1/2} \in \R^{n_\tot \times n_\tot}$, $h = L_{\ov x}^{-1} \ov\delta_\mu \in \R^{n_\tot}$, $c_s = H_{w,\ov x}^{-1/2} \ov\delta_\mu \in \R^{n_\tot}$ $\beta_x,\beta_s, \beta_{c_s} \in \R$, $\wt h = L_{\ov x}^{-1} A^\top \in \R^{n_\tot \times m}$, $\wt \beta_x, \wt \beta_s \in \R^{m}$, $\epsilon_x, \epsilon_s \in \R^{n_\tot}$.

The data structure supports the following functions:
\begin{itemize}
    \item \textsc{Initialize}$(x, s, \ov x, \ov s\in \R^{n_{\tot}}, \ov t \in \R_{>0})$: Initializes the data structure in 
    $\wt O(n \tau^{\omega-1} + n \tau m + n m^{\omega-1})$
    time, with initial value of the primal-dual pair $(x, s)$, its initial approximation $(\ov x, \ov s)$, and initial approximate timestamp $\ov t$.
    \item \textsc{Move}$()$: Performs robust central path step
    \begin{align} 
    x & \gets x + \ov t B_{\ov x}^{-1} \delta_\mu - \ov t B_{\ov x}^{-1} A^\top (A B_{\ov x}^{-1} A^\top)^{-1} A B_{\ov x}^{-1} \delta_\mu, \label{eqn:thm:treewidth:exactds:x-step} \\
    s & \gets s + \ov t \delta_\mu - \ov t^2 H_{w,\ov x} B_{\ov x}^{-1} \delta_\mu + \ov t^2 H_{w,\ov x} B_{\ov x}^{-1} A^\top (A B_{\ov x}^{-1} A^\top)^{-1} A B_{\ov x}^{-1} \delta_\mu \label{eqn:thm:treewidth:exactds:s-step}
    \end{align}
    in $O(m^\omega)$ time by updating its implicit representation.
    \item \textsc{Update}$(\delta_{\ov x}, \delta_{\ov s} \in \R^{n_{\tot}})$: Updates the approximation pair $(\ov x, \ov s)$ to $(\ov x^{\new} = \ov x + \delta_{\ov x} \in \R^{ n_{\tot} }, \ov s^{\new} = \ov s + \delta_{ \ov s} \in \R^{n_{\tot}} )$ in
    $
    \wt O((\tau^2 m+\tau m^2) (\|\delta_{\ov x}\|_0 + \|\delta_{\ov s}\|_0))
    $
    time, and output the changes in variables $\delta_{H_{w,\ov x}^{1/2}\wh x}$, $\delta_h$, $\delta_{\beta_x}$, $\delta_{\wt h}$, $\delta_{\wt \beta_x}$, $\delta_{\epsilon_x}$, $\delta_{H_{w,\ov x}^{-1/2}\wh s}$, $\delta_{\beta_s}$, $\delta_{\wt \beta_s}$, $\delta_{\epsilon_s}$.
    
    Furthermore, $h, \epsilon_x, \epsilon_s$ change in $O(\tau (\|\delta_{\ov x}\|_0 + \|\delta_{\ov s}\|_0))$ coordinates,
    $\wt h$ changes in $\wt O(\tau m (\|\delta_{\ov x}\|_0 + \|\delta_{\ov s}\|_0))$ coordinates,
    and $H_{\ov x}^{1/2} \wh x, H_{\ov x}^{-1/2} \wh s, c_s$ change in $O(\|\delta_{\ov x}\|_0 + \|\delta_{\ov s}\|_0)$  coordinates.
    
    \item \textsc{Output}$()$: Output $x$ and $s$ in $\wt O(n \tau m)$ time.

    \item \textsc{Query}$x(i\in [n])$: Output $x_i$ in $\wt O(\tau^2 m)$ time.
    This function is used by \textsc{ApproxDS}.
    \item \textsc{Query}$s(i\in [n])$: Output $s_i$ in $\wt O(\tau^2 m)$ time.
    This function is used by \textsc{ApproxDS}.
\end{itemize}
\end{theorem}

\begin{algorithm}[!ht]\caption{The \textsc{ExactDS} data structure used in Algorithm~\ref{alg:treewidth:cpm}. }\label{alg:treewidth:exactds-part1}
\begin{algorithmic}[1]
\State {\bf data structure} \textsc{ExactDS} \Comment{Theorem~\ref{thm:treewidth:exactds} 
}
\State {\bf members}
    \State \hspace{4mm} $\ov x , \ov s \in \R^{n_{\tot}}$, $\ov t\in \R_{+}$, $H_{w,\ov x}, B_{\ov x}, L_{\ov x} \in \R^{n_{\tot} \times n_{\tot}}$
    \State \hspace{4mm} $\wh x, \wh s, h, \epsilon_x, \epsilon_s, c_s\in \R^{n_{\tot}}$, $\wt h\in \R^{n_\tot \times m}$, $\beta_x, \beta_s, \beta_{c_s}\in \R$, $\wt \beta_x, \wt \beta_s \in \R^m$
    \State \hspace{4mm} $\wt u\in \R^{m\times m}$, $u \in \R^m$, $\ov\alpha \in \R, \ov\delta_\mu \in \R^{n}$
    \State \hspace{4mm} $k\in \bN$
\State {\bf end members}
\Procedure{Initialize}{$x,s, \ov x, \ov s \in \R^{n_\tot}, \ov t \in \R_+$}
    \State $\ov x \gets \ov x$, $\ov s \gets \ov s$, $\ov t \gets \ov t$, $\ov\alpha \gets 0$
    \State $\wh{x} \gets x$, $\wh{s} \gets s$, $\epsilon_x \gets 0$, $\epsilon_s \gets 0$, $\beta_x \gets 0$, $\beta_s \gets 0$, $\wt \beta_x \gets 0$, $\wt \beta_s \gets 0$, $\beta_{c_s} \gets 0$
    \State $H_{w,\ov x} \gets \nabla^2 \phi_w( \ov x )$, $B_{\ov x} \gets Q + \ov t H_{w,\ov x}$
    \State Compute lower Cholesky factor $L_{\ov x}$ where $L_{\ov x} L_{\ov x}^\top = B_{\ov x}$
    \State \textsc{Initialize$h$}($\ov x, \ov s, H_{w,\ov x}, L_{\ov x}$)
\EndProcedure
\Procedure{Initialize$h$}{$\ov x, \ov s \in \R^{n_{\tot}}, H_{w,\ov x}, L_{\ov x} \in \R^{n_{\tot} \times n_{\tot}}$}
    \For{$i \in [n]$}
        \State $( \ov\delta_{\mu})_i \gets - \frac{ \alpha \sinh( \frac{\lambda}{w_i} \gamma_i( \ov x, \ov s, \ov t ) ) }{ \gamma_i( \ov x, \ov s, \ov t ) } \cdot \mu_i( \ov x, \ov s, \ov t )$
        \State $\ov\alpha \gets \ov\alpha + w_i^{-1} \cosh^2( \frac{\lambda}{w_i} \gamma_i( \ov x, \ov s, \ov t ) )$
    \EndFor 
    \State $h \gets L_{\ov x}^{-1} \ov\delta_{\mu}$, $\wt h \gets L_{\ov x}^{-1} A^\top$, $c_s \gets H_{w,\ov x}^{-1/2} \ov\delta_\mu$
    \State $\wt u \gets \wt h^\top \wt h$, $u \gets \wt h^\top h$
\EndProcedure 
\Procedure{Move}{$ $}
    \State $\beta_x \gets \beta_x + \ov t \cdot (\ov\alpha)^{-1/2}$
    \State $\wt \beta_x \gets \wt \beta_x + \ov t \cdot (\ov\alpha)^{-1/2} \cdot \wt u^{-1} u$
    \State $\beta_{c_s} \gets \beta_{c_s} + \ov t \cdot (\ov\alpha)^{-1/2}$
    \State $\beta_s \gets \beta_s + \ov t^2 \cdot (\ov\alpha)^{-1/2}$
    \State $\wt \beta_s \gets \wt \beta_s + \ov t^2 \cdot (\ov\alpha)^{-1/2} \cdot \wt u^{-1} u$
    \State \Return $\beta_x, \beta_s, \beta_{c_s}, \wt \beta_x, \wt \beta_s$
\EndProcedure 
\Procedure{Output}{$ $}
    \State \Return $\wh x + H_{w,\ov x}^{-1/2} \cW^\top (h \beta_x - \wt h \wt \beta_x + \epsilon_x), \wh s + H_{w,\ov x}^{1/2} c_s \beta_{c_s} - H_{w,\ov x}^{1/2} \cW^{\top} (h \beta_s - \wt h \wt \beta_s + \epsilon_s)$
\EndProcedure
\Procedure{Query$x$}{$i\in [n]$}
    \State \Return $\wh x_i + H_{w,\ov x,(i,i)}^{-1/2} (\cW^\top(h \beta_x - \wt h \wt \beta_x + \epsilon_x))_i$
\EndProcedure
\Procedure{Query$s$}{$i\in [n]$}
    \State \Return $\wh s_i + H_{w,\ov x,(i,i)}^{1/2} c_{s,i} \beta_{c_s} - H_{w,\ov x,(i,i)}^{1/2} (\cW^\top (h \beta_s - \wt h \wt \beta_s + \epsilon_s))_i$
\EndProcedure
\State {\bf end data structure}
\end{algorithmic}
\end{algorithm}

\begin{algorithm}[!ht]\caption{Algorithm~\ref{alg:treewidth:exactds-part1} continued.
}\label{alg:treewidth:exactds-part2}
\begin{algorithmic}[1]
\State {\bf data structure} \textsc{ExactDS} \Comment{Theorem~\ref{thm:treewidth:exactds}}
\Procedure{Update}{$\delta_{\ov x}, \delta_{\ov s} \in \R^{n_{\tot}}$}
    \State $\Delta_{H_{w,\ov x}} \gets \nabla^2 \phi_w( \ov x + \delta_{\ov x}) - H_{w,\ov x}$
    \Comment{$\Delta_{H_{w,\ov x}}$ is non-zero only for diagonal blocks $(i,i)$ for which $\delta_{\ov x,i} \ne 0$}
    \State Compute $\Delta_{L_{\ov x}}$ where $(L_{\ov x}+\Delta_{L_{\ov x}})(L_{\ov x}+\Delta_{L_{\ov x}})^\top = B_{\ov x}+\ov t \Delta_{H_{w, \ov x}}$
    \State \textsc{Update$h$}$( \delta_{\ov x} , \delta_{\ov s}, \Delta_{H_{w,\ov x}}, \Delta_{L_{\ov x}})$
    \State \textsc{Update$\cW$}$(\Delta_{H_{w,\ov x}}, \Delta_{L_{\ov x}})$
    \State $\ov x \gets \ov x + \delta_{ \ov x}$, $\ov s \gets \ov s + \delta_{\ov s}$
    \State $H_{w,\ov x} \gets H_{w,\ov x} + \Delta_{H_{w,\ov x}}$, $B_{\ov x} \gets B_{\ov x} + \ov t \Delta_{H_{w,\ov x}}$, $L_{\ov x} \gets L_{\ov x} + \Delta_{L_{\ov x}}$
    \State \Return $\delta_{h}, \delta_{\wt h}, \delta_{\epsilon_x}, \delta_{\epsilon_s}, \delta_{H_{w,\ov x}^{1/2} \wh x}, \delta_{H_{w,\ov x}^{-1/2} \wh s}, \delta_{c_s}$
\EndProcedure 
\Procedure{Update$h$}{$\delta_{\ov x}$, $\delta_{\ov s} \in \R^{n_{\tot}}$, $\Delta_{H_{w,\ov x}}$, $\Delta_{L_{\ov x}} \in \R^{n_{\tot} \times n_{\tot}}$}
    \State $S \gets \{ i \in [n] ~|~ \delta_{\ov x,i} \ne 0 \mathrm{~or~} \delta_{\ov s,i} \ne 0 \}$
    \State $\delta_{ \ov\delta_{\mu} } \gets 0$
    \For{$i \in S$}
        \State Let $\gamma_i = \gamma_i(\ov x, \ov s, \ov t)$, $\gamma_i^{\new} = \gamma_i( \ov x + \delta_{\ov x}, \ov s + \delta_{\ov s}, \ov t )$, $\mu_i^{\new} = \mu_i ( \ov x + \delta_{\ov x}, \ov s + \delta_{\ov s}, \ov t )$
        \State $\ov\alpha \gets \ov\alpha - w_i^{-1} \cosh^2( \frac{\lambda}{w_i} \gamma_i ) + w_i^{-1} \cosh^2( \frac{\lambda}{ w_i } \gamma_i^{\new} )$
        \State $\delta_{\ov\delta_{\mu},i} \gets - \alpha \sinh( \frac{\lambda}{w_i} {\gamma}_i^{\new} ) \cdot \frac{ 1 }{ \gamma^{\new}_i } \cdot \mu^{\new}_i - \ov\delta_{\mu,i}$
    \EndFor
    \State $\delta_h \gets L_{\ov x}^{-1} \delta_{\ov\delta_\mu} -(L_{\ov x}+\Delta_{L_{\ov x}})^{-1} \Delta_{L_{\ov x}} (h + L_{\ov x}^{-1} \delta_{\ov\delta_\mu})$
    \State $\delta_{c_s} \gets \Delta_{H_{w,\ov x}^{-1/2}} (\ov\delta_\mu + \delta_{\ov\delta_\mu}) + H_{w,\ov x}^{-1/2} \delta_{\ov\delta_\mu}$
    \State $\delta_{\wt h} \gets -(L_{\ov x}+\Delta_{L_{\ov x}})^{-1} \Delta_{L_{\ov x}} \wt h$
    \State $\delta_{\hat s} \gets -\delta_{\ov\delta_\mu} \beta_{c_s} + \Delta_{H_{w,\ov x}} L_{\ov x}^{-\top} (h \beta_s - \wt h \wt \beta_s + \epsilon_s)$ \Comment{uses the values of $h,\wt h,\epsilon_s,L_{\ov x}$ before this call, and is supported on the blocks in $S$}
    \State $\delta_{\epsilon_x} \gets -\delta_h \beta_x + \delta_{\wt h} \wt \beta_x$
    \State $\delta_{\epsilon_s} \gets -\delta_h \beta_s + \delta_{\wt h} \wt \beta_s$
    \State $\delta_{\wt u} \gets \delta_{\wt h}^\top (\wt h + \delta_{\wt h}) + \wt h^\top \delta_{\wt h}$
    \State $\delta_{u} \gets \delta_{\wt h}^\top (h + \delta_{h}) + \wt h^\top \delta_h$
    \label{line:alg:treewidth:exactds-part2-updateh-breakpoint}
    \State $\delta_{H_{w,\ov x}^{1/2}\wh x} \gets (H_{w,\ov x}+\Delta_{H_{w,\ov x}})^{1/2}\wh x - H_{w,\ov x}^{1/2}\wh x$, $\delta_{H_{w,\ov x}^{-1/2}\wh s} \gets (H_{w,\ov x}+\Delta_{H_{w,\ov x}})^{-1/2}(\wh s+\delta_{\hat s}) - H_{w,\ov x}^{-1/2}\wh s$
    \State $\ov\delta_\mu \gets \ov\delta_\mu + \delta_{\ov\delta_\mu}$, $h\gets h+\delta_h$, $\wt h \gets \wt h+\delta_{\wt h}$, $c_s \gets c_s + \delta_{c_s}$, $\wh s \gets \wh s + \delta_{\hat s}$, $\epsilon_x\gets \epsilon_x + \delta_{\epsilon_x}$, $\epsilon_s \gets \epsilon_s + \delta_{\epsilon_s}$, $\wt u \gets \wt u + \delta_{\wt u}$, $u \gets u + \delta_{u}$
\EndProcedure 
\Procedure{Update$\cW$}{$\Delta_{H_{w,\ov x}}, \Delta_{L_{\ov x}} \in \R^{n_{\tot}}$}
    \State $\delta_{\epsilon_x} \gets \Delta_{L_{\ov x}}^{\top} L_{\ov x}^{-\top} (h \beta_x - \wt h \wt \beta_x + \epsilon_x)$
    \State $\delta_{\epsilon_s} \gets \Delta_{L_{\ov x}}^{\top} L_{\ov x}^{-\top} (h \beta_s - \wt h \wt \beta_s + \epsilon_s)$
    \label{line:alg:treewidth:exactds-part2-updateW-breakpoint}
    \State $\epsilon_x \gets \epsilon_x + \delta_{\epsilon_x}$, $\epsilon_s \gets \epsilon_s + \delta_{\epsilon_s}$
\EndProcedure 
\State {\bf end data structure}
\end{algorithmic}
\end{algorithm}

\begin{proof}[Proof of Theorem~\ref{thm:treewidth:exactds}]
By combining Lemma~\ref{lem:treewidth:exactds-correct} and~\ref{lem:treewidth:exactds-time}.
\end{proof}

\begin{lemma} \label{lem:treewidth:exactds-correct}
\textsc{ExactDS} correctly maintains an implicit representation of $(x,s)$, i.e., invariant
\begin{align*}
&x = \wh x + H_{w,\ov x}^{-1/2} \cW^\top ( h \beta_x - \wt h \wt \beta_x +\epsilon_x), \\
&s = \wh s + H_{w,\ov x}^{1/2} c_s \beta_{c_s}- H_{w,\ov x}^{1/2} \cW^\top ( h \beta_s - \wt h \wt \beta_s +\epsilon_s),\\
&h = L_{\ov x}^{-1} \ov\delta_\mu, \qquad c_s = H_{w,\ov x}^{-1/2} \ov\delta_\mu, \qquad \wt h=L_{\ov x}^{-1} A^\top, \\
&\wt u = \wt h^\top \wt h, \qquad u = \wt h^\top h, \\
& \ov\alpha = \sum_{i\in [n]} w_i^{-1} \cosh^2(\frac{\lambda}{w_i} \gamma_i(\ov x, \ov s, \ov t)),\\
& \ov\delta_\mu = \ov\alpha^{1/2} \delta_\mu(\ov x, \ov s, \ov t)
\end{align*}
always holds after every external call, and return values of the queries are correct.
\end{lemma}
\begin{proof}
\textsc{Initialize}:
By checking the definitions we see that all invariants are satisfied after \textsc{Initialize}.

\textsc{Move}:
By comparing the implicit representation~\eqref{eqn:thm:treewidth:exactds:x-rep}\eqref{eqn:thm:treewidth:exactds:s-rep} and the robust central path step~\eqref{eqn:thm:treewidth:exactds:x-step}\eqref{eqn:thm:treewidth:exactds:s-step}, we see that \textsc{Move} updates $(x,s)$ correctly.

\textsc{Update}:
We would like to prove that \textsc{Update} correctly updates the values of $h$, $c_s$, $\wt h$, $\wt u$, $u$, $\ov\alpha$, $\ov\delta_\mu$, while preserving the values of $(x,s)$.

First note that $H_{w,\ov x}$, $B_{\ov x}$, $L_{\ov x}$ are updated correctly.
The remaining updates are separated into two steps: \textsc{Update}$h$ and \textsc{Update}$\cW$.

\textbf{Step} \textsc{Update}$h$:
The first few lines of \textsc{Update}$h$ updates $\ov\alpha$ and $\ov\delta_\mu$ correctly.

We define $H_{w,\ov x}^{\new} := H_{w,\ov x} + \Delta_{H_{w,\ov x}}$, $B_{\ov x}^\new := B_{\ov x} + \Delta_{B_{\ov x}}$, $L_{\ov x}^\new := L_{\ov x} + \Delta_{L_{\ov x}}$, $\ov\delta_\mu^\new := \ov\delta_\mu + \delta_{\ov\delta_\mu}$, and so on.
Immediately after Algorithm~\ref{alg:treewidth:exactds-part2}, Line~\ref{line:alg:treewidth:exactds-part2-updateh-breakpoint}, we have
\begin{align*}
h + \delta_h &= L_{\ov x}^{-1} \ov\delta_\mu + L_{\ov x}^{-1} \delta_{\ov\delta_\mu} -(L_{\ov x}+\Delta_{L_{\ov x}})^{-1} \Delta_{L_{\ov x}} (L_{\ov x}^{-1} \ov\delta_\mu + L_{\ov x}^{-1} \delta_{\ov\delta_\mu})\\
&= (L_{\ov x}^{-1}-(L_{\ov x}+\Delta_{L_{\ov x}})^{-1} \Delta_{L_{\ov x}} L_{\ov x}^{-1}) \ov\delta_\mu^\new \\
&= (L_{\ov x}^{\new})^{-1} \ov\delta_\mu^\new,\\
c_s + \delta_{c_s} &= H_{w,\ov x}^{-1/2} \ov\delta_\mu + \Delta_{H_{w,\ov x}^{-1/2}} (\ov\delta_\mu + \delta_{\ov\delta_\mu}) + H_{w,\ov x}^{-1/2} \delta_{\ov\delta_\mu}\\
&= (H_{w,\ov x}^\new)^{-1/2} \ov\delta_\mu^\new, \\
\wt h + \delta_{\wt h} &= L_{\ov x}^{-1} A^\top -(L_{\ov x}+\Delta_{L_{\ov x}})^{-1} \Delta_{L_{\ov x}} A^\top  \\
&= (L_{\ov x}^{-1}-(L_{\ov x}+\Delta_{L_{\ov x}})^{-1} \Delta_{L_{\ov x}} L_{\ov x}^{-1}) A^\top \\
&= (L_{\ov x}^{\new})^{-1} A^\top.
\end{align*}
So $h, c_s, \wt h$ are updated correctly.
Also
\begin{align*}
    \wt u + \delta_{\wt u} &= \wt h^\top \wt h + \delta_{\wt h}^\top (\wt h + \delta_{\wt h}) + \wt h^\top \delta_{\wt h} = (\wt h+\delta_{\wt h})^\top (\wt h+\delta_{\wt h}), \\
    u + \delta_u &= \wt h^\top h + \delta_{\wt h}^\top (h + \delta_{h}) + \wt h^\top \delta_h
    = (\wt h + \delta_{\wt h})^\top (h + \delta_h).
\end{align*}
So $\wt u$ and $u$ are maintained correctly.
Furthermore, immediately after Algorithm~\ref{alg:treewidth:exactds-part2}, Line~\ref{line:alg:treewidth:exactds-part2-updateh-breakpoint}, we have
\begin{align*}
&~ (\wh x + L_{\ov x}^{-\top} (h^\new \beta_x - \wt h^\new \wt \beta_x + \epsilon_x^\new) )
- (\wh x + L_{\ov x}^{-\top} (h \beta_x - \wt h\wt \beta_x + \epsilon_x) ) \\
=&~ L_{\ov x}^{-\top} (\delta_h \beta_x - \delta_{\wt h} \wt \beta_x + \delta_{\epsilon_x}) \\
=&~0.
\end{align*}
Therefore, after \textsc{Update}$h$ finishes, we have
\begin{align*}
x = \hat x + L_{\ov x}^{-\top} (h \beta_x - \wt h \wt \beta_x + \epsilon_x).
\end{align*}
For $s$, note that $H_{w,\ov x}^{1/2}\cW^\top = H_{w,\ov x} L_{\ov x}^{-\top}$ and $H_{w,\ov x}^{1/2} c_s = \ov\delta_\mu$. Using $\delta_h \beta_s - \delta_{\wt h} \wt \beta_s + \delta_{\epsilon_s}=0$ and the definition of $\delta_{\wh s}$, we have
\begin{align*}
&~(\hat s^\new + (H_{w,\ov x}^\new)^{1/2} c_s^\new \beta_{c_s} - H_{w,\ov x}^\new L_{\ov x}^{-\top}(h^\new \beta_s-\wt h^\new \wt \beta_s+\epsilon_s^\new))\\
&-(\hat s + H_{w,\ov x}^{1/2} c_s \beta_{c_s} - H_{w,\ov x} L_{\ov x}^{-\top}(h\beta_s-\wt h \wt \beta_s+\epsilon_s)) \\
=&~ \delta_{\wh s} + \delta_{\ov\delta_\mu} \beta_{c_s} - \Delta_{H_{w,\ov x}} L_{\ov x}^{-\top}(h\beta_s-\wt h \wt \beta_s+\epsilon_s) - H_{w,\ov x}^\new L_{\ov x}^{-\top}(\delta_h \beta_s - \delta_{\wt h} \wt \beta_s + \delta_{\epsilon_s}) \\
=&~ 0.
\end{align*}
Therefore, after \textsc{Update}$h$ finishes, we have
\begin{align*}
s = \hat s + (H_{w,\ov x}^\new)^{1/2}c_s \beta_{c_s} - H_{w,\ov x}^\new L_{\ov x}^{-\top} (h \beta_s - \wt h \wt \beta_s + \epsilon_s).
\end{align*}
So $x$ and $s$ are both updated correctly. This proves the correctness of \textsc{Update}$h$.

\textbf{Step} $\textsc{Update}\cW$:
Define $\epsilon_x^\new := \epsilon_x + \delta_{\epsilon_x}$, $\epsilon_s^\new := \epsilon_s + \delta_{\epsilon_s}$.
Immediately after Algorithm~\ref{alg:treewidth:exactds-part2}, Line~\ref{line:alg:treewidth:exactds-part2-updateW-breakpoint}, we have
\begin{align*}
&~ (\wh x + (L_{\ov x}^\new)^{-\top} (h \beta_x - \wt h\wt \beta_x + \epsilon_x^\new) )
- (\wh x + L_{\ov x}^{-\top} (h \beta_x - \wt h\wt \beta_x + \epsilon_x) ) \\
=&~ ((L_{\ov x}^\new)^{-\top} - L_{\ov x}^{-\top}) (h \beta_x - \wt h\wt \beta_x + \epsilon_x)
+  (L_{\ov x}^\new)^{-\top} \delta_{\epsilon_x}\\
=&~0,\\
&~(\hat s + (H_{w,\ov x}^\new)^{1/2} c_s \beta_{c_s} - H_{w,\ov x}^\new (L_{\ov x}^\new)^{-\top}(h \beta_s-\wt h\wt \beta_s+\epsilon_s^\new))\\
&~ -(\hat s + (H_{w,\ov x}^\new)^{1/2} c_s \beta_{c_s} - H_{w,\ov x}^\new L_{\ov x}^{-\top}(h \beta_s-\wt h\wt \beta_s+\epsilon_s)) \\
=&~ H_{w,\ov x}^\new\big((-(L_{\ov x}^\new)^{-\top} + L_{\ov x}^{-\top}) (h \beta_s - \wt h \wt \beta_s + \epsilon_s) - (L_{\ov x}^\new)^{-\top} \delta_{\epsilon_s}\big) \\
=&~ 0.
\end{align*}
Therefore, after \textsc{Update}$\cW$ finishes, we have
\begin{align*}
x &= \hat x + (L_{\ov x}^\new)^{-\top} (h \beta_x - \wt h \wt \beta_x + \epsilon_x), \\
s &= \hat s + (H_{w,\ov x}^\new)^{1/2}c_s \beta_{c_s} - H_{w,\ov x}^\new (L_{\ov x}^\new)^{-\top} (h \beta_s - \wt h \wt \beta_s + \epsilon_s).
\end{align*}
So $x$ and $s$ are both updated correctly. This proves the correctness of \textsc{Update}$\cW$.
\end{proof}

\begin{lemma} \label{lem:treewidth:exactds-time}
We bound the running time of \textsc{ExactDS} as following.
\begin{enumerate}[label=(\roman*)]
    \item \textsc{ExactDS.Initialize} (Algorithm~\ref{alg:treewidth:exactds-part1}) runs in $\wt O(n \tau^{\omega-1} + n \tau m + n m^{\omega-1})$ time. \label{item:treewidth:exactds-init-time}
    \item \textsc{ExactDS.Move} (Algorithm~\ref{alg:treewidth:exactds-part1}) runs in $\wt O(m^\omega)$ time. \label{item:treewidth:exactds-move-time}
    \item \textsc{ExactDS.Output} (Algorithm~\ref{alg:treewidth:exactds-part1}) runs in $\wt O(n \tau m)$ time and correctly outputs $(x,s)$. \label{item:treewidth:exactds-output-time}
    \item \textsc{ExactDS.Query$x$} and \textsc{ExactDS.Query$s$} (Algorithm~\ref{alg:treewidth:exactds-part1}) runs in $\wt O(\tau^2 m)$ time and returns the correct answer. \label{item:treewidth:exactds-query-time}
    \item \textsc{ExactDS.Update} (Algorithm~\ref{alg:treewidth:exactds-part1}) runs in $\wt O((\tau^2 m+\tau m^2) (\|\delta_{\ov x}\|_0 + \|\delta_{\ov s}\|_0))$ time.
    Furthermore, $h, \epsilon_x, \epsilon_s$ change in $O(\tau (\|\delta_{\ov x}\|_0 + \|\delta_{\ov s}\|_0))$ coordinates,
    $\wt h$ changes in $\wt O(\tau m (\|\delta_{\ov x}\|_0 + \|\delta_{\ov s}\|_0))$ coordinates,
    and $H_{\ov x}^{1/2} \wh x, H_{\ov x}^{-1/2} \wh s, c_s$ change in $O(\|\delta_{\ov x}\|_0 + \|\delta_{\ov s}\|_0)$  coordinates.
    \label{item:treewidth:exactds-update-time}
\end{enumerate}
\end{lemma}
\begin{proof}
    \begin{enumerate}[label=(\roman*)]
    \item Computing $L_{\ov x}$ takes $\wt O(n \tau^{\omega-1})$ time by Lemma~\ref{lem:cholesky_fast}.
    Computing $h$ and $\wt h$ takes $\wt O(n \tau m)$ by Lemma~\ref{lem:mat_vec_mult_time}\ref{item:lem_mat_vec_mult_time_Linvv}.\footnote{Here we compute $\wt h$ by computing $\wt h_{*,i} = L_{\ov x}^{-1} (A_{i,*})^\top$ for $i\in [m]$ independently. Using fast rectangular matrix multiplication is possible to improve this running time and other similar places. We keep the current bounds for simplicity.}
    Computing $\wt u$ and $u$ takes $\Tmat(m,n,m) = \wt O(n m^{\omega-1})$ time.
    All other operations are cheap.
    \item Computing $\wt u^{-1}$ takes $\wt O(m^\omega)$ time. All other operations take $O(m^2)$ time.
    \item Running time is by Lemma~\ref{lem:mat_vec_mult_time}\ref{item:lem_mat_vec_mult_time_WTv}. Correctness is by Lemma~\ref{lem:treewidth:exactds-correct}.
    \item Running time is by Lemma~\ref{lem:mat_vec_mult_coord_time}\ref{item:mat_vec_mult_coord_time_Wtvi}. Correctness is by Lemma~\ref{lem:treewidth:exactds-correct}.
    \item Computing $\Delta_{L_{\ov x}}$ takes $\wt O(\tau^2\|\delta_{\ov x}\|_0)$ time by Lemma~\ref{lem:cholesky_update}, since each changed diagonal block of $\ov t\Delta_{H_{w,\ov x}}$ is a sum of $O(n_i)=O(1)$ rank-one updates and downdates. It is easy to see that $\nnz(\Delta_{H_{w,\ov x}}) = O(\|\delta_{\ov x}\|_0)$ and $\nnz(\Delta_{L_{\ov x}})= \wt O(\tau^2 \|\delta_{\ov x}\|_0)$. It remains to analyze \textsc{Update$h$} and \textsc{Update$\cW$}.
    For simplicity, we write $k = \|\delta_{\ov x}\|_0 + \|\delta_{\ov s}\|_0$ in this proof only.
    \begin{itemize}
        \item \textsc{Update$h$}: Updating $\ov \alpha$ and computing $\delta_{\ov \delta_\mu}$ takes $O(k)$ time.
        Also, $\|\delta_{\ov \delta_\mu}\|_0=O(k)$.

        Computing $\delta_h$ takes $\wt O(\tau^2 k)$ time by Lemma~\ref{lem:mat_vec_mult_time}\ref{item:lem_mat_vec_mult_time_Linvv_path}, applied to each of the $O(k)$ paths in the elimination tree.
        Also, $\delta_h$ is supported on $O(k)$ paths in the elimination tree, so $\|\delta_h\|_0 = \wt O(\tau k)$.
        Similarly we see that computing $\delta_{\wt h}$ take $\wt O(\tau^2 m k)$ time and $\nnz(\delta_{\wt h}) = \wt O(\tau m k)$.

        Computing $\delta_{c_s}$ takes $O(\tau^2 k)$ time. Computing $\delta_{\wh s}$ takes $\wt O(\tau^2 m k)$ time, since it only needs the restriction of $L_{\ov x}^{-\top}(h \beta_s - \wt h \wt \beta_s + \epsilon_s)$ to the blocks in $S$, computed as in \textsc{Update$\cW$} below. Also $\|\delta_{c_s}\|_0, \|\delta_{\wh s}\|_0 = O(k)$.
        
        Computing $\delta_{\epsilon_x}$ and $\delta_{\epsilon_s}$ takes $O(\tau m k)$ time after computing $\delta_h$ and $\delta_{\wt h}$. Furthermore, $\|\delta_{\epsilon_x}\|_0,\|\delta_{\epsilon_s}\|_0 = O(\tau k)$.

        Computing $\delta_{\wt u}$ takes $\Tmat(m, \tau k, m) = \wt O(\tau m^2 k)$ time.
        Computing $\delta_{u}$ takes $\wt O(\tau m k)$ time.
        \item \textsc{Update$\cW$}:
        To compute $\delta_{\epsilon_x}$ and $\delta_{\epsilon_s}$, we first compute $L_{\ov x}^{-\top}(h \beta_x - \wt h \wt \beta_x + \epsilon_x)$ and $L_{\ov x}^{-\top}(h \beta_s - \wt h \wt \beta_s + \epsilon_s)$, where $S\subseteq [n_\tot]$ is the row support of $\Delta_{L_{\ov x}}$, which can be decomposed into at most $O(\|\delta_{\ov x}\|_0)$ paths. 
        This takes $\wt O(\tau^2 m \|\delta_{\ov x}\|_0)$ time by Lemma~\ref{lem:mat_vec_mult_coord_time}\ref{item:mat_vec_mult_coord_time_LinvtvS} (the extra $m$ factor is due to $\wt h$).
    \end{itemize}
    Combining everything we finish the proof of running time of \textsc{ExactDS.Update}.\qedhere
    \end{enumerate}
\end{proof}

\subsubsection{\textsc{ApproxDS}} \label{sec:treewidth:cpm-ds:approxds}
In this section we present the data structure \textsc{ApproxDS}.
Given \textsc{BatchSketch}, a data structure maintaining a sketch of the primal-dual pair $(x,s)\in \R^{n_\tot} \times \R^{n_\tot}$, \textsc{ApproxDS} maintains a sparsely-changing $\ell_\infty$-approximation of $(x,s)$. This data structure is a slight variation of \textsc{ApproxDS} in \cite{gs22}.

\begin{algorithm}[!ht] \caption{The \textsc{ApproxDS} data structure used in Algorithm~\ref{alg:treewidth:cpm}.}\label{alg:treewidth:approxds-part1}
\begin{algorithmic}[1]
\small
\State {\bf data structure} \textsc{ApproxDS} \Comment{Theorem~\ref{thm:treewidth:approxds}}
\State {\bf private : members}
\State \hspace{4mm} $\epsilon_{\apx,x}, \epsilon_{\apx,s}\in \R$
\State \hspace{4mm} $\ell \in \bN$
\State \hspace{4mm} \textsc{BatchSketch} $\mathsf{bs}$ \Comment{This maintains a sketch of $H_{w,\ov x}^{1/2} x$ and $H_{w,\ov x}^{-1/2} s$. See Algorithm~\ref{alg:treewidth:batchsketch-part1},~\ref{alg:treewidth:batchsketch-part2},~\ref{alg:treewidth:batchsketch-part3}.}

\State \hspace{4mm} \textsc{ExactDS*} $\mathsf{exact}$ \Comment{This is a pointer to the \textsc{ExactDS} (Algorithm~\ref{alg:treewidth:exactds-part1}, \ref{alg:treewidth:exactds-part2}) we maintain in parallel to \textsc{ApproxDS}.}

\State \hspace{4mm} $\wt x, \wt s \in \R^{n_\tot}$
\Comment{$(\wt x, \wt s)$ is a sparsely-changing approximation of $(x, s)$. They have the same value as $(\ov x, \ov s)$, but for these local variables we use $(\wt x, \wt s)$ to avoid confusion.}

\State {\bf end members}

\Procedure{\textsc{Initialize}}{$x,s \in \R^{n_\tot},
h\in \R^{n_\tot}, \wt h \in \R^{n_\tot\times m}, \epsilon_x, \epsilon_s, H_{w,\ov x}^{1/2} \wh x, H_{w,\ov x}^{-1/2} \wh s, c_s \in \R^{n_\tot}, \beta_x, \beta_s, \beta_{c_s} \in \R, \wt \beta_x, \wt \beta_s \in \R^{m},
q\in \bN,
\textsc{ExactDS*}~\mathsf{exact}, \epsilon_{\apx,x}, \epsilon_{\apx,s}, \delta_{\apx}\in \R$}
    \State $\ell \gets 0$, $q \gets q$
    \State $\epsilon_{\apx,x} \gets \epsilon_{\apx,x}, \epsilon_{\apx,s} \gets \epsilon_{\apx,s}$
    
    \State $\mathsf{bs}.\textsc{Initialize}(x,h,\wt h,\epsilon_x, \epsilon_s, H_{w,\ov x}^{1/2} \wh x, H_{w,\ov x}^{-1/2} \wh s, c_s, \beta_x, \beta_s, \beta_{c_s}, \wt \beta_x, \wt \beta_s, \delta_{\apx}/q)$
    \Comment{Algorithm~\ref{alg:treewidth:batchsketch-part1}}
    \State $\wt x \gets x, \wt s \gets s$
    \State $\mathsf{exact} \gets \mathsf{exact}$
\EndProcedure
\Procedure{\textsc{Update}}{$\delta_{\ov x} \in \R^{n_\tot}, \delta_{h}\in \R^{n_\tot}, \delta_{\wt h} \in \R^{n_\tot \times m}, \delta_{\epsilon_x}, \delta_{\epsilon_s}, \delta_{H_{w,\ov x}^{1/2} \wh{x}}, \delta_{ H_{w,\ov x}^{-1/2} \wh{s} } , \delta_{c_s} \in \R^{n_\tot}$}
    \State {$\mathsf{bs}.\textsc{Update}(\delta_{\ov x}, \delta_{h}, \delta_{\wt h}, \delta_{\epsilon_x}, \delta_{\epsilon_s}, \delta_{H_{w,\ov x}^{1/2} \wh x}, \delta_{H_{w,\ov x}^{-1/2} \wh s}, \delta_{c_s})$}
    \Comment{Algorithm~\ref{alg:treewidth:batchsketch-part2}}
    \State $\ell \gets \ell+1$
\EndProcedure 
\Procedure{MoveAndQuery}{$\beta_x, \beta_s, \beta_{c_s}\in \R$, $\wt \beta_x, \wt \beta_s \in \R^m$}
    \State $\mathsf{bs}.\textsc{Move}(\beta_x, \beta_s, \beta_{c_s}, \wt \beta_x, \wt \beta_s)$
    \Comment{Algorithm~\ref{alg:treewidth:batchsketch-part2}. Do not update $\ell$ yet}
    \State $\delta_{\wt x} \gets \textsc{Query$x$}(\epsilon_{\apx,x}/(2\log_2 q+1))$ \Comment{Algorithm~\ref{alg:treewidth:approxds-part2}}
    \State $\delta_{\wt s} \gets \textsc{Query$s$}(\epsilon_{\apx,s}/(2\log_2 q+1))$ \Comment{Algorithm~\ref{alg:treewidth:approxds-part2}}
    \For{all $i\in [n]$ with $\delta_{\wt x,i}\ne 0$} \label{line:alg:treewidth:approxds-resync-s}
        \State $s_i \gets \mathsf{exact}.\textsc{Query}s(i)$, $\delta_{\wt s,i} \gets s_i-\wt s_i$
        \Comment{Whenever $\wt x_i$ changes, we set $\wt s_i$ to its exact value}
    \EndFor
    \State $\wt x \gets \wt x + \delta_{\wt x}$, $\wt s \gets \wt s + \delta_{\wt s}$
    \State \Return $(\delta_{\wt x}, \delta_{\wt s})$
\EndProcedure
\State {\bf end data structure}
\end{algorithmic}
\end{algorithm}

\begin{algorithm}[!ht] \caption{\textsc{ApproxDS} Algorithm~\ref{alg:treewidth:approxds-part1} continued.}\label{alg:treewidth:approxds-part2}
\begin{algorithmic}[1]
\State {\bf data structure} \textsc{ApproxDS} \Comment{Theorem~\ref{thm:treewidth:approxds}}
\State {\bf private:}
\Procedure{\textsc{Query$x$}}{$\epsilon \in \R$}
    \State $\cI \gets \emptyset$
    \For{$j=0\to \lfloor \log_2 (\ell+1)\rfloor$}
    \If{$(\ell+1) \bmod {2^j}=0$}
        \State $\ell' \gets \ell-2^j+1$
        \Comment{The window $[\ell',\ell]$ consists of the last $2^j$ steps}
        \For{all $i \in \mathsf{bs}.\textsc{Query$x$}(\ell', \epsilon)$}
        \Comment{Algorithm~\ref{alg:treewidth:batchsketch-part3}}
            \State Compute $x^{(\ell')}_i$ and $x^{(\ell+1)}_i$ by $\mathsf{exact}.\textsc{Query}x(i)$ at timestamp $\ell'$ and at the current time
            \Comment{Remark~\ref{rmk:snapshot}}
            \If{$\|H_{w,\ov x^{(\ell')}}^{1/2} x^{(\ell')}_i - H_{w,\ov x^{(\ell)}}^{1/2} x^{(\ell+1)}_i\|_2 \ge \epsilon$} \label{line:alg:treewidth:approxds-queryx-verify}
                \State $\cI \gets \cI \cup \{i\}$
                \Comment{Keep $i$ only if its exact change over the window is large}
            \EndIf
        \EndFor
    \EndIf
    \EndFor
    \State $\delta_{\wt x} \gets 0$
    \For{all $i\in \cI$} \label{line:alg:treewidth:approxds-queryx-xinI}
        \State $x_i \gets \mathsf{exact}.\textsc{Query}x(i)$
        \Comment{Algorithm~\ref{alg:treewidth:exactds-part1}}
        \If{$\|\wt x_i - x_i\|_{\wt x_i} > \epsilon$}  \label{line:alg:treewidth:approxds-queryx-checkxi}
            \State $\delta_{\wt x,i} \gets x_i-\wt x_i$
        \EndIf
    \EndFor
    \State \Return $\delta_{\wt x}$
\EndProcedure
\Procedure{\textsc{Query$s$}}{$\epsilon \in \R$}
    \State Same as \textsc{Query$x$}, except for replacing
    $x,\wt x,\cdots$ with $s,\wt s,\cdots$, replacing the test in Line~\ref{line:alg:treewidth:approxds-queryx-verify} with ``$\|H_{w,\ov x^{(\ell')}}^{-1/2} s^{(\ell')}_i - H_{w,\ov x^{(\ell)}}^{-1/2} s^{(\ell+1)}_i\|_2 \ge \epsilon$'', and replacing ``$\|\wt x_i - x_i\|_{\wt x_i}$'' in Line~\ref{line:alg:treewidth:approxds-queryx-checkxi} with ``$\|\wt s_i - s_i\|_{\wt x_i}^*$''.
\EndProcedure
\State {\bf end data structure}
\end{algorithmic}
\end{algorithm}

\begin{theorem}\label{thm:treewidth:approxds}
Given parameters $\epsilon_{\apx,x}, \epsilon_{\apx,s} \in (0, 1), \delta_\apx \in (0,\frac12]$, $\zeta_{x}, \zeta_{s}\in \R$ such that 
\begin{align*}
\|H_{w,\ov x^{(\ell)}}^{1/2} x^{(\ell)}-H_{w,\ov x^{(\ell)}}^{1/2} x^{(\ell+1)}\|_2 \le \zeta_{x}, ~~~ \|H_{w,\ov x^{(\ell)}}^{-1/2} s^{(\ell)}-H_{w,\ov x^{(\ell)}}^{-1/2} s^{(\ell+1)}\|_2 \le \zeta_{s}
\end{align*}
for all $\ell \in \{0,\ldots,q-1\}$, whatever the outputs of the data structure are,
data structure \textsc{ApproxDS} (Algorithm~\ref{alg:treewidth:approxds-part1} and Algorithm~\ref{alg:treewidth:approxds-part2}) supports the following operations:
\begin{itemize}
\item $\textsc{Initialize}(x,s \in \R^{n_\tot},
h\in \R^{n_\tot}, \wt h \in \R^{n_\tot\times m}, \epsilon_x, \epsilon_s, H_{w,\ov x}^{1/2} \wh x, H_{w,\ov x}^{-1/2} \wh s, c_s \in \R^{n_\tot}, \beta_x, \beta_s, \beta_{c_s} \in \R, \wt \beta_x, \wt \beta_s \in \R^{m},
q\in \bN,
\textsc{ExactDS*}~\mathsf{exact}, \epsilon_{\apx,x}, \epsilon_{\apx,s}, \delta_{\apx}\in \R)$: Initialize the data structure in $\wt O(n \tau^{\omega-1}+n\tau m)$ time.

\item $\textsc{MoveAndQuery}(\beta_x, \beta_s, \beta_{c_s}\in \R$, $\wt \beta_x, \wt \beta_s \in \R^m)$:
Update values of $\beta_x, \beta_s, \beta_{c_s}, \wt \beta_x, \wt \beta_s$ by calling $\textsc{BatchSketch}.\textsc{Move}$.
This effectively moves $(x^{(\ell)}, s^{(\ell)})$ to $(x^{(\ell+1)}, s^{(\ell+1)})$ while keeping $\ov x^{(\ell)}$ unchanged.

Then return the updates $(\delta_{\wt x},\delta_{\wt s})$ of the approximation, so that $\ov x^{(\ell+1)}=\ov x^{(\ell)}+\delta_{\wt x}$ and $\ov s^{(\ell+1)}=\ov s^{(\ell)}+\delta_{\wt s}$. Let $L_x^{(\ell)}, L_s^{(\ell)} \subseteq [n]$ be the sets of blocks in which $\delta_{\wt x}$ and $\delta_{\wt s}$ are nonzero.

Suppose that the inputs of all calls are determined by the state at initialization and by the earlier outputs of the data structure, as is the case in Algorithm~\ref{alg:treewidth:cpm}. Then with probability at least $1-\delta_\apx$ over the sketching matrix $\Phi$, the following holds for all $\ell\in\{0,\ldots,q-1\}$ simultaneously, even though the inputs depend on $\Phi$ through the earlier outputs. For all $i\in [n]$,
\begin{align*}
\|\ov x^{(\ell+1)}_i - x^{(\ell+1)}_i\|_{\ov x^{(\ell+1)}_i} \le \epsilon_{\apx,x}, \qquad
\|\ov s^{(\ell+1)}_i - s^{(\ell+1)}_i\|_{\ov x^{(\ell+1)}_i}^* \le w_i^{1/2}\epsilon_{\apx,s}.
\end{align*}
Moreover
\begin{align*}
\sum_{0\le \ell \le q-1} |L_x^{(\ell)}| & = \wt O(\epsilon_{\apx,x}^{-2} \zeta_x^2 q^2), \\
\sum_{0\le \ell \le q-1} |L_s^{(\ell)}| & = \wt O((\epsilon_{\apx,x}^{-2} \zeta_x^2+\epsilon_{\apx,s}^{-2} \zeta_s^2) q^2),
\end{align*}
and the total time cost over all queries is at most
\begin{align*}
\wt O\left( (1+\epsilon_{\apx,x}^{-2} \zeta_{x}^2 + \epsilon_{\apx,s}^{-2} \zeta_{s}^2)q^2 \tau^2 m\right).
\end{align*}

\item $\textsc{Update}(\delta_{\ov x} \in \R^{n_\tot}, \delta_{h}\in \R^{n_\tot}, \delta_{\wt h} \in \R^{n_\tot \times m}, \delta_{\epsilon_x}, \delta_{\epsilon_s}, \delta_{H_{w,\ov x}^{1/2} \wh{x}}, \delta_{ H_{w,\ov x}^{-1/2} \wh{s} } , \delta_{c_s} \in \R^{n_\tot})$: Update sketches of $H_{w,\ov x^{(\ell)}}^{1/2} x^{(\ell+1)}$ and $H_{w,\ov x^{(\ell)}}^{-1/2} s^{(\ell+1)}$ by calling \textsc{BatchSketch}.\textsc{Update}.
This effectively moves $\ov x^{(\ell)}$ to $\ov x^{(\ell+1)}$ while keeping $(x^{(\ell+1)}, s^{(\ell+1)})$ unchanged. Then advance timestamp $\ell$.

Each update costs
\begin{align*}
\wt O(\tau^2 m \|\delta_{\ov x}\|_0 + \nnz(\delta_h) + \nnz(\delta_{\wt h}) + \|\delta_{\epsilon_x}\|_0 + \|\delta_{\epsilon_s}\|_0 + \|\delta_{H_{w,\ov x}^{1/2}\wh x}\|_0 + \|\delta_{H_{w,\ov x}^{-1/2}\wh s}\|_0 + \|\delta_{c_s}\|_0)
\end{align*}
time.
\end{itemize}

\end{theorem}

\begin{proof}
The running times of \textsc{Initialize} and \textsc{Update} follow from Theorem~\ref{thm:treewidth:batchsketch}. For \textsc{MoveAndQuery} we follow the proof of \cite[Theorem 4.18]{gs22} and add three things. We take the weights $w_i$ into account. We use the exact test in Line~\ref{line:alg:treewidth:approxds-queryx-verify} of Algorithm~\ref{alg:treewidth:approxds-part2} and the loop in Line~\ref{line:alg:treewidth:approxds-resync-s} of Algorithm~\ref{alg:treewidth:approxds-part1}, which the data structure of~\cite{gs22} does not have. And we show that the guarantee holds although the inputs depend on $\Phi$.

\textbf{Notation.}
Let $\epsilon_x:=\epsilon_{\apx,x}/(2\log_2 q+1)$ and $\epsilon_s:=\epsilon_{\apx,s}/(2\log_2 q+1)$ be the thresholds that \textsc{MoveAndQuery} passes to \textsc{Query}$x$ and \textsc{Query}$s$. We call the call of \textsc{MoveAndQuery} that moves $(x^{(\ell)},s^{(\ell)})$ to $(x^{(\ell+1)},s^{(\ell+1)})$ step $\ell$, for $0\le \ell\le q-1$. We say that $\ov x_i$ is \emph{reset} at step $\ell$ if this step assigns $\delta_{\wt x,i}$, so that $\ov x_i^{(\ell+1)}=x_i^{(\ell+1)}$. Otherwise $\ov x_i^{(\ell+1)}=\ov x_i^{(\ell)}$. We use the same word for $\ov s_i$. The loop in Line~\ref{line:alg:treewidth:approxds-resync-s} resets $\ov s_i$ at every step that resets $\ov x_i$. Recall that $\ov x^{(0)}=x^{(0)}$ and $\ov s^{(0)}=s^{(0)}$.

An interval of steps $W=[\ell',\ell]$ is a \emph{dyadic window} if $\ell'=\ell-2^j+1$ and $2^j$ divides $\ell+1$ for some integer $j\ge 0$. These are the windows that \textsc{Query}$x$ and \textsc{Query}$s$ consider at step $\ell$. For every $j$ the dyadic windows of length $2^j$ are pairwise disjoint, so at most $q/2^j$ of them lie inside $[0,q-1]$. Every nonempty interval of steps inside $[0,q-1]$ is a disjoint union of at most $\max\{1,2\log_2 q\}$ dyadic windows. For a dyadic window $W=[\ell',\ell]$ let
\begin{align*}
v^x_W:=H_{w,\ov x^{(\ell')}}^{1/2} x^{(\ell')} - H_{w,\ov x^{(\ell)}}^{1/2} x^{(\ell+1)}, \qquad
v^s_W:=H_{w,\ov x^{(\ell')}}^{-1/2} s^{(\ell')} - H_{w,\ov x^{(\ell)}}^{-1/2} s^{(\ell+1)}.
\end{align*}
At step $\ell$, the call $\mathsf{bs}.\textsc{Query}x(\ell',\epsilon_x)$ compares sketches of the two vectors whose difference is $v^x_W$, and the test in Line~\ref{line:alg:treewidth:approxds-queryx-verify} reads $\|(v^x_W)_i\|_2\ge\epsilon_x$. The same holds for $s$ with $v^s_W$ and $\epsilon_s$. We call a run of the data structure \emph{accurate} if $\Phi$ is accurate, in the sense of Theorem~\ref{thm:treewidth:batchsketch}, for $v^x_W$ and $v^s_W$ for every dyadic window $W\subseteq[0,q-1]$.

\textbf{Step 1: the run is accurate with probability at least $1-\delta_\apx$.}
Consider the \emph{ideal} data structure, which differs from \textsc{ApproxDS} only in that every call $\mathsf{bs}.\textsc{Query}x(\ell',\epsilon_x)$ is replaced by the set $\{i\in[n]:\|(v^x_W)_i\|_2\ge\epsilon_x\}$, and every call $\mathsf{bs}.\textsc{Query}s(\ell',\epsilon_s)$ by the set $\{i\in[n]:\|(v^s_W)_i\|_2\ge\epsilon_s\}$. The ideal data structure does not use $\Phi$. Since the inputs are determined by the state at initialization and by the earlier outputs, the run with the ideal data structure, which we call the \emph{ideal run}, is determined by the state at initialization. So are its vectors $v^x_W$ and $v^s_W$. The matrix $\Phi$ is drawn in \textsc{Initialize}, independently of this state, so these vectors do not depend on $\Phi$. There are at most $4q$ of them, and the partition tree has $O(n)$ nodes. By Lemma~\ref{lem:jl_matrix} and the choice $r=\Theta(\log^3(n_\tot)\log(q/\delta_\apx))$ in $\mathsf{bs}.\textsc{Initialize}$, with a large enough constant, with probability at least $1-\delta_\apx$ the matrix $\Phi$ is accurate for all of them. Here we use $\delta_\apx\le\frac12$, which gives $\log(q/\delta_\apx)\ge\log 2$.

We claim that on this event the actual run coincides with the ideal run, so that the actual run is accurate. We show by induction on $\ell$ that steps $0,\ldots,\ell-1$ return the same outputs in both runs. Suppose that this holds for $\ell$. Then the inputs of step $\ell$ coincide in both runs, and so do $\ov x^{(t)},\ov s^{(t)}$ for $t\le\ell$ and $x^{(t)},s^{(t)}$ for $t\le\ell+1$. So the vectors $v^x_W$ and $v^s_W$ of the windows $W$ considered at step $\ell$ are those of the ideal run, and $\Phi$ is accurate for them. By Theorem~\ref{thm:treewidth:batchsketch}, the set returned by $\mathsf{bs}.\textsc{Query}x(\ell',\epsilon_x)$ contains $\{i:\|(v^x_W)_i\|_2\ge\epsilon_x\}$. After the test in Line~\ref{line:alg:treewidth:approxds-queryx-verify}, exactly this set is added to $\cI$, as in the ideal run. The same holds for $s$. The remaining operations of step $\ell$ do not use $\Phi$. So step $\ell$ returns the same output in both runs.

\textbf{Step 2: closeness in an accurate run.}
In an accurate run, the argument of Step 1 shows that at every step $\ell$ the set $\cI$ of \textsc{Query}$x$ is the set of all $i$ with $\|(v^x_W)_i\|_2\ge\epsilon_x$ for some dyadic window $W$ that ends at $\ell$. The same holds for $s$.

Fix $i\in[n]$ and a step $\ell$. Let $a\le\ell+1$ be the smallest index such that $\ov x_i$ is not reset at any of the steps $a,\ldots,\ell$. Then $\ov x_i^{(a)}=\cdots=\ov x_i^{(\ell+1)}=:y$, and $y=x_i^{(a)}$, by a reset at step $a-1$ or by initialization. If $a=\ell+1$ we are done, so let $a\le\ell$. Consider a dyadic window $W=[a',b']\subseteq[a,\ell]$ such that $i\notin\cI$ at step $b'$. Then $\|(v^x_W)_i\|_2<\epsilon_x$. Since $\ov x_i^{(a')}=\ov x_i^{(b')}=y$ and $H_{w,\ov x,(i,i)}=w_i\nabla^2\phi_i(\ov x_i)$, this reads
\begin{align}\label{eq:approxds-window-x}
\|x_i^{(a')}-x_i^{(b'+1)}\|_{y}=w_i^{-1/2}\|(v^x_W)_i\|_2<\epsilon_x .
\end{align}
Let $u$ be the last step in $[a,\ell]$ with $i\in\cI$, if there is one. Since $\ov x_i$ is not reset at step $u$, the test in Line~\ref{line:alg:treewidth:approxds-queryx-checkxi} gave $\|y-x_i^{(u+1)}\|_{y}\le\epsilon_x$. If there is no such step, let $u:=a-1$, so that $y=x_i^{(u+1)}$. We write the interval $[u+1,\ell]$, if it is not empty, as a disjoint union of $k$ dyadic windows. At their last steps we have $i\notin\cI$, so~\eqref{eq:approxds-window-x} applies to each of them, and the triangle inequality gives
\begin{align*}
\|\ov x_i^{(\ell+1)}-x_i^{(\ell+1)}\|_{\ov x_i^{(\ell+1)}}=\|y-x_i^{(\ell+1)}\|_{y}\le\|y-x_i^{(u+1)}\|_{y}+k\epsilon_x .
\end{align*}
If $u=a-1$, the first term is zero and $k\le\max\{1,2\log_2q\}$. If $u\ge a$, the first term is at most $\epsilon_x$, and either $u=\ell$ and $k=0$, or $q\ge 2$ and $k\le 2\log_2 q$. In all cases the right-hand side is at most $(2\log_2q+1)\epsilon_x=\epsilon_{\apx,x}$.

For $s$ we argue in the same way. Let $a\le\ell+1$ be the smallest index such that $\ov s_i$ is not reset at any of the steps $a,\ldots,\ell$. Then $\ov s_i^{(a)}=\cdots=\ov s_i^{(\ell+1)}=:z$ and $z=s_i^{(a)}$. A reset of $\ov x_i$ comes with a reset of $\ov s_i$, so $\ov x_i$ is not reset at these steps either, and $\ov x_i^{(a)}=\cdots=\ov x_i^{(\ell+1)}=:y$. This is the point where the loop in Line~\ref{line:alg:treewidth:approxds-resync-s} is needed: without it, $\ov x_i$ and with it the norm could change between two resets of $\ov s_i$. For a dyadic window $W=[a',b']\subseteq[a,\ell]$ such that $i$ is not in the set $\cI$ of \textsc{Query}$s$ at step $b'$, we get
\begin{align*}
\|s_i^{(a')}-s_i^{(b'+1)}\|_{y}^*=w_i^{1/2}\|(v^s_W)_i\|_2<w_i^{1/2}\epsilon_s .
\end{align*}
If $u$ is the last step in $[a,\ell]$ at which $i$ is in the set $\cI$ of \textsc{Query}$s$, then the test of \textsc{Query}$s$, which uses the norm of $\wt x_i=\ov x_i^{(u)}=y$, gave $\|z-s_i^{(u+1)}\|_{y}^*\le\epsilon_s\le w_i^{1/2}\epsilon_s$. The same case distinction as for $x$ gives
\begin{align*}
\|\ov s_i^{(\ell+1)}-s_i^{(\ell+1)}\|_{\ov x_i^{(\ell+1)}}^*=\|z-s_i^{(\ell+1)}\|_{y}^*\le(2\log_2q+1)w_i^{1/2}\epsilon_s=w_i^{1/2}\epsilon_{\apx,s}.
\end{align*}

\textbf{Step 3: the number of updates and the running time.}
The bound on $\sum_\ell|L_x^{(\ell)}|$ holds in every run. Let $i\in L_x^{(\ell)}$, and let $a\le\ell$ be the smallest index such that $\ov x_i$ is not reset at any of the steps $a,\ldots,\ell-1$. Then $\ov x_i^{(t)}=x_i^{(a)}$ for all $a\le t\le\ell$. Since $\ov x_i$ is reset at step $\ell$, the condition in Line~\ref{line:alg:treewidth:approxds-queryx-checkxi} was true. By the triangle inequality, the Cauchy-Schwarz inequality and $w_i\ge 1$,
\begin{align*}
\epsilon_x^2<\|x_i^{(a)}-x_i^{(\ell+1)}\|_{\ov x_i^{(\ell)}}^2
\le(\ell-a+1)\sum_{t=a}^{\ell}\|x_i^{(t)}-x_i^{(t+1)}\|_{\ov x_i^{(t)}}^2
\le q\sum_{t=a}^{\ell}\|(H_{w,\ov x^{(t)}}^{1/2}x^{(t)}-H_{w,\ov x^{(t)}}^{1/2}x^{(t+1)})_i\|_2^2 .
\end{align*}
For a fixed block $i$, the intervals $[a,\ell]$ that belong to different resets of $\ov x_i$ are disjoint. Summing over all resets of all blocks gives
\begin{align*}
\epsilon_x^2\sum_{0\le\ell\le q-1}|L_x^{(\ell)}|\le q\sum_{0\le t\le q-1}\|H_{w,\ov x^{(t)}}^{1/2}x^{(t)}-H_{w,\ov x^{(t)}}^{1/2}x^{(t+1)}\|_2^2\le q^2\zeta_x^2,
\end{align*}
which is the first bound, because $\epsilon_x^{-2}=(2\log_2q+1)^2\epsilon_{\apx,x}^{-2}$.

For the other claims we consider an accurate run. Let $W=[\ell',\ell]$ be a dyadic window of length $2^j$. By Theorem~\ref{thm:treewidth:batchsketch}, the hypothesis on $\zeta_x$, and $\|\ov x^{(t)}-\ov x^{(t+1)}\|_{2,0}=|L_x^{(t)}|$, the set returned by $\mathsf{bs}.\textsc{Query}x(\ell',\epsilon_x)$ at step $\ell$ has size
\begin{align*}
O\Big(\epsilon_x^{-2}4^j\zeta_x^2+\sum_{t\in W}|L_x^{(t)}|\Big).
\end{align*}
We sum this over all dyadic windows inside $[0,q-1]$. There are at most $q/2^j$ windows of length $2^j$, and $\sum_{j\le\log_2q}(q/2^j)4^j\le 2q^2$. Every step lies in at most $\log_2q+1$ of the windows. So \textsc{Query}$x$ receives $O(\epsilon_x^{-2}\zeta_x^2q^2+\log q\cdot\sum_\ell|L_x^{(\ell)}|)=\wt O(\epsilon_{\apx,x}^{-2}\zeta_x^2q^2)$ candidates from $\mathsf{bs}$ in total. In the same way, \textsc{Query}$s$ receives $\wt O((\epsilon_{\apx,x}^{-2}\zeta_x^2+\epsilon_{\apx,s}^{-2}\zeta_s^2)q^2)$ candidates in total, where the first term again comes from the changes of $\ov x$. Every block of $L_s^{(\ell)}$ is such a candidate at step $\ell$ or lies in $L_x^{(\ell)}$. This proves the bound on $\sum_\ell|L_s^{(\ell)}|$.

For the running time, Theorem~\ref{thm:treewidth:batchsketch} bounds the cost of a call of $\mathsf{bs}.\textsc{Query}x$ or $\mathsf{bs}.\textsc{Query}s$ by $\wt O(\tau^2m)$ times one plus the above bound on the size of its output, and there are $O(q\log q)$ such calls. Every candidate and every block of $L_x^{(\ell)}$ costs $O(1)$ calls of $\mathsf{exact}.\textsc{Query}x$ or $\mathsf{exact}.\textsc{Query}s$, at the current time or at an earlier timestamp. Each of them costs $\wt O(\tau^2m)$ by Theorem~\ref{thm:treewidth:exactds} and Remark~\ref{rmk:snapshot}. Summing up gives the stated total time.
\end{proof}

\subsubsection{\textsc{BatchSketch}} \label{sec:treewidth:cpm-ds:batchsketch}
In this section we present the data structure \textsc{BatchSketch}. It maintains a sketch of $H_{\ov x}^{1/2} x$ and $H_{\ov x}^{-1/2} s$. It is a variation of \textsc{BatchSketch} in \cite{gs22}.

We recall the following definition from \cite{gs22}. 
\begin{definition}[Partition tree]
A partition tree $(\cS,\chi)$ of $\R^n$ is a constant degree rooted tree $\cS = (V,E)$ and a labeling of the vertices $\chi: V\to 2^{[n]}$, such that
\begin{itemize}
    \item $\chi(\text{root}) = [n]$;
    \item if $v$ is a leaf of $\cS$, then $|\chi(v)|=1$;
    \item for any non-leaf node $v\in V$, the set $\{\chi(c): \text{$c$ is a child of $v$}\}$ is a partition of $\chi(v)$.
\end{itemize}
\end{definition}

\begin{algorithm}[!ht]\caption{The \textsc{BatchSketch} data structure used by Algorithm~\ref{alg:treewidth:approxds-part1} and \ref{alg:treewidth:approxds-part2}.}\label{alg:treewidth:batchsketch-part1}
\begin{algorithmic}[1]
\State {\bf data structure} \textsc{BatchSketch} \Comment{Theorem~\ref{thm:treewidth:batchsketch}}
\State {\bf members}
\State \hspace{4mm} $\Phi\in \R^{r\times n_\tot} $ \Comment{All sketches need to share the same sketching matrix}
\State \hspace{4mm} $\cS,\chi$ partition tree
\State \hspace{4mm} $\ell \in \bN$ \Comment{Current timestamp}
\State \hspace{4mm} \textsc{BalancedSketch} $\mathsf{sketch}\cW^\top h$, $\mathsf{sketch}\cW^\top \wt h$, $\mathsf{sketch}\cW^\top \epsilon_x$, $\mathsf{sketch}\cW^\top \epsilon_s$
\Comment{Algorithm \ref{alg:treewidth:bsketch-part1}}
\State \hspace{4mm} \textsc{VectorSketch} $\mathsf{sketch}H_{w,\ov x}^{1/2} \wh x$, $\mathsf{sketch}H_{w,\ov x}^{-1/2} \wh s$, $\mathsf{sketch} c_s$
\Comment{Algorithm \ref{alg:treewidth:vsketch}}
\State \hspace{4mm} $\beta_x, \beta_s, \beta_{c_s} \in \R$, $\wt \beta_x, \wt \beta_s \in \R^m$
\State \hspace{4mm} $(\mathsf{history}[t])_{t\ge 0}$
\Comment{Snapshot of data at timestamp $t$. See Remark~\ref{rmk:snapshot}.}
\State {\bf end members}
\Procedure{Initialize}{$\ov x\in \R^{n_\tot}, h\in \R^{n_\tot}, \wt h \in \R^{n_\tot\times m}, \epsilon_x, \epsilon_s, H_{w,\ov x}^{1/2} \wh x, H_{w,\ov x}^{-1/2} \wh s, c_s \in \R^{n_\tot}, \beta_x, \beta_s, \beta_{c_s} \in \R, \wt \beta_x, \wt \beta_s \in \R^{m}, \delta_\apx \in \R$}
    \State Construct partition tree $(\cS, \chi)$ as in Definition~\ref{defn:partition_tree_construction}
    \State $r \gets \Theta(\log^3(n_\tot)\log(1/\delta_\apx))$
    \State Initialize $\Phi\in \R^{r \times n_\tot}$ with iid $\cN(0, \frac 1{r})$
    \State $\beta_x \gets \beta_x$, $\beta_s \gets \beta_s$, $\beta_{c_s} \gets \beta_{c_s}$, $\wt \beta_x \gets \wt \beta_x$, $\wt \beta_s \gets \wt \beta_s$
    \State $\mathsf{sketch}\cW^\top h.\textsc{Initialize}( {\cal S} , {\chi}, \Phi, \ov x, h)$ 
    \Comment{Algorithm~\ref{alg:treewidth:bsketch-part1}}
    \State $\mathsf{sketch}\cW^\top \wt h.\textsc{Initialize}( {\cal S} , {\chi}, \Phi, \ov x, \wt h)$ 
    \Comment{Algorithm~\ref{alg:treewidth:bsketch-part1}}
    \State $\mathsf{sketch}\cW^\top \epsilon_x.\textsc{Initialize}( {\cal S} , {\chi}, \Phi, \ov x, \epsilon_x )$ 
    \Comment{Algorithm~\ref{alg:treewidth:bsketch-part1}}
    \State $\mathsf{sketch}\cW^\top \epsilon_s.\textsc{Initialize}( {\cal S} , {\chi}, \Phi, \ov x, \epsilon_s )$ 
    \Comment{Algorithm~\ref{alg:treewidth:bsketch-part1}}
    \State $\mathsf{sketch}H_{w,\ov x}^{1/2} \wh x.\textsc{Initialize}( {\cal S} , {\chi}, \Phi, H_{w,\ov x}^{1/2} \wh{x} )$
    \Comment{Algorithm~\ref{alg:treewidth:vsketch}}
    \State $\mathsf{sketch}H_{w,\ov x}^{-1/2} \wh s.\textsc{Initialize}( {\cal S} , {\chi}, \Phi, H_{w,\ov x}^{-1/2} \wh{s} )$
    \Comment{Algorithm~\ref{alg:treewidth:vsketch}}
    \State $\mathsf{sketch} c_s.\textsc{Initialize}( {\cal S} , {\chi}, \Phi, c_s )$
    \Comment{Algorithm~\ref{alg:treewidth:vsketch}}
    \State $\ell \gets 0$. Make snapshot $\mathsf{history}[\ell]$ \Comment{Remark~\ref{rmk:snapshot}}
\EndProcedure
\State {\bf end data structure}
\end{algorithmic}
\end{algorithm}

\begin{algorithm}[!ht]\caption{\textsc{BatchSketch} Algorithm~\ref{alg:treewidth:batchsketch-part1} continued.}\label{alg:treewidth:batchsketch-part2}
\begin{algorithmic}[1]
\State {\bf data structure} \textsc{BatchSketch} \Comment{Theorem~\ref{thm:treewidth:batchsketch}}
\Procedure{Move}{$\beta_x, \beta_s, \beta_{c_s}\in \R$, $\wt \beta_x, \wt \beta_s \in \R^m$}
    \State $\beta_x \gets \beta_x$, $\beta_s \gets \beta_s$, $\beta_{c_s} \gets \beta_{c_s}$, $\wt \beta_x \gets \wt \beta_x$, $\wt \beta_s \gets \wt \beta_s$
    \Comment{Do not update $\ell$ yet}
\EndProcedure
\Procedure{Update}{$\delta_{\ov x} \in \R^{n_\tot}, \delta_{h}\in \R^{n_\tot}, \delta_{\wt h} \in \R^{n_\tot \times m}, \delta_{\epsilon_x}, \delta_{\epsilon_s}, \delta_{H_{w,\ov x}^{1/2} \wh{x}}, \delta_{ H_{w,\ov x}^{-1/2} \wh{s} } , \delta_{c_s} \in \R^{n_\tot}$}
    \State $\mathsf{sketch}\cW^\top h.\textsc{Update}(\delta_{\ov x}, \delta_{h})$
    \Comment{Algorithm~\ref{alg:treewidth:bsketch-part2}}
    \State $\mathsf{sketch}\cW^\top \wt h.\textsc{Update}(\delta_{\ov x}, \delta_{\wt h})$
    \Comment{Algorithm~\ref{alg:treewidth:bsketch-part2}}
    \State $\mathsf{sketch}\cW^\top \epsilon_x.\textsc{Update}(\delta_{\ov x}, \delta_{\epsilon_x})$
    \Comment{Algorithm~\ref{alg:treewidth:bsketch-part2}}
    \State $\mathsf{sketch}\cW^\top \epsilon_s.\textsc{Update}(\delta_{\ov x}, \delta_{\epsilon_s})$
    \Comment{Algorithm~\ref{alg:treewidth:bsketch-part2}}
    \State $\mathsf{sketch}H_{w,\ov x}^{1/2} \wh x.\textsc{Update}( \delta_{ H_{w,\ov x}^{1/2} \wh{x} } )$
    \Comment{Algorithm~\ref{alg:treewidth:vsketch}}
    \State $\mathsf{sketch}H_{w,\ov x}^{-1/2} \wh s.\textsc{Update}(  \delta_{ H_{w,\ov x}^{-1/2} \wh{s}} )$ 
    \Comment{Algorithm~\ref{alg:treewidth:vsketch}}
    \State $\mathsf{sketch}c_s.\textsc{Update}( \delta_{ c_s } )$
    \Comment{Algorithm~\ref{alg:treewidth:vsketch}}
    \State $\ell \gets \ell+1$
    \State Make snapshot $\mathsf{history}[\ell]$ \Comment{Remark~\ref{rmk:snapshot}}
\EndProcedure
\State {\bf end data structure}
\end{algorithmic}
\end{algorithm}

\begin{algorithm}[!ht]\caption{\textsc{BatchSketch} Algorithm~\ref{alg:treewidth:batchsketch-part1},~\ref{alg:treewidth:batchsketch-part2} continued.}\label{alg:treewidth:batchsketch-part3}
\begin{algorithmic}[1]
\State {\bf data structure} \textsc{BatchSketch} \Comment{Theorem~\ref{thm:treewidth:batchsketch}}
\State {\bf private:}
    \Procedure{Query$x$Sketch}{$v\in \cS$}
        \Comment{Return the value of $\Phi_{\chi(v)} (H_{w,\ov x}^{1/2} x)_{\chi(v)}$}
        \State \Return $\mathsf{sketch} H_{w,\ov x}^{1/2} \wh x.\textsc{Query}(v) + \mathsf{sketch}\cW^\top h.\textsc{Query}(v) \cdot \beta_x - \mathsf{sketch}\cW^\top \wt h.\textsc{Query}(v) \cdot \wt \beta_x + \mathsf{sketch}\cW^\top \epsilon_x.\textsc{Query}(v)$
        \Comment{Algorithm~\ref{alg:treewidth:vsketch}, \ref{alg:treewidth:bsketch-part1}}
    \EndProcedure
    \Procedure{Query$s$Sketch}{$v\in \cS$}
        \Comment{Return the value of $\Phi_{\chi(v)} (H_{w,\ov x}^{-1/2} s)_{\chi(v)}$}
        \State \Return $\mathsf{sketch} H_{w,\ov x}^{-1/2} \wh s.\textsc{Query}(v) + \mathsf{sketch}c_s.\textsc{Query}(v) \cdot \beta_{c_s}  - \mathsf{sketch}\cW^\top h.\textsc{Query}(v) \cdot \beta_s + \mathsf{sketch}\cW^\top \wt h.\textsc{Query}(v) \cdot \wt \beta_s - \mathsf{sketch}\cW^\top \epsilon_s.\textsc{Query}(v)$
        \Comment{Algorithm~\ref{alg:treewidth:vsketch}, \ref{alg:treewidth:bsketch-part1}}
    \EndProcedure
\State {\bf public:}
\Procedure{Query$x$}{$\ell' \in \bN, \epsilon \in \R$}
    \State $L_0=\{\text{root}(\cS)\}$
    \State $S \gets \emptyset$
    \For{$d = 0 \to \infty$}
        \If{$L_d=\emptyset$}
        \State \Return $S$
        \EndIf
        \State $L_{d+1} \gets \emptyset$
        \For{$v \in L_d$}
            \If{$v$ is a leaf node}
            \State $S \gets S \cup \{v\}$
            \Else
            \For{$u$ child of $v$}
            \If{$\|\textsc{Query$x$Sketch}(u) - \mathsf{history}[\ell'].\textsc{Query$x$Sketch}(u)\|_2 > 0.9\epsilon$} \label{line:treewidth:batchsketch-queryx}
                \State $L_{d+1} \gets L_{d+1} \cup \{ u \}$
            \EndIf
            \EndFor
            \EndIf
        \EndFor 
    \EndFor 
\EndProcedure
\Procedure{Query$s$}{$\ell' \in \bN, \epsilon \in \R$}
    \State Same as \textsc{Query}$x$, except for replacing $\textsc{Query$x$Sketch}$ in Line~\ref{line:treewidth:batchsketch-queryx} with $\textsc{Query$s$Sketch}$.
\EndProcedure
\State {\bf end structure}
\end{algorithmic}
\end{algorithm}

\begin{theorem}\label{thm:treewidth:batchsketch}
Data structure \textsc{BatchSketch} (Algorithm~\ref{alg:treewidth:batchsketch-part1}, \ref{alg:treewidth:batchsketch-part3}) supports the following operations:
\begin{itemize}
\item $\textsc{Initialize}(\ov x\in \R^{n_\tot}, h\in \R^{n_\tot}, \wt h \in \R^{n_\tot\times m}, \epsilon_x, \epsilon_s, H_{w,\ov x}^{1/2} \wh x, H_{w,\ov x}^{-1/2} \wh s, c_s \in \R^{n_\tot}, \beta_x, \beta_s, \beta_{c_s} \in \R, \wt \beta_x, \wt \beta_s \in \R^{m}, \delta_\apx \in \R)$: Initialize the data structure in $\wt O(n\tau^{\omega-1} + n \tau m)$ time.

\item $\textsc{Move}(\beta_x, \beta_s, \beta_{c_s}\in \R, \wt \beta_x, \wt \beta_s \in \R^m)$: Update values of $\beta_x, \beta_s, \beta_{c_s}, \wt \beta_x, \wt \beta_s$ in $O(m)$ time. This effectively moves $(x^{(\ell)}, s^{(\ell)})$ to $(x^{(\ell+1)}, s^{(\ell+1)})$ while keeping $\ov x^{(\ell)}$ unchanged.

\item $\textsc{Update}(\delta_{\ov x} \in \R^{n_\tot}, \delta_{h}\in \R^{n_\tot}, \delta_{\wt h} \in \R^{n_\tot \times m}, \delta_{\epsilon_x}, \delta_{\epsilon_s}, \delta_{H_{w,\ov x}^{1/2} \wh{x}}, \delta_{ H_{w,\ov x}^{-1/2} \wh{s} } , \delta_{c_s} \in \R^{n_\tot})$:
Update sketches of $H_{w,\ov x^{(\ell)}}^{1/2} x^{(\ell+1)}$ and $H_{w,\ov x^{(\ell)}}^{-1/2} s^{(\ell+1)}$.
This effectively moves $\ov x^{(\ell)}$ to $\ov x^{(\ell+1)}$ while keeping $(x^{(\ell+1)}, s^{(\ell+1)})$ unchanged.
Then advance timestamp $\ell$.

Each update costs
\begin{align*}
\wt O(\tau^2 m \|\delta_{\ov x}\|_0 + \nnz(\delta_h) + \nnz(\delta_{\wt h}) + \|\delta_{\epsilon_x}\|_0 + \|\delta_{\epsilon_s}\|_0 + \|\delta_{H_{w,\ov x}^{1/2}\wh x}\|_0 + \|\delta_{H_{w,\ov x}^{-1/2}\wh s}\|_0 + \|\delta_{c_s}\|_0)
\end{align*}
time.

\item $\textsc{Query}x(\ell' \in \bN, \epsilon \in\R)$:
Given timestamp $\ell'$, return a set $S\subseteq [n]$. Let $\ell$ be the current timestamp and let
\begin{align*}
    v := H_{w,\ov x^{(\ell')}}^{1/2} x^{(\ell')} - H_{w,\ov x^{(\ell)}}^{1/2} x^{(\ell+1)} \in \R^{n_\tot}.
\end{align*}
We say that $\Phi$ is \emph{accurate} for a vector $v\in \R^{n_\tot}$ if $0.99\|v_{\chi(u)}\|_2 \le \|\Phi_{\chi(u)} v_{\chi(u)}\|_2 \le 1.01 \|v_{\chi(u)}\|_2$ for every node $u$ of the partition tree. If $\Phi$ is accurate for $v$, then
\begin{align*}
    S &\supseteq \{i \in [n]: \|v_i\|_2 \ge \epsilon\},
\end{align*}
and
\begin{align*}
    |S| & = O(\epsilon^{-2} (\ell-\ell'+1) \sum_{\ell' \le t \le \ell} \|H_{w,\ov x^{(t)}}^{1/2} x^{(t)} - H_{w,\ov x^{(t)}}^{1/2} x^{(t+1)}\|_2^2 + \sum_{\ell ' \le t \le \ell-1} \|\ov x^{(t)} - \ov x^{(t+1)}\|_{2,0}),
\end{align*}
and the query costs at most
\begin{align*}
\wt O\Big(\tau^2 m \cdot \Big(1 & +\epsilon^{-2} (\ell-\ell'+1) \sum_{\ell' \le t \le \ell} \|H_{w,\ov x^{(t)}}^{1/2} x^{(t)} - H_{w,\ov x^{(t)}}^{1/2} x^{(t+1)}\|_2^2 \\
& + \sum_{\ell ' \le t \le \ell-1} \|\ov x^{(t)} - \ov x^{(t+1)}\|_{2,0}\Big)\Big)
\end{align*}
running time.

If $v$ does not depend on $\Phi$, then $\Phi$ is accurate for $v$ with probability at least $1-\delta_\apx$, where $\delta_\apx\le\frac12$ is the parameter of \textsc{Initialize}. This statement does not cover a vector $v$ that depends on $\Phi$, for example through earlier answers of the data structure. Theorem~\ref{thm:treewidth:approxds} treats this case.

\item $\textsc{Query}s(\ell' \in \bN, \epsilon \in\R)$:
Given timestamp $\ell'$, return a set $S\subseteq [n]$. Let $\ell$ be the current timestamp and let
\begin{align*}
    v := H_{w,\ov x^{(\ell')}}^{-1/2} s^{(\ell')} - H_{w,\ov x^{(\ell)}}^{-1/2} s^{(\ell+1)} \in \R^{n_\tot}.
\end{align*}
If $\Phi$ is accurate for $v$, then
\begin{align*}
    S &\supseteq \{i \in [n]: \|v_i\|_2 \ge \epsilon\}
\end{align*}
and
\begin{align*}
    |S| & = O(\epsilon^{-2}  (\ell-\ell'+1)\sum_{\ell' \le t \le \ell} \|H_{w,\ov x^{(t)}}^{-1/2} s^{(t)} - H_{w,\ov x^{(t)}}^{-1/2} s^{(t+1)}\|_2^2 + \sum_{\ell' \le t \le \ell-1} \|\ov x^{(t)} - \ov x^{(t+1)}\|_{2,0}),
\end{align*}
and the query costs at most
\begin{align*}
\wt O\Big(\tau^2 m \cdot \Big(1 & +\epsilon^{-2} (\ell-\ell'+1) \sum_{\ell' \le t \le \ell} \|H_{w,\ov x^{(t)}}^{-1/2} s^{(t)} - H_{w,\ov x^{(t)}}^{-1/2} s^{(t+1)}\|_2^2 \\
& + \sum_{\ell ' \le t \le \ell-1} \|\ov x^{(t)} - \ov x^{(t+1)}\|_{2,0}\Big)\Big)
\end{align*}
running time.
If $v$ does not depend on $\Phi$, then $\Phi$ is accurate for $v$ with probability at least $1-\delta_\apx$.
\end{itemize}
\end{theorem}

\begin{proof}
The proof is essentially the same as proof of \cite[Theorem 4.21]{gs22}.
That proof derives the claims on $S$ and on the running time from the accuracy of $\Phi$ for $v$. It obtains the accuracy from the Johnson-Lindenstrauss lemma (Lemma~\ref{lem:jl_matrix}) and the choice of $r$, which applies to a vector $v$ that does not depend on $\Phi$.
For the running time claims, we plug in Lemma~\ref{lem:treewidth:vsketch} and~\ref{lem:treewidth:bsketch} when necessary. A call of \textsc{Query$x$Sketch} or \textsc{Query$s$Sketch} costs $\wt O(\tau^2 m)$, because of the sketch of $\cW^\top \wt h$ with $\wt h\in\R^{n_\tot\times m}$.
\end{proof}

\begin{remark}[Snapshot]\label{rmk:snapshot}
As in previous works, we use persistent data structures (e.g., \cite{driscoll1989making}) to keep a snapshot of the data structure after every update. This allows us to support query to historical data.
This incurs an $O(\log n_\tot)=\wt O(1)$ multiplicative factor in all running times, which we ignore in our analysis.
We keep \textsc{ExactDS} persistent in the same way, with one snapshot after \textsc{Initialize} (timestamp $0$) and one after the $\ell'$-th call of \textsc{Update} since the last \textsc{Initialize} (timestamp $\ell'$). So \textsc{ApproxDS} can call $\mathsf{exact}.\textsc{Query}x(i)$ and $\mathsf{exact}.\textsc{Query}s(i)$ at an earlier timestamp $\ell'$ of the current restart period. These calls return $x_i^{(\ell')}$ and $s_i^{(\ell')}$, and the snapshot also stores $\ov x_i^{(\ell')}$. Their cost is the cost of \textsc{Query}$x$ and \textsc{Query}$s$, up to the same logarithmic factor.
\end{remark}

\subsubsection{\textsc{VectorSketch}} \label{sec:treewidth:cpm-ds:vsketch}
\textsc{VectorSketch} is a data structure used to maintain sketches of sparsely-changing vectors.
It is a direct application of segment trees. For completeness, we include code (Algorithm~\ref{alg:treewidth:vsketch}) from \cite[Algorithm 10]{gs22}.

\begin{algorithm}[!ht]\caption{ \cite[Algorithm 10]{gs22}. Used in Algorithm~\ref{alg:treewidth:batchsketch-part1},~\ref{alg:treewidth:batchsketch-part2},~\ref{alg:treewidth:batchsketch-part3}.
}\label{alg:treewidth:vsketch}
\begin{algorithmic}[1]
\State {\bf data structure} \textsc{VectorSketch} \Comment{Lemma~\ref{lem:treewidth:vsketch}}
    \State {\bf private: members}
    \State \hspace{4mm} $\Phi \in \R^{r \times n_\tot}$
    \State \hspace{4mm} Partition tree $( {\cal S}, \chi )$
    \State \hspace{4mm} $x \in \R^{n_\tot}$
    \State \hspace{4mm} Segment tree $\cT$ on $[n]$ with values in $\R^r$
    \State {\bf end members}
    \Procedure{Initialize}{${\cal S}, \chi: \text{partition tree}, \Phi \in \R^{r \times n_{\tot}}, x \in \R^{n_{\tot}}$}
        \State $({\cal S}, \chi) \gets ( {\cal S} , \chi)$, $\Phi \gets \Phi$
        \State $x \gets x$
        \State Order leaves of $\cS$ (variable blocks) such that every node $\chi(v)$ corresponds to a contiguous interval $\subseteq [n]$.
        \State Build a segment tree $\cT$ on $[n]$ such that each segment tree interval $I\subseteq [n]$ maintains $\Phi_I x_I \in \R^r$.
    \EndProcedure
    
    \Procedure{Update}{$\delta_{x} \in \R^{n_{\tot}}$}
        \For{all $i\in [n_\tot]$ such that $\delta_{ x,i} \ne 0$}
            \State Let $j\in [n]$ be such that $i$ is in $j$-th block
            \State Update $\cT$ at $j$-th coordinate $\Phi_j x_j \gets \Phi_j x_j + \Phi_i \cdot \delta_{ x,i}$.
            \State $x_i \gets x_i + \delta_{ x,i}$
        \EndFor
    \EndProcedure
    \Procedure{Query}{$v \in V(\cS)$}
        \State Find interval $I$ corresponding to $\chi(v)$
        \State \Return range sum of $\cT$ on interval $I$
    \EndProcedure
\State {\bf end data structure}
\end{algorithmic}
\end{algorithm}

\begin{lemma}[{\cite[Lemma 4.23]{gs22}}]\label{lem:treewidth:vsketch}
Given a partition tree $(\cS,\chi)$ of $\R^{n}$, and a JL sketching matrix $\Phi \in \R^{r \times n_\tot}$, the data structure \textsc{VectorSketch} (Algorithm~\ref{alg:treewidth:vsketch}) maintains $\Phi_{\chi(v)} x_{\chi(v)}$ for all nodes $v$ in the partition tree implicitly through the following functions:
\begin{itemize}
    \item \textsc{Initialize}$(\cS, \chi,\Phi)$: Initializes the data structure in $O(r n_\tot)$ time.
    \item \textsc{Update}$(\delta_{x} \in \R^{n_\tot})$: Maintains the data structure for $x\leftarrow x + \delta_{x}$ in $O(r \|\delta_{x}\|_0\log n)$ time.
    \item \textsc{Query}$(v\in V(\cS))$: Outputs $\Phi_{\chi(v)}x_{\chi(v)}$ in $O(r \log n)$ time.
\end{itemize}
\end{lemma}

\subsubsection{\textsc{BalancedSketch}} \label{sec:treewidth:cpm-ds:bsketch}
In this section, we present data structure \textsc{BalancedSketch}. It is a data structure for maintaining a sketch of a vector of form $\cW^\top h$, where $\cW = L_{\ov x}^{-1} H_{w,\ov x}^{1/2}$ and $h\in \R^{n_\tot}$ is a sparsely-changing vector. This is a variation of \textsc{BlockBalancedSketch} in \cite{gs22}.

We use the following construction of a partition tree.
\begin{definition}[Construction of Partition Tree] \label{defn:partition_tree_construction}
We fix an ordering $\pi$ of $[n]$ using the heavy-light decomposition (Lemma~\ref{lem:heavy_light}).
Let $\cS$ be a complete binary tree with leaf set $[n]$ and ordering $\pi$.
Let $\chi$ map a node to the set of leaves in its subtree.
Then $(\cS,\chi)$ is a valid partition tree.
\end{definition}

\begin{algorithm}[!ht]\caption{The \textsc{BalancedSketch} data structure is used in Algorithm~\ref{alg:treewidth:batchsketch-part1},~\ref{alg:treewidth:batchsketch-part2},~\ref{alg:treewidth:batchsketch-part3}.
}  \label{alg:treewidth:bsketch-part1}
\begin{algorithmic}[1]
\State {\bf data structure} \textsc{BalancedSketch} \Comment{Lemma~\ref{lem:treewidth:bsketch}}
\State {\bf private: members}
    \State \hspace{4mm} $\Phi \in \R^{r \times n_\tot}$
    \State \hspace{4mm} Partition tree $(\cS, \chi)$ with balanced binary tree ${\cal B}$
    \State \hspace{4mm} $t \in \mathbb{N}$
    \State \hspace{4mm} $h \in \R^{n_\tot}$, $\ov x \in \R^{n_\tot}$, $H_{w,\ov x} \in \R^{n_\tot \times n_\tot}$
    \State \hspace{4mm} $\{ L[t] \in \R^{n_\tot \times n_\tot} \}_{t\geq 0}$
    \State \hspace{4mm} $\{ J_v \in \R^{r \times n_\tot} \}_{v\in \cS}$
    \State \hspace{4mm} $\{ Z_v \in \R^{r \times n_\tot} \}_{v\in \cB}$
    \State \hspace{4mm} $\{ y_v^{\triangledown} \in \R^r \}_{v\in \cB}$ 
    \State \hspace{4mm} $\{ t_v \in \mathbb{N} \}_{v\in \cB}$
\State {\bf end members}
\Procedure{Initialize}{${\cal S}, \chi: \text{partition tree}, \Phi \in \R^{r \times n_{\tot}}, \ov x \in \R^{n_{\tot}}, h \in \R^{n_\tot\times k}$} 
    \State $({\cal S}, \chi ) \gets ( {\cal S}, \chi )$, $\Phi \gets \Phi$
    \State $t \gets 0$, $h \gets h$
    \State $H_{w,\ov x} \gets \nabla^2 \phi( \ov x )$, $B_{\ov x} \gets Q + \ov t H_{w,\ov x}$
    \State Compute lower Cholesky factor $L_{\ov x}[t]$ of $B_{\ov x}$
    \For{all $v \in {\cal S}$}
        \State $J_v \gets \Phi_{ \chi(v) } H_{w,\ov x}^{1/2} $
    \EndFor 
    \For{all $v \in {\cal B}$}
        \State $Z_v \gets J_v L_{\ov x}[t]^{-\top}$
        \State $y_v^{\triangledown} \gets Z_v ( I - I_{\Lambda(v)} ) h$
        \State $t_v \gets t$
    \EndFor 
\EndProcedure
\Procedure{Query}{$v \in {\cal S}$} 
    \If{$v \in {\cal S} \backslash {\cal B}$}
        \State \Return $J_v \cdot L_{\ov x}[t]^{-\top} h$
    \EndIf
    \State $\Delta_{L_{\ov x}} \gets ( L_{\ov x}[t] - L_{\ov x}[t_v] ) \cdot I_{\Lambda(v)} $
    \State $\delta_{Z_v} \gets -( L_{\ov x}[t]^{-1} \cdot \Delta_{L_{\ov x}} \cdot Z_v^\top )^\top$
    \State $Z_v \gets Z_v + \delta_{Z_v}$
    \State $\delta_{y_v^{\triangledown}} \gets \delta_{Z_v} \cdot (I - I_{\Lambda(v)}) h$
    \State $y_v^{\triangledown} \gets y_v^{\triangledown} + \delta_{y_v^{\triangledown}}$
    \State $t_v \gets t$
    \State $y_v^{\triangleup} \gets Z_v \cdot I_{\Lambda(v)} \cdot h$
    \State \Return $y_v^{\triangleup} + y_v^{\triangledown}$
\EndProcedure 
\State {\bf end data structure}
\end{algorithmic}
\end{algorithm}

\begin{algorithm}[!ht]\caption{
\textsc{BalancedSketch} Algorithm~\ref{alg:treewidth:bsketch-part1} continued.
This is used in Algorithm~\ref{alg:treewidth:batchsketch-part1},~\ref{alg:treewidth:batchsketch-part2},~\ref{alg:treewidth:batchsketch-part3}.
} \label{alg:treewidth:bsketch-part2}
\begin{algorithmic}[1]
\State {\bf data structure} \textsc{BalancedSketch} 
\Procedure{Update}{$ \delta_{\ov x} \in \R^{n_{\tot}} , \delta_{h} \in \R^{n_\tot\times k}$}
    \For{$i \in [n]$ where $\delta_{\ov x,i} \ne 0$}
        \State \textsc{Update$\ov x$}$(\delta_{\ov x,i})$
    \EndFor 
    \For{all $\delta_{h,i} \ne 0$} \label{line:alg:treewidth:bsketch-part2_updateh_start}
        \State $v \gets \Lambda^{\circ}(i)$
        \For{all $u \in {\cal P}^{\cal B}(v)$}
            \State $y_u^{\triangledown} \gets y_u^\triangledown + Z_u \cdot I_{ \{ i\} } \cdot \delta_{h}$
        \EndFor
    \EndFor 
    \label{line:alg:treewidth:bsketch-part2_updateh_end}
    \State $h \gets h + \delta_{h}$
\EndProcedure 
\Procedure{Update$\ov x$}{$ \delta_{\ov x,i} \in \R^{n_i} $}   
    \State $t \gets t+1$ 
    \State $\ov x_i \gets \ov x_i + \delta_{\ov x,i}$
    \State $\Delta_{H_{w,\ov x}, (i,i)} \gets \nabla^2 \phi_i(\ov x_i) - H_{w,\ov x,(i,i)}$
    \State Compute $\Delta_{L_{\ov x}}$ such that $L_{\ov x}[t] \gets L_{\ov x}[t-1] + \Delta_{L_{\ov x}}$ is the lower Cholesky factor of $Q + \ov t (H_{w,\ov x} + \Delta_{ H_{w,\ov x}})$
    \State $S \gets {\cal P}^{\cal B} ( \Lambda^{\circ} ( \mathsf{low}^\cT(i) ) ) $
    \State $\textsc{Update$L$}(S, \Delta_{L_{\ov x}}) $
    \State $\textsc{Update$H$}(i, \Delta_{H_{w,\ov x}, (i,i)})$
\EndProcedure 
\State {\bf end data structure}
\end{algorithmic} 
\end{algorithm} 

\begin{algorithm}[!ht]\caption{
\textsc{BalancedSketch} Algorithm~\ref{alg:treewidth:bsketch-part1}, \ref{alg:treewidth:bsketch-part2} continued.
This is used in Algorithm~\ref{alg:treewidth:batchsketch-part1},~\ref{alg:treewidth:batchsketch-part2},~\ref{alg:treewidth:batchsketch-part3}.
}\label{alg:treewidth:bsketch-part3} 
\begin{algorithmic}[1]
\State {\bf data structure} \textsc{BalancedSketch} \Comment{Lemma~\ref{lem:treewidth:bsketch}}
\State {\bf private:}
\Procedure{Update$L$}{$S \subset {\cal B}$, $\Delta_{L_{\ov x}} \in \R^{n_\tot \times n_\tot}$} 
    \For{all $v \in S$}
        \State $\delta_{Z_v} \gets - ( L_{\ov x}[t-1]^{-1} ( L_{\ov x}[t-1] - L_{\ov x}[t_v] ) \cdot I_{\Lambda(v)} \cdot Z_v^\top )^\top$
        \label{line:alg:treewidth:bsketch-part3_updateZv1}
        \State $\delta_{Z_v}'\gets - ( L_{\ov x}[t]^{-1} \cdot \Delta_{L_{\ov x}} \cdot (Z_v+\delta_{Z_v})^\top )^\top$
        \label{line:alg:treewidth:bsketch-part3_updateZv2}
        \State $Z_v \gets Z_v +\delta_{Z_v}+\delta_{Z_v}'$
        \State $\delta_{y_v^{\triangledown}} \gets (\delta_{Z_v}+\delta_{Z_v}')(I - I_{\Lambda(v)}) h$
        \State $y_v^{\triangledown} \gets y_v^{\triangledown} + \delta_{y_v^{\triangledown}}$
        \State $t_v \gets t$
    \EndFor
\EndProcedure 
\State {\bf private:}
\Procedure{Update$H$}{$i\in [n], \Delta_{H_{w,\ov x}, (i,i)} \in \R^{n_i\times n_i}$} 
    \State Find $u$ such that $\chi(u) = \{i\}$
    \State $\Delta_{H_{w,\ov x}^{1/2}, (i,i)} \gets (H_{w,\ov x,(i,i)} + \Delta_{H_{w,\ov x},(i,i)})^{1/2} - H_{w,\ov x,(i,i)}^{1/2}$
    \State $\delta_{J_u} \gets \Phi_{i} \cdot \Delta_{H_{w,\ov x}^{1/2}, (i,i)}$
    \For{all $v \in {\cal P}^{\cS} (u)$}
        \State $J_v \gets J_v + \delta_{J_u} $
        \If{$v \in {\cal B}$}
            \State $\delta_{Z_v} \gets \delta_{J_v} \cdot L_{\ov x}[t_v]^{-\top}$
            \State $Z_v \gets Z_v + \delta_{Z_v}$
            \State $\delta_{y_v^{\triangledown}} \gets \delta_{Z_v} \cdot (I - I_{\Lambda(v)}) \cdot h$
            \State $y_v^{\triangledown} \gets y_v^{\triangledown} + \delta_{y_v^{\triangledown}}$
        \EndIf 
    \EndFor 
    \State $H_{w,\ov x} \gets H_{w,\ov x} + \Delta_{ H_{w,\ov x}, (i,i)}$
\EndProcedure
\State {\bf end data structure}
\end{algorithmic}
\end{algorithm}

\begin{lemma}
\label{lem:treewidth:bsketch}
Given an elimination tree $\cT$ with height $\eta$, a JL matrix $\Phi \in \R^{r \times n_\tot}$, and a partition tree $(\cS, \chi)$ constructed as in Definition~\ref{defn:partition_tree_construction}
with height $\wt O(1)$, the data structure \textsc{BalancedSketch} (Algorithm \ref{alg:treewidth:bsketch-part1}, \ref{alg:treewidth:bsketch-part2}, \ref{alg:treewidth:bsketch-part3}), maintains $\Phi_{\chi(v)} (\cW^\top h)_{\chi(v)}$ for each $v\in V(\cS)$ through the following operations
\begin{itemize}
    \item \textsc{Initialize}$((\cS, \chi): \text{partition tree}, \Phi\in \R^{n_\tot}, \ov x\in \R^{n_\tot}, h\in \R^{n_\tot\times k})$: Initializes the data structure in $\wt O(r(n\tau^{\omega-1} + n \tau k))$ time.
    \item \textsc{Update}$(\delta_{\ov x}\in \R^{n_\tot}, \delta_h \in \R^{n_\tot\times k})$: Updates all sketches in $\cS$ implicitly to reflect $(\cW, h)$ updating to $(\cW^\new, h^\new)$ in $\wt O(r (\tau^2 k\|\delta_{\ov x}\|_0 + \nnz(\delta_h)))$ time, that is, $\wt O(r\tau^2 k)$ per changed block of $\ov x$ and $\wt O(r)$ per changed entry of $h$.
    \item \textsc{Query}$(v\in \cS)$: Outputs $\Phi_{\chi(v)} (\cW^\top h)_{\chi(v)}$ in $\wt O(r \tau^2 k)$ time.
\end{itemize}
\end{lemma}
\begin{proof}
The proof is almost same as the proof of \cite[Lemma 4.24]{gs22}.
(In fact, our $\cW$ is simpler than the one used in \cite{gs22}.)

For \textsc{Initialize} running time, we note that computing $Z_v$ for all $v\in \cB$ takes $\wt O(r n \tau^{\omega-1})$ time by \cite[Lemma 8.3]{gs22}.
Because $Z_v$ is supported on the path from $v$ to the root in $\cT$, we know that $\nnz(Z) = O(r n \tau)$.
Therefore computing $y_v^\triangledown$ for all $v\in \cB$ takes $\wt O(r n \tau k)$ time.

Remaining claims follow from combining proof of \cite[Lemma 4.24]{gs22} and \cite[Lemma 8.3]{gs22}.
\end{proof}

\subsection{Analysis of \textsc{CentralPathMaintenance}} \label{sec:treewidth:cpm-analysis}
\begin{lemma}[Correctness of \textsc{CentralPathMaintenance}] \label{lem:treewidth:cpm-correct}
Algorithm~\ref{alg:treewidth:cpm} implicitly maintains the primal-dual solution pair $(x,s)$ via representation Eq.~\eqref{eqn:thm:treewidth:exactds:x-rep}\eqref{eqn:thm:treewidth:exactds:s-rep}.
It also explicitly maintains $(\ov x, \ov s)\in \R^{n_\tot} \times \R^{n_\tot}$ such that $\|\ov x_i-x_i\|_{\ov x_i} \le \ov \epsilon$ and $\|\ov s_i-s_i\|_{\ov x_i}^* \le t \ov \epsilon w_i$ for all $i\in [n]$ with probability at least $0.9$.
\end{lemma}
\begin{proof}
    We correctly maintain the implicit representation because of correctness of $\mathsf{exact}$.\textsc{Update} (Theorem~\ref{thm:treewidth:exactds}).

    We show that $\|\ov x_i-x_i\|_{\ov x_i} \le \ov \epsilon$ and $\|\ov s_i - s_i\|_{\ov x_i}^* \le t \ov \epsilon w_i$ for all $i\in [n]$ (c.f.~Algorithm~\ref{alg:centering}, Line~\ref{line:alg:centering:cpm-guarantee}).
    We first check the hypotheses of Theorem~\ref{thm:treewidth:approxds}.
    In every step we have
    \begin{align*}
    \|H_{w,\ov x}^{1/2} x^{(\ell+1)} - H_{w,\ov x}^{1/2} x^{(\ell)}\|_2 = \|\delta_x\|_{w,\ov x} \le \frac 98 \alpha \le \zeta_x
    \end{align*}
    where the first step from definition of $\|\cdot\|_{w,\ov x}$, the second step follows from Lemma~\ref{lem:DLY_Lemma_A.9}, the third step follows from definition of $\zeta_x$.
    Similarly,
    \begin{align*}
    \|H_{w,\ov x}^{-1/2} \delta_s\|_2 = \|\delta_s\|_{w,\ov x}^* \le \frac {17}8 \alpha \cdot t \le \zeta_s,
    \end{align*}
    where the last step uses $\zeta_s=3\alpha\ov t$ and $t\le(1+\epsilon_t)\ov t$.
    Lemma~\ref{lem:DLY_Lemma_A.9} only uses that the step is computed from $(\ov x,\ov s,\ov t)$ and that $|\ov t-t|\le\epsilon_t\ov t$, so these bounds hold whatever the outputs of $\mathsf{approx}$ are.
    Moreover, between two restarts the inputs of $\mathsf{approx}$ are determined by the state at the restart and by the earlier outputs of $\mathsf{approx}$, because $\mathsf{exact}$.\textsc{Move} and $\mathsf{exact}$.\textsc{Update} are deterministic and so is the sequence of values of $t$.

    A restart draws a new sketching matrix, which is independent of the state at the restart.
    Because \textsc{MultiplyAndMove} restarts whenever $\ell\ge q$, a period between two restarts has at most $q$ steps, which are among the steps $0,\ldots,q-1$ of Theorem~\ref{thm:treewidth:approxds}.
    So we can apply Theorem~\ref{thm:treewidth:approxds} to each period between two restarts, conditioned on the state at its start. With probability at least $1-\delta_{\apx}$, its guarantees hold throughout the period. There are at most $N$ periods, and $\delta_\apx=\frac1{20N}$. By a union bound, with probability at least $1-N\delta_\apx\ge 0.9$ the guarantees hold in all periods. On this event, after every step and for all $i\in[n]$,
    \begin{align*}
        \|\ov x_i - x_i\|_{\ov x_i} \le \epsilon_{\apx,x} = \ov \epsilon
    \end{align*}
    and
    \begin{align*}
    \|\ov s_i-s_i\|_{\ov x_i}^* \le w_i^{1/2} \epsilon_{\apx,s} = w_i^{1/2}(1-\epsilon_t)\ov \epsilon \cdot \ov t \le \ov\epsilon\cdot t\cdot w_i,
    \end{align*}
    where the last step uses $w_i\ge 1$ and $(1-\epsilon_t)\ov t\le t$. At the start of a period the approximation is exact.
    Note that the bound for $s$ is loose by a factor of $w_i^{1/2}$.
    When $w_i$s are large, we could improve running time by using a tighter choice of the thresholds, as did in \cite{gs22}. Here we use a loose bound for simplicity of presentation.
\end{proof}

\begin{lemma}\label{lem:treewidth:cpm-time}
We bound the running time of \textsc{CentralPathMaintenance} as following.
\begin{itemize}
    \item \textsc{CentralPathMaintenance.Initialize} takes $\wt O(n \tau^{\omega-1} + n \tau m + n m^{\omega-1})$ time.
    \item If \textsc{CentralPathMaintenance.MultiplyAndMove} is called $N$ times, then with probability at least $0.9$ the bounds of Lemma~\ref{lem:treewidth:cpm-correct} hold after every call and the total running time is 
    \begin{align*}
        \wt O((N n^{-1/2} + \epsilon_t^{-1}\log(t_{\max}/t_{\min})) \cdot n (\tau^2 m+\tau m^2)^{1/2} (\tau^{\omega-1} + \tau m + m^{\omega-1})^{1/2} + N m^{\omega}).
    \end{align*}
    When $m^{\omega-1}\le n\tau$, the term $N m^\omega$ is at most $N n^{1/2} (\tau^2 m+\tau m^2)^{1/2} (\tau^{\omega-1} + \tau m + m^{\omega-1})^{1/2}$ and can be dropped.
    \item \textsc{CentralPathMaintenance.Output} takes $\wt O(n \tau m)$ time.
\end{itemize}
\end{lemma}
\begin{proof}
    \textsc{Initialize} part:
    By Theorem~\ref{thm:treewidth:exactds} and~\ref{thm:treewidth:approxds}.

    \textsc{Output} part:
    By Theorem~\ref{thm:treewidth:exactds}.

    \textsc{MultiplyAndMove} part:
    Between two restarts, the total size of $|L_x|$ returned by $\mathsf{approx}$.\textsc{Query} is bounded by $\wt O(q^2 \zeta_x^2/\epsilon_{\apx,x}^2)$ by Theorem~\ref{thm:treewidth:approxds}.
    By plugging in $\zeta_x = 2\alpha$, $\epsilon_{\apx,x}=\ov \epsilon$, we have
    $\sum_{\ell\in [q]} |L_x^{(\ell)}| = \wt O(q^2).$
    Similarly, for $s$ the bound of Theorem~\ref{thm:treewidth:approxds} is $\wt O(q^2 (\zeta_x^2/\epsilon_{\apx,x}^2+\zeta_s^2/\epsilon_{\apx,s}^2))$. With $\zeta_s=3\alpha\ov t$ and $\epsilon_{\apx,s}=(1-\epsilon_t)\ov\epsilon\,\ov t$ we have 
    $\sum_{\ell\in [q]} |L_s^{(\ell)}| = \wt O(q^2).$
    These bounds, and the bound on the query time used below, hold on the event of the proof of Lemma~\ref{lem:treewidth:cpm-correct}, which has probability at least $0.9$. The running time bound of this lemma is for this event, on which the bounds of Lemma~\ref{lem:treewidth:cpm-correct} hold after every call.

    \textbf{Update time:}
    By Theorem~\ref{thm:treewidth:exactds}, an update changing $k$ blocks of $(\ov x,\ov s)$ costs $\wt O((\tau^2 m+\tau m^2)k)$ in $\mathsf{exact}$. In $\mathsf{approx}$ it costs $\wt O(\tau^2 m k)$. Lemma~\ref{lem:treewidth:bsketch}, read per changed block of $\ov x$, shows that each changed block of $\ov x$ costs $\wt O(\tau^2 m)$. Each of the $\wt O(\tau m k)$ changed entries of $h,\wt h,\epsilon_x,\epsilon_s$ costs $\wt O(1)$ (Algorithm~\ref{alg:treewidth:bsketch-part2}). Hence in a sequence of $q$ updates, total cost for update is $\wt O(q^2 (\tau^2 m+\tau m^2))$.
    So the amortized update cost per iteration is $\wt O(q (\tau^2 m+\tau m^2))$.
    The total update cost is 
    \begin{align*}
        \text{number of iterations} \cdot \text{time per iteration} 
        = \wt O(N q (\tau^2 m+\tau m^2)).
    \end{align*}

    \textbf{Move and query time:}
    By Theorem~\ref{thm:treewidth:exactds}, each call of $\mathsf{exact}$.\textsc{Move} costs $\wt O(m^\omega)$, because it solves a linear system with the $m\times m$ matrix $\wt u$. So the $N$ calls cost $\wt O(N m^\omega)$ in total.
    By Theorem~\ref{thm:treewidth:approxds}, the calls of $\mathsf{approx}$.\textsc{MoveAndQuery} between two restarts cost $\wt O(q^2\tau^2 m)$ in total, which is at most the update cost of the same period.

    \textbf{Init/restart time:}
    We restart the data structure whenever $\ell\ge q$ or $|\ov t-t| > \ov t \epsilon_t$, so there are $O(N/q + \log(t_{\max}/t_{\min}) \epsilon_t^{-1})$ restarts in total.
    By Theorem~\ref{thm:treewidth:exactds} and~\ref{thm:treewidth:approxds}, time cost per restart is $\wt O(n (\tau^{\omega-1} + \tau m + m^{\omega-1}))$.
    So the total initialization time is 
    \begin{align*}
        \text{number of restarts} \cdot \text{time per restart}
        = \wt O((N/q + \log(t_{\max}/t_{\min}) \epsilon_t^{-1}) \cdot n (\tau^{\omega-1} + \tau m + m^{\omega-1}) ).
    \end{align*}

    \textbf{Combine everything:}
    Overall running time is
    \begin{align*}
         \wt O(N q (\tau^2 m+\tau m^2) + (N/q + \log(t_{\max}/t_{\min}) \epsilon_t^{-1}) \cdot n (\tau^{\omega-1} + \tau m + m^{\omega-1}) + N m^{\omega}).
    \end{align*}
    The choice of $q$ that balances the first two terms is
    \begin{align*}
        q = n^{1/2} (\tau^2 m+\tau m^2)^{-1/2} (\tau^{\omega-1} + \tau m + m^{\omega-1})^{1/2},
    \end{align*}
    achieving overall running time
    \begin{align*}
    & ~ \wt O((N n^{-1/2} + \epsilon_t^{-1}\log(t_{\max}/t_{\min})) \cdot n (\tau^2 m+\tau m^2)^{1/2} (\tau^{\omega-1} + \tau m + m^{\omega-1})^{1/2} + N m^{\omega}).
    \end{align*}
    Finally, $(\tau^2 m+\tau m^2) (\tau^{\omega-1} + \tau m + m^{\omega-1}) \ge \tau m^2\cdot m^{\omega-1}=\tau m^{\omega+1}$. So when $m^{\omega-1}\le n\tau$, we have $m^{2\omega}\le n\tau m^{\omega+1}\le n (\tau^2 m+\tau m^2) (\tau^{\omega-1} + \tau m + m^{\omega-1})$, and hence $N m^\omega \le N n^{1/2} (\tau^2 m+\tau m^2)^{1/2} (\tau^{\omega-1} + \tau m + m^{\omega-1})^{1/2}$.
\end{proof}

\begin{proof}[Proof of Theorem~\ref{thm:treewidth:cpm}]
    Combining Lemma~\ref{lem:treewidth:cpm-correct} and~\ref{lem:treewidth:cpm-time}.
\end{proof}

\subsection{Proof of Main Statement}~\label{sec:treewidth:main}
\begin{proof}[Proof of Theorem~\ref{thm:treewidth-formal}]
Use \textsc{CentralPathMaintenance} (Algorithm~\ref{alg:treewidth:cpm}) as the maintenance data structure in Algorithm~\ref{alg:centering}.
Theorem~\ref{thm:robust_ipm} calls \textsc{MultiplyAndMove} $N=O(\sqrt \kappa \log n \log \frac{n \kappa R}{\epsilon r})$ times, with $t_{\max}/t_{\min} = 4\kappa/\epsilon^2$, and Algorithm~\ref{alg:centering} sets $\epsilon_t^{-1} = \wt O(\max_{i\in [n]}(1+\nu_i/w_i+\sqrt{\nu_i}))$. Since $\kappa\ge n$, plugging these into Theorem~\ref{thm:treewidth:cpm} gives the first stated running time, where the term $N m^\omega$ of Theorem~\ref{thm:treewidth:cpm} gives the term $\sqrt \kappa \log(\frac{n\kappa R}{r\epsilon}) \cdot m^{\omega}$.
By Lemma~\ref{lem:treewidth:cpm-time}, when $m^{\omega-1}\le n\tau$ this term is at most the first term, which gives the second stated running time.
In the simplified bound we have $m=\wt O(\tau^{\omega-2})$, so $m^{\omega-1} = \wt O(\tau^{(\omega-1)(\omega-2)}) = \wt O(n\tau)$, because $(\omega-1)(\omega-2)\le 2$ and $\tau\le n$. So the term is absorbed there as well, up to the logarithmic factors that $\wt O$ hides.
The accuracy guarantee holds whenever the closeness invariants of Theorem~\ref{thm:treewidth:cpm} hold, which happens with probability at least $0.9$. The running time bound of Theorem~\ref{thm:treewidth:cpm} holds on the same event. We stop a run that exceeds this bound, which does not happen on this event, so the bound also holds in expectation. Combining Theorem~\ref{thm:treewidth:cpm} and Theorem~\ref{thm:robust_ipm} finishes the proof.
\end{proof}

%% file: ds-rank.tex
\section{Algorithm for Low-Rank QP} \label{sec:rank}
In this section we present a nearly-linear time algorithm for solving low-rank QP with small number of linear constraints. We briefly describe the outline of this section.
\begin{itemize}
    \item In Section~\ref{sec:rank:intro}, we present the main statement of Section~\ref{sec:rank}.
    \item In Section~\ref{sec:rank:cpm-structure}, we present the main data structure \textsc{CentralPathMaintenance}.
    \item In Section~\ref{sec:rank:cpm-ds}, we present several data structures used in \textsc{CentralPathMaintenance}, including \textsc{ExactDS} (Section~\ref{sec:rank:cpm-ds:exactds}), \textsc{ApproxDS} (Section~\ref{sec:rank:cpm-ds:approxds}), \textsc{BatchSketch} (Section~\ref{sec:rank:cpm-ds:batchsketch}).
    \item In Section~\ref{sec:rank:cpm-analysis}, we prove correctness and running time of \textsc{CentralPathMaintenance} data structure.
    \item In Section~\ref{sec:rank:main}, we prove the main result (Theorem~\ref{thm:rank-formal}).
\end{itemize}

\subsection{Main Statement} \label{sec:rank:intro}
We consider programs of the form~\eqref{eqn:qp-general}, i.e.,
\begin{align*}
 \min_{x \in \R^n} ~ & ~ \frac{1}{2} x^\top Q x + c^\top x \\
\mathrm{s.t.~} & ~ A x = b  \\
& ~ x_i \in \cK_i \qquad \forall i\in [n]
\end{align*}
where $Q\in \cS^{n_\tot}$ is symmetric positive semidefinite, $c\in \R^{n_\tot}$, 
$A\in \R^{m\times n_\tot}$, $b\in \R^m$, $\cK_i \subset \R^{n_i}$ is a convex set.
For simplicity, we assume that $n_i = O(1)$ for all $i\in [n]$.

\begin{theorem} \label{thm:rank-formal}
Consider the convex program~\eqref{eqn:qp-general}.
Let $\phi_i: \cK_i \to \R$ be a $\nu_i$-self-concordant barrier for all $i\in [n]$.
Suppose the program satisfies the following properties:
\begin{itemize}
    \item Inner radius $r$: There exists $z\in \R^{n_\tot}$ such that $Az=b$ and $B(z,r) \subseteq \cK$.
    \item Outer radius $R$: $\cK \subseteq B(0,R)$ where $0\in \R^{n_\tot}$.
    \item Lipschitz constant $L$: $\|Q\|_{2\to 2}\le L$, $\|c\|_2 \le L$.
    \item Low rank: We are given a factorization $Q = UV^\top$ of the symmetric positive semidefinite matrix $Q$, where $U,V\in \R^{n_\tot \times k}$.
\end{itemize}
Let $(w_i)_{i\in [n]} \in \R_{\ge 1}^n$ and $\kappa = \sum_{i\in [n]} w_i \nu_i$.
Given any $0<\epsilon\le \frac 12$, we can find an approximate solution $x\in \cK$ satisfiying
\begin{align*}
\frac 12 x^\top Q x + c^\top x &\le \min_{Ax=b, x\in \cK} \left(\frac 12 x^\top Q x + c^\top x\right)+ \epsilon L R(R+1),\\
\|Ax-b\|_1 &\le 3\epsilon(R \|A\|_1 + \|b\|_1),
\end{align*}
in expected time
\begin{align*}
    \wt O\big( & (\sqrt{\kappa/n}\,\log(n\kappa R/(r\epsilon)) + \max_{i\in [n]}(1+\nu_i/w_i+\sqrt{\nu_i})\log(\kappa/\epsilon)) \cdot n (k+m)^{(\omega+1)/2} \\
    & + \sqrt{\kappa}\,\log(n\kappa R/(r\epsilon)) \cdot (k+m)^{\omega}\big).
\end{align*}
We assume that $k$ and $m$ are much smaller than $n$, precisely that $(k+m)^{\omega-1}\le n$. Under this assumption, which the bound above does not use, the last term is at most the first one, and the expected time is
\begin{align*}
    \wt O\big((\sqrt{\kappa/n}\,\log(n\kappa R/(r\epsilon)) + \max_{i\in [n]}(1+\nu_i/w_i+\sqrt{\nu_i})\log(\kappa/\epsilon)) \cdot n (k+m)^{(\omega+1)/2}\big).
\end{align*}
When in addition $\max_{i\in [n]} \nu_i = \wt O(1)$, $w_i=1$, the running time simplifies to
\begin{align*}
    \wt O(n (k+m)^{(\omega+1)/2} \log(R/(r\epsilon))).
\end{align*}
\end{theorem}

\begin{remark}[The assumption on $k$ and $m$]
\label{rem:rank-km-assumption}
The term $\sqrt{\kappa}\,\log(n\kappa R/(r\epsilon)) \cdot (k+m)^{\omega}$ in Theorem~\ref{thm:rank-formal} is the cost of inverting a $k\times k$ matrix and an $m\times m$ matrix in every iteration (Lemma~\ref{lem:rank:cpm-time}). The assumption $(k+m)^{\omega-1}\le n$ holds in our applications, where $k$ and $m$ are much smaller than $n$. The SVM programs of Section~\ref{sec:svm-formulations} have $m\le 2$, the linear SVM has $k=d$, and the Gaussian kernel application has $k=n^{o(1)}$ (Corollary~\ref{cor:gaussian_kernel_almost_linear}).
\end{remark}

\subsection{Algorithm Structure and Central Path Maintenance} \label{sec:rank:cpm-structure}
Similar to the low-treewidth case, our algorithm is based on the robust IPM. Details of the robust IPM will be given in Section~\ref{sec:robust_ipm}.
During the algorithm, we maintain a primal-dual solution pair $(x,s) \in \R^{n_\tot} \times \R^{n_\tot}$ on the robust central path. In addition, we maintain a sparsely-changing approximation $(\ov x, \ov s) \in \R^{n_\tot} \times \R^{n_\tot}$ to $(x,s)$.
In each iteration, we implicitly perform update
\begin{align*}
x&\gets x + \ov t B_{w,\ov x,\ov t}^{-1/2} (I-P_{w,\ov x,\ov t}) B_{w,\ov x,\ov t}^{-1/2} \delta_\mu \\
s&\gets s + \ov t \delta_\mu - \ov t^2 H_{w,\ov x} B_{w,\ov x,\ov t}^{-1/2} (I-P_{w,\ov x,\ov t}) B_{w,\ov x,\ov t}^{-1/2} \delta_\mu
\end{align*}
where
\begin{align*}
    H_{w,\ov x} &= \nabla^2 \phi_w(\ov x) \tag{see Eq.~\eqref{eqn:def:Hwx}} \\
    B_{w,\ov x,\ov t} &= Q + \ov t H_{w,\ov x} \tag{see Eq.~\eqref{eqn:def:B}}\\
    P_{w,\ov x,\ov t} &= B_{w,\ov x,\ov t}^{-1/2} A^\top (A B_{w,\ov x,\ov t}^{-1} A^\top)^{-1} A B_{w,\ov x,\ov t}^{-1/2} \tag{see Eq.~\eqref{eqn:def:P}}
\end{align*}
and explicitly maintain $(\ov x, \ov s)$ such that they remain close to $(x,s)$ in $\ell_\infty$-distance.

This task is handled by the \textsc{CentralPathMaintenance} data structure, which is our main data structure. The robust IPM algorithm (Algorithm~\ref{alg:robust-ipm},~\ref{alg:centering}) directly calls it in every iteration.

The \textsc{CentralPathMaintenance} data structure (Algorithm~\ref{alg:rank:cpm}) has two main sub data structures, \textsc{ExactDS} (Algorithm~\ref{alg:rank:exactds-part1},~\ref{alg:rank:exactds-part2}) and \textsc{ApproxDS} (Algorithm~\ref{alg:rank:approxds}).
\textsc{ExactDS} is used to maintain $(x,s)$, and \textsc{ApproxDS} is used to maintain $(\ov x, \ov s)$.

\begin{algorithm}[!ht]\caption{Main algorithm for low-rank QP.}\label{alg:rank:cpm}
\begin{algorithmic}[1]
\State {\bf data structure} \textsc{CentralPathMaintenance} \Comment{Theorem~\ref{thm:rank:cpm}}
\State {\bf private : members}
\State \hspace{4mm} \textsc{ExactDS} $\mathsf{exact}$ \Comment{Algorithm~\ref{alg:rank:exactds-part1}, \ref{alg:rank:exactds-part2}}
\State \hspace{4mm}  \textsc{ApproxDS} $\mathsf{approx}$ \Comment{ Algorithm~\ref{alg:rank:approxds}}
\State \hspace{4mm} $\ell\in \bN$
\State {\bf end members}
\Procedure{\textsc{Initialize}}{$x, s\in \R^{n_\tot}, t\in \R_+, \ov \epsilon \in (0, 1)$}
    \State $\mathsf{exact}.\textsc{Initialize}(x, s, x, s, t)$ \Comment{Algorithm~\ref{alg:rank:exactds-part1}}
    \State $\ell \leftarrow 0$
    \State $w \gets \nu_{\max}$, $N\gets \sqrt \kappa \log n \log \frac{n \kappa R}{\ov \epsilon r}$
    \State $q \gets n^{1/2} (k^2+m^2)^{-1/2} (k^{\omega-1}+m^{\omega-1})^{1/2}$
    \State $\epsilon_{\apx,x} \leftarrow \ov \epsilon, \zeta_x \leftarrow 2 \alpha, \delta_\apx \leftarrow \frac 1{20N}$
    \State $\epsilon_{\apx,s} \leftarrow (1-\epsilon_t)\cdot\ov \epsilon \cdot \ov t , \zeta_s \leftarrow 3\alpha \ov t$
    \State \begin{align*}
        & ~ \mathsf{approx}.\textsc{Initialize}(x,s,h, \wh h, \wt h, H_{w,\ov x}^{1/2}\wh x, H_{w,\ov x}^{-1/2}\wh s, \beta_x, \beta_s, \wh \beta_x, \wh \beta_s, \wt \beta_x, \wt \beta_s, q, \& \mathsf{exact}, \\
        & ~ \epsilon_{\apx,x}, \epsilon_{\apx,s}, \delta_{\apx})
    \end{align*}
    \State \Comment{Algorithm~\ref{alg:rank:approxds}.Parameters from $x$ to $\wt \beta_s$ come from $\mathsf{exact}$. $\& \mathsf{exact}$ is pointer to $\mathsf{exact}$}
\EndProcedure
\Procedure{\textsc{MultiplyAndMove}}{$t\in \R_+$}
    \State $\ell\gets \ell + 1$
    \If{$|\ov t - t| > \ov t \cdot \epsilon_t$ or $\ell \ge q$}
        \State $x, s \gets \mathsf{exact}.\textsc{Output}()$ \Comment{Algorithm~\ref{alg:rank:exactds-part2}}
        \State \textsc{Initialize}$(x,s,t,\ov \epsilon)$
    \EndIf
    \State $\beta_x, \beta_s, \wh \beta_x, \wh \beta_s, \wt \beta_x, \wt \beta_s \gets \mathsf{exact}.\textsc{Move}()$ \Comment{Algorithm~\ref{alg:rank:exactds-part1}}
    \State $\delta_{\ov x}, \delta_{\ov s } \gets \mathsf{approx}.\textsc{MoveAndQuery}(\beta_x, \beta_s, \wh \beta_x, \wh \beta_s, \wt \beta_x, \wt \beta_s)$ \Comment{Algorithm~\ref{alg:rank:approxds}}
    \State $\delta_h, \delta_{\wh h}, \delta_{\wt h}, \delta_{H_{w,\ov x}^{1/2}\wh x}, \delta_{H_{w,\ov x}^{-1/2}\wh s}\gets \mathsf{exact}.\textsc{Update}(\delta_{\ov x}, \delta_{\ov s})$ \Comment{Algorithm~\ref{alg:rank:exactds-part2}}
    \State $\mathsf{approx}.\textsc{Update}(\delta_{\ov x}, \delta_h, \delta_{\wh h}, \delta_{\wt h}, \delta_{H_{w,\ov x}^{1/2}\wh x}, \delta_{H_{w,\ov x}^{-1/2}\wh s})$ \Comment{Algorithm~\ref{alg:rank:approxds}}
\EndProcedure
\Procedure{\textsc{Output}}{$ $}
    \State \Return $\mathsf{exact}.\textsc{Output}()$ \Comment{Algorithm~\ref{alg:rank:exactds-part2}}
\EndProcedure
\State {\bf end data structure}
\end{algorithmic}
\end{algorithm}
\begin{theorem} \label{thm:rank:cpm}
Data structure \textsc{CentralPathMaintenance} (Algorithm~\ref{alg:rank:cpm}) implicitly maintains the central path primal-dual solution pair $(x,s) \in \R^{n_\tot} \times \R^{n_\tot}$ and explicitly maintains its approximation $(\ov x, \ov s) \in \R^{n_\tot} \times \R^{n_\tot}$ using the following functions:
\begin{itemize}
    \item \textsc{Initialize}$(x\in \R^{n_\tot}, s\in \R^{n_\tot}, t_0\in \R_{>0}, \epsilon\in (0,1))$: Initializes the data structure with initial primal-dual solution pair $(x,s)\in \R^{n_\tot} \times \R^{n_\tot}$, initial central path timestamp $t_0\in \R_{>0}$ in $\wt O(n(k^{\omega-1}+m^{\omega-1}))$ time.
    \item \textsc{MultiplyAndMove}$(t\in \R_{>0})$: It implicitly maintains
    \begin{align*}
    x&\gets x + \ov t B_{w,\ov x,\ov t}^{-1/2} (I-P_{w,\ov x,\ov t}) B_{w,\ov x,\ov t}^{-1/2} \delta_\mu(\ov x, \ov s, \ov t) \\
    s&\gets s + \ov t \delta_\mu - \ov t^2 H_{w,\ov x} B_{w,\ov x,\ov t}^{-1/2} (I-P_{w,\ov x,\ov t}) B_{w,\ov x,\ov t}^{-1/2} \delta_\mu(\ov x, \ov s, \ov t)
    \end{align*}
    where $H_{w,\ov x}$, $B_{w,\ov x, \ov t}$, $P_{w,\ov x, \ov t}$ are defined in Eq.~\eqref{eqn:def:Hwx}\eqref{eqn:def:B}\eqref{eqn:def:P} respectively, and $\ov t$ is some timestamp satisfying $|\ov t-t| \le \epsilon_t \cdot \ov t$.
    
    It also explicitly maintains $(\ov x, \ov s) \in \R^{n_{\tot} \times n_{\tot}}$ such that $\|\ov x_i-x_i\|_{\ov x_i} \le \ov\epsilon$ and $\|\ov s_i-s_i\|^*_{\ov x_i} \le t \ov\epsilon w_i$ for all $i\in [n]$ with probability at least $0.9$.
    
    Assuming the function is called at most $N$ times and $t$ decreases from $t_{\max}$ to $t_{\min}$, on the same event of probability at least $0.9$ the total running time is 
    \begin{align*}
        \wt O((N n^{-1/2} + \epsilon_t^{-1}\log(t_{\max}/t_{\min})) \cdot n (k^{(\omega+1)/2} + m^{(\omega+1)/2}) + N (k^{\omega}+m^{\omega})).
    \end{align*}
    
    \item \textsc{Output}: Computes $(x,s) \in \R^{n_{\tot}} \times \R^{n_{\tot}}$ exactly and outputs them in 
    $\wt O(n(k+m))$
    time.

\end{itemize}
\end{theorem}

\subsection{Data Structures Used in \textsc{CentralPathMaintenance}} \label{sec:rank:cpm-ds}
In this section we present several data structures used in \textsc{CentralPathMaintenance}, including:
\begin{itemize}
    \item \textsc{ExactDS} (Section~\ref{sec:rank:cpm-ds:exactds}): This data structure maintains an implicit representation of the primal-dual solution pair $(x,s)$. This is directly used by \textsc{CentralPathMaintenance}.
    \item \textsc{ApproxDS} (Section~\ref{sec:rank:cpm-ds:approxds}): This data structure explicitly maintains an approximation $(\ov x,\ov s)$ of $(x,s)$. This data structure is directly used by \textsc{CentralPathMaintenance}.
    \item \textsc{BatchSketch} (Section~\ref{sec:rank:cpm-ds:batchsketch}): This data structure maintains a sketch of $(x,s)$. This data structure is used by \textsc{ApproxDS}.
\end{itemize}


\subsubsection{\textsc{ExactDS}} \label{sec:rank:cpm-ds:exactds}
In this section we present the data structure \textsc{ExactDS}. It maintains an implicit representation of the primal-dual solution pair $(x,s)$ by maintaining several sparsely-changing vectors (see Eq.~\eqref{eqn:thm:rank:exactds:x-rep}\eqref{eqn:thm:rank:exactds:s-rep}).

\begin{theorem} \label{thm:rank:exactds}
Data structure \textsc{ExactDS} (Algorithm~\ref{alg:rank:exactds-part1}, \ref{alg:rank:exactds-part2}) implicitly maintains the primal-dual pair $(x,s) \in \R^{n_{\tot}} \times \R^{ n_{\tot} }$, computable via the expression
\begin{align}
x &= \wh x + H_{w,\ov x}^{-1/2} h \beta_x + H_{w,\ov x}^{-1/2} \wh h \wh \beta_x + H_{w,\ov x}^{-1/2} \wt h \wt \beta_x, \label{eqn:thm:rank:exactds:x-rep} \\
s &= \wh s + H_{w,\ov x}^{1/2} h\beta_s + H_{w,\ov x}^{1/2} \wh h \wh \beta_s + H_{w,\ov x}^{1/2} \wt h \wt \beta_s, \label{eqn:thm:rank:exactds:s-rep}
\end{align}
where
$\wh x,\wh s \in \R^{n_\tot}$, $h = H_{w,\ov x}^{-1/2} \ov\delta_\mu \in \R^{n_\tot}$, $\wh h = H_{w,\ov x}^{-1/2} U^\top \in \R^{n_\tot \times k}$ (in this subsection we write $Q = U^\top V$ with $U, V\in \R^{k\times n_\tot}$, the transposes of the factors in Theorem~\ref{thm:rank-formal}), $\wt h = H_{w,\ov x}^{-1/2} A^\top \in \R^{n_\tot\times m}$, $\beta_x, \beta_s\in \R$, $\wh \beta_x, \wh \beta_s\in \R^{k}$, $\wt \beta_x, \wt \beta_s \in \R^m$.

The data structure supports the following functions:
\begin{itemize}
    \item \textsc{Initialize}$(x, s, \ov x, \ov s\in \R^{n_{\tot}}, \ov t \in \R_{>0})$: Initializes the data structure in 
    $\wt O(n(k^{\omega-1}+m^{\omega-1}))$
    time, with initial value of the primal-dual pair $(x, s)$, its initial approximation $(\ov x, \ov s)$, and initial approximate timestamp $\ov t$.
    \item \textsc{Move}$()$: Performs robust central path step
    \begin{align} 
    x & \gets x + \ov t B_{\ov x}^{-1} \delta_\mu - \ov t B_{\ov x}^{-1} A^\top (A B_{\ov x}^{-1} A^\top)^{-1} A B_{\ov x}^{-1} \delta_\mu, \label{eqn:thm:rank:exactds:x-step} \\
    s & \gets s + \ov t \delta_\mu - \ov t^2 H_{w,\ov x} B_{\ov x}^{-1} \delta_\mu + \ov t^2 H_{w,\ov x} B_{\ov x}^{-1} A^\top (A B_{\ov x}^{-1} A^\top)^{-1} A B_{\ov x}^{-1} \delta_\mu \label{eqn:thm:rank:exactds:s-step}
    \end{align}
    in $O(k^\omega+m^\omega)$ time by updating its implicit representation.
    \item \textsc{Update}$(\delta_{\ov x}, \delta_{\ov s} \in \R^{n_{\tot}})$: Updates the approximation pair $(\ov x, \ov s)$ to $(\ov x^{\new} = \ov x + \delta_{\ov x} \in \R^{ n_{\tot} }, \ov s^{\new} = \ov s + \delta_{ \ov s} \in \R^{n_{\tot}} )$ in
    $
    \wt O((k^2+m^2) (\|\delta_{\ov x}\|_0 + \|\delta_{\ov s}\|_0))
    $
    time, and output the changes in variables $h, \wh h, \wt h, H_{w,\ov x}^{1/2} \wh x, H_{w,\ov x}^{-1/2} \wh s$.
    
    Furthermore, $h, H_{w,\ov x}^{1/2} \wh x, H_{w,\ov x}^{-1/2} \wh s$ changes in $O(\|\delta_{\ov x}\|_0 + \|\delta_{\ov s}\|_0)$ coordinates,
    $\wh h$ changes in $O(k(\|\delta_{\ov x}\|_0 + \|\delta_{\ov s}\|_0))$ coordinates,
    $\wt h$ changes in $O(m(\|\delta_{\ov x}\|_0 + \|\delta_{\ov s}\|_0))$ coordinates.
    
    \item \textsc{Output}$()$: Output $x$ and $s$ in $\wt O(n(k+m))$ time.

    \item \textsc{Query}$x(i\in [n])$: Output $x_i$ in $\wt O(k+m)$ time.
    This function is used by \textsc{ApproxDS}.
    \item \textsc{Query}$s(i\in [n])$: Output $s_i$ in $\wt O(k+m)$ time.
    This function is used by \textsc{ApproxDS}.
\end{itemize}
\end{theorem}

\begin{algorithm}[!ht]\caption{ This is used in Algorithm~\ref{alg:rank:cpm}. }\label{alg:rank:exactds-part1}
\begin{algorithmic}[1]
\State {\bf data structure} \textsc{ExactDS} \Comment{Theorem~\ref{thm:rank:exactds} 
}
\State {\bf members}
    \State \hspace{4mm} $\ov x , \ov s \in \R^{n_{\tot}}$, $\ov t\in \R_{+}$, $H_{w,\ov x} \in \R^{n_{\tot} \times n_{\tot}}$
    \State \hspace{4mm} $\wh x, \wh s, \in \R^{n_\tot}$, $\wh h \in \R^{n_\tot \times k}$, $\wt h\in \R^{n_\tot\times m}$, $\beta_x, \beta_s\in \R$, $\wh \beta_x, \wh \beta_s \in \R^k$, $\wt \beta_x, \wt \beta_s \in \R^m$
    \State \hspace{4mm} $u_1, u_2 \in \R^{k\times m}, u_3 \in \R^{m\times m}, u_4 \in \R^m, u_5\in \R^k, u_6\in \R^{k\times k}$
    \State \hspace{4mm} $\ov\alpha \in \R, \ov\delta_\mu \in \R^{n}$
    \State \hspace{4mm} $K\in \bN$
\State {\bf end members}
\Procedure{Initialize}{$x,s, \ov x, \ov s \in \R^{n_\tot}, \ov t \in \R_+$}
    \State $\ov x \gets \ov x$, $\ov s \gets \ov s$, $\ov t \gets \ov t$, $\ov\alpha \gets 0$
    \State $\wh{x} \gets x$, $\wh{s} \gets s$, $\beta_x\gets 0$, $\beta_s \gets 0$, $\wh \beta_x \gets 0$, $\wh \beta_s\gets 0$, $\wt \beta_x \gets 0$, $\wt \beta_s \gets 0$
    \State $H_{w,\ov x} \gets \nabla^2 \phi_w( \ov x )$
    \State \textsc{Initialize$h$}($\ov x, \ov s, H_{w,\ov x}$)
\EndProcedure
\Procedure{Initialize$h$}{$\ov x, \ov s \in \R^{n_{\tot}}, H_{w,\ov x} \in \R^{n_{\tot} \times n_{\tot}}$}
    \For{$i \in [n]$}
        \State $( \ov\delta_{\mu})_i \gets - \frac{ \alpha \sinh( \frac{\lambda}{w_i} \gamma_i( \ov x, \ov s, \ov t ) ) }{ \gamma_i( \ov x, \ov s, \ov t ) } \cdot \mu_i( \ov x, \ov s, \ov t )$
        \State $\ov\alpha \gets \ov\alpha + w_i^{-1} \cosh^2( \frac{\lambda}{w_i} \gamma_i( \ov x, \ov s, \ov t ) )$
    \EndFor 
    \State $h \gets H_{w,\ov x}^{-1/2} \ov\delta_\mu, \wh h \gets H_{w,\ov x}^{-1/2} U^\top, \wt h \gets H_{w,\ov x}^{-1/2} A^\top$
    \State $u_1 \gets UH_{w,\ov x}^{-1} A^\top, u_2 \gets VH_{w,\ov x}^{-1} A^\top, u_3 \gets A H_{w,\ov x}^{-1} A^\top$
    \State $u_4 \gets A H_{w,\ov x}^{-1} \ov \delta_\mu, u_5 \gets V H_{w,\ov x}^{-1} \ov \delta_\mu, u_6 \gets V H_{w,\ov x}^{-1} U^\top$
\EndProcedure 
\Procedure{Move}{$ $}
    \State $v_0 \gets I+\ov t^{-1} u_6 \in \R^{k\times k}$
    \State $v_1\gets \ov t^{-1} u_3 - \ov t^{-2} u_1^\top v_0^{-1} u_2 \in \R^{m\times m}$
    \State $v_2 \gets \ov t^{-1} u_4 - \ov t^{-2} u_1^\top v_0^{-1} u_5 \in \R^{m}$
    \State $\beta_x \gets \beta_x + (\ov\alpha)^{-1/2}$
    \State $\wh \beta_x \gets \wh \beta_x - (\ov\alpha)^{-1/2} \cdot \ov t^{-1} v_0^{-1} u_5 + (\ov\alpha)^{-1/2} \cdot \ov t^{-1} v_0^{-1} u_2 v_1^{-1} v_2$
    \State $\wt \beta_x \gets \wt \beta_x - (\ov\alpha)^{-1/2} \cdot v_1^{-1} v_2$
    \State $\beta_s \gets \beta_s$
    \State $\wh \beta_s \gets \wh \beta_s + (\ov\alpha)^{-1/2} \cdot v_0^{-1} u_5 - (\ov\alpha)^{-1/2} \cdot v_0^{-1} u_2 v_1^{-1} v_2$
    \State $\wt \beta_s \gets \wt \beta_s + (\ov\alpha)^{-1/2} \cdot \ov t v_1^{-1} v_2$
    \State \Return $\beta_x, \beta_s, \wh \beta_x, \wh \beta_s, \wt \beta_x, \wt \beta_s$
\EndProcedure 
\State {\bf end data structure}
\end{algorithmic}
\end{algorithm}

\begin{algorithm}[!ht]\caption{Algorithm~\ref{alg:rank:exactds-part1} continued.
}\label{alg:rank:exactds-part2}
\begin{algorithmic}[1]
\State {\bf data structure} \textsc{ExactDS} \Comment{Theorem~\ref{thm:rank:exactds}}
\Procedure{Output}{$ $}
    \State \Return $\wh x + H_{w,\ov x}^{-1/2} h \beta_x + H_{w,\ov x}^{-1/2} \wh h \wh \beta_x + H_{w,\ov x}^{-1/2} \wt h \wt \beta_x, \wh s + H_{w,\ov x}^{1/2} h\beta_s + H_{w,\ov x}^{1/2} \wh h \wh \beta_s + H_{w,\ov x}^{1/2} \wt h \wt \beta_s$
\EndProcedure
\Procedure{Query$x$}{$i\in [n]$}
    \State \Return $\wh x_i + H_{w,\ov x}^{-1/2} h_{i,*} \beta_x + H_{w,\ov x}^{-1/2} \wh h_{i,*} \wh \beta_x + H_{w,\ov x}^{-1/2} \wt h_{i,*} \wt \beta_x$
\EndProcedure
\Procedure{Query$s$}{$i\in [n]$}
    \State \Return $\wh s_i + H_{w,\ov x}^{1/2} h_{i,*}\beta_s + H_{w,\ov x}^{1/2} \wh h_{i,*} \wh \beta_s + H_{w,\ov x}^{1/2} \wt h_{i,*} \wt \beta_s$
\EndProcedure
\Procedure{Update}{$\delta_{\ov x}, \delta_{\ov s} \in \R^{n_{\tot}}$}
    \State $\Delta_{H_{w,\ov x}} \gets \nabla^2 \phi_w( \ov x + \delta_{\ov x}) - H_{w,\ov x}$
    \Comment{$\Delta_{H_{w,\ov x}}$ is non-zero only for diagonal blocks $(i,i)$ for which $\delta_{\ov x,i} \ne 0$}

    \State $S \gets \{ i \in [n] ~|~ \delta_{\ov x,i} \ne 0 \mathrm{~or~} \delta_{\ov s,i} \ne 0 \}$
    \State $\delta_{ \ov\delta_{\mu} } \gets 0$
    \For{$i \in S$}
        \State Let $\gamma_i = \gamma_i(\ov x, \ov s, \ov t)$, $\gamma_i^{\new} = \gamma_i( \ov x + \delta_{\ov x}, \ov s + \delta_{\ov s}, \ov t )$, $\mu_i^{\new} = \mu_i ( \ov x + \delta_{\ov x}, \ov s + \delta_{\ov s}, \ov t )$
        \State $\ov\alpha \gets \ov\alpha - w_i^{-1} \cosh^2( \frac{\lambda}{w_i} \gamma_i ) + w_i^{-1} \cosh^2( \frac{\lambda}{ w_i } \gamma_i^{\new} )$
        \State $\delta_{\ov\delta_{\mu},i} \gets - \alpha \sinh( \frac{\lambda}{w_i} {\gamma}_i^{\new} ) \cdot \frac{ 1 }{ \gamma^{\new}_i } \cdot \mu^{\new}_i - \ov\delta_{\mu,i}$
    \EndFor
    \State $\delta_h \gets \Delta_{H_{w,\ov x}^{-1/2}} (\ov \delta_\mu + \delta_{\ov \delta_\mu}) + H_{w,\ov x}^{-1/2} \delta_{\ov \delta_\mu}$
    \State $\delta_{\wh h} \gets \Delta_{H_{w,\ov x}^{-1/2}} U^\top$
    \State $\delta_{\wt h} \gets \Delta_{H_{w,\ov x}^{-1/2}} A^\top$
    \State $\delta_{\wh x} \gets H_{w,\ov x}^{-1/2}(h \beta_x + \wh h \wh \beta_x + \wt h \wt \beta_x) - (H_{w,\ov x}+\Delta_{H_{w,\ov x}})^{-1/2}((h+\delta_h) \beta_x + (\wh h+\delta_{\wh h}) \wh \beta_x + (\wt h+\delta_{\wt h}) \wt \beta_x)$ \Comment{supported on the blocks in $S$}
    \State $\delta_{\wh s} \gets H_{w,\ov x}^{1/2}(h \beta_s + \wh h \wh \beta_s + \wt h \wt \beta_s) - (H_{w,\ov x}+\Delta_{H_{w,\ov x}})^{1/2}((h+\delta_h) \beta_s + (\wh h+\delta_{\wh h}) \wh \beta_s + (\wt h+\delta_{\wt h}) \wt \beta_s)$
    \State $\delta_{H_{w,\ov x}^{1/2}\wh x} \gets (H_{w,\ov x}+\Delta_{H_{w,\ov x}})^{1/2}(\wh x+\delta_{\wh x}) - H_{w,\ov x}^{1/2}\wh x$, $\delta_{H_{w,\ov x}^{-1/2}\wh s} \gets (H_{w,\ov x}+\Delta_{H_{w,\ov x}})^{-1/2}(\wh s+\delta_{\wh s}) - H_{w,\ov x}^{-1/2}\wh s$

    \State $h \gets h + \delta_h$, $\wh h \gets \wh h + \delta_{\wh h}$, $\wt h \gets \wt h + \delta_{\wt h}$, $\wh x \gets \wh x + \delta_{\wh x}$, $\wh s \gets \wh s + \delta_{\wh s}$

    \State $u_1 \gets u_1 + U \Delta_{H_{w,\ov x}^{-1}} A^\top$, $u_2 \gets u_2 + V \Delta_{H_{w,\ov x}^{-1}} A^\top$, $u_3 \gets u_3 + A \Delta_{H_{w,\ov x}^{-1}} A^\top$
    \State $u_4 \gets u_4 + A (\Delta_{H_{w,\ov x}^{-1}} (\ov \delta_\mu + \delta_{\ov \delta_\mu}) + H_{w,\ov x}^{-1} \delta_{\ov \delta_\mu})$
    \State $u_5 \gets u_5 + V (\Delta_{H_{w,\ov x}^{-1}} (\ov \delta_\mu + \delta_{\ov \delta_\mu}) + H_{w,\ov x}^{-1} \delta_{\ov \delta_\mu})$
    \State $u_6 \gets u_6 + V \Delta_{H_{w,\ov x}^{-1}} U^\top$
    
    \State $\ov x \gets \ov x + \delta_{ \ov x}$, $\ov s \gets \ov s + \delta_{\ov s}$
    \State $H_{w,\ov x} \gets H_{w,\ov x} + \Delta_{H_{w,\ov x}}$
    \State \Return $\delta_h, \delta_{\wh h}, \delta_{\wt h}, \delta_{H_{w,\ov x}^{1/2}\wh x}, \delta_{H_{w,\ov x}^{-1/2}\wh s}$
\EndProcedure
\State {\bf end data structure}
\end{algorithmic}
\end{algorithm}

\begin{proof}[Proof of Theorem~\ref{thm:rank:exactds}]
By combining Lemma~\ref{lem:rank:exactds-correct} and~\ref{lem:rank:exactds-time}.
\end{proof}

\begin{lemma} \label{lem:rank:exactds-correct}
\textsc{ExactDS} correctly maintains an implicit representation of $(x,s)$, i.e., invariant
\begin{align*}
&x = \wh x + H_{w,\ov x}^{-1/2} h \beta_x + H_{w,\ov x}^{-1/2} \wh h \wh \beta_x + H_{w,\ov x}^{-1/2} \wt h \wt \beta_x,\\
&s = \wh s + H_{w,\ov x}^{1/2} h\beta_s + H_{w,\ov x}^{1/2} \wh h \wh \beta_s + H_{w,\ov x}^{1/2} \wt h \wt \beta_s,\\
&h = H_{w,\ov x}^{-1/2} \ov\delta_\mu \in \R^{n_\tot}, \wh h=H_{w,\ov x}^{-1/2} U^\top \in \R^{n_\tot \times k}, \wt h = H_{w,\ov x}^{-1/2} A^\top \in \R^{n_\tot \times m}, \\
&u_1 = UH_{w,\ov x}^{-1} A^\top \in \R^{k\times m},
u_2 = VH_{w,\ov x}^{-1} A^\top \in \R^{k\times m},
u_3 = A H_{w,\ov x}^{-1} A^\top \in \R^{m\times m},\\
&u_4 = A H_{w,\ov x}^{-1} \ov \delta_\mu \in \R^{m},
u_5 = V H_{w,\ov x}^{-1} \ov \delta_\mu \in \R^{k},
u_6 = V H_{w,\ov x}^{-1} U^\top \in \R^{k\times k},
\\
& \ov\alpha = \sum_{i\in [n]} w_i^{-1} \cosh^2(\frac{\lambda}{w_i} \gamma_i(\ov x, \ov s, \ov t)),\\
& \ov\delta_\mu = \ov\alpha^{1/2} \delta_\mu(\ov x, \ov s, \ov t)
\end{align*}
always holds after every external call, and return values of the queries are correct.
\end{lemma}
\begin{proof}
\textsc{Initialize}:
By checking the definitions we see that all invariants are satisfied after \textsc{Initialize}.

\textsc{Move}:
By the invariants, we have
\begin{align*}
v_0 &= I + \ov t^{-1} V H_{w,\ov x}^{-1} U^\top,\\
v_1 &=  \ov t^{-1} A H_{w,\ov x}^{-1} A^\top - \ov t^{-2} A H_{w,\ov x}^{-1} U^\top (I + \ov t^{-1} V H_{w,\ov x}^{-1} U^\top)^{-1} V H_{w,\ov x}^{-1} A^\top \\
&= A B_{\ov x}^{-1} A^\top \\
v_2 &=  \ov t^{-1} A H_{w,\ov x}^{-1} \ov \delta_\mu - \ov t^{-2} A H_{w,\ov x}^{-1} U^\top (I + \ov t^{-1} V H_{w,\ov x}^{-1} U^\top)^{-1} V H_{w,\ov x}^{-1} \ov \delta_\mu \\
&= A B_{\ov x}^{-1} \ov \delta_\mu.
\end{align*}
By implicit representation~\eqref{eqn:thm:rank:exactds:x-rep},
\begin{align*}
\delta_x =&~ H_{w,\ov x}^{-1/2} h \delta_{\beta_x} + H_{w,\ov x}^{-1/2} \wh h \delta_{\wh \beta_x} + H_{w,\ov x}^{-1/2} \wt h \delta_{\wt \beta_x} \\
=&~ H_{w,\ov x}^{-1} \ov \delta_\mu \cdot (\ov \alpha)^{-1/2} \\
&~+ H_{w,\ov x}^{-1} U^\top \cdot (\ov \alpha)^{-1/2} \ov t^{-1} v_0^{-1} (-u_5 + u_2 v_1^{-1} v_2)\\
&~-H_{w,\ov x}^{-1} A^\top \cdot (\ov \alpha)^{-1/2} v_1^{-1} v_2 \\
=&~ H_{w,\ov x}^{-1} \delta_\mu \\
&~+ H_{w,\ov x}^{-1} U^{\top} \ov t^{-1} (I + \ov t^{-1} V H_{w,\ov x}^{-1} U^\top)^{-1} (-V H_{w,\ov x}^{-1} \delta_\mu + V H_{w,\ov x}^{-1} A^\top (A B_{\ov x}^{-1} A^\top)^{-1} A B_{\ov x}^{-1} \delta_\mu  )\\
&~- H_{w,\ov x}^{-1} A^\top (A B_{\ov x}^{-1} A^\top)^{-1} A B_{\ov x}^{-1} \delta_\mu \\
=&~ \ov t \cdot (\ov t^{-1} H_{w,\ov x}^{-1} - \ov t^{-2} H_{w,\ov x}^{-1}U^\top (I + \ov t^{-1} V H_{w,\ov x}^{-1} U^\top)^{-1} V H_{w,\ov x}^{-1}) \delta_\mu \\
&~- \ov t (\ov t^{-1} H_{w,\ov x}^{-1} - \ov t^{-2} H_{w,\ov x}^{-1}U^\top (I + \ov t^{-1} V H_{w,\ov x}^{-1} U^\top)^{-1} V H_{w,\ov x}^{-1}) A^\top (A B_{\ov x}^{-1} A^\top)^{-1} A B_{\ov x}^{-1} \delta_\mu \\
=&~ \ov t B_{\ov x}^{-1} \delta_\mu - \ov t B_{\ov x}^{-1} A^\top (A B_{\ov x}^{-1} A^\top)^{-1} A B_{\ov x}^{-1} \delta_\mu.
\end{align*}
Comparing with the robust central path step~\eqref{eqn:thm:rank:exactds:x-step}, we see that $x$ is updated correctly.

For $s$, from implicit representation~\ref{eqn:thm:rank:exactds:s-rep} we have
\begin{align*}
\delta_s =&~ H_{w,\ov x}^{1/2} h \delta_{\beta_s} + H_{w,\ov x}^{1/2} \wh h \delta_{\wh \beta_s} + H_{w,\ov x}^{1/2} \wt h \delta_{\wt \beta_s} \\
=&~ -U^\top \cdot (\ov \alpha)^{-1/2} \cdot v_0^{-1} (-u_5 + u_2 v_1^{-1} v_2) + A^\top \cdot (\ov \alpha)^{-1/2} \cdot \ov t v_1^{-1} v_2 \\
=&~ \ov t \delta_\mu - \ov t^2 H_{w,\ov x} B_{\ov x}^{-1} \delta_\mu + \ov t^2 H_{w,\ov x} B_{\ov x}^{-1} A^\top (A B_{\ov x}^{-1} A^\top)^{-1} A B_{\ov x}^{-1} \delta_\mu.
\end{align*}
Comparing with robust central path step~\eqref{eqn:thm:rank:exactds:s-step}, we see that $s$ is updated correctly.

\textsc{Update}:
We would like to prove that \textsc{Update} correctly updates the values of $\wh x, \wh s, h, \wh h, \wt h$, $u_1,u_2,u_3,u_4,u_5,u_6$, $\ov\alpha$, $\ov\delta_\mu$, while preserving the values of $(x,s)$.
In fact, by checking the definitions, it is easy to see that $h, \wh h, \wt h$, $u_1,u_2,u_3,u_4,u_5,u_6$, $\ov\alpha$, $\ov\delta_\mu$ are updated correctly.
Furthermore
\begin{align*}
\delta_x &= \delta_{\wh x} + (H_{w,\ov x}^\new)^{-1/2}(h^\new \beta_x + \wh h^\new \wh \beta_x + \wt h^\new \wt \beta_x) - H_{w,\ov x}^{-1/2}(h \beta_x + \wh h \wh \beta_x + \wt h \wt \beta_x)=0,\\
\delta_s &= \delta_{\wh s} + (H_{w,\ov x}^\new)^{1/2}(h^\new \beta_s + \wh h^\new \wh \beta_s + \wt h^\new \wt \beta_s) - H_{w,\ov x}^{1/2}(h \beta_s + \wh h \wh \beta_s + \wt h \wt \beta_s)=0,
\end{align*}
where $H_{w,\ov x}^\new = H_{w,\ov x}+\Delta_{H_{w,\ov x}}$, $h^\new = h+\delta_h$, and so on. The returned vectors $\delta_{H_{w,\ov x}^{1/2}\wh x}$ and $\delta_{H_{w,\ov x}^{-1/2}\wh s}$ are correct by their definition.
So values of $(x,s)$ are preserved.
\end{proof}

\begin{lemma} \label{lem:rank:exactds-time}
We bound the running time of \textsc{ExactDS} as following.
\begin{enumerate}[label=(\roman*)]
    \item \textsc{ExactDS.Initialize} (Algorithm~\ref{alg:rank:exactds-part1}) runs in $\wt O(n(k^{\omega-1}+m^{\omega-1}))$ time. \label{item:rank:exactds-init-time}
    \item \textsc{ExactDS.Move} (Algorithm~\ref{alg:rank:exactds-part1}) runs in $\wt O(k^\omega+m^\omega)$ time. \label{item:rank:exactds-move-time}
    \item \textsc{ExactDS.Output} (Algorithm~\ref{alg:rank:exactds-part2}) runs in $\wt O(n(k+m))$ time and correctly outputs $(x,s)$. \label{item:rank:exactds-output-time}
    \item \textsc{ExactDS.Query$x$} and \textsc{ExactDS.Query$s$} (Algorithm~\ref{alg:rank:exactds-part2}) runs in $\wt O(k+m)$ time and returns the correct answer. \label{item:rank:exactds-query-time}
    \item \textsc{ExactDS.Update} (Algorithm~\ref{alg:rank:exactds-part2}) runs in $\wt O((k^2+m^2) (\|\delta_{\ov x}\|_0 + \|\delta_{\ov s}\|_0))$ time.
    Furthermore,
    $\|\delta_h\|_0,\|\delta_{\wh x}\|_0,\|\delta_{\wh s}\|_0 = O(\|\delta_{\ov x}\|_0 + \|\delta_{\ov s}\|_0)$,
    $\nnz(\delta_{\wh h}) = O(k(\|\delta_{\ov x}\|_0 + \|\delta_{\ov s}\|_0))$,
    $\nnz(\delta_{\wt h}) = O(m(\|\delta_{\ov x}\|_0 + \|\delta_{\ov s}\|_0))$.
    \label{item:rank:exactds-update-time}
\end{enumerate}
\end{lemma}
\begin{proof}
    \begin{enumerate}[label=(\roman*)]
    \item \textsc{ExactDS.Initialize}:
    Computing $u_1$ and $u_2$ takes $\Tmat(k,n,m) = \wt O(n (k^{\omega-1} + m^{\omega-1}))$ time.
    Computing $u_3$ takes $\Tmat(m,n,m) = \wt O(n m^{\omega-1})$ time.
    Computing $u_4$ takes $O(n m)$ time.
    Computing $u_5$ takes $O(n k)$ time.
    Computing $u_6$ takes $\Tmat(k,n,k) = \wt O(n k^{\omega-1})$ time.
    All other computations are cheaper.
    \item \textsc{ExactDS.Move}: Computing $v_0^{-1}$ takes $\wt O(k^\omega)$ time.
    Computing $v_1^{-1}$ takes $\wt O(m^\omega)$ time.
    All other computations are cheaper.
    \item \textsc{ExactDS.Output}: Takes $\wt O(n(k+m))$ time.
    \item \textsc{ExactDS.Query$x$} and \textsc{ExactDS.Query$s$}: Takes $\wt O(k+m)$ time.
    \item \textsc{ExactDS.Update}: 
    For simplicity, write $t=\|\delta_{\ov x}\|_0 + \|\delta_{\ov s}\|_0$.
    Computing $\delta_h$ takes $\wt O(t)$ time.
    Computing $\delta_{\wh h}$ takes $\wt O(t k)$ time.
    Computing $\delta_{\wt h}$ takes $\wt O(t m)$ time.
    Computing $\delta_{\wh x}$, $\delta_{\wh s}$, $\delta_{H_{w,\ov x}^{1/2}\wh x}$ and $\delta_{H_{w,\ov x}^{-1/2}\wh s}$ takes $\wt O(t (k+m))$ time, since all of them are supported on the blocks in $S$.
    The sparsity statements follow directly.
    Computing $u_1$ and $u_2$ takes $\wt O(tkm)$ time.
    Computing $u_3$ takes $\wt O(t m^2)$ time.
    Computing $u_4$ takes $\wt O(t m)$ time.
    Computing $u_5$ takes $\wt O(t k)$ time.
    Computing $u_6$ takes $\wt O(t k^2)$ time. \qedhere
    \end{enumerate}
\end{proof}

\subsubsection{\textsc{ApproxDS}} \label{sec:rank:cpm-ds:approxds}
In this section we present the data structure \textsc{ApproxDS}.
Given \textsc{BatchSketch}, a data structure maintaining a sketch of the primal-dual pair $(x,s)\in \R^{n_\tot} \times \R^{n_\tot}$, \textsc{ApproxDS} maintains a sparsely-changing $\ell_\infty$-approximation of $(x,s)$. 

\begin{algorithm}[!ht] \caption{This is used in Algorithm~\ref{alg:rank:cpm}.}\label{alg:rank:approxds}
\begin{algorithmic}[1]
\small
\State {\bf data structure} \textsc{ApproxDS} \Comment{Theorem~\ref{thm:rank:approxds}}
\State {\bf private : members}
\State \hspace{4mm} $\epsilon_{\apx,x}, \epsilon_{\apx,s}\in \R$
\State \hspace{4mm} $\ell \in \bN$
\State \hspace{4mm} \textsc{BatchSketch} $\mathsf{bs}$ \Comment{This maintains a sketch of $H_{w,\ov x}^{1/2} x$ and $H_{w,\ov x}^{-1/2} s$. See Algorithm~\ref{alg:rank:batchsketch-part1} and \ref{alg:rank:batchsketch-part2}.}

\State \hspace{4mm} \textsc{ExactDS*} $\mathsf{exact}$ \Comment{This is a pointer to the \textsc{ExactDS} (Algorithm~\ref{alg:rank:exactds-part1}, \ref{alg:rank:exactds-part2}) we maintain in parallel to \textsc{ApproxDS}.}

\State \hspace{4mm} $\wt x, \wt s \in \R^{n_\tot}$
\Comment{$(\wt x, \wt s)$ is a sparsely-changing approximation of $(x, s)$. They have the same value as $(\ov x, \ov s)$, but for these local variables we use $(\wt x, \wt s)$ to avoid confusion.}

\State {\bf end members}

\Procedure{\textsc{Initialize}}{$x,s \in \R^{n_\tot},
h\in \R^{n_\tot}, \wh h\in \R^{n_\tot\times k}, \wt h\in \R^{n_\tot \times m},
H_{w,\ov x}^{1/2}\wh x, H_{w,\ov x}^{-1/2}\wh s \in \R^{n_\tot},
\beta_x,\beta_s\in \R, \wh \beta_x, \wh \beta_s \in \R^k, \wt \beta_x, \wt \beta_s \in \R^m,
q\in \bN,
\textsc{ExactDS*}~\mathsf{exact}, \epsilon_{\apx,x}, \epsilon_{\apx,s}, \delta_{\apx}\in \R$}
    \State $\ell \gets 0$, $q \gets q$
    \State $\epsilon_{\apx,x} \gets \epsilon_{\apx,x}, \epsilon_{\apx,s} \gets \epsilon_{\apx,s}$
    
    \State $\mathsf{bs}.\textsc{Initialize}(x,h,\wh h,\wt h,H_{w,\ov x}^{1/2}\wh x, H_{w,\ov x}^{-1/2}\wh s,\beta_x,\beta_s,\wh \beta_x,\wh \beta_s,\wt \beta_x,\wt \beta_s, \delta_{\apx}/q)$
    \Comment{Algorithm~\ref{alg:rank:batchsketch-part1}}
    \State $\wt x \gets x, \wt s \gets s$
    \State $\mathsf{exact} \gets \mathsf{exact}$
\EndProcedure
\Procedure{\textsc{Update}}{$\delta_{\ov x} \in \R^{n_\tot}, \delta_h \in \R^{n_\tot}, \delta_{\wh h} \in \R^{n_\tot \times k}, \delta_{\wt h}\in \R^{n_\tot \times m}, \delta_{H_{w,\ov x}^{1/2}\wh x}, \delta_{H_{w,\ov x}^{-1/2}\wh s} \in \R^{n_\tot}$}
    \State {$\mathsf{bs}.\textsc{Update}(\delta_{\ov x}, \delta_h, \delta_{\wh h}, \delta_{\wt h}, \delta_{H_{w,\ov x}^{1/2}\wh x}, \delta_{H_{w,\ov x}^{-1/2}\wh s} )$}
    \Comment{Algorithm~\ref{alg:rank:batchsketch-part1}}
    \State $\ell \gets \ell+1$
\EndProcedure 
\Procedure{MoveAndQuery}{$\beta_x,\beta_s\in \R, \wh \beta_x, \wh \beta_s \in \R^k, \wt \beta_x, \wt \beta_s \in \R^m$}
    \State $\mathsf{bs}.\textsc{Move}(\beta_x, \beta_s, \wh \beta_x, \wh \beta_s, \wt \beta_x, \wt \beta_s)$
    \Comment{Algorithm~\ref{alg:rank:batchsketch-part1}. Do not update $\ell$ yet}
    \State $\delta_{\wt x} \gets \textsc{Query$x$}(\epsilon_{\apx,x}/(2\log_2 q+1))$ \Comment{Algorithm~\ref{alg:rank:approxds}}
    \State $\delta_{\wt s} \gets \textsc{Query$s$}(\epsilon_{\apx,s}/(2\log_2 q+1))$ \Comment{Algorithm~\ref{alg:rank:approxds}}
    \For{all $i\in [n]$ with $\delta_{\wt x,i}\ne 0$}
        \State $s_i \gets \mathsf{exact}.\textsc{Query}s(i)$, $\delta_{\wt s,i} \gets s_i-\wt s_i$
        \Comment{Whenever $\wt x_i$ changes, we set $\wt s_i$ to its exact value}
    \EndFor
    \State $\wt x \gets \wt x + \delta_{\wt x}$, $\wt s \gets \wt s + \delta_{\wt s}$
    \State \Return $(\delta_{\wt x}, \delta_{\wt s})$
\EndProcedure
\Procedure{\textsc{Query$x$}}{$\epsilon \in \R$}
   \State Same as Algorithm~\ref{alg:treewidth:approxds-part2}, \textsc{Query$x$}, with $\mathsf{bs}$ and $\mathsf{exact}$ of this section.
\EndProcedure
\Procedure{\textsc{Query$s$}}{$\epsilon \in \R$}
    \State Same as Algorithm~\ref{alg:treewidth:approxds-part2}, \textsc{Query$s$}, with $\mathsf{bs}$ and $\mathsf{exact}$ of this section.
\EndProcedure
\State {\bf end data structure}
\end{algorithmic}
\end{algorithm}

\begin{theorem}\label{thm:rank:approxds}
Given parameters $\epsilon_{\apx,x}, \epsilon_{\apx,s} \in (0, 1), \delta_\apx \in (0,\frac12]$, $\zeta_{x}, \zeta_{s}\in \R$ such that 
\begin{align*}
\|H_{w,\ov x^{(\ell)}}^{1/2} x^{(\ell)}-H_{w,\ov x^{(\ell)}}^{1/2} x^{(\ell+1)}\|_2 \le \zeta_{x}, ~~~ \|H_{w,\ov x^{(\ell)}}^{-1/2} s^{(\ell)}-H_{w,\ov x^{(\ell)}}^{-1/2} s^{(\ell+1)}\|_2 \le \zeta_{s}
\end{align*}
for all $\ell \in \{0,\ldots,q-1\}$, whatever the outputs of the data structure are,
data structure \textsc{ApproxDS} (Algorithm~\ref{alg:rank:approxds}) supports the following operations:
\begin{itemize}
\item $\textsc{Initialize}(x,s \in \R^{n_\tot},
h\in \R^{n_\tot}, \wh h\in \R^{n_\tot\times k}, \wt h\in \R^{n_\tot \times m},
H_{w,\ov x}^{1/2}\wh x, H_{w,\ov x}^{-1/2}\wh s \in \R^{n_\tot},
\beta_x,\beta_s\in \R, \wh \beta_x, \wh \beta_s \in \R^k, \wt \beta_x, \wt \beta_s \in \R^m,
q\in \bN,
\textsc{ExactDS*}~\mathsf{exact}, \epsilon_{\apx,x}, \epsilon_{\apx,s}, \delta_{\apx}\in \R)$: Initialize the data structure in $\wt O(n(k+m))$ time.

\item $\textsc{MoveAndQuery}(\beta_x,\beta_s\in \R, \wh \beta_x, \wh \beta_s \in \R^k, \wt \beta_x, \wt \beta_s \in \R^m)$:
Update values of $\beta_x, \beta_s, \wh \beta_x, \wh \beta_s, \wt \beta_x, \wt \beta_s$ by calling $\textsc{BatchSketch}.\textsc{Move}$.
This effectively moves $(x^{(\ell)}, s^{(\ell)})$ to $(x^{(\ell+1)}, s^{(\ell+1)})$ while keeping $\ov x^{(\ell)}$ unchanged.

Then return the updates $(\delta_{\wt x},\delta_{\wt s})$ of the approximation, so that $\ov x^{(\ell+1)}=\ov x^{(\ell)}+\delta_{\wt x}$ and $\ov s^{(\ell+1)}=\ov s^{(\ell)}+\delta_{\wt s}$. Let $L_x^{(\ell)}, L_s^{(\ell)} \subseteq [n]$ be the sets of blocks in which $\delta_{\wt x}$ and $\delta_{\wt s}$ are nonzero.

Suppose that the inputs of all calls are determined by the state at initialization and by the earlier outputs of the data structure, as is the case in Algorithm~\ref{alg:rank:cpm}. Then with probability at least $1-\delta_\apx$ over the sketching matrix $\Phi$, the following holds for all $\ell\in\{0,\ldots,q-1\}$ simultaneously, even though the inputs depend on $\Phi$ through the earlier outputs. For all $i\in [n]$,
\begin{align*}
\|\ov x^{(\ell+1)}_i - x^{(\ell+1)}_i\|_{\ov x^{(\ell+1)}_i} \le \epsilon_{\apx,x}, \qquad
\|\ov s^{(\ell+1)}_i - s^{(\ell+1)}_i\|_{\ov x^{(\ell+1)}_i}^* \le w_i^{1/2}\epsilon_{\apx,s}.
\end{align*}
Moreover
\begin{align*}
\sum_{0\le \ell \le q-1} |L_x^{(\ell)}| & = \wt O(\epsilon_{\apx,x}^{-2} \zeta_x^2 q^2), \\
\sum_{0\le \ell \le q-1} |L_s^{(\ell)}| & = \wt O((\epsilon_{\apx,x}^{-2} \zeta_x^2+\epsilon_{\apx,s}^{-2} \zeta_s^2) q^2),
\end{align*}
and the total time cost over all queries is at most
\begin{align*}
\wt O\left( (1+\epsilon_{\apx,x}^{-2} \zeta_{x}^2 + \epsilon_{\apx,s}^{-2} \zeta_{s}^2)q^2 (k+m) \right).
\end{align*}

\item $\textsc{Update}(\delta_{\ov x} \in \R^{n_\tot}, \delta_h \in \R^{n_\tot}, \delta_{\wh h} \in \R^{n_\tot \times k}, \delta_{\wt h}\in \R^{n_\tot \times m}, \delta_{H_{w,\ov x}^{1/2}\wh x}, \delta_{H_{w,\ov x}^{-1/2}\wh s} \in \R^{n_\tot})$: Update sketches of $H_{w,\ov x^{(\ell)}}^{1/2} x^{(\ell+1)}$ and $H_{w,\ov x^{(\ell)}}^{-1/2} s^{(\ell+1)}$ by calling \textsc{BatchSketch}.\textsc{Update}.
This effectively moves $\ov x^{(\ell)}$ to $\ov x^{(\ell+1)}$ while keeping $(x^{(\ell+1)}, s^{(\ell+1)})$ unchanged. Then advance timestamp $\ell$.

Each update costs
\begin{align*}
\wt O(\|\delta_h\|_0 + \nnz(\delta_{\wh h}) + \nnz(\delta_{\wt h}) + \|\delta_{H_{w,\ov x}^{1/2} \wh x}\|_0 + \|\delta_{H_{w,\ov x}^{-1/2} \wh s}\|_0 )
\end{align*}
time.
\end{itemize}

\end{theorem}

\begin{proof}
The proof is the same as the proof of Theorem~\ref{thm:treewidth:approxds}. That proof uses the data structures only through the guarantees of \textsc{BatchSketch} and the cost of the exact queries. Here we use Theorem~\ref{thm:rank:batchsketch} in place of Theorem~\ref{thm:treewidth:batchsketch}, and Theorem~\ref{thm:rank:exactds} with Remark~\ref{rmk:snapshot} for the exact queries, which cost $\wt O(k+m)$ each.
\end{proof}

\subsubsection{\textsc{BatchSketch}} \label{sec:rank:cpm-ds:batchsketch}
In this section we present the data structure \textsc{BatchSketch}. It maintains a sketch of $H_{\ov x}^{1/2} x$ and $H_{\ov x}^{-1/2} s$. It is a variation of \textsc{BatchSketch} in \cite{gs22}.

\begin{algorithm}[!ht]\caption{This is used by Algorithm~\ref{alg:rank:approxds}.}\label{alg:rank:batchsketch-part1}
\begin{algorithmic}[1]
\State {\bf data structure} \textsc{BatchSketch} \Comment{Theorem~\ref{thm:rank:batchsketch}}
\State {\bf members}
\State \hspace{4mm} $\Phi\in \R^{r\times n_\tot} $ \Comment{All sketches need to share the same sketching matrix}
\State \hspace{4mm} $\cS,\chi$ partition tree
\State \hspace{4mm} $\ell \in \bN$ \Comment{Current timestamp}
\State \hspace{4mm} \textsc{VectorSketch} $\mathsf{sketch}H_{w,\ov x}^{1/2} \wh x$, $\mathsf{sketch}H_{w,\ov x}^{-1/2} \wh s$, $\mathsf{sketch} h$, $\mathsf{sketch} \wh h$, $\mathsf{sketch} \wt h$
\Comment{Algorithm \ref{alg:treewidth:vsketch}}
\State \hspace{4mm} $\beta_x, \beta_s \in \R, \wh \beta_x, \wh \beta_s \in \R^k, \wt \beta_x, \wt \beta_s \in \R^m$
\State \hspace{4mm} $(\mathsf{history}[t])_{t\ge 0}$
\Comment{Snapshot of data at timestamp $t$. See Remark~\ref{rmk:snapshot}.}
\State {\bf end members}
\Procedure{Initialize}{$\ov x\in \R^{n_\tot}, h\in \R^{n_\tot}, \wh h \in \R^{n_\tot \times k}, \wt h \in \R^{n_\tot\times m}, H_{w,\ov x}^{1/2} \wh x, H_{w,\ov x}^{-1/2} \wh s \in \R^{n_\tot}, \beta_x, \beta_s \in \R, \wh \beta_x, \wh \beta_s \in \R^k, \wt \beta_x, \wt \beta_s \in \R^{m}, \delta_\apx \in \R$}
    \State Construct any valid partition tree $(\cS, \chi)$
    \State $r \gets \Theta(\log^3(n_\tot)\log(1/\delta_\apx))$
    \State Initialize $\Phi\in \R^{r \times n_\tot}$ with iid $\cN(0, \frac 1{r})$
    \State $\beta_x \gets \beta_x, \beta_s \gets \beta_s, \wh \beta_x \gets \wh \beta_x, \wh \beta_s \gets \wh \beta_s, \wt \beta_x \gets \wt \beta_x, \wt \beta_s \gets \wt \beta_s$
    \State $\mathsf{sketch}H_{w,\ov x}^{1/2} \wh x.\textsc{Initialize}( {\cal S} , {\chi}, \Phi, H_{w,\ov x}^{1/2} \wh{x} )$
    \Comment{Algorithm~\ref{alg:treewidth:vsketch}}
    \State $\mathsf{sketch}H_{w,\ov x}^{-1/2} \wh s.\textsc{Initialize}( {\cal S} , {\chi}, \Phi, H_{w,\ov x}^{-1/2} \wh{s} )$
    \Comment{Algorithm~\ref{alg:treewidth:vsketch}}
    \State $\mathsf{sketch} h.\textsc{Initialize}( {\cal S} , {\chi}, \Phi, h )$
    \Comment{Algorithm~\ref{alg:treewidth:vsketch}}
    \State $\mathsf{sketch} \wh h.\textsc{Initialize}( {\cal S} , {\chi}, \Phi, \wh h )$
    \Comment{Algorithm~\ref{alg:treewidth:vsketch}. Here we construct one sketch for $\wh h_{*,i}$ for every $i\in [k]$.}
    \State $\mathsf{sketch} \wt h.\textsc{Initialize}( {\cal S} , {\chi}, \Phi, \wt h)$
    \Comment{Algorithm~\ref{alg:treewidth:vsketch}. Here we construct one sketch for $\wt h_{*,i}$ for every $i\in [m]$.}
    \State $\ell \gets 0$
    \State Make snapshot $\mathsf{history}[\ell]$ \Comment{Remark~\ref{rmk:snapshot}}
\EndProcedure
\Procedure{Move}{$\beta_x, \beta_s\in \R, \wh \beta_x, \wh \beta_s \in \R^k, \wt \beta_x, \wt \beta_s \in \R^m$}
    \State $\beta_x \gets \beta_x, \beta_s \gets \beta_s, \wh \beta_x \gets \wh \beta_x, \wh \beta_s \gets \wh \beta_s, \wt \beta_x \gets \wt \beta_x, \wt \beta_s \gets \wt \beta_s$
    \Comment{Do not update $\ell$ yet}
\EndProcedure
\Procedure{Update}{$\delta_{\ov x} \in \R^{n_\tot}, \delta_h \in \R^{n_\tot}, \delta_{\wh h} \in \R^{n_\tot \times k}, \delta_{\wt h}\in \R^{n_\tot \times m}, \delta_{H_{w,\ov x}^{1/2}\wh x}, \delta_{H_{w,\ov x}^{-1/2}\wh s} \in \R^{n_\tot}$}
    \State $\mathsf{sketch}H_{w,\ov x}^{1/2} \wh x.\textsc{Update}( \delta_{ H_{w,\ov x}^{1/2} \wh{x} } )$
    \Comment{Algorithm~\ref{alg:treewidth:vsketch}}
    \State $\mathsf{sketch}H_{w,\ov x}^{-1/2} \wh s.\textsc{Update}(  \delta_{ H_{w,\ov x}^{-1/2} \wh{s}} )$ 
    \Comment{Algorithm~\ref{alg:treewidth:vsketch}}
    \State $\mathsf{sketch}h.\textsc{Update}( \delta_{ h } )$
    \Comment{Algorithm~\ref{alg:treewidth:vsketch}}
    \State $\mathsf{sketch}\wh h.\textsc{Update}( \delta_{ \wh h } )$
    \Comment{Algorithm~\ref{alg:treewidth:vsketch}}
    \State $\mathsf{sketch}\wt h.\textsc{Update}( \delta_{ \wt h } )$
    \Comment{Algorithm~\ref{alg:treewidth:vsketch}}
    \State $\ell \gets \ell+1$
    \State Make snapshot $\mathsf{history}[\ell]$ \Comment{Remark~\ref{rmk:snapshot}}
\EndProcedure
\State {\bf end data structure}
\end{algorithmic}
\end{algorithm}

\begin{algorithm}[!ht]\caption{\textsc{BatchSketch} Algorithm~\ref{alg:rank:batchsketch-part1} continued. This is used by Algorithm~\ref{alg:rank:approxds}.}\label{alg:rank:batchsketch-part2}
\begin{algorithmic}[1]
\State {\bf data structure} \textsc{BatchSketch} \Comment{Theorem~\ref{thm:rank:batchsketch}}
\State {\bf private:}
    \Procedure{Query$x$Sketch}{$v\in \cS$}
        \Comment{Return the value of $\Phi_{\chi(v)} (H_{w,\ov x}^{1/2} x)_{\chi(v)}$}
        \State \Return $\mathsf{sketch} H_{w,\ov x}^{1/2} \wh x.\textsc{Query}(v) + \mathsf{sketch} h.\textsc{Query}(v) \cdot \beta_x + \mathsf{sketch} \wh h.\textsc{Query}(v) \cdot \wh \beta_x + \mathsf{sketch} \wt h.\textsc{Query}(v) \cdot \wt \beta_x$
        \Comment{Algorithm~\ref{alg:treewidth:vsketch}}
    \EndProcedure
    \Procedure{Query$s$Sketch}{$v\in \cS$}
        \Comment{Return the value of $\Phi_{\chi(v)} (H_{w,\ov x}^{-1/2} s)_{\chi(v)}$}
        \State \Return $\mathsf{sketch} H_{w,\ov x}^{-1/2} \wh s.\textsc{Query}(v) + \mathsf{sketch} h.\textsc{Query}(v) \cdot \beta_s + \mathsf{sketch} \wh h.\textsc{Query}(v) \cdot \wh \beta_s + \mathsf{sketch} \wt h.\textsc{Query}(v) \cdot \wt \beta_s$
        \Comment{Algorithm~\ref{alg:treewidth:vsketch}}
    \EndProcedure
\State {\bf public:}
\Procedure{Query$x$}{$\ell' \in \bN, \epsilon \in \R$}
    \State Same as Algorithm~\ref{alg:treewidth:batchsketch-part2}, \textsc{Query$x$}, using \textsc{Query$x$Sketch} defined here instead of the one in Algorithm~\ref{alg:treewidth:batchsketch-part2}.
\EndProcedure
\Procedure{Query$s$}{$\ell' \in \bN, \epsilon \in \R$}
    \State Same as Algorithm~\ref{alg:treewidth:batchsketch-part2}, \textsc{Query$s$}, using \textsc{Query$s$Sketch} defined here instead of the one in Algorithm~\ref{alg:treewidth:batchsketch-part2}.
\EndProcedure
\State {\bf end structure}
\end{algorithmic}
\end{algorithm}

\begin{theorem}\label{thm:rank:batchsketch}
Data structure \textsc{BatchSketch} (Algorithm~\ref{alg:rank:batchsketch-part1}, \ref{alg:rank:batchsketch-part2}) supports the following operations:
\begin{itemize}
\item $\textsc{Initialize}(\ov x\in \R^{n_\tot}, h\in \R^{n_\tot}, \wh h \in \R^{n_\tot \times k}, \wt h \in \R^{n_\tot\times m}, H_{w,\ov x}^{1/2} \wh x, H_{w,\ov x}^{-1/2} \wh s \in \R^{n_\tot}, \beta_x, \beta_s \in \R, \wh \beta_x, \wh \beta_s \in \R^k, \wt \beta_x, \wt \beta_s \in \R^{m}, \delta_\apx \in \R)$: Initialize the data structure in $\wt O(n (k+m))$ time.

\item $\textsc{Move}(\beta_x, \beta_s\in \R, \wh \beta_x, \wh \beta_s \in \R^k, \wt \beta_x, \wt \beta_s \in \R^m)$: Update values of $\beta_x, \beta_s, \wh \beta_x, \wh \beta_s, \wt \beta_x, \wt \beta_s$ in $O(k+m)$ time. This effectively moves $(x^{(\ell)}, s^{(\ell)})$ to $(x^{(\ell+1)}, s^{(\ell+1)})$ while keeping $\ov x^{(\ell)}$ unchanged.

\item $\textsc{Update}(\delta_{\ov x} \in \R^{n_\tot}, \delta_h \in \R^{n_\tot}, \delta_{\wh h} \in \R^{n_\tot \times k}, \delta_{\wt h}\in \R^{n_\tot \times m}, \delta_{H_{w,\ov x}^{1/2}\wh x}, \delta_{H_{w,\ov x}^{-1/2}\wh s} \in \R^{n_\tot})$:
Update sketches of $H_{w,\ov x^{(\ell)}}^{1/2} x^{(\ell+1)}$ and $H_{w,\ov x^{(\ell)}}^{-1/2} s^{(\ell+1)}$.
This effectively moves $\ov x^{(\ell)}$ to $\ov x^{(\ell+1)}$ while keeping $(x^{(\ell+1)}, s^{(\ell+1)})$ unchanged.
Then advance timestamp $\ell$.

Each update costs
\begin{align*}
\wt O(\|\delta_h\|_0 + \nnz(\delta_{\wh h}) + \nnz(\delta_{\wt h}) + \|\delta_{H_{w,\ov x}^{1/2} \wh x}\|_0 + \|\delta_{H_{w,\ov x}^{-1/2} \wh s}\|_0 ).
\end{align*}

\item $\textsc{Query}x(\ell' \in \bN, \epsilon \in\R)$:
Given timestamp $\ell'$, return a set $S\subseteq [n]$. Let $\ell$ be the current timestamp and let
\begin{align*}
    v := H_{w,\ov x^{(\ell')}}^{1/2} x^{(\ell')} - H_{w,\ov x^{(\ell)}}^{1/2} x^{(\ell+1)} \in \R^{n_\tot}.
\end{align*}
If $\Phi$ is accurate for $v$ (as defined in Theorem~\ref{thm:treewidth:batchsketch}), then
\begin{align*}
    S &\supseteq \{i \in [n]: \|v_i\|_2 \ge \epsilon\},
\end{align*}
and
\begin{align*}
    |S| & = O(\epsilon^{-2} (\ell-\ell'+1) \sum_{\ell' \le t \le \ell} \|H_{w,\ov x^{(t)}}^{1/2} x^{(t)} - H_{w,\ov x^{(t)}}^{1/2} x^{(t+1)}\|_2^2 + \sum_{\ell ' \le t \le \ell-1} \|\ov x^{(t)} - \ov x^{(t+1)}\|_{2,0}),
\end{align*}
and the query costs at most
\begin{align*}
\wt O\Big((k+m) \cdot \Big(1 & +\epsilon^{-2} (\ell-\ell'+1) \sum_{\ell' \le t \le \ell} \|H_{w,\ov x^{(t)}}^{1/2} x^{(t)} - H_{w,\ov x^{(t)}}^{1/2} x^{(t+1)}\|_2^2 \\
& + \sum_{\ell ' \le t \le \ell-1} \|\ov x^{(t)} - \ov x^{(t+1)}\|_{2,0}\Big)\Big)
\end{align*}
running time.

If $v$ does not depend on $\Phi$, then $\Phi$ is accurate for $v$ with probability at least $1-\delta_\apx$, where $\delta_\apx\le\frac12$ is the parameter of \textsc{Initialize}. This statement does not cover a vector $v$ that depends on $\Phi$, for example through earlier answers of the data structure. Theorem~\ref{thm:rank:approxds} treats this case.

\item $\textsc{Query}s(\ell' \in \bN, \epsilon \in\R)$:
Given timestamp $\ell'$, return a set $S\subseteq [n]$. Let $\ell$ be the current timestamp and let
\begin{align*}
    v := H_{w,\ov x^{(\ell')}}^{-1/2} s^{(\ell')} - H_{w,\ov x^{(\ell)}}^{-1/2} s^{(\ell+1)} \in \R^{n_\tot}.
\end{align*}
If $\Phi$ is accurate for $v$, then
\begin{align*}
    S &\supseteq \{i \in [n]: \|v_i\|_2 \ge \epsilon\}
\end{align*}
and
\begin{align*}
    |S| & = O(\epsilon^{-2}  (\ell-\ell'+1)\sum_{\ell' \le t \le \ell} \|H_{w,\ov x^{(t)}}^{-1/2} s^{(t)} - H_{w,\ov x^{(t)}}^{-1/2} s^{(t+1)}\|_2^2 + \sum_{\ell' \le t \le \ell-1} \|\ov x^{(t)} - \ov x^{(t+1)}\|_{2,0}),
\end{align*}
and the query costs at most
\begin{align*}
\wt O\Big((k+m) \cdot \Big(1 & +\epsilon^{-2} (\ell-\ell'+1) \sum_{\ell' \le t \le \ell} \|H_{w,\ov x^{(t)}}^{-1/2} s^{(t)} - H_{w,\ov x^{(t)}}^{-1/2} s^{(t+1)}\|_2^2 \\
& + \sum_{\ell ' \le t \le \ell-1} \|\ov x^{(t)} - \ov x^{(t+1)}\|_{2,0}\Big)\Big)
\end{align*}
running time.
If $v$ does not depend on $\Phi$, then $\Phi$ is accurate for $v$ with probability at least $1-\delta_\apx$.
\end{itemize}
\end{theorem}

\begin{proof}
The proof is essentially the same as proof of \cite[Theorem 4.21]{gs22}.
That proof derives the claims on $S$ and on the running time from the accuracy of $\Phi$ for $v$. It obtains the accuracy from the Johnson-Lindenstrauss lemma (Lemma~\ref{lem:jl_matrix}) and the choice of $r$, which applies to a vector $v$ that does not depend on $\Phi$.
\end{proof}

\subsection{Analysis of \textsc{CentralPathMaintenance}} \label{sec:rank:cpm-analysis}
\begin{lemma}[Correctness of \textsc{CentralPathMaintenance}] \label{lem:rank:cpm-correct}
Algorithm~\ref{alg:rank:cpm} implicitly maintains the primal-dual solution pair $(x,s)$ via representation Eq.~\eqref{eqn:thm:rank:exactds:x-rep}\eqref{eqn:thm:rank:exactds:s-rep}.
It also explicitly maintains $(\ov x, \ov s)\in \R^{n_\tot} \times \R^{n_\tot}$ such that $\|\ov x_i-x_i\|_{\ov x_i} \le \ov \epsilon$ and $\|\ov s_i-s_i\|_{\ov x_i}^* \le t \ov \epsilon w_i$ for all $i\in [n]$ with probability at least $0.9$.
\end{lemma}
\begin{proof}
    Same as proof of Lemma~\ref{lem:treewidth:cpm-correct}, with Theorem~\ref{thm:rank:exactds} and Theorem~\ref{thm:rank:approxds} in place of Theorem~\ref{thm:treewidth:exactds} and Theorem~\ref{thm:treewidth:approxds}.
\end{proof}

\begin{lemma}\label{lem:rank:cpm-time}
We bound the running time of \textsc{CentralPathMaintenance} as following.
\begin{itemize}
    \item \textsc{CentralPathMaintenance.Initialize} takes $\wt O(n (k^{\omega-1}+m^{\omega-1}))$ time.
    \item If \textsc{CentralPathMaintenance.MultiplyAndMove} is called $N$ times, then with probability at least $0.9$ the bounds of Lemma~\ref{lem:rank:cpm-correct} hold after every call and the total running time is 
    \begin{align*}
        \wt O((N n^{-1/2} + \epsilon_t^{-1}\log(t_{\max}/t_{\min})) \cdot n (k+m)^{(\omega+1)/2} + N (k^{\omega}+m^{\omega})).
    \end{align*}
    When $(k+m)^{\omega-1}\le n$, the term $N (k^{\omega}+m^{\omega})$ is at most $N n^{1/2} (k+m)^{(\omega+1)/2}$ and can be dropped.
    \item \textsc{CentralPathMaintenance.Output} takes $\wt O(n(k+m))$ time.
\end{itemize}
\end{lemma}
\begin{proof}
    \textsc{Initialize} part:
    By Theorem~\ref{thm:rank:exactds} and~\ref{thm:rank:approxds}.

    \textsc{Output} part:
    By Theorem~\ref{thm:rank:exactds}.

    \textsc{MultiplyAndMove} part:
    Between two restarts, the total size of $|L_x|$ returned by $\mathsf{approx}$.\textsc{Query} is bounded by $\wt O(q^2 \zeta_x^2/\epsilon_{\apx,x}^2)$ by Theorem~\ref{thm:rank:approxds}.
    By plugging in $\zeta_x = 2\alpha$, $\epsilon_{\apx,x}=\ov \epsilon$, we have
    $\sum_{\ell\in [q]} |L_x^{(\ell)}| = \wt O(q^2).$
    Similarly, for $s$ the bound of Theorem~\ref{thm:rank:approxds} is $\wt O(q^2 (\zeta_x^2/\epsilon_{\apx,x}^2+\zeta_s^2/\epsilon_{\apx,s}^2))$. With $\zeta_s=3\alpha\ov t$ and $\epsilon_{\apx,s}=(1-\epsilon_t)\ov\epsilon\,\ov t$ we have 
    $\sum_{\ell\in [q]} |L_s^{(\ell)}| = \wt O(q^2).$
    These bounds, and the bound on the query time used below, hold on the event of the proof of Lemma~\ref{lem:rank:cpm-correct}, which has probability at least $0.9$. The running time bound of this lemma is for this event, on which the bounds of Lemma~\ref{lem:rank:cpm-correct} hold after every call.

    \textbf{Update time:}
    By Theorem~\ref{thm:rank:exactds} and~\ref{thm:rank:approxds}, in a sequence of $q$ updates, total cost for update is $\wt O(q^2 (k^2+m^2))$.
    So the amortized update cost per iteration is $\wt O(q (k^2+m^2))$.
    The total update cost is 
    \begin{align*}
        \text{number of iterations} \cdot \text{time per iteration} 
        = \wt O(N q (k^2+m^2)).
    \end{align*}

    \textbf{Move and query time:}
    By Theorem~\ref{thm:rank:exactds}, each call of $\mathsf{exact}$.\textsc{Move} costs $\wt O(k^\omega+m^\omega)$, because it inverts the $k\times k$ matrix $v_0$ and the $m\times m$ matrix $v_1$. So the $N$ calls cost $\wt O(N(k^\omega+m^\omega))$ in total.
    By Theorem~\ref{thm:rank:approxds}, the calls of $\mathsf{approx}$.\textsc{MoveAndQuery} between two restarts cost $\wt O(q^2(k+m))$ in total, which is at most the update cost of the same period.

    \textbf{Init/restart time:}
    We restart the data structure whenever $\ell\ge q$ or $|\ov t-t| > \ov t \epsilon_t$, so there are $O(N/q + \log(t_{\max}/t_{\min}) \epsilon_t^{-1})$ restarts in total.
    By Theorem~\ref{thm:rank:exactds} and~\ref{thm:rank:approxds}, time cost per restart is $\wt O(n (k^{\omega-1} + m^{\omega-1}))$.
    So the total initialization time is 
    \begin{align*}
        \text{number of restarts} \cdot \text{time per restart}
        = \wt O((N/q + \log(t_{\max}/t_{\min}) \epsilon_t^{-1}) \cdot n (k^{\omega-1}+m^{\omega-1}) ).
    \end{align*}

    \textbf{Combine everything:}
    Overall running time is
    \begin{align*}
         \wt O(N q (k^2+m^2) + (N/q + \log(t_{\max}/t_{\min}) \epsilon_t^{-1}) \cdot n (k^{\omega-1}+m^{\omega-1}) + N (k^{\omega}+m^{\omega})).
    \end{align*}
    The choice of $q$ that balances the first two terms is
    \begin{align*}
        q = n^{1/2} (k^2+m^2)^{-1/2} (k^{\omega-1}+m^{\omega-1})^{1/2},
    \end{align*}
    achieving overall running time
    \begin{align*}
    &~\wt O((N n^{-1/2} + \epsilon_t^{-1}\log(t_{\max}/t_{\min})) \cdot n (k^2+m^2)^{1/2} (k^{\omega-1}+m^{\omega-1})^{1/2} + N (k^{\omega}+m^{\omega})) \\
    =&~ \wt O((N n^{-1/2} + \epsilon_t^{-1}\log(t_{\max}/t_{\min})) \cdot n (k+m)^{(\omega+1)/2} + N (k^{\omega}+m^{\omega})).
    \end{align*}
    Finally, when $(k+m)^{\omega-1}\le n$, we have $k^\omega+m^\omega\le (k+m)^{\omega} = (k+m)^{(\omega-1)/2}(k+m)^{(\omega+1)/2}\le n^{1/2}(k+m)^{(\omega+1)/2}$, and hence $N(k^\omega+m^\omega)\le N n^{1/2}(k+m)^{(\omega+1)/2}$.
\end{proof}

\begin{proof}[Proof of Theorem~\ref{thm:rank:cpm}]
    Combining Lemma~\ref{lem:rank:cpm-correct} and~\ref{lem:rank:cpm-time}.
\end{proof}

\subsection{Proof of Main Statement}~\label{sec:rank:main}
\begin{proof}[Proof of Theorem~\ref{thm:rank-formal}]
Use \textsc{CentralPathMaintenance} (Algorithm~\ref{alg:rank:cpm}) as the maintenance data structure in Algorithm~\ref{alg:centering}.
Theorem~\ref{thm:robust_ipm} calls \textsc{MultiplyAndMove} $N=O(\sqrt \kappa \log n \log \frac{n \kappa R}{\epsilon r})$ times, with $t_{\max}/t_{\min} = 4\kappa/\epsilon^2$, and Algorithm~\ref{alg:centering} sets $\epsilon_t^{-1} = \wt O(\max_{i\in [n]}(1+\nu_i/w_i+\sqrt{\nu_i}))$. Plugging these into Theorem~\ref{thm:rank:cpm} gives the first stated running time, where the term $N(k^\omega+m^\omega)$ of Theorem~\ref{thm:rank:cpm} gives the term $\sqrt{\kappa}\,\log(n\kappa R/(r\epsilon)) \cdot (k+m)^{\omega}$.
When $(k+m)^{\omega-1}\le n$, we have $(k+m)^{\omega}\le n^{1/2}(k+m)^{(\omega+1)/2}$, so this term is at most $\sqrt{\kappa/n}\,\log(n\kappa R/(r\epsilon))\cdot n(k+m)^{(\omega+1)/2}$, which gives the second stated running time.
Combining Theorem~\ref{thm:rank:cpm} (with $\delta_\apx \gets \frac{1}{20N}$ in Algorithm~\ref{alg:rank:cpm}) and Theorem~\ref{thm:robust_ipm} gives the guarantee with probability at least $0.9$. The running time bound of Theorem~\ref{thm:rank:cpm} holds on the same event. At the end we compute $(x,s)$ exactly by \textsc{Output} and check in $\wt O(n(k+m))$ time that $x\in\cK$ and $\|\mu_i(x,s,t)\|_{x_i}^*\le w_i$ for all $i$. Since $Ax=b$ and $-Qx+A^\top y+s=c$ hold exactly after every step, a passing check certifies the output by Lemmas~\ref{lem:approx-sol-from-potential} and~\ref{lem:init-point}. On failure, or if a run exceeds the time bound, we restart, so the expected number of runs is $O(1)$.
\end{proof}

%% file: ipm.tex
\section{Robust IPM Analysis} \label{sec:robust_ipm}
In this section we present a robust IPM algorithm for quadratic programming.
The algorithm is a modification of previous robust IPM algorithms for linear programming \cite{lsz19,lv21}.

Convention: Variables are in $n$ blocks of dimension $n_i$ ($i\in [n]$).
Total dimension is $n_\tot = \sum_{i\in [n]} n_i$.
We write $x = (x_1,\ldots, x_n) \in \R^{n_\tot}$ where $x_i \in \R^{n_i}$.
We consider programs of the following form:
\begin{align} \label{eqn:qp-general}
 \min_{x \in \R^n} ~ & ~ \frac{1}{2} x^\top Q x + c^\top x \\
\mathrm{s.t.~} & ~ A x = b  \nonumber \\
& ~ x_i \in \cK_i \qquad \forall i\in [n] \nonumber
\end{align}
where $Q\in \cS^{n_\tot}$ is symmetric positive semidefinite, $c\in \R^{n_\tot}$, 
$A\in \R^{m\times n_\tot}$, $b\in \R^m$, $\cK_i \subset \R^{n_i}$ is a convex set.
Let $\cK = \prod_{i\in [n]} \cK_i$.

\begin{theorem} \label{thm:robust_ipm}
Consider the convex program~\eqref{eqn:qp-general}.
Let $\phi_i: \cK_i \to \R$ be a $\nu_i$-self-concordant barrier for all $i\in [n]$.
Suppose the program satisfies the following properties:
\begin{itemize}
    \item Inner radius $r$: There exists $z\in \R^{n_\tot}$ such that $Az=b$ and $B(z,r) \subseteq \cK$.
    \item Outer radius $R$: $\cK \subseteq B(0,R)$ where $0\in \R^{n_\tot}$.
    \item Lipschitz constant $L$: $\|Q\|_{2\to 2}\le L$, $\|c\|_2 \le L$.
\end{itemize}

Let $(w_i)_{i\in [n]} \in \R^n_{\ge 1}$ and $\kappa = \sum_{i\in [n]} w_i \nu_i$.
For any $0<\epsilon\le \frac 12$, Algorithm~\ref{alg:robust-ipm} outputs an approximate solution $x$ in $O(\sqrt \kappa \log n \log \frac{n \kappa R}{\epsilon r})$ steps, satisfying
\begin{align*}
\frac 12 x^\top Q x + c^\top x &\le \min_{Ax=b, x\in \cK} \left(\frac 12 x^\top Q x + c^\top x\right)+ \epsilon L R(R+1),\\
\|Ax-b\|_1 &\le 3\epsilon(R \|A\|_1 + \|b\|_1),\\
x&\in \cK.
\end{align*}
\end{theorem}

\begin{algorithm}[!ht]\caption{Our main algorithm}\label{alg:robust-ipm}
\begin{algorithmic}[1]
\Procedure{RobustQPIPM}{$Q \in \cS^{n_\tot}, c\in \R^{n_\tot}, A \in \R^{m\times n_\tot}, b\in \R^m, (\phi_i: \cK_i \to \R)_{i\in [n]}, w \in \R^n$}
    \State {\color{blue} /* Initial point reduction */}
    \State $\epsilon \gets \min\{\epsilon, \frac 1\lambda\}$, where $\lambda = 64 \log(256 (n+1) (1+\sum_{i\in [n]} w_i))$ is the parameter of Algorithm~\ref{alg:centering} for the augmented program
    \State $\rho \gets \frac{1}{LR(R+1)}$, $x^{(0)} \gets \arg \min_x \sum_{i\in [n]} w_i \phi_i(x_i)$, $s^{(0)} \gets \epsilon \rho (c + Qx^{(0)})$
    \State $\ov x \gets \begin{bmatrix}
        x^{(0)} \\ 1
    \end{bmatrix}$, $\ov s \gets \begin{bmatrix}
        s^{(0)} \\ 1
    \end{bmatrix}$, $\ov Q \gets \begin{bmatrix}
        \epsilon \rho Q & 0 \\ 0 & 0
    \end{bmatrix}$, $\ov A \gets \begin{bmatrix}
        A \mid b-Ax^{(0)}
    \end{bmatrix}$
    \State $\ov w \gets \begin{bmatrix}
        w \\ 1
    \end{bmatrix}$, $\ov \phi_i = \phi_i \forall i\in [n]$, $\ov \phi_{n+1}(x) := -\log x$
    \State $(x,s)\gets \textsc{Centering}(\ov Q, \ov A, (\ov \phi_i)_{i\in [n+1]}, \ov w, \ov x, \ov s, t_{\mathrm{start}}=1, t_{\mathrm{end}} = \frac{\epsilon^2}{4\kappa})$
    \State \Return $(x_{1:n},s_{1:n})$
\EndProcedure
\end{algorithmic}
\end{algorithm}

\begin{algorithm}[!ht]\caption{Subroutine used by Algorithm~\ref{alg:robust-ipm}} \label{alg:centering}
\begin{algorithmic}[1]
\Procedure{Centering}{$Q \in \cS^{n_\tot}, A \in \R^{m\times n_\tot}, (\phi_i: \cK_i \to \R)_{i\in [n]}, w \in \R^n, x \in \R^{n_\tot}, s \in \R^{n_\tot}, t_{\mathrm{start}} \in \R_{>0}, t_{\mathrm{end}} \in \R_{>0}$}
    \State {\color{blue}/* Parameters */}
    \State $\lambda = 64 \log(256 n \sum_{i\in [n]} w_i)$, $\ov \epsilon = \frac 1{1440\lambda}$, $\alpha = \frac{\ov \epsilon}2$
    \State $\epsilon_t = \frac{\ov \epsilon}4 (\min_{i\in [n]} \min\{\frac{w_i}{w_i+\nu_i}, \frac{1}{1+\sqrt{\nu_i}}\})$, $h = \frac{\alpha}{64 \sqrt{ \kappa }}$
    \State {\color{blue}/* Definitions */}
    \State $\phi_w(x) := \sum_{i\in [n]} w_i \phi_i(x_i)$
    \State $\mu_i(x,s,t) := s_i/t + w_i \nabla \phi_i(x_i)$, $\forall i\in [n]$ \Comment{Eq.~\eqref{eqn:def:mu_i}}
    \State $\gamma_i(x,s,t) \gets \|\mu_i^t(x,s)\|_{x_i}^*$, $\forall i\in [n]$ \Comment{Eq.~\eqref{eqn:def:gamma_i}}
    \State $c_i(x,s,t) := \frac{\sinh(\frac{\lambda}{w_i} \gamma_i(x,s,t))}{\gamma_i(x,s,t) \sqrt{\sum_{j\in [n]} w_j^{-1} \cosh^2 (\frac{\lambda}{w_j} \gamma_j(x,s,t))}}$, $\forall i\in [n]$ \Comment{Eq.~\eqref{eqn:def:c_i}}
    \State $H_{w,x} := \nabla^2 \phi_w(x)$ \Comment{Eq.~\eqref{eqn:def:Hwx}}
    \State $B_{w,x,t} := Q + t H_{w,x}$ \Comment{Eq.~\eqref{eqn:def:B}}
    \State $P_{w,x,t} := B_{w,x,t}^{-1/2} A^\top (A B_{w,x,t}^{-1} A^\top)^{-1} A B_{w,x,t}^{-1/2}$ \Comment{Eq.~\eqref{eqn:def:P}}
    \State {\color{blue}/* Main loop */}
    \State $\ov t \gets t \gets t_{\mathrm{start}}$, $\ov x \gets x$, $\ov s \gets s$
    \While{$t > t_{\mathrm{end}}$}
    \State Maintain $\ov x, \ov s, \ov t$ such that $\|\ov x_i - x_i\|_{\ov x_i} \le \ov \epsilon$, $\|\ov s_i - s_i\|_{\ov x_i}^* \le t \ov \epsilon w_i$ and $|\ov t-t| \le \epsilon_t \ov t$ \label{line:alg:centering:cpm-guarantee}
    \State $\delta_{\mu,i} \gets -\alpha \cdot c_i(\ov x, \ov s, \ov t) \cdot \mu_i(\ov x, \ov s, \ov t)$, $\forall i\in [n]$ \Comment{Eq.~\eqref{eqn:def:delta_mu}}
    \State  \label{line:centering:delta_x_assump}
    Pick $\delta_x$ and $\delta_s$ such that $A \delta_x = 0$, $\delta_s - Q \delta_x \in \mathrm{Range}(A^\top)$ and 
    \begin{align*}
    \|\delta_x -  \ov t B_{w,\ov x,\ov t}^{-1/2} (I-P_{w,\ov x,\ov t}) B_{w,\ov x,\ov t}^{-1/2} \delta_\mu\|_{w,\ov x} &\le \ov \epsilon\alpha,\\
    \|\ov t^{-1} \delta_s - (\delta_\mu - \ov t H_{w,\ov x} B_{w,\ov x,\ov t}^{-1/2} (I-P_{w,\ov x,\ov t}) B_{w,\ov x,\ov t}^{-1/2} \delta_\mu )\|_{w,\ov x}^* &\le \ov \epsilon\alpha.
    \end{align*}
    \State $t \gets \max\{(1-h)t, t_{\mathrm{end}}\}$, $x\gets x + \delta_x$, $s \gets s + \delta_s$
    \EndWhile
    \State \Return $(x,s)$
\EndProcedure
\end{algorithmic}
\end{algorithm}

\subsection{Preliminaries} \label{sec:robust_ipm:prelim}
Previous works on linear programming (e.g.~\cite{lsz19}, \cite{lv21}) use the following path:
\begin{align*}
s/t + \nabla \phi_w(x) = & ~ \mu, \\
A x = & ~ b, \\
A^\top y + s = & ~ c
\end{align*}
where
$
\phi_w(x) := \sum_{i=1}^n w_i \phi_i(x_i).
$

For quadratic programming, we modify the above central path as following:
\begin{align*}
s/t + \nabla \phi_w(x) = & ~ \mu, \\
A x = & ~ b, \\
-Q x + A^\top y + s = & ~ c.
\end{align*}

We make the following definitions.
\begin{definition}
For each $i \in [n]$, we define the $i$-th coordinate error
\begin{align} \label{eqn:def:mu_i}
 \mu_i (x,s,t):= \frac{s_i}{t} + w_i \nabla \phi_i(x_i)
\end{align}

We define $\mu_i$'s norm as 
\begin{align} \label{eqn:def:gamma_i}
\gamma_i(x,s,t) := \| \mu_i(x,s,t) \|_{x_i}^* .
\end{align}

We define the soft-max function by
\begin{align} \label{eqn:def:Psi}
\Psi_{\lambda}(r):= \sum_{i=1}^n \cosh ( \lambda \frac{r_i}{w_i} )
\end{align}
for some $\lambda > 0$ and the potential function is the soft-max of the norm of the error of each coordinate
\begin{align} \label{eqn:def:Phi}
 \Phi(x,s,t) = \Psi_{\lambda} ( \gamma(x,s,t) )
\end{align}

We choose the step direction $\delta_\mu$ as
\begin{align} \label{eqn:def:delta_mu}
\delta_{\mu,i} := - \alpha \cdot c_i (x,s,t) \cdot \mu_i(x,s,t)
\end{align}
where
\begin{align} \label{eqn:def:c_i}
    c_i(x,s,t) := \frac{\sinh(\frac{\lambda}{w_i} \gamma_i(x,s,t))}{\gamma_i(x,s,t) \sqrt{\sum_{j\in [n]} w_j^{-1} \cosh^2 (\frac{\lambda}{w_j} \gamma_j(x,s,t))}}
\end{align}

\end{definition}

We define induced norms as following. Note that we include the weight vector $w$ in the subscript to avoid confusion.
\begin{definition}\label{def:w_x_norm}
For each block $\cK_i$, we define
\begin{align*}
\| v \|_{x_i} := & ~ \| v \|_{ \nabla^2 \phi_i(x_i) }, \\
\| v \|_{x_i}^* := & ~ \| v \|_{ ( \nabla^2 \phi_i(x_i) )^{-1} }
\end{align*}
for $v \in \R^{n_i}$.

For the whole domain $\cK = \prod_{i=1}^n \cK_i$, we define
\begin{align*}
 \| v \|_{w,x} := & ~ \| v \|_{\nabla^2 \phi_w(x)} = ( \sum_{i=1}^n w_i \| v_i \|_{x_i}^2 )^{1/2}, \\
 \| v \|_{w,x}^* := & ~ \| v \|_{ (\nabla^2 \phi_w(x) )^{-1} } = ( \sum_{i=1}^n w_i^{-1} ( \| v_i \|_{x_i}^* )^2 )^{1/2}
\end{align*}
for $v \in \R^{n_\tot}$.
\end{definition}

The Hessian matrices of the barrier functions appear a lot in the computation.
\begin{definition}\label{def:H}
We define matrices $H_{x,i}\in \R^{n_i \times n_i}$ and $H_{w,x} \in \R^{n_\tot \times n_\tot}$ as
\begin{align}
H_{x,i} &:= \nabla^2 \phi_i(x_i), \label{eqn:def:Hxi} \\
H_{w,x} &:= \nabla^2 \phi_w(x). \label{eqn:def:Hwx}
\end{align}
\end{definition}
From the definition, we see that
\begin{align*}
H_{w,x,(i,i)} &= w_i H_{x,i}.
\end{align*}

The following equations are immediate from definition.
\begin{claim}\label{cla:w_x_norm}
Let $H_{w,x} \in \R^{n_\tot \times n_\tot}$ be defined as Definition~\ref{def:H}. 
For $v\in \R^{n_\tot}$, we have
\begin{align*}
\|v\|_{w,x} &= \|H_{w,x}^{1/2} v\|_2, \\
\|v\|_{w,x}^* &= \|H_{w,x}^{-1/2} v\|_2.\\
\end{align*}
\end{claim}

\begin{claim}\label{cla:x_i_norm}
For each $i \in [n]$, let $H_{x,i}$ be defined as Definition~\ref{def:H}.
For $v\in \R^{n_i}$, $i\in [n]$, we have
\begin{align*}
\|v\|_{x_i} &= \|H_{x,i}^{1/2} v\|_2, \\
\|v\|_{x_i}^* &= \|H_{x,i}^{-1/2} v\|_2. \\
\end{align*}
\end{claim}

We define matrices $B$ and $P$ used in the algorithm. 
\begin{definition}\label{def:B_P}
Let $A,Q$ denote two fixed matrices. 
Let $H_{w,x} \in \R^{n_\tot \times n_\tot}$ be defined as Definition~\ref{def:H}. 
    We define matrix $B_{w,x,t} \in \R^{n_\tot \times n_\tot}$ as
    \begin{align} \label{eqn:def:B}
    B_{w,x,t} := Q + t \cdot H_{w,x}
    \end{align}
    We define projection matrix $P_{w,x,t} \in \R^{n_\tot \times n_\tot}$ as
    \begin{align} \label{eqn:def:P}
    P_{w,x,t} \gets B_{w,x,t}^{-1/2} A^\top (A B_{w,x,t}^{-1} A^\top)^{-1} A B_{w,x,t}^{-1/2}.
    \end{align}
\end{definition}

\subsection{Deriving the Central Path Step} \label{sec:robust_ipm:derive}
In this section we explain how to derive the central path step.

We follow the central path
\begin{align*}
s/t + \nabla \phi_w(x) &= \mu \\
Ax &= b \\
-Qx +  A^\top y + s &= c \\
\end{align*}

We perform gradient descent on $\mu$ with step $\delta_\mu$. Then Newton step gives
\begin{align}
\frac{1}{t} \delta_s + \nabla^2 \phi_w(x) \delta_x   = & ~ \delta_{\mu} \label{eq:cpm_3} \\
A \delta_x  = & ~ 0 \label{eq:cpm_2} \\
- Q \delta_x + A^\top \delta_y + \delta_s = & ~ 0 \label{eq:cpm_1} 
\end{align}
where $\delta_x$ (resp.~$\delta_y$, $\delta_s$) is the step taken by $x$ (resp.~$y$, $s$).

For simplicity, we define $H \in \R^{n_\tot \times n_\tot}$ to represent $\nabla^2 \phi_w(x)$.\footnote{In this section, and in this section only, we omit the subscript in $H$, $B$, $P$ for simplicity.}

From Eq.~\eqref{eq:cpm_3} we get
\begin{align} \label{eq:cpm-step1}
\delta_s = t \delta_\mu - t H\delta_x.
\end{align}
Plug the above equation into Eq.~\eqref{eq:cpm_1} we get
\begin{align} \label{eq:cpm-step2}
-Q \delta_x + A^\top \delta_y + t  \delta_\mu - t H \delta_x = 0.
\end{align}
Let $B=Q+t H$,  
multiply by $A B^{-1}$ we get
\begin{align*}
-A \delta_x + A B^{-1} A^\top \delta_y + t AB^{-1} \delta_\mu = 0.
\end{align*}
Using Eq.~\eqref{eq:cpm_2} we get 
\begin{align*}
A B^{-1} A^\top \delta_y + t AB^{-1} \delta_\mu = 0.
\end{align*}
Solve for $\delta_y$ (assuming that $AB^{-1} A$ is invertible), we get
\begin{align*}
\delta_y = -t (AB^{-1} A^\top)^{-1} A B^{-1} \delta_\mu.
\end{align*}
Plug into Eq.~\eqref{eq:cpm-step2} we get
\begin{align*}
-B \delta_x - t A^\top (AB^{-1} A^\top)^{-1} A B^{-1} \delta_\mu + t \delta_\mu = 0.
\end{align*}
Solve for $\delta_x$ we get
\begin{align*}
\delta_x &= t B^{-1} \delta_\mu - t B^{-1} A^\top  (AB^{-1} A^\top)^{-1} A B^{-1}\delta_\mu \\
&= t B^{-1/2} (I-P) B^{-1/2} \delta_\mu
\end{align*}
where $P = B^{-1/2}  A^\top  (AB^{-1} A^\top)^{-1} A B^{-1/2}$ is the projection matrix.
Solve for $\delta_s$ in Eq.~\eqref{eq:cpm-step1} we get
\begin{align*}
\delta_s = t \delta_\mu - t^2 H B^{-1/2} (I-P) B^{-1/2} \delta_\mu.
\end{align*}

In summary, we have
\begin{align*}
\delta_x &= t B^{-1/2} (I-P) B^{-1/2} \delta_\mu,\\
\delta_y &= -t (AB^{-1} A^\top)^{-1} A B^{-1} \delta_\mu,\\
\delta_s &= t \delta_\mu - t^2 H B^{-1/2} (I-P) B^{-1/2} \delta_\mu,\\
P &= B^{-1/2}  A^\top  (AB^{-1} A^\top)^{-1} A B^{-1/2}.
\end{align*}

These equations will guide the design of our actual algorithm.

\subsection{Bounding Movement of Potential Function}
The goal of this section is to bound the movement of potential function during the robust IPM algorithm.

In robust IPM, we do not need to follow the ideal central path exactly over the entire algorithm. Instead, we only use an approximate version.
For convenience of analysis we state two assumptions (see Algorithm~\ref{alg:centering}, Line~\ref{line:centering:delta_x_assump}).
 
\begin{assumption} \label{ass:ipm_approx_step_direction}
    We make the following assumptions on $\delta_x \in \R^{n_\tot}$ and $\delta_s \in \R^{n_\tot}$.
    \begin{align*}
    \|\delta_x -  \ov t B_{w,\ov x,\ov t}^{-1/2} (I-P_{w,\ov x,\ov t}) B_{w,\ov x,\ov t}^{-1/2} \delta_\mu\|_{w,\ov x} &\le \ov \epsilon\alpha,\\
    \|\ov t^{-1} \delta_s - (\delta_\mu - \ov t H_{w,\ov x} B_{w,\ov x,\ov t}^{-1/2} (I-P_{w,\ov x,\ov t}) B_{w,\ov x,\ov t}^{-1/2} \delta_\mu )\|_{w,\ov x}^* &\le \ov \epsilon\alpha.
    \end{align*}
\end{assumption}

The following lemma bounds the movement of potential function $\Psi$ assuming bound on $\delta_\gamma$.
\begin{lemma}[{\cite[Lemma A.5]{y20}}]

For any $r \in \R^{n}$, and $w \in \R^{n}_{\geq 1}$. Let $\alpha$ and $\lambda$ denote the parameters that are satisfying $0 \leq \alpha \leq \frac{1}{8 \lambda}$.

Let $\epsilon_r \in \R^{n}$ denote a vector satisfying 
\begin{align*}
( \sum_{i=1}^n w_i^{-1} \epsilon_{r,i}^2 )^{1/2} \leq \alpha / 8.
\end{align*}
Suppose that vector $\ov{r} \in \R^{n}$ is satisfying the following property
\begin{align*}
|r_i - \ov{r}_i| \leq \frac{w_i}{8 \lambda}, ~~~ \forall i \in [n]
\end{align*}

We define vector $\delta_r \in \R^{n}$ as follows:
\begin{align*}
 \delta_{r,i} := \frac{ - \alpha \cdot \sinh( \frac{\lambda}{w_i} \ov{r}_i ) }{ \sqrt{ \sum_{j=1}^n w_j^{-1} \cosh^2( \frac{\lambda}{w_j} \ov{r}_j ) } } + \epsilon_{r,i}.
\end{align*}
Then, we have that
\begin{align*}
\Psi_{\lambda} (r + \delta_r) \leq \Psi_{\lambda}(r) - \frac{\alpha \lambda}{2} ( \sum_{i=1}^n w_i^{-1} \cosh^2 ( \lambda \frac{r_i}{w_i} ) )^{1/2} + \alpha \lambda ( \sum_{i=1}^n w_i^{-1} )^{1/2}
\end{align*}
\end{lemma}

The following lemma bounds the norm of $\delta_\mu.$
\begin{lemma}[Bounding norm of $\delta_\mu$] \label{lem:ipm_delta_mu_length}
\begin{align*}
\|\delta_\mu(\ov x, \ov s, \ov t)\|_{w,\ov x}^*  \le \alpha.
\end{align*}
\end{lemma}
\begin{proof}
\begin{align*}
(\|\delta_\mu(\ov x, \ov s, \ov t)\|_{w,\ov x}^*)^2 
= & ~ \sum_{i=1}^n w_i^{-1} ( \| \delta_{\mu,i} (\ov x, \ov s, \ov t) \|_{\ov{x}_i}^* )^2 \\
= &~ \alpha^2 \sum_{i\in [n]} w_i^{-1}  c_i^2(\ov x, \ov s, \ov t) \cdot (\|  \mu_i(\ov x, \ov s, \ov t)\|_{\ov{x}_i}^*)^2 \\
= & ~ \alpha^2 \sum_{i\in [n]} w_i^{-1}  c_i^2(\ov x, \ov s, \ov t) \cdot \|H_{\ov x,i}^{-1/2} \mu_i(\ov x, \ov s, \ov t)\|_2^2  \\
= & ~ \alpha^2 \sum_{i\in [n]} w_i^{-1} c_i^2(\ov x, \ov s, \ov t) \cdot \gamma_i^2 (\ov{x}, \ov{s}, \ov{t}) \\
= & ~ \alpha^2 \sum_{i\in [n]} \frac{w_i^{-1} \sinh^2(\frac{\lambda}{w_i} \gamma_i(\ov x, \ov s, \ov t))}{ \gamma_i^2(\ov x, \ov s, \ov t) \cdot {\sum_{j\in [n]} w_j^{-1} \cosh^2(\frac{\lambda}{w_j} \gamma_j(\ov x, \ov s, \ov t))}} \cdot \gamma_i^2(\ov x, \ov s, \ov t) \\
= & ~ \alpha^2 \frac{ \sum_{j\in [n]} w_j^{-1} \sinh^2 ( \frac{\lambda}{w_j} \gamma_j (\ov{x}, \ov{s}, \ov{t} ) ) }{ \sum_{j\in [n]} w_j^{-1} \cosh^2 ( \frac{\lambda}{w_j} \gamma_j (\ov{x}, \ov{s}, \ov{t} ) ) } \\
\le & ~ \alpha^2.
\end{align*}
where the first step follows from Definition~\ref{def:w_x_norm}, the second step follows from $\delta_{\mu, i} (\ov{x}, \ov{s}, \ov{t}) = -\alpha \cdot c_i(\ov{x}, \ov{s}, \ov{t})  \cdot \mu_i( \ov{x}, \ov{s}, \ov{t} )$, the third step follows from norm of $\ov{x}_i$ (see Definition~\ref{def:w_x_norm}), the forth step follows from $\gamma_i (\ov{x}, \ov{s}, \ov{t}) = \| H_{\ov{x}, i }^{-1/2} \mu_i( \ov{x}, \ov{s}, \ov{t} ) \|_2$ (see~Eq.~\eqref{eqn:def:gamma_i}), the fifth step follows from $ c_i(\ov x,\ov s,\ov t)^2 = \frac{\sinh^2(\frac{\lambda}{w_i} \gamma_i(\ov x,\ov s,\ov t))}{\gamma_i^2(\ov x,\ov s,\ov t) \sum_{j\in [n]} w_j^{-1} \cosh^2 (\frac{\lambda}{w_j} \gamma_j(\ov x,\ov s,\ov t))}$ (see Eq.~\eqref{eqn:def:c_i}), the sixth step follows from canceling the term $\gamma_i^2(\ov{x}, \ov{s}, \ov{t})$, and the last step follows from $\cosh^2(x) \geq \sinh^2(x)$ for all $x$.
\end{proof}

The following lemma bounds the norm of $\delta_x$ and $\delta_s$.
\begin{lemma} \label{lem:DLY_Lemma_A.9}
For each $i \in [n]$, we define
$\alpha_i := \|\delta_{x,i}\|_{\ov x_i}.$
Then, we have
\begin{align}
\|\delta_x\|_{w,\ov x} = ( \sum_{i\in [n]} w_i \alpha_i^2 )^{1/2} \le \frac{9}{8} \alpha.
\end{align}
In particular, we have
$
\alpha_i \le \frac 98 \alpha.
$
Similarly, for $\delta_s$, we have
\begin{align}
\|\delta_s\|_{w,\ov x}^* = \sqrt{\sum_{i\in [n]} w_i^{-1} (\|\delta_{s,i}\|_{\ov x_i}^*)^2} \le \frac {17}8 \alpha \cdot t.
\end{align}
 
\end{lemma}
\begin{proof}

For $\delta_x$, we have
\begin{align*}
\|\delta_x\|_{w,\ov x}
&\le \|\ov t H_{w,\ov x}^{1/2} B_{w,\ov x,\ov t}^{-1/2} (I-P_{w,\ov x,\ov t}) B_{w,\ov x,\ov t}^{-1/2} \delta_\mu\|_2 + \ov \epsilon \alpha \\
& \le \| \ov t^{1/2} (I-P_{w,\ov x,\ov t}) B_{w,\ov x,\ov t}^{-1/2} \delta_\mu\|_{2} + \ov \epsilon\alpha\\
& \le \| \ov t^{1/2} B_{w,\ov x,\ov t}^{-1/2} \delta_\mu\|_{2} + \ov \epsilon \alpha \\
& \le \| H_{w,\ov x}^{-1/2} \delta_\mu\|_{2} + \ov \epsilon \alpha \\
& \le \alpha + \ov \epsilon \alpha\\
&\le \frac 98 \alpha.
\end{align*}
First step follows from Assumption~\ref{ass:ipm_approx_step_direction}. 
Second step is because $\ov t H_{w,\ov x} \preceq B_{w,\ov x,\ov t}$. 
Third step is because $P_{w,\ov x,\ov t}$ is a projection matrix.
Fourth step is because $\ov t H_{w,\ov x} \preceq B_{w,\ov x,\ov t}$.
Fifth step is by Lemma~\ref{lem:ipm_delta_mu_length}. 
Sixth step is because $\ov \epsilon \le \frac 18$.

For $\delta_s$, we have
\begin{align*}
\|\delta_s\|_{w,\ov x}^* & \le \|\ov t \delta_\mu\|_{w,\ov x}^* + \|\ov t^2 H_{w,\ov x} B_{w,\ov x,\ov t}^{-1/2} (I-P_{w,\ov x,\ov t}) B_{w,\ov x,\ov t}^{-1/2} \delta_\mu\|_{w,\ov x}^* + \ov \epsilon \alpha \ov t \\
& \le \alpha \ov t + \alpha \ov t + \ov \epsilon \alpha \ov t \\
& \le \frac {17}8 \alpha \cdot t.
\end{align*}
First step is by triangle inequality and the assumption that
\begin{align*}
\delta_s \approx \ov t \delta_\mu - \ov t^2 H_{w,\ov x} B_{w,\ov x,\ov t}^{-1/2} (I-P_{w,\ov x,\ov t}) B_{w,\ov x,\ov t}^{-1/2} \delta_\mu.
\end{align*}
Second step is by same analysis as the analysis for $\delta_x$.
Third step is by $\ov t \le \frac{33}{32} t$ and $\ov \epsilon \le \frac{1}{32}$.
\end{proof}

The following lemma shows that $\mu^{\new}$ is close to $\mu + \delta_\mu$ under an approximate step.
\begin{lemma}[Variation of~{\cite[Lemma A.9]{y20}}] \label{lem:DLY_Lemma_A.10}
For each $i \in [n]$, we define
\begin{align*}
\beta_i := \| \epsilon_{\mu,i} \|_{x_i}^*
\end{align*}
For each $i \in [n]$, let
\begin{align*}
\mu_i (x^{\new}, s^{\new}, t) = \mu_i(x,s, t) + \delta_{\mu,i} + \epsilon_{\mu,i}.
\end{align*}
Then, we have
\begin{align*}
( \sum_{i=1}^n w_i^{-1} \beta_i^2 )^{1/2} \leq 15 \ov{\epsilon} \alpha.
\end{align*}

\end{lemma}
\begin{proof}

The proof is similar as \cite[Lemma A.9]{y20}, except for changing the definitions of $\epsilon_1$ and $\epsilon_2$:
\begin{align*}
\epsilon_1 &:= H_{w,\ov x}^{1/2} \delta_x - \ov{t} \cdot H_{w, \ov{x} }^{1/2} B_{w,\ov x,\ov t}^{-1/2} (I-P_{w, \ov{x}, \ov{t} } ) B_{w,\ov x,\ov t}^{-1/2} \delta_\mu, \\
\epsilon_2 &:= \ov t^{-1} H_{w,\ov x}^{-1/2} \delta_s - H_{w,\ov x}^{-1/2} (\delta_\mu - \ov t H_{w,\ov x} B_{w,\ov x,\ov t}^{-1/2} (I-P_{w,\ov x,\ov t}) B_{w,\ov x,\ov 
t}^{-1/2} \delta_\mu ).
\end{align*}

One key step in the proof of \cite{y20} is the following property:
\begin{align*}
\delta_{\mu,i} = \ov{t}^{-1} \cdot \delta_{s,i} + H_{w,\ov x} \delta_{x,i} - H_{w,\ov x}^{1/2} (\epsilon_1+\epsilon_2) .
\end{align*}
Under our new definition of $\epsilon_1$ and $\epsilon_2$, the above property still holds.
Remaining parts of the proof are similar and we omit the details here.
\end{proof}

The following lemma shows that error $\mu(\ov x, \ov s, \ov t)$ on the robust central path is close to error $\mu(x,s,t)$ on the ideal central path. Furthermore, norms of errors $\gamma_i(x,s,t)$ and $\gamma_i(\ov x, \ov s, \ov t)$ are also close to each other.
\begin{lemma}[{\cite[Lemma A.10]{y20}}] \label{lem:DLY_Lemma_A.11}
Assume that $\gamma_i(x,s,t) \leq w_i$ for all $i$, and that $\ov x, \ov s, \ov t$ satisfy the invariants of Algorithm~\ref{alg:centering}, Line~\ref{line:alg:centering:cpm-guarantee}. For all $i \in [n]$, we have
\begin{align*}
\| \mu_i(x,s, t) - \mu_i ( \ov{x}, \ov{s} , \ov{t} ) \|_{x_i}^* \leq 3 \ov{\epsilon} w_i.
\end{align*}
Furthermore, we have that
\begin{align*}
| \gamma_i(x,s, t) - \gamma_i ( \ov{x}, \ov{s} , \ov{t}) | \leq 5 \ov{\epsilon} w_i.
\end{align*}
 
\end{lemma}
\begin{proof}
We write out the proof, because its last step is where the choice of $\epsilon_t$ in Algorithm~\ref{alg:centering} is used. Fix $i\in[n]$. We use two facts about the $\nu_i$-self-concordant barrier $\phi_i$. First, $\|\nabla\phi_i(x_i)\|_{x_i}^*\le\sqrt{\nu_i}$, by the definition of a $\nu_i$-self-concordant barrier. Second, for every $y$ with $\vartheta:=\|y-\ov x_i\|_{\ov x_i}<1$ we have~\cite[Theorem 4.1.6]{n98}
\begin{align}\label{eq:hessian_sandwich}
(1-\vartheta)^2\,\nabla^2\phi_i(\ov x_i)\preceq\nabla^2\phi_i(y)\preceq(1-\vartheta)^{-2}\,\nabla^2\phi_i(\ov x_i).
\end{align}
The invariants give $\|x_i-\ov x_i\|_{\ov x_i}\le\ov\epsilon$, and $\ov\epsilon\le\frac1{1440}$ because $\lambda\ge 1$. Taking $y=x_i$ in~\eqref{eq:hessian_sandwich} and inverting the matrices, we get for every $v\in\R^{n_i}$
\begin{align}\label{eq:dual_norm_compare}
(1-\ov\epsilon)\|v\|_{\ov x_i}^*\le\|v\|_{x_i}^*\le(1-\ov\epsilon)^{-1}\|v\|_{\ov x_i}^*.
\end{align}
Next, let $x^{(u)}:=\ov x_i+u(x_i-\ov x_i)$ for $u\in[0,1]$, so that $\|x^{(u)}-\ov x_i\|_{\ov x_i}\le u\ov\epsilon$. By~\eqref{eq:hessian_sandwich} with $y=x^{(u)}$, every $v\in\R^{n_i}$ satisfies $\|\nabla^2\phi_i(x^{(u)})v\|_{\ov x_i}^*\le(1-u\ov\epsilon)^{-2}\|v\|_{\ov x_i}$. Hence
\begin{align}\label{eq:grad_phi_change}
\|\nabla\phi_i(x_i)-\nabla\phi_i(\ov x_i)\|_{\ov x_i}^*
=\Big\|\int_0^1\nabla^2\phi_i(x^{(u)})(x_i-\ov x_i)\,\d u\Big\|_{\ov x_i}^*
\le\int_0^1\frac{\ov\epsilon}{(1-u\ov\epsilon)^2}\,\d u=\frac{\ov\epsilon}{1-\ov\epsilon}.
\end{align}

For the first claim we write
\begin{align*}
\mu_i(x,s,t)-\mu_i(\ov x,\ov s,\ov t)=\Big(\frac1t-\frac1{\ov t}\Big)s_i+\frac{s_i-\ov s_i}{\ov t}+w_i\big(\nabla\phi_i(x_i)-\nabla\phi_i(\ov x_i)\big)
\end{align*}
and bound the three terms. For the first term, $|\ov t-t|\le\epsilon_t\ov t$ gives $|1-t/\ov t|\le\epsilon_t$. Moreover $s_i/t=\mu_i(x,s,t)-w_i\nabla\phi_i(x_i)$, so
\begin{align*}
\Big\|\Big(\frac1t-\frac1{\ov t}\Big)s_i\Big\|_{x_i}^*
\le\epsilon_t\Big\|\frac{s_i}{t}\Big\|_{x_i}^*
\le\epsilon_t\big(\gamma_i(x,s,t)+w_i\|\nabla\phi_i(x_i)\|_{x_i}^*\big)
\le\epsilon_t w_i(1+\sqrt{\nu_i})
\le\frac{\ov\epsilon}{4}w_i,
\end{align*}
where the last step uses $\epsilon_t\le\frac{\ov\epsilon}{4(1+\sqrt{\nu_i})}$. For the second term, the invariants give $\|s_i-\ov s_i\|_{\ov x_i}^*\le t\ov\epsilon w_i$ and $t\le(1+\epsilon_t)\ov t$, so $\|(s_i-\ov s_i)/\ov t\|_{\ov x_i}^*\le(1+\epsilon_t)\ov\epsilon w_i$. For the third term we use~\eqref{eq:grad_phi_change}. We convert the last two bounds to the norm $\|\cdot\|_{x_i}^*$ by~\eqref{eq:dual_norm_compare} and get
\begin{align}\label{eq:mu_distance_sharp}
\|\mu_i(x,s,t)-\mu_i(\ov x,\ov s,\ov t)\|_{x_i}^*
\le\Big(\frac14+\frac{1+\epsilon_t}{1-\ov\epsilon}+\frac{1}{(1-\ov\epsilon)^2}\Big)\ov\epsilon w_i
\le\frac52\ov\epsilon w_i,
\end{align}
because $\epsilon_t\le\ov\epsilon\le\frac1{1440}$. This proves the first claim.

For the second claim, let $v:=\mu_i(x,s,t)$, so that $\gamma_i(x,s,t)=\|v\|_{x_i}^*\le w_i$. Then
\begin{align*}
|\gamma_i(x,s,t)-\gamma_i(\ov x,\ov s,\ov t)|
\le & ~ \big|\|v\|_{x_i}^*-\|v\|_{\ov x_i}^*\big|+\|v-\mu_i(\ov x,\ov s,\ov t)\|_{\ov x_i}^* \\
\le & ~ \frac{\ov\epsilon}{1-\ov\epsilon}w_i+\frac{1}{1-\ov\epsilon}\cdot\frac52\ov\epsilon w_i \\
\le & ~ 5\ov\epsilon w_i,
\end{align*}
where the second step uses~\eqref{eq:dual_norm_compare} twice and~\eqref{eq:mu_distance_sharp}.
\end{proof}
\begin{remark}[The choice of $\epsilon_t$]
\label{rem:eps_t}
The robust interior point method of~\cite{dly21} uses $\epsilon_t=\frac{\ov\epsilon}{4}\min_{i\in[n]}\frac{w_i}{w_i+\nu_i}$. Its proof of the statement of Lemma~\ref{lem:DLY_Lemma_A.11} bounds the norm of $s_i/t$ by $2(w_i+\nu_i)$. The identity $s_i/t=\mu_i(x,s,t)-w_i\nabla\phi_i(x_i)$ only gives the bound $w_i(1+\sqrt{\nu_i})$ in the norm $\|\cdot\|_{x_i}^*$, and this bound is attained. It exceeds $2(w_i+\nu_i)$ when $w_i(\sqrt{\nu_i}-1)>2\nu_i$, that is, when $w_i$ is large compared with $\sqrt{\nu_i}$. This is why $\epsilon_t$ in Algorithm~\ref{alg:centering} also has the term $\frac{1}{1+\sqrt{\nu_i}}$, and why the general running times of Theorems~\ref{thm:treewidth-formal} and~\ref{thm:rank-formal} have the factor $\max_{i\in[n]}(\nu_i/w_i+\sqrt{\nu_i})$. When $w_i=1$ for all $i$, as in our applications, the two choices of $\epsilon_t$ agree.
\end{remark}

The following lemma bounds the change of $\gamma$ under one robust IPM step.
\begin{lemma}[{\cite[Lemma A.12]{y20}}] \label{lem:DLY_Lemma_A.13}
Assume $\Phi(x,s,t) \leq \cosh(\lambda)$. For all $i \in [n]$, we define
\begin{align*}
\epsilon_{r,i} := \gamma_i(x^{\new}, s^{\new}) - \gamma_i (x,s,t) + \alpha \cdot c_i( \ov{x}, \ov{s}, \ov{t} ) \cdot \gamma_i ( \ov{x} , \ov{s}, \ov{t} ).
\end{align*}
Then, we have
\begin{align*}
( \sum_{i=1}^n w_i^{-1} \epsilon_{r,i}^2 )^{1/2} \leq 90 \cdot \ov{\epsilon} \cdot \lambda \alpha + 4 \cdot \max_{i\in [n]} ( w_i^{-1} \gamma_i(x,s, t) ) \cdot \alpha.
\end{align*}
\end{lemma}
\begin{proof}
    The proof is similar to the proof of \cite[Lemma A.12]{y20}. By replacing corresponding references in \cite{y20} by our versions (Lemma~\ref{lem:DLY_Lemma_A.9},~\ref{lem:DLY_Lemma_A.10},~\ref{lem:DLY_Lemma_A.11}) we get proof of this lemma.
\end{proof}

Finally, the following theorem bounds the movement of potential function $\Phi$ under one robust IPM step.
\begin{theorem}[Variation of {\cite[Theorem A.15]{y20}}] \label{thm:DLY_Thm_A.16}

Assume $\Phi(x,s,t) \le \cosh(\lambda/64)$. Then for any $0\le h \le \frac{\alpha}{64 \sqrt{\sum_{i\in [n]} w_i \nu_i}}$, we have
\begin{align*}
\Phi(x^\new, s^\new, t^\new) \le (1-\frac{\alpha \lambda}{4\sqrt{\sum_{i\in [n]} w_i}} ) \cdot \Phi(x,s,t) + \alpha \lambda \sqrt{\sum_{i\in [n]} w_i^{-1}}.
\end{align*}
In particular, for any $\cosh(\lambda/128) \le \Phi(x,s,t) \le \cosh(\lambda/64)$, we have that
\begin{align*}
\Phi(x^\new, s^\new, t^\new) \le \Phi(x,s,t).
\end{align*}
\end{theorem}
\begin{proof}
Similar to the proof of \cite[Theorem A.15]{y20}, but replacing lemmas with the corresponding QP versions.
\end{proof}

\subsection{Initial Point Reduction}
In this section, we propose an initial point reduction scheme for quadratic programming.
Our scheme is closer to \cite{lsz19} rather than \cite{y20,lv21}. The reason is that \cite{lv21}'s initial point reduction requires an efficient algorithm for finding the optimal solution to an unconstrained program, which may be difficult in quadratic programming.

\begin{lemma}[{\cite[Theorem 4.1.7 and Lemma 4.2.4]{n98}}] \label{lem:LSZ_Lemma_D.1}
Let $\phi$ be a $\nu$-self-concordant barrier. Then for any $x, y \in \mathrm{dom} (\phi)$, we have
\begin{align*}
\langle \nabla \phi(x), y - x \rangle \leq & ~ \nu, \\
\langle \nabla \phi(y) - \nabla \phi(x), y - x \rangle \geq & ~ \frac{\| y - x \|_x^2}{ 1 + \| y - x \|_x }.
\end{align*}
Let $x^* = \arg \min_x \phi(x)$. For any $x \in \R^n$ such that $\| x - x^* \|_{x^*} \leq 1$, we have that $x \in \mathrm{dom} (\phi)$.

\end{lemma}

\begin{lemma}[QP version of {\cite[Lemma D.2]{lsz19}}] \label{lem:init-point}
Work under the setting of Theorem~\ref{thm:robust_ipm}.
Let $x^{(0)} = \arg \min_x \sum_{i\in [n]} w_i \phi_i(x_i)$.
Let $\rho = \frac 1{LR(R+1)}$.
For any $0<\epsilon\le \frac 12$, the modified program
\begin{align*}
\min_{ \ov{A} \ov{x} = \ov{b} , \ov{x} \in \cK\times \R_{\ge 0} } \left(\frac{1}{2} \ov{x}^\top \ov{Q} \ov{x} + \ov c^\top \ov x\right)
\end{align*}
with
\begin{align*}
\ov Q = \begin{bmatrix}
    \epsilon \rho  Q & 0 \\ 0 & 0
\end{bmatrix}, \qquad \ov{A} = [ A ~|~ b - A x^{(0)} ], \qquad \ov{b} = b, \qquad \ov{c} = \begin{bmatrix} \epsilon \rho c \\ 1 \end{bmatrix}
\end{align*}
satisfies the following:
\begin{itemize}
 \item $\ov{x} = \begin{bmatrix} x^{(0)} \\ 1 \end{bmatrix}$, $\ov{y} = 0\in \R^m$ and $\ov{s} = \begin{bmatrix}  \epsilon\rho (c + Q x^{(0)}) \\ 1 \end{bmatrix}$ are feasible primal dual vectors with $\| \ov{s} + \nabla \ov{\phi}_w ( \ov{x} ) \|_{ \ov{x} }^* \leq \epsilon$ where $\ov{\phi}_w ( \ov{x} ) = \sum_{i=1}^n w_i \phi_i( \ov{x}_i ) - \log ( \ov{x}_{n+1} ) $.
 \item For any $\ov{x}\in \cK \times \R_{\ge 0}$ satisfying $\ov{A} \ov{x} = \ov{b}$ and
 \begin{align}
 \frac{1}{2} \ov{x}^\top \ov{Q} \ov{x} + \ov c^\top \ov x \leq \min_{ \ov{A} \ov{x} = \ov{b}, \ov{x} \in \cK \times \R_{\ge 0}} \left(\frac{1}{2} \ov{x}^\top \ov{Q} \ov{x} + \ov c^\top \ov x\right)+ \epsilon^2, \label{eqn:lem:init-point-second-cond}
 \end{align}
 the vector $\ov{x}_{1:n}$ ( $\ov{x}_{1:n}$ is the first $n$ coordinates of $\ov{x}$ ) is an approximate solution to the original convex program in the following sense:
 \begin{align*}
 \frac{1}{2} \ov{x}_{1:n}^\top Q \ov{x}_{1:n} + c^\top \ov{x}_{1:n} &\leq \min_{A x = b, x \in \cK} \left(\frac{1}{2} x^\top Q x + c^\top x\right) + \epsilon \rho^{-1}, \\
 \| A \ov{x}_{1:n} - b \|_1 &\leq  3 \epsilon \cdot ( R \| A \|_1 + \| b \|_1  ), \\
 \ov{x}_{1:n} &\in \cK.
 \end{align*}
\end{itemize}
 \end{lemma}
\begin{proof}

\textbf{First bullet point:}
Direct computation shows that $(\ov{x}, \ov{y}, \ov{s})$ is feasible.

Let us compute $\| \ov{s} + \nabla \ov{\phi}_w ( \ov{x} ) \|_{ \ov{x} }^*$.
We have
\begin{align*}
 \| \ov{s} + \nabla \ov{\phi}_w( \ov{x} ) \|_{ \ov{x} }^* = \| \epsilon \rho (c + Q x^{(0)}) \|_{ \nabla^2 \phi_w(x^{(0)})^{-1} } 
\end{align*}

Lemma~\ref{lem:LSZ_Lemma_D.1} says that for all $x \in \R^n$ with $ \| x - x^{(0)} \|_{ w,x^{(0)} } \leq 1$, we have $x\in \cK$, because $x^{(0)} = \arg\min_x \phi_w(x)$.
Therefore for any $v$ such that $v^\top \nabla^2 \phi_w( x^{(0)} ) v \leq 1$, we have $x^{(0)} \pm v \in \cK$ and hence $\| x^{(0)} \pm v \|_2 \leq R$. This implies $\| v \|_2 \leq R$ for any $v^\top \nabla^2 \phi_w(x^{(0)}) v \leq 1$. Hence $( \nabla^2 \phi_w( x^{(0)} ) )^{-1} \preceq R^2 \cdot I$. So we have
\begin{align*}
\| \ov{s} + \nabla \ov{\phi}_w( \ov{x} ) \|_{ \ov{x} } ^*
& =\| \epsilon \rho (c + Q x^{(0)}) \|_{ \nabla^2 \phi_w(x^{(0)})^{-1} } \\
& \le \epsilon \rho R \| c + Q x^{(0)} \|_2\\
&\le \epsilon \rho R (\|c\|_2 + \|Q\|_{2\to 2} \|x^{(0)}\|_2) \\
& \le \epsilon \rho R (L + LR) \\
& \le \epsilon.
\end{align*}

\textbf{Second bullet point:}
We define
\begin{align}
\OPT &:= \min_{A x = b, x \in \cK} \left(\frac{1}{2} x^\top Q x + c^\top x\right), \label{eqn:lem:init-point-original-program}\\
\ov{\OPT} &:= \min_{ \ov{A} \ov{x} = \ov{b}, \ov{x} \in \cK \times \R_{\ge 0}} \left(\frac{1}{2} \ov{x}^\top \ov{Q} \ov{x} + \ov{c}^\top \ov{x} \right). \label{eqn:lem:init-point-modified-program}
\end{align}

For any feasible $x$ in the original problem~\eqref{eqn:lem:init-point-original-program}, $\ov{x} = \begin{bmatrix} x \\ 0 \end{bmatrix}$ is feasible in the modified problem~\eqref{eqn:lem:init-point-modified-program}.
Therefore we have
\begin{align*}
\ov{\OPT} \leq \epsilon \rho (\frac 12 x^\top Qx + c^\top x) = \epsilon\rho \cdot \OPT.
\end{align*}
Given a feasible $\ov{x}$ satisfying~\eqref{eqn:lem:init-point-second-cond}, we write $\ov{x} = \begin{bmatrix} \ov{x}_{1:n} \\ \tau \end{bmatrix}$ for some $\tau \geq 0$.
Then we have
\begin{align*}
\epsilon\rho(\frac 12 \ov x_{1:n}^\top Q \ov x_{1:n} + c^\top \ov x_{1:n}) + \tau
\le \ov \OPT + \epsilon^2 \le \epsilon \rho \cdot \OPT + \epsilon^2. 
\end{align*}
Therefore
\begin{align*}
\frac 12 \ov x_{1:n}^\top Q \ov x_{1:n} + c^\top \ov x_{1:n} \le \OPT + \epsilon \rho^{-1}.
\end{align*}
We have
\begin{align*}
\tau \le -\epsilon\rho(\frac 12 \ov x_{1:n}^\top Q \ov x_{1:n} + c^\top \ov x_{1:n}) + \epsilon \rho \cdot \OPT + \epsilon^2 \le 3\epsilon
\end{align*}
because
$
\left|\frac 12 x^\top Q x + c^\top x\right| \le LR(R+1)
$
for all $x\in \cK$.

Note that $\ov x$ satisfies $A \ov x_{1:n} + (b-A x^{(0)}) \tau = b$.
So
\begin{align*}
\|A \ov x_{1:n} - b\|_1 \le \|b-A x^{(0)}\|_1 \cdot \tau.
\end{align*}
This finishes the proof.
\end{proof}

The following lemma is a generalization of~\cite[Lemma D.3]{lsz19} to quadratic program, and with weight vector $w$.
\begin{lemma}[QP version of{~\cite[Lemma D.3]{lsz19}}] \label{lem:approx-sol-from-potential}
Work under the setting of Theorem~\ref{thm:robust_ipm}.
Suppose we have $\frac{s_i} t + w_i \nabla \phi_i(x_i) = \mu_i$ for all $i\in [n]$, $- Q x + A^\top y + s = c$ and $Ax=b$.
If $\|\mu_i\|_{x_i}^* \le w_i$ for all $i\in [n]$, then we have
\begin{align*}
\frac 12 x^\top Q x + c^\top x \le \frac 12 x^{* \top} Q x^* + c^\top x^* + 4t\kappa
\end{align*}
where $x^* = \arg \min_{Ax = b, x\in \cK} \left(\frac 12 x^\top Q x + c^\top x\right)$.
\end{lemma}
\begin{proof}
    Let $x_\alpha = (1-\alpha) x + \alpha x^*$ for some $\alpha$ to be chosen.
    By Lemma~\ref{lem:LSZ_Lemma_D.1}, we have $\langle \nabla \phi_w(x_\alpha), x^*-x_\alpha\rangle \le \kappa$.
    (Note that $\phi_w$ is a $\kappa$-self-concordant barrier for $\cK$.)
    Therefore we have
    \begin{align*}
    \frac{\kappa \alpha}{1-\alpha} &\ge \langle \nabla \phi_w(x_\alpha), x_\alpha-x\rangle\\
    &= \langle \nabla \phi_w(x_\alpha) -\nabla \phi_w(x), x_\alpha-x\rangle + \langle \mu-\frac st, x_\alpha-x\rangle \\
    &\ge \sum_{i\in [n]} w_i \frac{\|x_{\alpha,i}-x_i\|_{x_i}^2}{1+\|x_{\alpha,i} - x_i\|_{x_i}} + \langle \mu, x_\alpha-x\rangle - \frac 1t \langle c-A^\top y + Qx,x_\alpha-x\rangle\\
    &\ge \sum_{i\in [n]} w_i \frac{\alpha^2 \|x_i^*-x_i\|_{x_i}^2}{1+\alpha \|x_i^* - x_i\|_{x_i}} -\alpha \sum_{i\in [n]} \|\mu_i\|_{x_i}^* \|x_i^*-x_i\|_{x_i} - \frac {\alpha}t \langle c + Qx,x^*-x\rangle.
    \end{align*}
    First step is because $\langle \nabla \phi_w(x_\alpha), x^*-x_\alpha\rangle \le \kappa$.
    Second step is because $\mu = \frac st+\nabla \phi_w(x)$.
    Third step is by Lemma~\ref{lem:LSZ_Lemma_D.1} and $c=-Qx+A^\top y+s$.
    Fourth step is by Cauchy-Schwarz and $A x_\alpha=Ax$.
    
    So we get
    \begin{align*}
    &~\frac 1t(x^\top Q x + c^\top x)\\
    \le&~ \frac 1t(x^\top Q x^* + c^\top x^*) + \frac{\kappa}{1-\alpha}
    + \sum_{i\in [n]} \|\mu_i\|_{x_i}^* \|x_i^*-x_i\|_{x_i}
    - \sum_{i\in [n]} w_i \frac{\alpha \|x_i^*-x_i\|_{x_i}^2}{1+\alpha \|x_i^* - x_i\|_{x_i}} \\
    \le&~ \frac 1t(\frac12 x^\top Q x + \frac 12 x^{*\top} Q x^* + c^\top x^*) + \frac{\kappa}{1-\alpha}
    + \sum_{i\in [n]} w_i \|x_i^*-x_i\|_{x_i}
    - \sum_{i\in [n]} w_i \frac{\alpha \|x_i^*-x_i\|_{x_i}^2}{1+\alpha \|x_i^* - x_i\|_{x_i}}\\
    =&~ \frac 1t(\frac12 x^\top Q x + \frac 12 x^{*\top} Q x^* + c^\top x^*) + \frac{\kappa}{1-\alpha}
    + \sum_{i\in [n]} w_i \frac{\|x_i^*-x_i\|_{x_i}}{1+\alpha \|x_i^* - x_i\|_{x_i}} \\
    \le &~ \frac 1t(\frac12 x^\top Q x + \frac 12 x^{*\top} Q x^* + c^\top x^*) + \frac{\kappa}{1-\alpha} + \sum_{i\in [n]} \frac{w_i}\alpha \\
    \le &~ \frac 1t(\frac12 x^\top Q x + \frac 12 x^{*\top} Q x^* + c^\top x^*) + \frac{\kappa}{\alpha(1-\alpha)}.
    \end{align*}
    First step is by rearranging terms in the previous inequality.
    Second step is by AM-GM inequality and $\|\mu_i\|_{x_i}^*\le w_i$.
    Third step is by merging the last two terms.
    Fourth step is by bounding the last term.
    Fifth step is by $\sum_{i\in [n]} w_i \le \sum_{i\in [n]} w_i \nu_i = \kappa$.

    Finally,
    \begin{align*}
    \frac 12 x^\top Q x + c^\top x &\le \frac 12 x^{*\top} Q x^* + c^\top x^* + \frac{\kappa t}{\alpha(1-\alpha)}\\
    &\le \frac 12 x^{*\top} Q x^* + c^\top x^* + 4\kappa t.
    \end{align*}
    First step is by rearranging terms in the previous inequality.
    Second step is by taking $\alpha=\frac 12$.
    This finishes the proof.
\end{proof}

\subsection{Proof of Theorem~\ref{thm:robust_ipm}}
In this section we combine everything and prove Theorem~\ref{thm:robust_ipm}.

\begin{proof}[Proof of Theorem~\ref{thm:robust_ipm}]

Lemma~\ref{lem:init-point} shows that the initial $x$ and $s$ satisfies
\begin{align*}
    \| \mu \|_{w,x}^* \le \epsilon.
\end{align*}
This implies $w_i^{-1} \|\mu_i\|_{x_i}^* \le \epsilon$ because $w_i \ge 1$.

Because Algorithm~\ref{alg:robust-ipm} first replaces $\epsilon$ by $\min\{\epsilon, \frac 1\lambda\}$ (with $\lambda$ as in Algorithm~\ref{alg:centering} for the augmented program), we have $\epsilon \le \frac 1{\lambda}$ and
\begin{align*}
\Phi(x,s,t) = \sum_{i\in [n]} \cosh(\lambda w_i^{-1} \|\mu_i\|_{x_i}^*) \le n \cosh(1) \le \cosh(\lambda/64)
\end{align*}
for the initial $x$ and $s$, by the choice of $\lambda$.

Using Theorem~\ref{thm:DLY_Thm_A.16}, we see that
\begin{align*}
    \Phi(x,s,t) \le \cosh(\lambda/64)
\end{align*}
during the entire algorithm.

So at the end of the algorithm, we have
$w_i^{-1} \|\mu_i\|_{x_i}^* \le \frac 1{64}$
for all $i\in [n]$.
In particular, $\|\mu_i\|_{x_i}^* \le w_i$ for all $i\in [n]$.

Therefore, applying Lemma~\ref{lem:approx-sol-from-potential} we get
\begin{align*}
\frac 12 x^\top Q x + c^\top x &\le \frac 12 x^{*\top} Q x^* + c^\top x^* + 4t\kappa \\
& \le \frac 12 x^{*\top} Q x^* + c^\top x^* + \epsilon^2
\end{align*}
where we used the stop condition for $t$ at the end.

So Lemma~\ref{lem:init-point} shows how to get an approximate solution for the original quadratic program with error $\epsilon LR(R+1)$.

The number of iterations is because we decrease $t$ by a factor of $1-h$ every iteration, and the choice $h = \frac{\alpha}{64 \sqrt{\kappa}}$.
\end{proof}

%% file: kernel.tex
\section{Gaussian Kernel SVM: Almost-Linear Time Algorithm and Hardness}
\label{sec:gaussian_kernel} 

In this section, we provide both algorithm and hardness for Gaussian kernel SVM problem. For the algorithm, we utilize a result due to~\cite{aa22} in conjunction with our low-rank QP solver to obtain an algorithm for the box-constrained program, with $\alpha\in[0,R]^n$, that runs in expected time $n^{1+o(1)}$ for accuracy $\epsilon\ge 1/\poly(n)$ and $R\le\poly(n)$ (Theorem~\ref{thm:gaussian_kernel} and Corollary~\ref{cor:gaussian_kernel_almost_linear}). For much smaller $\epsilon$ its time bound can be superpolynomial, for example at $d=\Theta(\log n)$ and $\epsilon=\exp(-\Theta(\log^2 n))$ (Remark~\ref{rem:gaussian_kernel_svm_parameters}). For the hardness, we build upon the framework outlined in~\cite{backurs2017fine} and improve their results in terms of dependence on dimension $d$. 

We start by proving a simple lemma that shows that if $K=UV^\top$ for low-rank $U, V$, then the quadratic objective $K\circ (yy^\top)$ also admits such a factorization via a simple scaling.

\begin{lemma}
Let $U, V\in \R^{n\times k}$ and $y\in \R^n$. Then, there exists a pair of matrices $\wt U, \wt V\in \R^{n\times k}$ such that
\begin{align*}
    \wt U\wt V^\top = & ~ (UV^\top)\circ (yy^\top)
\end{align*}
moreover, $\wt U, \wt V$ can be computed in time $O(nk)$.
\end{lemma}

\begin{proof}
The proof relies on the following identity for Hadamard product: for any matrix $A$ and conforming vectors $x, y$ (all real), one has
\begin{align*}
    A\circ (yx^\top) = & ~ D_y A D_x
\end{align*}
where $D_y, D_x\in \R^{n\times n}$ are diagonal matrices that put $y, x$ on their diagonals. Thus, we can simply compute $\wt U, \wt V$ as follows:
\begin{align*}
    \wt U = & ~ D_y U, \\
    \wt V = & ~ D_y V,
\end{align*}
consequently, 
\begin{align*}
    \wt U \wt V^\top = & ~ D_y U V^\top D_y \\
    = & ~ (yy^\top) \circ (UV^\top) \\
    = & ~ (UV^\top) \circ (yy^\top),
\end{align*}
as desired. Moreover, the diagonal scaling of $U, V$ can be indeed performed in $O(nk)$ time, as advertised.
\end{proof}

Throughout this section, we will let $B$ denote the squared radius of the dataset.

\subsection{Almost-Linear Time Algorithm for Gaussian Kernel SVM}

We state a result due to~\cite{aa22}, in which they present an optimal-degree polynomial approximation to the function $e^{-x}$ and consequentially, this produces an efficient approximate scheme to the Batch Gaussian Kernel Density Estimation problem.

We start by introducing a notion that captures the minimum degree polynomial that well-approximates $e^{-x}$:

\begin{definition}
Let $f: [0, B]\rightarrow \R$, we let $q_{B;\epsilon}(f)\in \mathbb{N}$ denote the minimum degree of a non-constant polynomial $p(x)$ such that
\begin{align*}
    \sup_{x\in [0, B]} |p(x)-f(x)|\leq \epsilon
\end{align*}
\end{definition}

Utilizing the Chebyshev polynomial machinery together with the orthgonal polynomial families,~\cite{aa22} provides the following characterization on $q_{B;\epsilon}(f)$:

\begin{theorem}[Theorem 1.2 of~\cite{aa22}]
Let $B\geq 1$ and $\epsilon\in (0,1)$. Then
\begin{align*}
    q_{B;\epsilon}(e^{-x}) = & ~ \Theta(\max\{\sqrt{B\log(1/\epsilon)}, \frac{\log(1/\epsilon)}{\log(B^{-1} \log(1/\epsilon))} \})
\end{align*}
\end{theorem}

\begin{theorem}[Corollary 1.7 of~\cite{aa22}]
\label{thm:aa22_main}
Let $x_1,\ldots,x_n\in \R^d$ be a dataset with squared radius $B$ and $\epsilon\in (0,1)$. Let $q=q_{B;\epsilon}(e^{-x})$. Let $K\in \R^{n\times n}$ be the Gaussian kernel matrix formed by $x_1,\ldots,x_n$. Finally, let $k=\binom{2d+2q}{2q}$. Then, there exists a deterministic algorithm that computes a pair of matrices $U, V\in \R^{n\times k}$ such that for any vector $v\in \R^n$, 
\begin{align*}
    \|Kv-UV^\top v\|_\infty \leq & ~ \epsilon \|v\|_1.
\end{align*}
Moreover, matrices $U, V$ can be computed in time $O(nkd)$.
\end{theorem}

Even though $\ell_\infty$ error in terms of $\ell_1$ norm of vector $v$ seems quite weak, it can be conveniently translated into more standard guarantees, e.g., spectral norm error. The following lemma provides a conversion of errors that come in handy later when integrating the kernel approximation to our low-rank QP solver.

\begin{lemma}
\label{lem:linf_to_spectral}
Let $K\in \R^{n\times n}$ be a PSD kernel matrix and $\epsilon\in (0,1)$ be a parameter. Let $\wt K\in \R^{n\times n}$ be an approximation to $K$ with the guarantee that for any $v\in \R^n$,
\begin{align*}
    \|Kv-\wt Kv\|_\infty \leq & ~ \epsilon \|v\|_1,
\end{align*}
then
\begin{align*}
    |v^\top Kv-v^\top \wt Kv| \leq \epsilon \|v\|_1^2 \leq \epsilon n \|v\|_2^2.
\end{align*}
\end{lemma}

\begin{proof}
The proof is a simple application of H{\"o}lder's inequality:
\begin{align*}
    |v^\top (Kv-\wt Kv)| = & ~ |\langle v, Kv-\wt Kv\rangle |\\
    \leq & ~ \|v\|_1 \|Kv-\wt Kv\|_\infty \\
    \leq & ~ \epsilon \|v\|_1^2 \\
    \leq & ~ \epsilon n \|v\|_2^2,
\end{align*}
where the second step is by H{\"o}lder's inequality, and the last step is by Cauchy-Schwarz. This completes the proof.
\end{proof}

We can now combine the Gaussian kernel low-rank decomposition with our low-rank QP solver to provide an almost-linear time algorithm for Gaussian kernel SVM. We restate the kernel SVM formulation here.

\begin{definition}[Restatement of Definition~\ref{def:kernel_svm}]
Given a data matrix $X \in \R^{n \times d}$ and labels $y \in \R^n$. 
Let $Q \in \R^{n \times n}$ denote a matrix where $Q_{i,j} = \mathsf{K}(x_i,x_j) \cdot y_i y_j $ for $i,j\in [n]$. The hard-margin kernel SVM problem with bias asks to solve the following program.
\begin{align*}
    \max_{\alpha\in \R^n}~&~{\bf 1}_n^\top \alpha -\frac{1}{2} \alpha^\top Q \alpha \\
    {\rm s.t.~} & ~ \alpha^\top y = 0 \\
    & ~ \alpha \geq 0.
\end{align*}
\end{definition}

\begin{theorem}
\label{thm:gaussian_kernel}
Let ${\sf K}(x_i,x_j)=\exp(-\|x_i-x_j\|_2^2)$ be the Gaussian kernel, let $K\in\R^{n\times n}$ be its kernel matrix on $x_1,\ldots,x_n\in\R^d$, and let $Q=D_yKD_y$, where $D_y$ is the diagonal matrix with the labels $y\in\{\pm 1\}^n$ on its diagonal. Suppose that $\|x_i-x_j\|_2^2\le B$ for all $i,j\in[n]$, where $B\ge 1$, and that both labels occur among $y_1,\ldots,y_n$. (If all labels are equal, the only feasible point below is $0$.) Let $\epsilon\in(0,1)$ and $R>0$, and let ${\rm OPT}_R$ be the optimal value of the box-constrained program
\begin{align}\label{eqn:gaussian_box}
    \max_{\alpha\in \R^n}~&~{\bf 1}_n^\top \alpha -\frac{1}{2} \alpha^\top Q \alpha \\
    {\rm s.t.~} & ~ \alpha^\top y = 0, \notag\\
    & ~ \alpha\in[0,R]^n. \notag
\end{align}
Let $\epsilon_1:=\min\{\frac12, c_0\epsilon/(n^2R^2)\}$, where $c_0>0$ is an absolute constant fixed in the proof, let $q=q_{B;\epsilon_1}(e^{-x})$, and let $k=\binom{2d+2q}{2q}$. Then there is a randomized algorithm whose output $\wh \alpha\in[0,R]^n$ always satisfies
\begin{align*}
    {\bf 1}_n^\top \wh \alpha-\frac{1}{2} \wh \alpha^\top Q \wh \alpha \geq & ~ {\rm OPT}_R-\epsilon , \\
    |\wh \alpha^\top y| \leq & ~ 3\epsilon,
\end{align*}
and whose expected running time is
\begin{align*}
    \wt O\big((nk^{(\omega+1)/2}+n^{1/2}k^{\omega})\log(n(1+R)/\epsilon)\big).
\end{align*}
When $k$ is much smaller than $n$, precisely when $(2k+1)^{\omega-1}\le n$, the expected running time is
\begin{align*}
    \wt O\big(nk^{(\omega+1)/2}\log(n(1+R)/\epsilon)\big).
\end{align*}
With $R=C$, program~\eqref{eqn:gaussian_box} is the soft-margin program of Definition~\ref{def:C_SVC}. If some optimal solution $\alpha^*$ of the hard-margin program~\eqref{eqn:kernel_svm} satisfies $\|\alpha^*\|_\infty\le R$, then ${\rm OPT}_R={\rm OPT}$, the optimal value of program~\eqref{eqn:kernel_svm}.
\end{theorem}

\begin{proof}
Let $n_+$ and $n_-$ be the numbers of indices $i$ with $y_i=1$ and with $y_i=-1$. Both are at least $1$, so $n\ge 2$. We take $c_0=1/8$, so that $\epsilon_1\le\frac12$ and $\epsilon_1\le \epsilon/(8n^2R^2)$. Write
\begin{align*}
f(\alpha):={\bf 1}_n^\top\alpha-\frac12\alpha^\top Q\alpha, \qquad \mathcal{F}_R:=\{\alpha\in[0,R]^n:\alpha^\top y=0\}.
\end{align*}
The set $\mathcal{F}_R$ is compact and contains $0$, so ${\rm OPT}_R=\max_{\alpha\in\mathcal{F}_R}f(\alpha)$ is attained. Every $\alpha\in[0,R]^n$ satisfies
\begin{align}\label{eq:box_norm}
\|\alpha\|_2^2\le nR^2 .
\end{align}

\textbf{Step 1: a positive semidefinite low-rank proxy for $Q$.} Theorem~\ref{thm:aa22_main} uses the squared radius only through the fact that every squared distance $\|x_i-x_j\|_2^2$ lies in $[0,B]$, so it applies with our $B\ge 1$. We invoke it with accuracy $\epsilon_1\in(0,1)$ and obtain, in time $O(nkd)$, matrices $U,V\in\R^{n\times k}$ such that $\wt K:=UV^\top$ satisfies $\|Kv-\wt Kv\|_\infty\le\epsilon_1\|v\|_1$ for all $v\in\R^n$. Let $\wt Q:=D_y\wt KD_y=(D_yU)(D_yV)^\top$. Lemma~\ref{lem:linf_to_spectral} with $v=D_yx$, together with $\|D_yx\|_2=\|x\|_2$ for $y\in\{\pm1\}^n$, gives
\begin{align}\label{eq:exact_approx_Q_gap}
|x^\top(Q-\wt Q)x|=|(D_yx)^\top(K-\wt K)(D_yx)|\le \epsilon_1 n\|x\|_2^2 \qquad \text{for all } x\in\R^n.
\end{align}
Theorem~\ref{thm:rank-formal} needs a symmetric positive semidefinite matrix, but the matrix $\wt Q$ need not be symmetric or positive semidefinite. Let $W\in\R^{n\times k'}$, with $k'\le\min\{n,2k\}$, have orthonormal columns that span the column space of the $n\times 2k$ matrix $[D_yU,\ D_yV]$. The symmetric matrix $S:=\frac12(\wt Q+\wt Q^\top)$ has its column space inside that of $W$, so $S=WMW^\top$ with $M:=W^\top SW\in\R^{k'\times k'}$, and $x^\top Sx=x^\top\wt Qx$ for all $x\in\R^n$. For $u\in\R^{k'}$ let $x=Wu$, so that $\|x\|_2=\|u\|_2$. Since $Q$ is positive semidefinite, \eqref{eq:exact_approx_Q_gap} gives
\begin{align*}
u^\top Mu=x^\top\wt Qx\ge x^\top Qx-\epsilon_1n\|x\|_2^2\ge-\epsilon_1n\|u\|_2^2 ,
\end{align*}
so every eigenvalue of $M$ is at least $-\epsilon_1n$. Let $M_+:=M+\epsilon_1nI_{k'}$. Then $M_+$ is positive semidefinite and $\|M-M_+\|=\epsilon_1n$. Set $\wt Q_+:=WM_+W^\top=S+\epsilon_1nWW^\top$. For every $x\in\R^n$,
\begin{align*}
|x^\top(Q-\wt Q_+)x|
\le & ~ |x^\top(Q-S)x|+|(W^\top x)^\top(M-M_+)(W^\top x)| \\
\le & ~ \epsilon_1n\|x\|_2^2+\epsilon_1n\|W^\top x\|_2^2 \\
\le & ~ 2\epsilon_1n\|x\|_2^2 .
\end{align*}
Since $Q-\wt Q_+$ is symmetric, this gives
\begin{align}\label{eq:psd_proxy_gap}
\|Q-\wt Q_+\|\le 2\epsilon_1 n .
\end{align}
The matrix $\wt Q_+$ is symmetric positive semidefinite and has the factorization $\wt Q_+=FG^\top$ with $F:=W$ and $G:=WM_+$ in $\R^{n\times k'}$. Computing $W$ (by orthogonalizing the $2k$ columns of $[D_yU,\ D_yV]$), the matrix $M=\frac12 W^\top(D_yU)(D_yV)^\top W+\frac12 W^\top (D_yV)(D_yU)^\top W$, and $G$ takes $O(nk^{\omega-1})$ time, because $k'\le\min\{n,2k\}$. No eigendecomposition is needed.

\textbf{Step 2: the norm of $\wt Q_+$.} The Gaussian kernel matrix $K$ is positive semidefinite with $K_{i,i}=1$, so $Q=D_yKD_y$ is positive semidefinite and
\begin{align*}
\|Q\|\le\tr[Q]=\tr[K]=n .
\end{align*}
By~\eqref{eq:psd_proxy_gap} and $\epsilon_1\le\frac12$,
\begin{align*}
\|\wt Q_+\|\le\|Q\|+\|Q-\wt Q_+\|\le(1+2\epsilon_1)n\le 2n.
\end{align*}

\textbf{Step 3: solving the proxy program.} We apply Theorem~\ref{thm:rank-formal} to the program
\begin{align}\label{eqn:gaussian_box_proxy}
\min_{\alpha\in\R^n}~&~\frac12\alpha^\top\wt Q_+\alpha-{\bf 1}_n^\top\alpha \\
{\rm s.t.~} &~ y^\top\alpha=0, \notag\\
&~\alpha_i\in\cK_i:=[0,R] \qquad\text{for all } i\in[n], \notag
\end{align}
that is, with $m=1$, $A=y^\top$, $b=0$ and linear term $-{\bf 1}_n$. On each $\cK_i$ we use the barrier $\phi_i(a)=-\log a-\log(R-a)$, which is a $2$-self-concordant barrier for the interval $[0,R]$, so $\nu_i=2$. We take $w_i=1$, so $\kappa=2n$. The hypotheses of Theorem~\ref{thm:rank-formal} hold with the following parameters.
\begin{itemize}
\item Inner radius $r:=R/(2n)$. Let $z_i:=R/(2n_+)$ if $y_i=1$ and $z_i:=R/(2n_-)$ if $y_i=-1$. Then $y^\top z=n_+\cdot\frac{R}{2n_+}-n_-\cdot\frac{R}{2n_-}=0$, and $z_i\in[\frac R{2n},\frac R2]$ for all $i$. Every point of $B(z,r)$ has its $i$-th coordinate in $[z_i-\frac R{2n},z_i+\frac R{2n}]\subseteq[0,R]$, so $B(z,r)\subseteq\cK=[0,R]^n$.
\item Outer radius $R_o:=R\sqrt n$, by~\eqref{eq:box_norm}.
\item Lipschitz constant $L:=2n$, because $\|\wt Q_+\|\le 2n$ by Step 2 and $\|{\bf 1}_n\|_2=\sqrt n\le 2n$.
\item Low rank: we are given $\wt Q_+=FG^\top$ with $F,G\in\R^{n\times k'}$ and $k'\le 2k$.
\end{itemize}
We use the accuracy parameter
\begin{align*}
\epsilon_{\rm QP}:=\min\Big\{\frac12,\ \frac{\epsilon}{4nR_o(R_o+1)}\Big\}.
\end{align*}
Theorem~\ref{thm:rank-formal} returns some $\wh\alpha\in[0,R]^n$. Let $\wt f(\alpha):={\bf 1}_n^\top\alpha-\frac12\alpha^\top\wt Q_+\alpha$. The first guarantee of Theorem~\ref{thm:rank-formal}, multiplied by $-1$, reads
\begin{align}\label{eq:gaussian_proxy_value}
\wt f(\wh\alpha)\ge\max_{\alpha\in\mathcal{F}_R}\wt f(\alpha)-\epsilon_{\rm QP}LR_o(R_o+1)\ge \max_{\alpha\in\mathcal{F}_R}\wt f(\alpha)-\frac{\epsilon}{2}.
\end{align}
Its second guarantee, with $\|A\|_1=\|y\|_1=n$ (the entrywise $\ell_1$ norm) and $b=0$, reads
\begin{align*}
|\wh\alpha^\top y|\le 3\epsilon_{\rm QP}R_on\le\frac{3\epsilon}{4(R_o+1)}\le 3\epsilon .
\end{align*}

\textbf{Step 4: from the proxy objective to the true objective.} For every $\alpha\in[0,R]^n$, the bounds~\eqref{eq:psd_proxy_gap}, \eqref{eq:box_norm} and $\epsilon_1\le\epsilon/(8n^2R^2)$ give
\begin{align}\label{eq:gaussian_objective_gap}
|f(\alpha)-\wt f(\alpha)|=\frac12|\alpha^\top(\wt Q_+-Q)\alpha|\le\epsilon_1n\|\alpha\|_2^2\le\epsilon_1n^2R^2\le\frac\epsilon8 .
\end{align}
Let $\alpha^R\in\mathcal{F}_R$ attain ${\rm OPT}_R$. Then
\begin{align*}
f(\wh\alpha)\ge & ~ \wt f(\wh\alpha)-\frac\epsilon8 \\
\ge & ~ \max_{\alpha\in\mathcal{F}_R}\wt f(\alpha)-\frac{5\epsilon}8 \\
\ge & ~ \wt f(\alpha^R)-\frac{5\epsilon}8 \\
\ge & ~ f(\alpha^R)-\frac{3\epsilon}4 \\
\ge & ~ {\rm OPT}_R-\epsilon ,
\end{align*}
where the first and fourth steps use~\eqref{eq:gaussian_objective_gap}, and the second step uses~\eqref{eq:gaussian_proxy_value}.

\textbf{Step 5: running time.} Step 1 takes $O(nkd+nk^{\omega-1})$ time. Since $q\ge1$, we have $k=\binom{2d+2q}{2q}\ge\binom{2d+2}{2}\ge d^2$, so $d\le k^{1/2}\le k^{(\omega-1)/2}$ and $nkd\le nk^{(\omega+1)/2}$. Also $nk^{\omega-1}\le nk^{(\omega+1)/2}$ because $\omega\le 3$. In Step 3, the first running time bound of Theorem~\ref{thm:rank-formal}, which does not use the assumption $(k+m)^{\omega-1}\le n$, with $\kappa=2n$, $\nu_i=2$, $w_i=1$, $m=1$ and rank $k'\le 2k$ gives expected time
\begin{align*}
\wt O\Big( & \big(\sqrt2\log(n\kappa R_o/(r\epsilon_{\rm QP}))+(3+\sqrt2)\log(\kappa/\epsilon_{\rm QP})\big)\cdot n(2k+1)^{(\omega+1)/2} \\
& +\sqrt{2n}\,\log(n\kappa R_o/(r\epsilon_{\rm QP}))\cdot(2k+1)^{\omega}\Big).
\end{align*}
Here $R_o/r=2n^{3/2}$. Moreover $R_o(R_o+1)\le nR(R+1)\le n(1+R)^2$, so $1/\epsilon_{\rm QP}\le\max\{2,4n^2(1+R)^2/\epsilon\}=4n^2(1+R)^2/\epsilon$. Hence both logarithms are $O(\log(n(1+R)/\epsilon))$, and the total expected time is
\begin{align*}
\wt O\big((nk^{(\omega+1)/2}+n^{1/2}k^{\omega})\log(n(1+R)/\epsilon)\big).
\end{align*}
If $(2k+1)^{\omega-1}\le n$, then $(2k+1)^{\omega}\le n^{1/2}(2k+1)^{(\omega+1)/2}$, so $n^{1/2}k^{\omega}=O(nk^{(\omega+1)/2})$ and the total expected time is $\wt O\big(nk^{(\omega+1)/2}\log(n(1+R)/\epsilon)\big)$.

\textbf{The last two claims.} With $R=C$, program~\eqref{eqn:gaussian_box} is the program of Definition~\ref{def:C_SVC}, whose matrix $(yy^\top)\circ K$ equals $D_yKD_y=Q$. Now let $\alpha^*$ be an optimal solution of program~\eqref{eqn:kernel_svm} with $\|\alpha^*\|_\infty\le R$. Then $\alpha^*\in\mathcal{F}_R$, so ${\rm OPT}_R\ge f(\alpha^*)={\rm OPT}$. Conversely, $\mathcal{F}_R$ lies inside the feasible set of program~\eqref{eqn:kernel_svm}, so ${\rm OPT}_R\le{\rm OPT}$.
\end{proof}

\begin{corollary}[Almost-linear time]
\label{cor:gaussian_kernel_almost_linear}
In the setting of Theorem~\ref{thm:gaussian_kernel}, suppose in addition that $d=O(\log n)$, $B=o(\frac{\log n}{\log\log n})$, $\epsilon\ge n^{-O(1)}$ and $R\le n^{O(1)}$. Then $k=n^{o(1)}$, and the algorithm of Theorem~\ref{thm:gaussian_kernel} runs in expected time $n^{1+o(1)}$. In particular, in expected time $n^{1+o(1)}$ one can find $\wh\alpha$ with $|\wh\alpha^\top y|\le 3\epsilon$ whose objective value is at least the optimum minus $\epsilon$, both for the soft-margin program of Definition~\ref{def:C_SVC} with $C\le n^{O(1)}$ and for the hard-margin program~\eqref{eqn:kernel_svm} when some optimal solution satisfies $\|\alpha^*\|_\infty\le R$ for a known $R\le n^{O(1)}$.
\end{corollary}

\begin{proof}
All asymptotic statements are as $n\to\infty$, and logarithms are natural. Fix constants $a,b\ge 1$ with $d\le a\log n$, $\epsilon\ge n^{-b}$ and $R\le n^{b}$, and let $\beta(n):=B/\log n$, so that $\beta(n)\to 0$. Assume $n\ge 8$.

\textbf{The degree.} With $c_0=1/8$ we have $c_0\epsilon/(n^2R^2)\ge n^{-(3b+2)}/8\ge n^{-(3b+3)}$ and $\frac12\ge n^{-(3b+3)}$, so $\epsilon_1\ge\delta:=n^{-(3b+3)}$. A polynomial that approximates $e^{-x}$ on $[0,B]$ within $\delta$ also does so within $\epsilon_1$, so $q=q_{B;\epsilon_1}(e^{-x})\le q_{B;\delta}(e^{-x})$. Let $\ell:=\log(1/\delta)=(3b+3)\log n$. Theorem 1.2 of~\cite{aa22}, quoted above, gives
\begin{align*}
q_{B;\delta}(e^{-x})=O\Big(\max\Big\{\sqrt{B\ell},\ \frac{\ell}{\log(\ell/B)}\Big\}\Big).
\end{align*}
Here $\sqrt{B\ell}=\sqrt{(3b+3)\beta(n)}\,\log n=o(\log n)$. Also $\ell/B=(3b+3)/\beta(n)\to\infty$, so $\log(\ell/B)\to\infty$ and $\ell/\log(\ell/B)=o(\log n)$. Hence $\gamma(n):=q/\log n$ satisfies $\gamma(n)\to 0$, and $\gamma(n)>0$ because $q\ge 1$.

\textbf{The rank.} The bound $\binom{N}{j}\le(eN/j)^j$ gives $k=\binom{2d+2q}{2q}\le(e(1+d/q))^{2q}$, so
\begin{align*}
\log k\le 2q\big(1+\log(1+d/q)\big)\le 2q\big(1+\log(1+a\log n/q)\big).
\end{align*}
Substituting $q=\gamma(n)\log n$ gives
\begin{align*}
\log k\le 2\gamma(n)\big(1+\log(1+a/\gamma(n))\big)\cdot\log n=o(\log n),
\end{align*}
since $x(1+\log(1+a/x))\to 0$ as $x\to 0^+$. So $k=n^{o(1)}$.

\textbf{The time.} Since $1+R\le 2n^b$ and $\epsilon\ge n^{-b}$, we have $\log(n(1+R)/\epsilon)=O(\log n)$. Theorem~\ref{thm:gaussian_kernel} then gives expected time $\wt O((n\cdot n^{o(1)}+n^{1/2}\cdot n^{o(1)})\cdot\log n)=n^{1+o(1)}$. The last sentence of the corollary follows from the last two claims of Theorem~\ref{thm:gaussian_kernel}, with $R=C$ for the soft-margin program.
\end{proof}

\begin{remark}
\label{rem:gaussian_kernel_svm_parameters}
We discuss the parameters of Theorem~\ref{thm:gaussian_kernel} and Corollary~\ref{cor:gaussian_kernel_almost_linear}.
\begin{itemize}
    \item \textit{Bounded radius.} Let $d=O(\log n)$, $B=\Theta(1)$, $\epsilon=n^{-\Theta(1)}$ and $1\le R\le n^{O(1)}$. Then $\log(1/\epsilon_1)=\Theta(\log n)$, and the theorem of~\cite{aa22} quoted above gives the degree $q=\Theta(\frac{\log n}{\log\log n})$. The rank is then
    \begin{align*}
        k = \binom{2d+2q}{2q} \leq (e(1+d/q))^{2q} = 2^{O(\log n \log\log\log n/\log\log n)} = n^{o(1)},
    \end{align*}
    and the algorithm runs in expected time $n^{1+o(1)}$, as Corollary~\ref{cor:gaussian_kernel_almost_linear} states for every $B=o(\frac{\log n}{\log\log n})$.
    \item \textit{Small $\epsilon$.} Corollary~\ref{cor:gaussian_kernel_almost_linear} covers $\epsilon\ge n^{-O(1)}$. For much smaller $\epsilon$ and $d=\Theta(\log n)$, the rank $k$, and with it the time bound of Theorem~\ref{thm:gaussian_kernel}, can become superpolynomial. For example, let $\epsilon=\exp(-\Theta(\log^2 n))$, $d=\Theta(\log n)$, $B=\Theta(1)$ and $1\le R\le n^{O(1)}$. Then $\log(1/\epsilon_1)=\Theta(\log^2 n)$, so $q=\Theta(\frac{\log^2 n}{\log\log n})$. The rank $k=\binom{2d+2q}{2d}$ lies between $(q/d)^{2d}$ and $(e(1+q/d))^{2d}$, so $k=n^{\Theta(\log\log n)}$. We do not know an algorithm that runs in time $n^{1+o(1)}\poly\log(1/\epsilon)$ in this regime. For kernel density estimation, \cite{aa22} note that other standard techniques, such as the fast multipole method and core-sets, also lead to time $n^{1+\Omega(1)}$ when $d=\Omega(\log n)$ and the error is below $n^{-\Theta(1)}$.
    \item \textit{The box.} For the soft-margin program the box is part of the problem, with $R=C$. For the hard-margin program, it needs a known bound $R\ge\|\alpha^*\|_\infty$ for some optimal solution $\alpha^*$. For instance, if the data are separable in feature space with margin $\gamma$, then $\|\alpha^*\|_1=1/\gamma^2$, as for linear SVMs (Theorem~\ref{thm:main-rank-informal}). So $R=1/\gamma_0^2$ works for any known $\gamma_0\le\gamma$, and $R\le n^{O(1)}$ when $\gamma_0\ge n^{-O(1)}$. The inner radius needed by Theorem~\ref{thm:rank-formal} is automatic once both labels occur (see the proof of Theorem~\ref{thm:gaussian_kernel}).
\end{itemize}

The runtime we obtain can be viewed as a consequence of the ``phase transition'' phenomenon observed in~\cite{aa22}, in which they prove that if $B=\omega(\log n)$, then quadratic time in $n$ is essentially needed to approximate the Gaussian kernel assuming \seth. 
\end{remark}

\subsection{Hardness of Gaussian Kernel SVM with Large Radius}

In this section, we show that for $d=\Theta(\log n)$, any algorithm that solves the program associated to hard-margin Gaussian kernel SVM without the bias term up to error $\exp(-\Theta(\log^2 n))$ would require $\Omega(n^{2-o(1)})$ time on some inputs with $B=\Theta(\log^2 n)$ (Theorem~\ref{thm:svm_hard_no_bias}). This shows that, at error $\exp(-\Theta(\log^2 n))$, some upper bound on $B$ is needed. A gap remains between the $B=o(\frac{\log n}{\log\log n})$ of Corollary~\ref{cor:gaussian_kernel_almost_linear} and $B=\Theta(\log^2 n)$. A second gap concerns the accuracy, since we prove almost-linear time only for $\epsilon\ge 1/\poly(n)$. At error $\exp(-\Theta(\log^2 n))$ and $d=\Theta(\log n)$, our time bound is superpolynomial even for $B=\Theta(1)$ (Remark~\ref{rem:gaussian_kernel_svm_parameters}). To prove the hardness result, we need to introduce the approximate Hamming nearest neighbor problem.

\begin{definition}
For $\delta>0$ and $n, d\in \mathbb{N}$, let $\{a_1,\ldots,a_n\}, \{b_1,\ldots,b_n\}\subseteq \{0,1\}^d$ be two sets of vectors, and let $t\in \{0,1,\ldots,d\}$ be a threshold. The \emph{$(1+\delta)$-Approximate Hamming Nearest Neighbor} problem asks to distinguish the following two cases:
\begin{itemize}
    \item If there exists some $a_i$ and $b_j$ such that $\|a_i-b_j\|_1\leq t$, output ``Yes'';
    \item If for any $i,j\in [n]$, we have $\|a_i-b_j\|_1>(1+\delta)\cdot t$, output ``No''.
\end{itemize}
\end{definition}

Note that the algorithm can output any value if the datasets fall in neither of these two cases. We will utilize the following hardness result due to Rubinstein.

\begin{theorem}[\cite{r18}]
\label{thm:r18}
Assuming \seth, for every $q>0$, there exists $\delta>0$ and $C>0$ such that $(1+\delta)$-Approximate Hamming Nearest Neighbor in dimension $d=C\log n$ requires time $\Omega(n^{2-q})$.
\end{theorem}
\begin{remark}
Rubinstein proves Theorem~\ref{thm:r18} by a reduction from the Orthogonal Vectors problem. He states the theorem as holding unless both \seth{} and the Orthogonal Vectors hypothesis fail~\cite{r18}. The theorem therefore needs only the Orthogonal Vectors hypothesis, which \seth{} implies, and so do the lower bounds below that use it (Remark~\ref{rem:eth-2026}).
\end{remark}

A final ingredient is a rewriting of the dual SVM into its primal form, without resorting to optimize over an infinite-dimensional hyperplane.

\begin{lemma}
Consider the dual hard-margin kernel SVM defined as
\begin{align*}
    \max_{\alpha\in \R^n} & ~ {\bf 1}^\top \alpha-\frac{1}{2} \sum_{i,j\in [n]\times [n]} \alpha_i \alpha_j y_iy_j {\sf K}(w_i, w_j), \\
    \textnormal{s.t.} & ~ \alpha^\top y = 0, \\
    & ~ \alpha \geq 0.
\end{align*}
The primal program can be written as
\begin{align*}
    \min_{\alpha\in \R^n} & ~ \frac{1}{2}\sum_{i,j\in [n]\times [n]} \alpha_i\alpha_jy_iy_j {\sf K}(w_i,w_j), \\
    \textnormal{s.t.} & ~ y_i f(w_i) \geq 1, \\
    & ~ \alpha \geq 0,
\end{align*}
where $f:\R^d\rightarrow \R$ is defined as
\begin{align*}
    f(w) = & ~ \sum_{j=1}^n \alpha_j y_j {\sf K}(w_j, w)-b.
\end{align*}
Moreover, the primal and dual program has no duality gap and the optimal solution $\alpha$ to both programs are the same.
\end{lemma}

\begin{proof}
Recall that the primal hard-margin SVM is the following program:
\begin{align*}
    \min_{v} & ~ \frac{1}{2}\|v\|_2^2, \\
    \text{s.t.} & ~ y_i (v^\top \phi(w_i)-b) \geq 1,
\end{align*}
where $b\in \R$ is the bias term and $\phi:\R^d\rightarrow \R^K$ is the feature mapping corresponding to the kernel in the sense that ${\sf K}(w_i, w_j)=\phi(w_i)^\top \phi(w_j)$. The optimal weight $v=\sum_{i=1}^n \alpha_i y_i \phi(w_i)$ where $\alpha\in \R^n$ is the optimal solution to the dual program. Consequently, 
\begin{align*}
    \|v\|_2^2 = & ~ (\sum_{i=1}^n \alpha_i y_i \phi(w_i))^2 \\
    = & ~ \sum_{i,j\in [n]\times [n]}\alpha_i\alpha_j y_iy_j \phi(w_i)^\top \phi(w_j) \\
    = & ~ \sum_{i,j\in [n]\times [n]} \alpha_i\alpha_jy_iy_j {\sf K}(w_i, w_j) \\
    = & ~ \alpha^\top Q \alpha,
\end{align*}
where the matrix $Q$ is the usual
\begin{align*}
    Q = & ~ (yy^\top) \circ K,
\end{align*}
the constraint can be rewritten as
\begin{align*}
    y_i(v^\top \phi(w_i)-b) = & ~ y_i ((\sum_{i=1}^n \alpha_i y_i \phi(w_i))^\top \phi(w_i)-b) \\
    = & ~ y_i(\sum_{j=1}^n \alpha_j y_j \phi(w_j)^\top \phi(w_i))-y_ib \\
    = & ~ y_i (\sum_{j=1}^n \alpha_j y_j {\sf K}(w_i, w_j))-y_ib \\
    = & ~ y_if(w_i),
\end{align*}
where $f:\R^d\rightarrow \R$ is defined as
\begin{align*}
    f(w) = & ~ \sum_{j=1}^n \alpha_j y_j {\sf K}(w_j, w)-b.
\end{align*}
Thus, we can alternatively write the primal as
\begin{align*}
    \min_{\alpha\in \R^n} & ~ \frac{1}{2} \alpha^\top Q \alpha, \\
    \text{s.t.} & ~ y_i f(w_i) \geq 1.
\end{align*}
For the strong duality and optimal solution, see, e.g.,~\cite{mmrts01}.
\end{proof}

We will now prove the almost-quadratic lower bound for Gaussian kernel SVM. Our proof strategy is similar to that of~\cite{backurs2017fine} with different set of parameters. It is also worth noting that the construction of~\cite{backurs2017fine} without the bias term requires dimension $d=\omega(\log n)$ and squared dataset radius $B=\Theta(d\log n)=\omega(\log^2 n)$. Their construction with the bias term appends $\Theta(\log^3 n)$ entries to each vector, so it has $d=\Theta(\log^3 n)$ and $B=\Theta(\log^4 n)$. Without the bias term, we improve both parameters, to $d=\Theta(\log n)$ and $B=\Theta(\log^2 n)$. With the bias term, we improve the dimension to $d=\Theta(\log n)$, at the price of a larger squared radius $B=\Theta(\log^6 n)$ (Corollary~\ref{cor:svm_hard_bias}).

\begin{theorem}
\label{thm:svm_hard_no_bias}
Assuming \seth, for every $q>0$, there exists a hard-margin Gaussian kernel SVM without the bias term, i.e., the program of Definition~\ref{def:kernel_svm} with the constraint $\alpha^\top y=0$ removed, with $d=\Theta(\log n)$ and error $\epsilon=\exp(-\Theta(\log^2 n))$ for inputs whose squared radius is at most $B=\Theta( \log^2 n)$ requiring time $\Omega(n^{2-q})$ to solve. 
\end{theorem}

\begin{proof}
Let $l=\sqrt{2(c'\delta)^{-1}\log n}$, where $c'>0$ is an absolute constant with $c'\delta\le 1/6$. We will provide a reduction from $(1+\delta)$-Approximate Hamming Nearest Neighbor to Gaussian kernel SVM. Let $A:=\{a_1,\ldots,a_n\}, B:=\{b_1,\ldots,b_n\}\subseteq \{0,1\}^d$ be the datasets, we assign label 1 to all vectors $a_i$ and label $-1$ to all vectors $b_j$, moreover, we scale both $A$ and $B$ by $l$, this results in two datasets with points in $\{0,l\}^d$. The squared radius of this dataset is then
\begin{align*}
    B = & ~ \max \{\max_{i,j}\|la_i-la_j\|_2^2,  \max_{i,j}\|lb_i-lb_j\|_2^2,\max_{i,j}\|la_i-lb_j\|_2^2\} \\
    \leq & ~ l^2 d \\
    = & ~ \Theta(\delta^{-1} \log^2 n).
\end{align*}
To simplify the notation, we will implictly assume $A$ and $B$ are scaled by $l$ without explicitly writing out $la_i$, $lb_j$. Now consider the following three programs:
\begin{itemize}
    \item Classifying $A$:
    \begin{align}\label{eq:prog_A}
        \min_{\alpha\in \R^n_{\geq 0}} & ~ \frac{1}{2} \sum_{i,j\in [n]\times [n]} \alpha_i\alpha_j{\sf K}(a_i, a_j), \notag\\
        \text{s.t.} & ~ \sum_{j=1}^n \alpha_j {\sf K}(a_i,a_j)\geq 1, \qquad\forall i\in [n]
    \end{align}
    \item Classifying $B$:
    \begin{align}\label{eq:prog_B}
        \min_{\beta\in \R^n_{\geq 0}} & ~ \frac{1}{2} \sum_{i,j\in [n]\times [n]} \beta_i\beta_j{\sf K}(b_i, b_j), \notag\\
        \text{s.t.} & ~ -\sum_{j=1}^n \beta_j {\sf K}(b_i,b_j)\leq -1, \qquad\forall i\in [n]
    \end{align}
    \item Classifying both $A$ and $B$:
    \begin{align}\label{eq:prog_AB}
        \min_{\alpha,\beta\in \R_{\geq 0}^n} & ~ \frac{1}{2}\sum_{i,j\in [n]\times [n]} \alpha_i\alpha_j{\sf K}(a_i, a_j)+\frac{1}{2} \sum_{i,j\in [n]\times [n]} \beta_i\beta_j{\sf K}(b_i, b_j)-\sum_{i,j\in [n]\times [n]} \alpha_i\beta_j {\sf K}(a_i, b_j),\notag \\
        \text{s.t.} & ~ \sum_{j=1}^n \alpha_j {\sf K}(a_i,a_j)-\sum_{j=1}^n \beta_j {\sf K}(a_i,b_j) \geq 1, \qquad \forall i\in [n],\notag \\
        & ~ \sum_{j=1}^n \alpha_j {\sf K}(b_i,a_j)-\sum_{j=1}^n \beta_j {\sf K}(b_i,b_j) \leq -1, \qquad \forall i\in [n]
    \end{align}
\end{itemize}
We will prove that the optimal solution $\alpha^*_i$'s and $\beta^*_i$'s are both lower and upper bounded. Use ${\rm Val}(A), {\rm Val}(B)$ and ${\rm Val}(A,B)$ to denote the value of program~\eqref{eq:prog_A},~\eqref{eq:prog_B} and~\eqref{eq:prog_AB} respectively, then note that
\begin{align*}
    {\rm Val}(A) \leq & ~ \frac{n^2}{2}
\end{align*}
by plugging in $\alpha={\bf 1}$ and setting all vectors to be the same. On the other hand, 
\begin{align*}
    {\rm Val}(A) \geq & ~ \frac{1}{2}\sum_{i=1}^n (\alpha_i^*)^2 {\sf K}(a_i,a_i) \\
    = & ~ \frac{1}{2}\sum_{i=1}^n (\alpha_i^*)^2.
\end{align*}
Combining these two, we can conclude that for any $\alpha_i^*$, it must be $\alpha_i^*\leq n$. For the lower bound, consider the inequality constraint for the $i$-th point:
\begin{align*}
    \alpha_i^*+\sum_{j\neq i} \alpha_j^* {\sf K}(a_i,a_j) \geq & ~ 1,
\end{align*}
to estimate ${\sf K}(a_i, a_j)$, note that $\|a_i-a_j\|_2^2=\|a_i-a_j\|_1\geq 1$ for $j\neq i$,\footnote{We without loss of generality that during preprocess, we have remove duplicates in $A$ and $B$.} and 
\begin{align*}
    {\sf K}(a_i, a_j) = & ~ \exp(-l^2 \|a_i-a_j\|_2^2) \\
    = & ~ \exp(-2(c'\delta)^{-1}\log n\|a_i-a_j\|_1) \\
    \leq & ~ \exp(-2(c'\delta)^{-1}\log n) \\
    \leq & ~ n^{-10}/100,
\end{align*}
combining with $\alpha_j^*\leq n$, we have
\begin{align*}
    \alpha_i^* \geq & ~ 1-\sum_{j\neq i}\alpha_j^* {\sf K}(a_i,a_j) \\
    \geq & ~ 1-n\cdot n\cdot n^{-10}/100 \\
    \geq & ~ 1/2.
\end{align*}
This lower bound on $\alpha_i^*$ is helpful when we attempt to lower bound ${\rm Val}(A,B)$ with ${\rm Val}(A)+{\rm Val}(B)$. Following the outline of~\cite{backurs2017fine}, we consider the three dual programs:
\begin{itemize}
    \item Dual of classifying $A$:
    \begin{align}\label{eq:prog_dual_A}
        \max_{\alpha\in \R^n_{\geq 0}} & ~ \sum_{i=1}^n \alpha_i-\frac{1}{2}\sum_{i,j} \alpha_i\alpha_j {\sf K}(a_i,a_j)
    \end{align}
    \item  Dual of classifying $B$:
    \begin{align}\label{eq:prog_dual_B}
        \max_{\beta\in \R^n_{\geq 0}} & ~ \sum_{i=1}^n \beta_i-\frac{1}{2}\sum_{i,j} \beta_i\beta_j {\sf K}(b_i,b_j)
    \end{align}
    \item Dual of classifying $A$ and $B$:
    \begin{align}\label{eq:prog_dual_AB}
        \max_{\alpha,\beta\in \R^n_{\geq 0}} & ~ \sum_{i=1}^n \alpha_i+\sum_{i=1}^n \beta_i-\frac{1}{2}\sum_{i,j} \alpha_i\alpha_j {\sf K}(a_i,a_j)-\frac{1}{2}\sum_{i,j} \beta_i\beta_j {\sf K}(b_i,b_j)+\sum_{i,j}\alpha_i\beta_j {\sf K}(a_i,b_j)
    \end{align}
\end{itemize}
as the SVM program exhibits strong duality, we know that the optimal value of the primal equals to the dual, so we can alternatively bound ${\rm Val}(A,B)$ using the dual program. Plug in $\alpha^*,\beta^*$ to program~\eqref{eq:prog_dual_AB}, we have
\begin{align*}
    {\rm Val}(A,B) \geq & ~ \sum_{i=1}^n \alpha_i^*+\sum_{i=1}^n \beta_i^*-\frac{1}{2}\sum_{i,j} \alpha_i^*\alpha_j^* {\sf K}(a_i,a_j)-\frac{1}{2}\sum_{i,j} \beta_i^*\beta_j^* {\sf K}(b_i,b_j)+\sum_{i,j}\alpha_i^*\beta_j^* {\sf K}(a_i,b_j) \\
    = & ~ {\rm Val}(A)+{\rm Val}(B)+\sum_{i,j}\alpha_i^*\beta_j^* {\sf K}(a_i,b_j),
\end{align*}
to bound the third term, we consider the pair $(a_{i_0},b_{j_0})$ such that $\|a_{i_0}-b_{j_0}\|_1\leq t$, and note that
\begin{align*}
    \sum_{i,j}\alpha_i^*\beta_j^* {\sf K}(a_i,b_j) \geq & ~ \alpha_{i_0}^*\beta_{j_0}^* {\sf K}(a_{i_0},b_{j_0}) \\
    \geq & ~ \frac{1}{4} \exp(-2(c'\delta)^{-1}\log n\cdot t).
\end{align*}
To wrap up, we have
\begin{align*}
    {\rm Val}(A,B) \geq & ~ {\rm Val}(A)+{\rm Val}(B)+\frac{1}{4}\exp(-2(c'\delta)^{-1}\log n\cdot t)
\end{align*}
We now prove the ``No'' case, where for any $a_i, b_j$, $\|a_i-b_j\|_1>(1+\delta)t$. Every such distance exceeds $(1+\delta)t$ and is an integer, hence is at least $t+1$. We have
\begin{align*}
    {\sf K}(a_i,b_j) = & ~ \exp(-l^2 \|a_i-b_j\|_2^2) \\
    \leq & ~ \exp(-2(c'\delta)^{-1} \log n\cdot (t+1)),
\end{align*}
we let $m:=\exp(-2(c'\delta)^{-1}\log n\cdot (t+1))$, set $\alpha'=\alpha^*+10n^2m$ and $\beta'=\beta^*+10n^2m$, we let $V$ to denote the value when evaluating program~\eqref{eq:prog_AB} with $\alpha',\beta'$. We will essentially show that 
\begin{align*}
    {\rm Val}(A,B) \leq & ~ V
\end{align*}
and 
\begin{align*}
    V \leq & ~ {\rm Val}(A)+{\rm Val}(B)+400n^6m
\end{align*}
chaining these two gives us a certificate for the ``No'' case. To prove the first assertion, we show that $\alpha',\beta'$ are feasible solution to program~\eqref{eq:prog_AB} since
\begin{align*}
    \sum_{j=1}^n \alpha'_j {\sf K}(a_i,a_j) = & ~ \sum_{j=1}^n (\alpha^*_j{\sf K}(a_i,a_j)+10n^2 m{\sf K}(a_i,a_j)) \\
    = & ~ \alpha^*_i+10n^2m+\sum_{j\neq i} (\alpha^*_j+10n^2m){\sf K}(a_i,a_j) \\
    \geq & ~ \alpha_i^*+\sum_{j\neq i} \alpha_j^*{\sf K}(a_i,a_j)+10n^2m \\
    = & ~ 10n^2 m+\sum_{j=1}^n \alpha^*_j{\sf K}(a_i,a_j) \\
    \geq & ~ 10n^2m+1
\end{align*}
where we use $\alpha_i^*$ satisfy the inequality constraint of program~\eqref{eq:prog_A}. We compute an upper bound on $\sum_{j=1}^n \beta_j'{\sf K}(a_i,b_j)$:
\begin{align*}
    \sum_{j=1}^n \beta_j'{\sf K}(a_i,b_j) \leq & ~ \sum_{j=1}^n 2n m \\
    = & ~ 2n^2 m,
\end{align*}
where we use the fact that $m=\exp(-2(c'\delta)^{-1}\log n\cdot (t+1))\leq n^{-10}/10$ therefore $\beta^*+10n^2m\leq 2n$. Thus, it must be the case that
\begin{align*}
    \sum_{j=1}\alpha_j' {\sf K}(a_i,a_j)-\sum_{j=1}^n \beta'_j{\sf K}(a_i,b_j) \geq & ~ 8n^2m+1 \\
    \geq & ~ 1,
\end{align*}
as desired. The other linear constraint follows by a symmetric argument. This indeed shows that $\alpha',\beta'$ are feasible solutions to program~\eqref{eq:prog_AB} and ${\rm Val}(A,B)\leq V$.

To prove an upper bound on $V$, we note that
\begin{align*}
    V = & ~ \frac{1}{2}\sum_{i,j}\alpha_i'\alpha_j'{\sf K}(a_i,a_j)+\frac{1}{2}\sum_{i,j}\beta_i'\beta_j'{\sf K}(b_i,b_j)-\sum_{i,j}\alpha'\beta_j'{\sf K}(a_i,b_j) \\
    \leq & ~ \frac{1}{2}\sum_{i,j}\alpha_i'\alpha_j'{\sf K}(a_i,a_j)+\frac{1}{2}\sum_{i,j}\beta_i'\beta_j'{\sf K}(b_i,b_j),
\end{align*}
we bound the first quantity, as the second follows similarly:
\begin{align*}
    \frac{1}{2}\sum_{i,j}\alpha_i'\alpha_j'{\sf K}(a_i,a_j) = & ~ \frac{1}{2}\sum_{i,j}(\alpha_i^*\alpha_j^*+10n^2m(\alpha_i^*+\alpha_j^*)+100n^4m^2){\sf K}(a_i,a_j) \\
    \leq & ~ {\rm Val}(A)+\sum_{i,j}10n^3m{\sf K}(a_i,a_j)+\sum_{i,j}100n^4m^2 {\sf K}(a_i,a_j) \\
    \leq & ~ {\rm Val}(A)+10n^5m+100n^6m^2 \\
    \leq & ~ {\rm Val}(A)+200n^6m,
\end{align*}
we can thus conclude
\begin{align*}
    V \leq & ~ {\rm Val}(A)+{\rm Val}(B)+400n^6m.
\end{align*}
Chaining these two, we obtain the following threshold for the ``No'' case:
\begin{align*}
    {\rm Val}(A,B)\leq & ~ {\rm Val}(A)+{\rm Val}(B)+400n^6m.
\end{align*}
Finally, we observe that 
\begin{align*}
    400n^6\exp(-2(c'\delta)^{-1}\log n\cdot (t+1)) \le & ~ \frac{1}{8}\exp(-2(c'\delta)^{-1}\log n\cdot t),
\end{align*}
since $\exp(-2(c'\delta)^{-1}\log n)\le n^{-12}$, so we can distinguish these two cases. 

Note that when one considers solving the program with additive error, we need to make sure that the error is smaller than the distinguishing threshold. Since three programs are solved and the gap is at least $\frac 18\exp(-2(c'\delta)^{-1}\log n\cdot t)$, it suffices that
\begin{align*}
    \epsilon \leq & ~ \frac{1}{64}\exp(-2(c'\delta)^{-1}d\log n) \\
    \leq & ~ \frac{1}{64}\exp(-2(c'\delta)^{-1}t\log n),
\end{align*}
where we use $t\leq d$. Since $d=\Theta(\log n)$, the error $\exp(-\Theta(\log^2 n))$ suffices. This concludes the proof.
\end{proof}

\begin{remark}
Our proof can be interpreted as using a stronger complexity theoretical tool in place of the one used by~\cite{backurs2017fine}, to obtain a better dependence on dimension $d$ and $B$. We also note that the construction due to~\cite{backurs2017fine} has the relation that $B=\Theta(d\log n)$, this is because in order to lower bound ${\rm Val}(A,B)$, one has to lower bound the optimal values of $\alpha_i^*$'s and $\beta_j^*$'s. To do so, one needs to further scale up $a_i$'s and $b_j$'s so that within datasets $A$ and $B$, the radius is at least $\Theta(\log n)$. This is in contrast to the Batch Gaussian KDE studied in~\cite{aa22}, where they show the almost-quadratic lower bound can be achieved for both $d, B=\Theta(\log n)$.
\end{remark}
\begin{remark}[The hypothesis after October 2026]
\label{rem:eth-2026}
If the proofs announced in October 2026 are correct (see the note on hardness assumptions in the introduction), Theorem~\ref{thm:svm_hard_no_bias} has a false hypothesis and, as stated, gives no lower bound.

The lower bound itself survives under a weaker hypothesis. The proof above uses \seth{} only through Theorem~\ref{thm:r18}. Rubinstein derives that theorem by a reduction from the Orthogonal Vectors problem~\cite{r18}. Hence Theorem~\ref{thm:svm_hard_no_bias} remains true when \seth{} is replaced by the Orthogonal Vectors hypothesis, which states that for every $\delta>0$ there is $c>0$ such that no algorithm decides in time $O(N^{2-\delta})$ whether two given sets of $N$ vectors in $\{0,1\}^{c\log N}$ contain an orthogonal pair. The same holds for Corollaries~\ref{cor:svm_hard_bias} and~\ref{cor:svm_soft_bias}, which are derived from the same hard instances. The lower bound of~\cite{backurs2017fine} that we cite reduces from an exact bichromatic Hamming closest-pair problem. That problem is at least as hard as the approximate problem of Theorem~\ref{thm:r18}, so this lower bound also holds under the Orthogonal Vectors hypothesis. The lower bound of~\cite{aa22} that we cite also uses \seth{} only through Theorem~\ref{thm:r18}~\cite[Theorem~5.2]{aa22}, so it also holds under the Orthogonal Vectors hypothesis. \seth{} implies the Orthogonal Vectors hypothesis~\cite{williams05}, but the converse is not known. Neither announcement claims to refute the Orthogonal Vectors hypothesis. The Anthropic paper states that its main theorem ``does not refute the Orthogonal Vectors hypothesis''~\cite{anthropic26eth}. The announcements change the evidence for this hypothesis. Its main support was that \seth{} implies it, and it now stands as a conjecture in its own right. The algorithms of this paper, including Theorems~\ref{thm:main-rank-informal} and~\ref{thm:gaussian_kernel}, do not depend on any of these hypotheses.
\end{remark}

Similar to~\cite{backurs2017fine}, we obtain hardness results for hard-margin kernel SVM with bias, and soft-margin kernel SVM with bias. 

\begin{corollary}
\label{cor:svm_hard_bias}
Assuming \seth, for every $q>0$, there exists a hard-margin Gaussian kernel SVM with the bias term with $d=\Theta(\log n)$ and error $\epsilon=\exp(-\Theta(\log^2 n))$ for inputs whose squared radius is at most $B=\Theta( \log^6 n)$ requiring time $\Omega(n^{2-q})$ to solve. 
\end{corollary}

\begin{proof}[Proof sketch]
The proof is similar to~\cite{backurs2017fine}. Given a hard instance of Theorem~\ref{thm:svm_hard_no_bias}, except we append $\Theta(\log n)$ entries with magnitude $\Theta(\log^2 n)$ instead of distributing the values across $\Theta(\log^3 n)$ entries. Rest of the proof follows exactly the same as~\cite{backurs2017fine}.
\end{proof}

\begin{corollary}
\label{cor:svm_soft_bias}
Assuming \seth, for every $q>0$, there exists a soft-margin Gaussian kernel SVM with the bias term with $d=\Theta(\log n)$ and error $\epsilon=\exp(-\Theta(\log^2 n))$ for inputs whose squared radius is at most $B=\Theta( \log^6 n)$ requiring time $\Omega(n^{2-q})$ to solve. 
\end{corollary}
\begin{proof}[Proof sketch]
As in~\cite{backurs2017fine}, this follows from the hard-margin instance with bias of Corollary~\ref{cor:svm_hard_bias}. We do not write out the details, so this corollary, like Corollary~\ref{cor:svm_hard_bias}, rests on a proof sketch.
\end{proof}
\begin{remark}
Corollary~\ref{cor:svm_hard_bias} modifies the hard instance of Theorem~\ref{thm:svm_hard_no_bias}. Corollary~\ref{cor:svm_soft_bias} follows from the hard-margin instance with bias, as in~\cite{backurs2017fine}. Both corollaries therefore use \seth{} only through Theorem~\ref{thm:r18}, and by Remark~\ref{rem:eth-2026} they also hold under the Orthogonal Vectors hypothesis.
\end{remark}

\begin{remark}
Compared to the construction of~\cite{backurs2017fine} in which they distribute a total mass of $\Theta(\log^3 n)$ across $\Theta(\log^3 n)$ entries so that they ensure after the reduction, the vectors take values in $\{-1,0,1\}$, we instead distribute the mass across $\Theta(\log n)$ entries so that each entry has magnitude $\Theta(\log^2 n)$. To make the reduction work, the total mass of $\Theta(\log^3 n)$ is needed, and for~\cite{backurs2017fine}, it is fine to append an extra $\Theta(\log^3 n)$ entries, as their hardness result for hard-margin SVM with bias is stated for $d=O(\log^3 n)$. For us, we need to restrict $d=\Theta(\log n)$ at the price of each entry has a larger magnitude of $\Theta(\log^2 n)$. This blows up the squared radius from $\log^2 n$ to $\log^6 n$. We note that the~\cite{backurs2017fine} construction has squared radius $\log^4 n$.
\end{remark}